 \documentclass[10pt,reqno]{article}

\usepackage{amscd,amssymb}
\usepackage[mathscr]{eucal}
\usepackage[all]{xy}
\usepackage{mathrsfs}
\usepackage{xypic}
\usepackage{amsfonts}
\usepackage{amsmath}
\usepackage{amsthm}
\usepackage{latexsym}
\usepackage{eufrak}

\theoremstyle{plain}
\newtheorem{theorem}[equation]{Theorem}
\newtheorem{corollary}[equation]{Corollary}
\newtheorem{proposition}[equation]{Proposition}
\newtheorem{lemma}[equation]{Lemma}

\theoremstyle{definition}
\newtheorem{definition}[equation]{Definition}

\newtheorem{claim}[equation]{Claim}

\theoremstyle{remark}
\newtheorem{notation}[equation]{Notation}
\newtheorem{remark}[equation]{Remark}

\makeatletter
\renewcommand{\subsection}{\@startsection{subsection}{2}{0pt}{-3ex
plus -1ex minus -0.2ex}{-2mm plus -0pt minus
-2pt}{\normalfont\bfseries}} \makeatother

\numberwithin{equation}{subsection}

\newcommand{\beq}{\begin{equation}\label}
\newcommand{\eeq}{\end{equation}}
\newcommand{\idot}{_{\:\raisebox{1pt}{\text{\circle*{1.5}}}}}
\newcommand{\hdot}{^{\:\raisebox{3pt}{\text{\circle*{1.5}}}}}
\newcommand{\hhdot}{{\:\raisebox{3pt}{\text{\circle*{1.5}}}}}

\renewcommand{\proof}[1][Proof.\quad]{\smallskip\noindent{\em #1}}
\def\endproof{\hfill\ensuremath{\square}\par\medskip}

\newcommand{\nc}{\newcommand}

\CompileMatrices 
\def\op{\operatorname}
\newcommand{\tR}{R}
\newcommand{\Gm}{{\mathbb{G}}_{\mathbf{m}}}
\newcommand{\dis}{\displaystyle} 

\newcommand{\tU}{U}
\newcommand{\tB}{B}
\newcommand{\scra}{{\mathscr{A}}}
\newcommand{\scrc}{{\mathscr{C}}}

\newcommand{\Co}{{\mathsf{Comp}}^G}
\newcommand{\fff}{{\mathbb{F}}}

\newcommand{\adh}{{\mathtt {Ad}}_{_{\sf hopf}}}
\newcommand{\K}{{{\mathcal{K}}}} 
\newcommand{\Mod}{\op{{Mod}}}

\newcommand{\pone}{{\mathsf{P}}}

\newcommand{\fcal}{{\mathcal{F}}}
\newcommand{\hht}{{\mathtt{ht}}}
\newcommand{\lo}{{\,\overset{L}{\otimes}\,}}
 
\newcommand{\ev}{{\mathsf{ev}}}
\newcommand{\mt}{{\mathsf{m}}}
\def\Mt{\widetilde{\mathcal{M}}}
\newcommand{\St}{{\mathfrak{L}}}
\newcommand{\ste}{{\mathfrak{L}}\red}

\nc{\U}{{\mathsf U}} 
\newcommand{\opp}{^{^{_{\sf{op}}}}}

\def\pb{$\bullet\quad$\parbox[t]{142mm}}
\def\pbox{\parbox[t]{142mm}}
 
\newcommand{\limind}{{\textsl{lim}^{\,}\textsl{ind}^{\,}}}
\newcommand{\limproj}{{\textsl{lim}^{\,}\textsl{proj}^{\,}}}
\newcommand{\DD}{D_{\mathsf{proj}}}

\newcommand{\alg}{^{\op{alg}}}

\renewcommand{\SS}{{\mathsf{S}}}
\newcommand{\uSS}{{\underline{\mathsf{S}}}}

\newcommand{\tK}{{\mathsf{R}}}
\newcommand{\EE}{{\mathsf{E}}}
\newcommand{\uK}{{\underline{\mathsf{R}}}}

\newcommand{\La}{{\mathsf{\Lambda}}}
\newcommand{\CA}{{{\boldsymbol{\Gamma}}[\NN]}}
\newcommand{\mcac}{\op{Mod}^{G\times\C^*\!}(\CA)}
\newcommand{\uaa}{\underline{\bb}}
\newcommand{\mca}{\op{Mod}^{G\times\Gm}(\CA)}
\newcommand{\mcaf}{\op{Mod}^{G\times\Gm}_{\tt{thin}}(\CA)}

\newcommand{\ka}{{\tK\bb}}

\newcommand{\agr}{{\op{Mod}^{B\times\Gm}_f(\La)}}
\newcommand{\ag}{{\op{Mod}^B_f(\La)}}

\newcommand{\dgcoh}{\op{DG}_{\text{coherent}}}
\newcommand{\dcoh}{D_{\text{coherent}}}

\newcommand{\wh}{\widehat}
\newcommand{\block}{{\textsf{block}}}
\newcommand{\Ugh}{{\mathcal{U}}_k{\wh{\g}}}
\newcommand{\bi}{{\mathbf{i}}}

\nc{\Rhom}{\op{ RHom}\bul}
\nc{\End}{\op{ End}} 
\renewcommand{\AA}{{\mathsf A}}

\renewcommand{\aa}{{\mathsf a}} 
\nc{\tilb}{\til{\mathsf b}}
\nc{\al}{\alpha}
\nc{\fU}{{\mathfrak{U}}}
\nc{\fB}{{\mathfrak{B}}}
\nc{\Ind}{\op{ Ind}}
\nc{\lotimes}{{\stackrel{_L}{\otimes}}}
\nc{\Coind}{\op{ Coind}} 
\nc\blangle{{\boldsymbol{\langle}}}
\nc\brangle{{\boldsymbol{\rangle}}}
\nc{\Ker}{\op{ Ker}}
\nc{\im}{\op{ Im}} 
\nc{\Coker}{\op{ Coker}}
\nc{\dirlim}{\underset{\rightarrow}{\op{ lim}}}
\nc{\invlim}{\underset{\leftarrow}{\op{ lim}}} 
\nc{\CN}{{\mathcal N}}
\nc{\fA}{{\mathfrak{A}}}
\nc{\fC}{{Z}}
\nc{\ff}{{\mathscr{F}}}
\nc{\lth}{{l} t} 
\nc{\Tor}{\op{ Tor}\idot}
\nc{\qis}{\,\stackrel{\tt{qis}}{\,\sim\,}\,}
\nc{\qisto}{\stackrel{\tt{qis}}{\too}}
\nc{\mult}{{\tt{mult}}}
\nc{\Wm}{{\mathsf{W}}}
\nc{\barU}{{\widehat{\mathcal U}\g}}
\nc{\hb}{{\widehat{\Ub}}}
\nc{\uhb}{{\k_{_\B}}}
\nc{\zz}{{\mathfrak{z}}} 
\nc{\til}{\widetilde}
\nc{\fI}{{\mathfrak{I}}}
\nc{\bbI}{{\mathbf{I}}}
\nc{\red}{'}
\nc{\X}{{\mathbb X}} \nc{\Y}{{\mathbb Y}} 
\nc{\YY}{{\fr\Y}}

\nc{\ten}{{\otimes}} 
\nc{\tenl}{\overset{{\mathtt L}}\ten}
\nc{\map}{\longrightarrow}

\nc{\triv}{\mathsf{triv}} 
\nc{\bs}{\bigskip\\} 
\nc{\ms}{\smallskip\\}
\nc{\tP}{{P'}}

 \def\beq{\begin{equation}}

 \def\eeq{\end{equation}}
\newcommand{\mto}{\longmapsto}
\newcommand{\s}{t}

\newcommand{\p}{{\mathcal P}}
\newcommand{\iso}{{\;\stackrel{_\sim}{\longrightarrow}\;}}
\newcommand{\cd}{\!\cdot\!}

\newcommand{\on}{{\mathcal{O}_{_{\widetilde{\mathcal{N}}}}}}
\newcommand{\B}{{\mathsf{B}}}

\def\vi{{\sf{(i)}\,\,}}
\def\vii{{\sf{(ii)}\,\,}} 
\def\viii{{\sf{(iii)}\,\,}}
\def\iv{{\sf{(iv)}\,\,}}
\def\bplus{{\mbox{$\bigoplus$}}}

\newcommand{\longeq}{{\;{=\!\!\!=\!\!\!=\!\!\!=\!\!\!=\!\!\!=\!\!\!=}\;}}
\newcommand{\fh}{{\mathfrak{h}}}
\nc{\vp}{\varpi} 
\nc{\li}{l}
\nc{\supp}{\op{{supp}^{\,}}}
 \nc{\starb}{\star}
\nc{\starg}{\star} 
\nc{\stara}{\star_{\!\!_{{\mathtt{act}}}}}

\def\DGMod{{\op{DGM}}}

\newcommand{\one}{{\mathbf{1}_{_{\Gr}}}}

\newcommand{\id}{\op{{id}}}
 \newcommand{\Id}{\op{{Id}}}

\newcommand{\m}{{\mathfrak{m}}}
\newcommand{\fZ}{{\mathfrak{Z}}}
\newcommand{\into}{\;\hookrightarrow\;}
\newcommand{\onto}{\,\twoheadrightarrow\,}

\newcommand{\R}{{\mathcal{R}}} 

\newcommand{\too}{\;\longrightarrow\;} 
\newcommand{\g}{{\mathfrak{g}}}

\newcommand{\ind}{{\protect{\op{Ind}}}}

\newcommand{\indf}{{\protect{\op{Ind}}}}
\newcommand{\rind}{{\protect{\op{RInd}}}}

\newcommand{\res}{{{\op{Res}}}}
\newcommand{\sym}{{\op{{Sym}}}}
 \newcommand{\fr}{{}^{\phi_{\!}}}
\newcommand{\D}{\op{D}}

\newcommand{\fF}{{\mathfrak{F}}}
\newcommand{\bff}{{\mathfrak{F}}\red}
\newcommand{\umod}{{\op{Rep}}({\mathsf u})}
\newcommand{\Umod}{{\op{Rep}}({\mathsf U})}
\newcommand{\modt}{{\op{Rep}}}

 \newcommand{\mmod}{\mbox{-}\op{ {mod}}}
 \newcommand{\bimod}{\mbox{-}{\op{bimod}}}
\newcommand{\rep}{{\op{Rep}}} 

 \newcommand{\bb}{{{\mathsf b}}}
\newcommand{\n}{{\mathfrak n}}
 \newcommand{\NN}{{\widetilde{{\mathcal N}}}} 
  
 \newcommand{\oo}{{\mathcal{O}}}
\newcommand{\dd}{{\mathcal{D}}}
\newcommand{\modB}{{D_{_{\mathsf{triv}}}(\B)}}

\newcommand{\eps}{_\epsilon}
\newcommand{\Rt}{\op{{R}}} 
\newcommand{\RHom}{\op{{RHom}}}
\newcommand{\REnd}{\op{{REnd}}}
\newcommand{\RInd}{\op{{RInd}}} 
\newcommand{\dg}{$\op{{dg}}$}
\newcommand{\bF}{{\mathbb{F}}} 
\newcommand{\st}{\ltimes}

\newcommand{\CC}{{\mathcal C}}

\newcommand{\bfF}{{\boldsymbol{\mathfrak{F}}}}

\newcommand{\mix}{^{^{\sf {mix}}}} 

\newcommand{\tg}{{\mathfrak t}}
 
\newcommand{\Gd}{{G^\vee}}
\newcommand{\Bd}{{B^\vee}}
\newcommand{\Td}{{T^\vee}}
\newcommand{\cU}{{\mathcal{U}}}
\newcommand{\tb}{{\mathsf{p}}}

 \newcommand{\tp}{P}
\newcommand{\Hoch}{H\!H}
\def\C{{\mathbb{C}}}
  
 \def\q{\Y}

 \def\Z{{\mathbb{Z}}} 
 \def\I{{\mathbf I}} 
 \def\J{{{\mathbf I}\mbox{-}{\mathsf{mon}}}}
\def\baf{{\mathcal B}} 
\def\baft{{\widetilde{\mathcal B}}}
 
\def\FI{{\mathcal B}}
\def\FIt{\widetilde{\mathcal B}}
\def\OO{{\mathcal O} } 
\def\K{{\mathcal K}} 
\def\k{{\Bbbk}} 
 
\def\Oe{{\mathcal N}^{^{\tt{reg}}}}
\def\A{{\mathcal{A}}}

\def\GG{{\mathsf{G}}}

\def\g{\mathfrak{g}}
\def\N{{\mathcal N}} 
\def\Ug{{{\mathcal{U}}\g}}
\def\Ub{{{\mathcal{U}}{\mathfrak{b}}}} 
\def\Uh{{{\mathcal{U}}{\mathfrak{t}}}}
\def\Un{{{\mathcal{U}}{\mathfrak{n}}}} 
\def\yh{_\Y^{\Ub}}
\def\u{{{\mathsf u}}} 

\def\cat{{\mathsf{block}}({\mathsf{U}})}
\def\catm{{\mathsf{block}}^{^{\mathsf{mix}}}({\mathsf{U}})}

\def\cohgm{{{{\mathcal{C}}oh}^{G\times\Gm}({\mathcal N})}}

\def\hu{\widehat{{\mathcal{U}}\mathfrak{g}}}
 
\def\Coh{{{\mathcal{C}}oh}}
\def\b{\mathfrak{b}}

\newcommand{\tooo}{{\;{-\!\!\!-\!\!\!-\!\!\!-\!\!\!\longrightarrow}\;}}
 
\def\PI{{{{\mathcal{P}}\!{erv}}}_{_\I}({\mathsf {Gr}})}
\def\PJ{{{{\mathcal{P}}\!{erv}}}_{_\J}({\mathsf {Gr}})} 
\def\Pmix{{{{\mathcal{P}}\!{erv}}}^{^{\mathsf{mix}}}({\mathsf {Gr}})}
\def\Perv{{{{\mathcal{P}}\!{erv}}}}
\def\pe{{\mathsf{Perv}}}
\def\M{{\mathcal M}}
\def\PO{{{\mathcal{P}}\!{erv}}_{_{G^{\!_{^{\tiny{\tiny\vee}}\!}}\!(\OO)}}
({\mathsf {Gr}})}
\def\BI{{{\mathcal{P}}\!{erv}}_{_\I}({\mathcal{B}})}
 
\def\DI{D^b_{_\I}}
 
\def\W{{\mathcal{W}}}
\def\bW{\overline{{\mathcal{W}}}}
\newcommand{\GO}{_{_{G^{\!_{^{\tiny{\tiny\vee}}\!}}\!({\mathcal O})}}}
 \def\bs{\overline{\sigma}}

\def\F{{\mathcal F}} 
\def\L{{\mathcal L}}

\def\zz{{\mathcal{Z}}}

\def\a{{\mathfrak{z}}}
\def\b{{\mathfrak b}} 
\def\la{\lambda}

\def\g{{\mathfrak{g}}} 
\def\h{{\mathfrak{t}}}

\def\e{{e}} \def\f{{f}}

  \def\ad{\op{{ad}^{\,}}}
\def\rk{{{\sl {rk}}^{\,}}} \def\Hom{\op{{Hom}}}
\def\Ext{\op{{Ext}}} 

 \def\gee{{\g}^e} \def\emod{{\rep\hdot(G^e)}}
 
\newcommand{\ab}{{\hskip 8mm}}

\newcommand{\Gr}{{\mathsf {Gr}}} \newcommand{\gr}{\op{{gr}^{}}}
\newcommand{\gd}{\check{\mathfrak{g}}}
\newcommand{\hd}{{\check{\mathfrak{t}}}} 
\newcommand{\sll}{{\mathfrak{s}\mathfrak{l}}_2}
\newcommand{\Lie}{{\mbox{\textsl{Lie}}^{\,}}}

\newcommand{\waf}{W_{\mbox{\tiny aff}}}

 \def\triv{{\mathsf {triv}}}
\newcommand{\sB}{{\mathsf{B}}} 

\newcommand{\BL}{{\pmb{\mathcal{F}}\boldsymbol\ell}}

\newcommand{\sset}{\subset} \newcommand{\Spec}{\op{{Spec}}}
\newcommand{\greg}{{\mathfrak{g}}^{^{\mathsf {reg}}}}
\newcommand{\z}{{\mathcal G}} \newcommand{\Ad}{\op{{Ad}_{\,}}}

\newcommand{\cO}{\fB_{\!_\A}} 
\newcommand{\cB}{{\mathsf{B}}_{\!_\A}}

\def\ccirc{{{}_{^{\,^\circ}}}}

\def\square{\hbox{\vrule\vbox{\hrule\phantom{o}\hrule}\vrule}}

\begin{document}

\setlength{\parindent}{0pt}
\setlength{\parskip}{3pt plus 5pt minus 0pt}

\centerline{\Large {\bf Quantum Groups, the loop Grassmannian,}}
\vskip 4pt

\centerline{\Large {\bf and the Springer resolution.}}

\vskip 10mm
\centerline{\large {\sc Sergey Arkhipov, Roman Bezrukavnikov, and Victor Ginzburg}}
\vskip 1cm


\begin{abstract}{\footnotesize
\noindent
We establish  equivalences of the following
 three triangulated  categories:
$$
D_\text{quantum}(\g)\enspace\longleftrightarrow\enspace
D^G_\text{coherent}(\NN)\enspace\longleftrightarrow\enspace
D_\text{perverse}(\Gr).
$$
Here, $D_\text{quantum}(\g)$ is the derived category
of the principal block of finite dimensional representations
of the quantized enveloping algebra
(at an odd root of
unity)  of
a complex semisimple Lie algebra $\g$;
the category $D^G_\text{coherent}(\NN)$ is defined 
in terms of 
coherent sheaves on the cotangent bundle
on the (finite dimensional) flag manifold for $G$ (= semisimple group
with Lie algebra~$\g$), and the  category 
$D_\text{perverse}(\Gr)$ is the derived category of perverse
sheaves on the Grassmannian $\Gr$ associated with
the loop group $LG^\vee$, where $G^\vee$ is the
Langlands dual  group, smooth along the Schubert stratification.

The equivalence between $D_\text{quantum}(\g)$ and $D^G_\text{coherent}(\NN)$ is  
 an `enhancement' of the known expression (due to Ginzburg-Kumar) for 
quantum group cohomology in terms of nilpotent variety. 
The equivalence between $D_\text{perverse}(\Gr)$ and $D^G_\text{coherent}(\NN)$
 can be viewed as  a `categorification'
 of the
isomorphism between two  completely different  geometric realizations
of the (fundamental polynomial representation
of  the)
affine Hecke algebra that has played a key role in the proof of
the Deligne-Langlands-Lusztig conjecture.
One realization is in terms
of locally constant functions on the flag manifold of a $p$-adic
reductive group, while the other  is in terms of 
equivariant $K$-theory of a complex (Steinberg) variety 
for the dual group.

The composite of the two  equivalences above
yields an equivalence between
{\em abelian} categories of
quantum group representations and perverse sheaves. A similar 
equivalence at an even root of unity can be deduced, following
 Lusztig program, from  earlier
deep
results of Kazhdan-Lusztig and Kashiwara-Tanisaki.
Our approach is independent of these results
and is totally different (it does not rely on representation theory
of Kac-Moody algebras). It also gives way to proving Humphreys'
conjectures on tilting $U_q(\g)$-modules, as will be explained in a 
separate paper.
}
\end{abstract}
\medskip

\centerline{\sf Table of Contents}
\smallskip
\def\hp{\hphantom{x}}
$\hspace{30mm}$ {\footnotesize \parbox[t]{115mm}{

\hp${}_{}$\!
\hp\!1.{ $\;\,$} {\tt Introduction} \newline
\hp\qquad{\bf  I. Algebraic part}\newline
\hp2.{ $\;\,$} {\tt Various quantum algebras}\newline
\hp3.{ $\;\,$} {\tt Algebraic category equivalences}\newline
\hp4.{ $\;\,$} {\tt Proof of Induction theorem}\newline
\hp5.{ $\;\,$} {\tt Proof of Quantum group "formality" theorem}\newline
\hp\qquad{\bf II. Geometric part}\newline
\hp6.{ $\;\,$} {\tt The loop Grassmannian and the Principal nilpotent}\newline
\hp7.{ $\;\,$} {\tt Self-extensions of the Regular sheaf}\newline
\hp8.{ $\;\,$} {\tt Wakimoto sheaves}\newline
\hp9.{ $\;\,$} {\tt Geometric Equivalence theorems}\newline
10. {$\;\,$}  {\tt Quantum group cohomology and the  loop Grassmannian}\newline
}}

\bigskip

\section{Introduction}\label{sec_intro}
\subsection{Main players.}
 Most of the contents of this paper may be roughly
summed-up in the following  diagram of category
equivalences
(where `$Q$' stands for "quantum", and `$P$' stands for
"perverse"):
\begin{equation}\label{sum-up}
\xymatrix{
{D^b\catm}\quad
\ar@{=>}[rr]^<>(.5){Q}
\ar@/^2pc/[rrrr]|-{\,{\mbox{\footnotesize{ Lusztig multiplicity }\,}}}
&&\quad D^b\Coh^{G\times\C^*}(\NN)\quad
\ar@{<=}[rr]^<>(.5){P}_{\;\;
\;\;{\,\mbox{ \tiny{\cite{CG}-conjecture}}}}&&
\quad
{D^b\Pmix}}
\end{equation}
In this diagram, 
 $G$ is a  connected complex
semisimple group of adjoint type with Lie algebra $\g$.
We  fix a Borel
subgroup
 $B\subset G$, write   $\b=\Lie B\sset\g$ for
 the corresponding 
 Borel subalgebra,  and  $\n$ for
 the nilradical of 
$\b$. Let 
$\NN:= G\times_B\,\n$ be the Springer resolution,
 and
 $\Coh^{G\times\C^*}(\NN)$ the abelian
category of $G\times\C^*$-equivariant coherent
sheaves on $\NN$,
where the group $G$ acts on $\NN$ by conjugation
and $\C^*$ acts by dilations along the fibers.
Further, let
 $\U$ be the quantized universal enveloping
algebra of  $\g$ specialized at a root of unity.
The category $\catm$ on the left of (\ref{sum-up}) 
 stands for a mixed version, see \cite[Definition 4.3.1]{BGS} or
Sect. \ref{mixed_categories} below, of
the abelian category  of
finite-dimensional $\U$-modules in the linkage class
of the trivial 1-dimensional module. 
Finally, we write  $D^b\CC$ 
 for the
 bounded derived category of an abelian category $\CC$.

\ab Forgetting part of the structure
one may consider, instead of  $\catm,$
the category $\cat$ of actual (non-mixed)  $\U$-modules
as well. 
Forgetting the mixed structure on the left of diagram   (\ref{sum-up})
 corresponds to forgetting
the $\C^*$-equivariance in the middle term of  (\ref{sum-up}),
i.e., to
 replacing
$G\times\C^*$-equivariant sheaves on $\NN$
by $G$-equivariant ones.
Although this sort of simplification may
look rather attractive,
the resulting triangulated category  $D^G_\text{coherent}(\NN)$
that will have to replace the middle term
in the diagram above will no longer be
the derived category of the
corresponding abelian category $\Coh^G(\NN)$ 
and, in effect, of any abelian category.
This subtlety is rather technical;
the reader may ignore it at first reading.

\ab 
Finally, let  $\Gd$ denote the complex connected and simply-connected
semisimple
group  dual to $G$ in the sense
of Langlands. We write $\Gr$ for the loop Grassmannian
of $\Gd$. 
The Grassmannian has a standard stratification by
Iwahori (= affine Borel) orbits. The strata, usually called
Schubert cells, are isomorphic to finite dimensional
affine-linear spaces.
We let $\Perv(\Gr)$ denote the  abelian category of 
perverse sheaves on $\Gr$ which are constructible with respect
to this stratification,
and we write $\Pmix$ for its mixed counterpart,
the category of mixed $\ell$-adic perverse sheaves,
see [BBD].

\ab The main result of the paper says that all
three  categories in
 (\ref{sum-up}) are equivalent
as triangulated categories. 
Furthermore, we show that the composite 
equivalence
$P^{-1}\ccirc Q$ is compatible with the natural $t$-structures
on the categories on the LHS and RHS of (\ref{sum-up}),
hence induces  equivalences of abelian categories:
\begin{equation}\label{abelian_eq}
\cat\simeq\Perv(\Gr)\,.
\end{equation}
This yields,  in particular,
 the conjecture formulated in [GK,\S4.3],
relating quantum group cohomology to perverse sheaves on
the loop Grassmannian. 

\ab The equivalence in
 (\ref{abelian_eq})  also provides  character formulas,
 conjectured by Lusztig [L3] and
referred to as "Lusztig multiplicity" formulas in (\ref{sum-up}),
for simple $\U$-modules in the principal block in terms
of intersection homology sheaves on the loop Grassmannian.
Very similar character formulas  have been proved earlier
by combining several known deep
results due to Kazhdan-Lusztig [KL2], Kashiwara-Tanisaki [KT].

\subsection{Relation to results by Kazhdan-Lusztig and
Kashiwara-Tanisaki.}  For each negative rational number $k$(=`level'),
Kashiwara and Tanisaki consider an abelian 
category $\Mod^G_k({\mathscr{D}})$ of $G$-equivariant holonomic
modules over ${\mathscr{D}}$, a sheaf of twisted
differential operators
on the affine flag variety, 
with  a certain  monodromy determined by
(the denominator of) $k$.
This category is equivalent
 via 
the Riemann-Hilbert correspondence
to  $\Perv_k(\Gr^G)$, a category  of monodromic
perverse sheaves. The latter
category 
is defined similarly to the category
$\Perv(\Gr)$ considered in \eqref{abelian_eq},
with the following two differences:

\pb{The Grassmannian $\Gr^G$ stands for the loop Grassmannian for the
group $G$ rather than for the Langlands dual group $G^\vee$; and}

\pb{The objects of $\Perv_k(\Gr^G)$
are perverse  sheaves not on the
Grassmannian $\Gr^G$ itself but on the
total space of a $\C^*$-bundle (so-called {\em determinant bundle})
on the
Grassmannian, with monodromy along the fibers (determined by the rational number $k$).}

\ab
Further,
let $\wh{\g}$ be the affine Lie algebra
associated with $\g$ and $\rep(\Ugh)$ the category
of $\g$-integrable highest weight $\wh{\g}$-modules of level $k-h$,
where $h$ denotes the dual Coxeter number of the Lie algebra
$\g$. Kazhdan and Lusztig used a `fusion type' product
to make  $\rep(\Ugh)$ a tensor category.
On the other hand, let $\U$ be the quantized
enveloping algebra
with parameter $q:=\exp(\pi\sqrt{-1}/d\cdot k)$,
where   $d=1$ if $\g$ is a simple Lie algebra
of types  $\mathbf{A}, \mathbf{D}, \mathbf{E}$, $d=2$ for types
$\mathbf{B},\mathbf{C}, \mathbf{F}$, and $d=3$ for  type
$\mathbf{G}$. 
Let $\rep(\U)$ be the tensor category of finite dimensional
$\U$-modules.
In [KL2] the authors  have established
an  equivalence of  tensor categories
$\rep(\U)\stackrel{\sim}\longleftrightarrow \rep(\Ugh)$.
The subcategory $\block(\U)\sset \rep(\U)$ goes
under the equivalence to the corresponding
principal block $\block(\Ugh)\sset \rep(\Ugh)$.

\ab Each of the categories $\Perv_k(\Gr^G),\,
\block(\Ugh)$ and  $\Mod^G_k({\mathscr{D}}),$
 comes equipped with
collections $\{\Delta_\mu\}_{\mu\in\Y},$ resp.,
$\{\nabla_\mu\}_{\mu\in\Y},$ of so-colled {\em standard},  resp.,
{\em costandard} objects, all labelled by the same partially ordered set
$\Y$. In each case, one has
\begin{equation}\label{highest_weight}
\Ext^i(\Delta_\la,\nabla_\mu)=\begin{cases} \C & \text{if}\enspace \la=\mu\\
0 & \text{if}\enspace \la\neq\mu.\end{cases}
\end{equation}
This is essentially well-known:
in the case of category 
 $\block(\Ugh)$ a proof can be found
 e.g. [KL2]; in the case of $\Perv_k(\Gr^G)$,  isomorphism
\eqref{highest_weight} follows  from a similar formula for the Ext-groups
in  $D^b(\Gr^G)$, a larger triangulated category containing
the abelian category  $\Perv_k(\Gr^G)$ as a subcategory,
and a result of \cite[Corollary 3.3.2]{BGS} saying
that the Ext-groups in the two categories are the same.
By the equivalence  $\Mod^G_k({\mathscr{D}})\cong \Perv_k(\Gr^G)$, the isomorphism in
\eqref{highest_weight} holds also for the category  $\Mod^G_k({\mathscr{D}})$.
To sum-up, the categories  $\block(\Ugh),\,\Perv_k(\Gr^G)$ and $\Mod^G_k({\mathscr{D}})$
are
{\em highest weight categories} in the
terminology of \cite{CPS}.

\ab
Kashiwara-Tanisaki consider the 
global sections functor
$\Gamma: \Mod^G_k({\mathscr{D}})\to\block(\Ugh),\,
{\mathscr M}\mapsto \Gamma({\mathscr M})$.
One of the main results of \cite{KT} says that this functor
  provides a bijection
$$\{\text{standard/costandard objects in  $\Mod^G_k({\mathscr{D}})$}\}
\enspace\longleftrightarrow
\enspace\{\text{standard/costandard objects in  $\block(\Ugh)$}\}
$$
which is compatible with the labelling of the objects involved
by the set $\Y$.
More recently, Beilinson-Drinfeld \cite{BeDr}
proved that $\Gamma$ is an exact functor,
cf. also \cite{FG}.
Now, by an elementary general result
(proved using  \eqref{highest_weight} and `devissage', cf. Lemma \ref{abstract_nonsense}),
 any exact functor
between  highest weight categories
that gives  bijections (compatible with labelling)
 both of the sets of (isomorphism classes of) 
standard and costandard
objects, respectively,
 must be an equivalence.
It follows that the category
 $\block(\Ugh)$ is equivalent to $\Mod^G_k({\mathscr{D}})$.
Thus,
one obtains the following equivalences:
\begin{equation}\label{KLKT}
\xymatrix{
\Perv_k(\Gr^G)\enspace
\ar@{=}[rr]^<>(0.5){\text{Rieman-Hilbert}}&&
\enspace\Mod^G_k({\mathscr{D}})\enspace\ar@{=}[r]^<>(0.5){[\op{KT}]}&
\enspace\block(\Ugh)\enspace\ar@{=}[r]^<>(0.5){[\op{KL}]}&
\enspace\block(\U)}
\eeq

\ab In this paper we consider the special case where $q$ is an odd root
of unity of order prime to $3$.
In that case, the corresponding  rational number
$k$, such that $q=\exp(\pi\sqrt{-1}/d\cdot k)$,
has a small denominator.
Compairing the composite
equivalence in \eqref{KLKT} with the one in \eqref{abelian_eq},
we get $\Perv_k(\Gr^G)\cong \block(\U)\cong \Perv(\Gr)$.
Although we do not know how to construct a direct
equivalence $\Perv_k(\Gr^G)\iso\Perv(\Gr)$
by geometric means,  the results
of Lusztig \cite{L5} imply that the
character formulas for simple objects
in these categories are identical. This explains the relation of our results
with those of [KL2] and [KT].

\subsection{Outline of our strategy.}\label{strategy}
The  construction of both  equivalences in (\ref{sum-up})
 is carried out according to the following 
rather general pattern. Let $D$ denote any of the
three triangulated categories in (\ref{sum-up}).
In each case, we  find an appropriate object
$P\in D$, and form  the  differential graded (\dg-) algebra
$\RHom\hdot_{\!_{D}}(P, P)$. Then, the assignment
$F: M\mapsto  \RHom\hdot_{\!_{D}}(P, M)\,$ gives a functor
 from the category $D$ to the derived category of \dg-modules
over $\RHom\hdot_{\!_{D}}(P, P)$.
We show, as a first step, that the functor $F$ 
 is an equivalence. We express this by saying  that
the category $D$ is `governed' by the \dg-algebra 
$\RHom\hdot_{\!_{D}}(P, P)$. The second step 
consists of proving that the \dg-algebra $\RHom\hdot_{\!_{D}}(P, P)$
is {\it formal}, that is,
quasi-isomorphic to $\Ext\hdot_{\!_{D}}(P, P)$, the corresponding
Ext-algebra under the 
Yoneda product. The formality implies  
  that
the category $D$ is `governed' by the algebra
$\Ext\hdot_{\!_{D}}(P, P)$, considered as a graded algebra with
trivial differential. The third step consists of an 
explicit calculation of this Ext-algebra.
An  exciting outcome of the calculation (Theorem \ref{ranee}) is 
 that the Ext-algebras turn out to be the same
for all three categories in question. Thus,  all three categories
are  `governed' by the same algebra, and we are done.

\subsection{The functor $Q$:} $D^b\catm\to
D^b\Coh^{G\times\C^*}(\NN)$
giving the first equivalence in (\ref{sum-up})
is a refinement of a very naive functor
introduced in [GK]. Specifically,
let $\bb$ be  the "Borel part" of 
the "small" quantum group $\u\sset \U$,
and $H\hdot(\bb,\C)$ the cohomology 
 algebra of $\bb$ with trivial coefficients.
Since $\bb\subset\U$, any $\U$-module may
be viewed as a $\bb$-module, 
by restriction, and the cohomology $H\hdot(\bb, M|_\bb)$
has a canonical  graded $H\hdot(\bb,\C)$-module structure.
The following  functor has been
considered in [GK]:
\begin{equation}\label{naive}
Q_{_{\sf{naive}}}: \cat\too H\hdot(\bb,\C)\mmod
,\quad
M\mto H\hdot(\bb, M|_\bb)\,.
\end{equation}

\ab 
Now, we have fixed a Borel
subgroup $B\sset G$ with Lie algebra $\b$.
 According to [GK] one has a natural
$\Ad B$-equivariant (degree doubling) algebra isomorphism
$H^{2\hhdot}(\bb,\C)\simeq \C\hdot[\n]$,
where the group $B$ acts on  $\n$,
 the nilradical of 
$\b$, by the adjoint action.
This puts, for any $\U$-module $M$,
the  structure of a $B$-equivariant
graded $H\hdot(\bb,\C)$-module, hence
 $\C[\n]$-module,
on $H\hdot(\bb, M|_\bb)$.
The module   $H\hdot(\bb, M|_\bb)$ is finitely generated,
provided $\dim M<\infty$,
hence,  gives rise to an object
of $\Coh^{B\times\C^*}(\n)$,
the category
a $B\times\C^*$-equivariant coherent sheaves on $\n$.
Further, inducing sheaves from the vector space $\n$ up
to the Springer resolution
$\NN=G\times_{_B}\n,$ we obtain from (\ref{naive}) 
the following composite functor:
\begin{equation}\label{naive2}
Q_{_{\sf{naive}}}:
\cat\too 
\Coh^{B\times\C^*}(\n)
\;\stackrel{_{\tt induction}}\iso\;
\Coh^{G\times\C^*}(\NN)\,,
\end{equation}
where the second arrow denotes the obvious equivalence,
whose inverse is given by restricting to the fiber
$\n=\{1\}\times_{_B}\n\into G\times_{_B}\n=\NN$.

\ab The functor $Q_{_{\sf{naive}}}$
may be viewed as a "naive" analogue of the functor
$Q$ in (\ref{sum-up}). In order to construct $Q$ itself,
 one has  to `lift'  considerations above
to the level of  derived categories.
To this end, we will prove  in \S\ref{sec_formality}
that the dg-algebra 
 $\RHom_\bb(\C,\C)$ is {\em formal}, that is, we
will construct an $\Ad B$-equivariant (degree doubling) dg-algebra map
\begin{equation}\label{intro_formal}
\C\hdot[\n]\too\RHom^{2\hhdot}_\bb(\C,\C),
\end{equation}
where $\C\hdot[\n]=\sym(\n^*[-2])$ is viewed as a dg-algebra  
(generated by the space $\n^*$ of linear functions placed in
degree $2$) and equipped with zero 
differential. The map in \eqref{intro_formal}
will be shown to induce the above mentioned
isomorphism of cohomology $\C\hdot[\n]\iso\Ext^{2\hhdot}_\bb(\C,\C)=H\hdot(\bb,\C)$
proved in \cite{GK}, in particular, it is a quasi-isomorphism.

\ab The main idea of our approach to constructing  quasi-isomorphism
 \eqref{intro_formal} is as follows.
Recall first  a well-known
result due to Gerstenhaber saying that any
associative algebra $\aa$ and a 1-st order deformation of
 $\aa$ parametrized by a vector space $V$,
give rise to a canonical   linear map 
$V\to\Hoch^2(\aa)$,
the second Hochschild cohomology group of $\aa$. 
The Hochschild cohomology being a commutative
algebra, the latter map extends to a unique
degree doubling
 algebra homomorphism
 $\sym\hdot(V[-2])\to
\Hoch^{2\hhdot}(\aa)$.
We show in [BG] that
  any extension of the 1-st order deformation 
to a deformation of infinite order
provides a canonical lift of  the homomorphism
$\sym\hdot(V[-2])\to \Hoch^{2\hhdot}(\aa)$
 {\em to the 
\dg-level}, i.e. to a dg-algebra homomorphism
$\sym\hdot(V[-2])\to \RHom_{\aa\bimod}(\aa,\aa),$
(where the graded algebra $\sym(V[-2])$ is viewed as a dg-algebra
 with zero 
differential) and such that  
the induced  map on cohomology
is the Gerstenhaber map mentioned above, see Theorem
\ref{BG_formality2} for a precise statement.

\ab Our crucial observation is that
the  De Concini-Kac version (without divided powers) of the quantum Borel
algebra provides  a formal
(infinite order) deformation
of the algebra $\bb$, with $V=\n^*$ being the parameter space.
Further, the algebra $\bb$ has
a natural Hopf algebra structure,
hence,  the  Hochschild cohomology algebra maps naturally to 
the algebra $H^{2\hhdot}(\bb,\C)$, the cohomology with 
trivial coefficients. Adapting the general
 construction of the dg-algebra homomorphism
$\sym(V[-2])\to \RHom_{\aa\bimod}(\aa,\aa)$
 to the Hopf algebra $\aa:=\bb$,
yields  the desired  dg-algebra map
\eqref{intro_formal}.

\ab It is worth mentioning perhaps that we actually need 
a stronger, $\Ub$-equivariant version,
of  quasi-isomorphism \eqref{intro_formal}. The construction of such an  equivariant
quasi-isomorphism
exploits the existence of {\em Steinberg representation},
and also a Hopf-adjoint action of the Lusztig
version (with  divided powers) of the quantum  Borel algebra on the
De Concini-Kac  version (without divided powers) of the same algebra.
 We refer to \S\ref{sec_formality} for details.

\ab One may compose a quasi-inverse of
the equivalence $Q$ on left of  (\ref{sum-up})
 with the forgetful functor $\catm\to\cat$.
This way, we obtain  the  following result
involving no mixed categories (see Theorems \ref{proposition_roma1} and
\ref{formality2}):
\begin{corollary}\label{roma_functor} There exists
 a triangulated functor $F: D^b\Coh^{G\times\C^*}(\NN)\map
D^b\cat$
such that 
the image of $F$ generates $D^b\cat$
as a triangulated category and we have:
\vskip 1pt

\vi $\dis F(\oo_{\NN}(\la))=\RInd_{_\B}^{^\U}(l\la),
\quad\text{and}\quad F(z^i\otimes{\mathscr{F}})=F({\mathscr{F}})[i],
\quad \forall\la\in\Y,\,
i\in \Z ,\, {\mathscr{F}}\in D^b\Coh^{G\times\C^*}(\NN).$
\vskip 1pt

\vii Write $i: \n=\{1\}\times_{_B}\n\into G\times_{_B}\n=\NN$
for the natural imbedding. Then, cf. \eqref{naive}-\eqref{naive2}, we have
$$ \op{R\Gamma}(\n, i^*{\mathscr{F}})=\RHom_{_{\bb\mmod}}(\C_\bb, F({\mathscr{F}})),
\quad \forall {\mathscr{F}}\in D^b\Coh^{G\times\C^*}(\NN).
$$

\viii The functor $F$
induces, for any ${\mathscr{F}},{\mathscr{F}}'
\in D^b\Coh^{G\times\C^*}(\NN),$ canonical isomorphisms
$$\bigoplus\nolimits_{i\in\Z}\,\Hom\hdot_{D^b\Coh^{G\times\C^*}(\NN)}
({\mathscr{F}},z^i\otimes{\mathscr{F}}')
\iso
\Hom\hdot_{D^b\cat}(F({\mathscr{F}}), F({\mathscr{F}}')).$$
\end{corollary}

Here, given a  $\C^*$-equivariant sheaf (or complex of sheaves) ${\mathscr{F}}$,
we write $z^i\otimes{\mathscr{F}}$ for the sheaf (or complex of sheaves)
obtained by twisting the $\C^*$-equivariant structure
by means of the character $z\mapsto z^i$, and
let ${\mathscr{F}}[k]$ denote the homological shift of 
${\mathscr{F}}$ by $k$ in the derived category.

\subsection{The functor $P$:}\label{functor_intro_P}
$D^b\Pmix\too D^b\Coh^{G\times\C^*}(\NN).$
The
point of departure in constructing the functor
on the right of (\ref{sum-up}) is the fundamental result of
 geometric Langlands theory
 saying that there is an equivalence
${\p} : \rep(G) \iso \PO,$
between the tensor category of
finite dimensional rational representations of the
group $G$ and the  tensor category of $\Gd(\oo)$-equivariant
perverse sheaves on the loop Grassmannian
equipped with a convolution-type monoidal structure:
${\mathcal{M}}_1,{\mathcal{M}}_2 \mto {\mathcal{M}}_1\star{\mathcal{M}}_2,$
see [G2],[MV] and also [Ga].
In particular, write $\one={\p}(\C)\in \PO$ for the
sky-scrapper sheaf at the base point of $\Gr$
that corresponds to the trivial
one-dimensional $G$-module,
and write $\R={\p}(\C[G])$
for the ind-object in $\PO$
corresponding to the regular $G$-representation.
The standard algebra structure on the coordinate
ring $\C[G],$ by pointwise multiplication, makes $\R$ 
 a ring-object in $\PO$.
It is easy
to see that this 
gives  a canonical commutative
graded algebra structure on the 
Ext-group $\Ext\hdot_{_{D^b(\Gr)}}(\one,\,\R),$
and that the $G$-action on $\C[G]$ by right translations
gives a $G$-action on the Ext-algebra.
Furthermore, for any perverse sheaf ${\mathcal{M}}$ on $\Gr$,
the Ext-group  $\Ext\hdot_{_{D^b(\Gr)}}(\one,\,{\mathcal{M}}\star\R)$
has the natural structure of a
$G$-equivariant finitely-generated graded
$\Ext\hdot_{_{D^b(\Gr)}}(\one,\,\R)$-module, via the Yoneda product.

\ab A crucial Ext-calculation, carried out in section \ref{sec_extensions},
provides a canonical
$G$-equivariant (degree doubling) algebra isomorphism
\begin{equation}\label{roma_ext}
\Ext^{2\hhdot}_{_{D^b(\Gr)}}(\one,\,\R)\;
\simeq\;\C\hdot[\N]\quad\left(\text{and}\quad
\Ext^{\tt{odd}}_{_{D^b(\Gr)}}(\one,\,\R)=0\right),
\end{equation}
where $\N$ is the nilpotent variety in $\g$.
The  homomorphism
$\C\hdot[\N]\to\Ext^{2\hhdot}(\one,\R)$
in \eqref{roma_ext} is induced by a morphism of
algebraic varieties
$\,
\Spec\left(\Ext\hdot(\one,\,\R)\right)
\longrightarrow \N.
$
The latter is constructed by means of
{\em Tannakian formalism} as follows. 

\ab In general, let $Y$ be an affine algebraic
variety. 
Constructing a map $Y\to \N$ is equivalent to producing
a family $\big\{\Phi_V: \oo_Y\otimes V\to \oo_Y\otimes V\big\}_{V\in \rep(G)}$,
(endomorphisms  of the trivial vector bundle with fiber $V$)
such that, for any $V, V'\in\rep(G),$ one has
$\dis\Phi_{_{V\otimes V'}} \simeq 
\Phi_{_V}\otimes\Id_{_{V'}} + \Id_{_V}\otimes
\Phi_{_{V'}}.$

\ab In the special case $Y=\Spec(\Ext\hdot(\one,\R))$, 
the geometric Satake isomorphism, cf. \S\ref{sec_loop}, provides a
canonical  isomorphism
$\Gamma(Y,\oo_Y\otimes V)=\Ext(\one,\R)\otimes V 
\cong\Ext\hdot(\one,\R\star \p V)$.
To construct  a nilpotent endomorphism $\Phi_{_V}: \Ext(\one,\R\star \p V)\to
\Ext(\one,\R\star \p V),$
 consider the first Chern class  $c\in H^2(\Gr,\C)$
of the standard determinant line bundle on the loop Grassmannian, see [G2].
Cup-product with $c$ induces a 
morphism\footnote{For any variety $X$ and $M\in D^b(X),$
 in the derived category,
there is a natural map $H^p(X,\C) \map \Ext^p_{_{D^b(X)}}(M,M)=
\Hom_{_{D^b(X)}}(M,M[p]),$ cf. e.g., \cite[(8.3.17)]{CG}.
We apply this to $X=\Gr$, $c\in H^2(\Gr,\C),$ and $M=\p{V}$.}
$c:\p( V) \to\p( V)[2]$.
We let $\Phi_{_V}$ be the  map
$\mbox{Id}_{\R}\star c:\,\Ext\hdot(\one,\,\R\star \p V) \longrightarrow 
\Ext^{\hhdot+2}(\one,\,\R\star \p V)$,
 obtained
by applying the functor $\Ext\hdot(\one,\, \R\star (-))$ to the morphism
above.
  For further ramifications of this
construction
see \S\S\ref{fiber_fun1},\ref{another}.

\ab 
Using \eqref{roma_ext}, we may 
view a $G$-equivariant  graded
$\Ext\hdot_{_{D^b(\Gr)}}(\one,\,\R)$-module as a $\C\hdot[\N]$-module,
equivalently,
as a $G\times\C^*$-equivariant sheaf on $\N$.
This way we obtain
a functor:
\begin{equation}\label{P_naive}
P_{_{\sf{naive}}}:\;
\Perv(\Gr) \too \Coh^{G\times\C^*}(\N)
,\quad
{\mathcal{M}}\mto \Ext\hdot_{_{D^b(\Gr)}}(\one,\,{\mathcal{M}}\star\R)\,.
\end{equation}
The functor thus obtained 
may be viewed as a "naive" analogue of the functor
$P$ in (\ref{sum-up}).
The actual construction of the equivalence $P$
is more involved:
one has to replace $\N$ by the Springer resolution $\NN$,
and to make everything work on the level of derived categories.
This is made possible by the technique of {\em weights}
of mixed $\ell$-adic sheaves combined with known results
on the purity of intersection cohomology for flag varieties,
due to \cite{KL1}.

\subsection{Relation to affine Hecke algebras.}
 One of the motivations for the present work was 
an attempt to understand an old mystery
surrounding the existence of two completely different realisations of the
affine Hecke algebra. The first realisation is in terms
of locally constant functions on the flag manifold of a $p$-adic
reductive group, while the other is in terms of
equivariant $K$-theory of a complex variety 
(Steinberg variety) acted on
by the Langlands dual complex reductive group,
see [KL3], [CG]. The existence of the two realisations indicates 
 a possible link between {\it perverse} sheaves on the
affine flag manifold, on one hand, and {\it coherent} sheaves on
the Steinberg variety  
(over $\C$) for the Langlands dual
group, on the other hand. Specifically,
it has been conjectured in [CG,~p.15]
that there should be a functor
\begin{equation}\label{aff_heck_fun}
F_{_{\tt{Hecke\; alg}}}: \;
D^b\Perv^{^{\mathsf{mix}}}(\text{\it affine flag manifold})
 \too D^b{{{\mathcal C}\mbox{\it oh}}}^{G\times\C^*}(\text{\it 
Steinberg variety})\,.
\end{equation}
  For the {\it finite} Hecke algebra,
a functor of this kind has been constructed
by Tanisaki,
see [Ta], by means of ${\mathcal{D}}$-modules:
each perverse sheaf gives rise to a
${\mathcal{D}}$-module, and taking  associated graded module
of that ${\mathcal{D}}$-module with respect to a certain filtration
yields a coherent sheaf on the Steinberg variety.
Tanisaki's construction does not extend, however, to the whole 
{\it affine} Hecke algebra; also, it by no means explains the
appearance of the Langlands dual group. 
Our equivalence $P^{-1}\ccirc Q$ in (\ref{sum-up})
provides a "correct" construction of a counterpart 
of the functor (\ref{aff_heck_fun})
for the fundamental polynomial representation
of the affine Hecke algebra instead of the algebra
itself.\footnote{cf. also [AB] for an alternative approach
which is, in a sense,
`Koszul dual' to ours.}
A complete construction of (\ref{aff_heck_fun})
in the algebra case will be
carried out in a forthcoming paper.
Here we mention only that
 replacing the module by the algebra 
results, geometrically, in replacing
the loop Grassmannian by the
affine flag manifold, on one hand,
and replacing
the Springer resolution by the Steinberg variety, on the other hand.
In addition to that, handling  the algebra case involves
an important extra-ingredient: the  geometric
construction of the center of the affine Hecke algebra
by means of nearby cycles, due to Gaitsgory [Ga].

\ab To conclude the Introduction,
the following remark is worth mentioning.
None of the equivalences $P$ and $Q$ taken separately,
as opposed to the composite
$P^{-1}\ccirc Q$ in (\ref{sum-up}),
is compatible with the natural $t$-structures. In other words,
the abelian subcategory
$Q(\catm)=P(\Pmix)$ 
of the triangulated category
 $D^b\Coh^{G\times\C^*}(\NN)$  does not coincide
with $\Coh^{G\times\C^*}(\NN)$. The "exotic" $t$-structure on
$D^b\Coh^{G\times\C^*}(\NN)$ arising, via $Q$
(equivalently, via $P$), from the
natural  $t$-structure on $Q(\catm)$
is, in effect, closely related to the 
{\it perverse coherent $t$-structure} studied in
[Be2]. Specifically, it will be shown
in a subsequent paper that the functor $P: D^b\Pmix
\too D^b\Coh^{G\times\C^*}(\NN)$ takes indecomposable {\it tilting},
resp.   {\it simple}, objects
of $\Pmix$
into {\it simple}, resp.  {\it tilting}, (with respect to perverse coherent $t$-structure)
objects of $\Coh^{G\times\C^*}(\NN)$. This, combined with the results
of the present paper, implies
that the {\it tilting} $\U$-modules in the category  $\cat$ go, under the
equivalence $Q: D^b\catm\too D^b\Coh^{G\times\C^*}(\NN)$ to
  {\it simple}
perverse coherent sheaves on $\NN$. Moreover, an
additional argument based on results of \cite{AB} shows that
the parameters labelling the tilting objects in $\cat$
and the simple perverse coherent sheaves on $\NN$ correspond
to each other, see \cite{B5}. 
It follows, in particular, that the  support of the quantum
group  cohomology of a tilting  $\U$-module agrees with
the one  conjectured by Humphreys.

\subsection{Organization of the paper.}
In \S 2 we recall basic constructions regarding various
versions of quantum groups that will be used later in the
paper. The main result of this section is
Theorem \ref{GB_thm} which is closely related to the
 De Concini-Kac-Procesi results \cite{DKP} on {\it quantum coadjoint action}.
 In 
\S\ref{sec_block} we introduce basic categories of $\U$-modules
and state two main results of the algebraic part of the paper.
 Section \ref{sec_induction} is devoted to the
proof of the first result, {\it Induction theorem},
saying that the derived category of $\U$-modules
in the principal block is equivalent to
an appropriate derived category of 
modules over the Borel part of $\U$.
The proof exploits the techniques of {\it wall-crossing functors}.
In \S\ref{sec_formality}, we prove the second main result
saying that the
\dg-algebra
of derived endomorphisms of the trivial 1-dimensional
module over $\bb$ (=  Borel part of the "small" quantum group)
is formal, i.e., is
quasi-isomorphic to its cohomology algebra.
In \S\ref{sec_loop} we review the (known) relation between finite
dimensional
representations of a semisimple group and
perverse sheaves on $\Gr$, the loop Grassmannian for the Langlands
dual group. We remind also the role of the principal 
nilpotent element in describing the cohomology
of $\Gr$. In \S\ref{sec_extensions} we prove an algebra isomorphism
that generalizes isomorphism (\ref{roma_ext}).
Section \ref{sec_wakimoto} is devoted to the basics of the theory
of {\it Wakimoto perverse sheaves} on the affine flag manifold,
due to Mirkovi\'c (unpublished). The classes of these sheaves 
in the Grothendieck group correspond,
under the standard isomorphism with the affine Hecke algebra,
to base elements of an important large commutative subalgebra
in  the affine Hecke algebra that has been introduced by Bernstein.
The main results of the paper are proved in \S\ref{sec_equivalence}
where the functors $Q, P$ are
constructed and
the category equivalences (\ref{sum-up}) are established.
The arguments there use both  algebraic 
and geometric results obtained in all the previous sections.
In \S \ref{Gr} we prove Ginzburg-Kumar conjecture [GK, \S4.3]
relating  quantum
group  cohomology to perverse sheaves.
\medskip

{\footnotesize{
{\bf Aknowledgements.} We are especially indebted to Ivan 
Mirkovi\'c 
 who suggested one of the key
ideas of the paper ("cohomological localization to the cotangent bundle")
to one of us (R.B.) back in 1999. We also thank him for critical reading
of the maniscript, and for the permission to use his unpublished results
on Wakimoto modules.
We are  also grateful to 
M. Finkelberg for many 
useful discussions, and
to H. H. Andersen for 
pointing out several inaccuracies in the
original draft of the paper.
Finally, we would like to thank
V. Drinfeld whose question has led us,
indirectly, to a construction
of bi-functor in \S\ref{bi_functor}
that is a key element in our proof
of the main result of Sect.~\ref{sec_formality}.
}}

\bigskip

\centerline{\bf PART I$\,$:$\;$ Algebra}

\section{Various quantum algebras}\label{sec_algebras}
\subsection{}
Let $\k$ be an algebraically closed field of characteristic zero,
and set  $\otimes=\otimes_\k$. We write $\k[X]$ for the
coordinate ring of an algebraic variety $X$.

\ab Given  a   $\k$-algebra $A$ with 
an augmentation $\epsilon: A\to\k$, let $A\eps$ 
denote its kernel.
 Thus, $A\eps$ is a two-sided ideal of $A$, called
the {\it augmentation ideal}, and $\k_{_A} :=
A/A\eps$ is a 1-dimensional
$A$-module.

\begin{definition}\label{normal} Given an associative algebra $A$
and a subalgebra $\aa\sset A$ with augmentation
$\aa\to\k$, we say that $\aa$ is a
 {\it normal}  subalgebra if one has $A\cdot\aa\eps=\aa\eps\cdot A$.
We then write $(\aa):=A\cdot\aa\eps\sset A$ for this
 two-sided ideal.
\end{definition}

\ab Given a   $\k$-algebra $A$,  we write either $A\mmod$ or
$\Mod(A)$ for the abelian category
of left $A$-modules. The notation $\rep(A)$  is
reserved for the tensor category
of finite-dimensional modules over a Hopf algebra $A$,
unless specified otherwise
(this convention will be altered slightly in \ref{rep}).
In case $A$ is a Hopf algebra, we always
assume that the augmentation $\epsilon: A\to\k$
coincides with the counit.

\subsection{} Let $\h$ be a finite dimensional $\k$-vector space,
$\h^*$ the dual space, and write $\langle-,-\rangle: \h^*\times\h\to\k$
for the canonical pairing. Let
 $R\subset \h^*$  be a finite reduced root system. From now on we
fix the set 
$R_+\sset R$  of positive roots of our  root system, and write $\{\alpha_i\}_{i\in I}$ 
for the corresponding set
of simple roots (labelled by a finite set $I$).
Let $\check\alpha$ denote the coroot corresponding
to a root $\alpha\in R$, so that  $a_{ij}=\langle \check\alpha_i,
\alpha_j\rangle$ is the  Cartan
matrix.

\ab 
Let $W$ be the Weyl group of our root system, acting naturally
on the lattices $\X$ and $\Y,$ see \eqref{XY}.
There is a unique $W$-invariant inner product
$(-,-)_{_{\!\Y}}: \Y \times \Y\map{\mathbb Q}$,
normalized so that $(\alpha_i,\alpha_i)_{_{\!\Y}}=2d_i
\,,\,\forall i\in I,$
where the integers $d_i\geq 1$ are mutually prime.
It is known further that  $d_i\in\{1,2,3\}$ and that 
$a_{ij}=(\alpha_i,\alpha_j)_{_{\!\Y}}/(\alpha_i,\alpha_i)_{_{\!\Y}}$.
In particular the
 matrix
$\|d_i\cdot a_{i,j}\|$ is symmetric.

{\renewcommand{\arraystretch}{1.3}
\beq\label{XY}
\begin{tabular}{|cc|}
\hline
$\X=\{\mu\in\h^*\mid
 \langle\mu,\check\alpha_i\rangle\in\Z\,,\,
 \forall i\in I\}$ &
weight lattice\\
$\X^{++}=\{\mu\in\X\mid\langle\mu,\check\alpha_i\rangle\ge 0\,,\,
 \forall i\in I\}$ & dominant Weyl chamber\\
$\Y=\sum_{i\in I}\,\Z\cd\alpha_i \,\subset\,\X$ &
root lattice\\
$\Y^\vee=\Hom(\Y,\Z)\sset \h$ & coweight lattice\\
$\Y^{++}=\Y\cap \X^{++}$&\\
\hline
\end{tabular}
\eeq
}

\ab
Let  $\g=\n\oplus\h\oplus\overline{\n}$ be a 
semisimple Lie algebra over $\k$ with a fixed
triangular decomposition, such that $R$ is
the root system of $(\g,\h)$, and such that
$\n$ is spanned by root vectors for $R_+$.

\subsection{}
Let $\k(q)$ be the field of rational functions in the variable $q$.
We write
$\tU_q=\tU_q(\g)$ for the
Drinfeld-Jimbo quantized enveloping algebra of $\g$.
Thus, $\tU_q$ is 
  a
$\k(q)$-algebra with  generators  $E_i, F_i,\ i\in
I$,
and
$K_\mu, \mu\in \Y^\vee$, and with the following  defining  relations:
$$
 K_{\mu_1}\cdot K_{\mu_2}=K_{\mu_1+ \mu_2}, $$ $$
 K_\mu\cdot E_i\cdot K_\mu^{-1}=q^{\langle\mu,\alpha_i\rangle}\cdot E_i
,\quad K_\mu
\cdot F_i\cdot K_\mu^{-1}=q^{-\langle\mu,\alpha_i\rangle}\cdot F_i  $$ $$
E_i\cdot F_j-F_j\cdot E_i=\delta_{i,j}\cdot\frac {K_i-K_i^{-1}}{q^{d_i}-
q^{-d_i}}\;,\quad \text{ where } K_i=K_{d_i\cdot\check\alpha_i}, $$ 
and some $q$-analogues of the Serre relations, see e.g. [L2].


\ab We will freely use Lusztig results on
 quantum groups at roots of
unity, see [L2] and also  [AP], pp.579-580.

\ab Fix an {\it odd} positive integer $l$ which is greater than the
Coxeter number of the root system $R$,
and which is moreover prime to 3 if our root system
has factors of type ${\mathbf{G_2}}$.
Fix $\zeta\in\k^\times$, a primitive $l$-th
root of unity, and let $\A\subset \k(q)$ be the local ring at $\zeta$
 and $\m\subset \A$ the maximal ideal in $\A$.

\begin{remark} One may alternatively take $\A=\k[q,q^{-1}]$ as is done in
[L2], [AP]; our choice of $\A$ leads to the same theory.
We alert the reader that the variable `$q$' that we are using here
was denoted by `$v$' in
[L2], [AP].
 $\quad\lozenge$\end{remark}
\subsection{$\A$-forms of $\tU_q$.}
Let $\U_\A$ be the Lusztig's integral form of
$\tU_q$, the $\A$-subalgebra in $\tU_q$
generated by  divided powers $E^{(n)}_i=
E^n_i/[n]_{d_i}!\,,\,F^{(n)}_i=
F^n_i/[n]_{d_i}!\,,\,
i\in I, n\ge 1$ (where $[m]_d!
:=\prod_{s=1}^m {\frac{q^{d\cdot s}-q^{-d\cdot s}}{q^d-q^{-d}}},\,$)
and also various
divided powers $\big[\!{\footnotesize
{\begin{array}{c}K_\mu, m\\ n\end{array}}}
\!\big]_{d_i}$, as defined in [L2].
We will also use a different $\A$-form of
$\tU_q$, {\it without divided powers}, introduced by  De Concini-Kac, see [DK].
This is an  $\A$-subalgebra $\fU_\A\subset\tU_q$
generated by the elements
 $E_i\,,\,F_i\,,\,\frac{K_i-K_i^{-1}}{q^{d_i}-
q^{-d_i}},\,
i\in I,$ and $K_\mu,\,\mu \in\Y^\vee.$
We set $\U:=\U_\A/\m\cd\U_\A$, 
the specialization of  $\U_\A$
at  $q=\zeta$. Further, the
 elements $\{K^\li_i\}_{i\in I}$ are known to be central in 
the algebra $\fU_\A/\m\cdot \fU_\A,$ see [DK, Corollary 3.1].
Put
$\fU:=\fU_\A/\bigl(\m\cdot \fU_\A + \sum_{i\in I}\,(K_i^\li-1)\cdot \fU_\A\bigr)$.
Thus, $\U$ and $\fU$ are
 $\k$-algebras, which are known as,
respectively, the Lusztig
and the De Concini-Kac quantum algebras at a root of unity.

\ab The algebra $\tU_q$ has a Hopf algebra structure over $\k(q)$.
It is known that both $\U_\A$ and  $\fU_\A$ are Hopf
$\A$-subalgebras in $\tU_q$. Therefore,  $\U$ and $\fU$ are Hopf
algebras over $\k$.

\ab By definition, one has  $\fU_\A\subset\U_\A$. Hence,
the imbedding of $\A$-forms
induces, after the specialization at $\zeta$,
a canonical (not necessarily injective) Hopf algebra
homomorphism  $\fU\to \U$. The image of this homomorphism
is a Hopf subalgebra $\u\subset\U$, first introduced by
Lusztig, and referred to as the {\it small quantum group}.
Equivalently, $\u$ is the subalgebra in $\U$ generated by
 the elements
 $E_i\,,\,F_i\,,\,\frac{K_i-K_i^{-1}}{q^{d_i}-
q^{-d_i}},\,
i\in I,$ and $K_\mu,\,\mu \in\Y^\vee.$

\subsection{}
The algebra $\tU_q$ has a triangular
decomposition  $\tU_q = \tU_q^+ \otimes_{_{\k(q)}}
\tU_q^\circ \otimes_{_{\k(q)}}\tU_q^-,$ where
$\tU_q^+,\tU_q^\circ$ and $ \tU_q^-$ are
 the $\k(q)$-subalgebras generated by the set
$\{E_i\}_{i\in I},$ the  set
$\{K_\mu\}_{\mu\in\Y^\vee}$, and the  set
$\{F_i\}_{i\in I}$, respectively.
Given any subring $A\subset \tU_q$, we set
$A^\pm:=A \cap \tU_q^\pm,$ and $A^\circ:=A \cap \tU_q^\circ$.
With this notation, both Lusztig
and De Concini-Kac $\A$-forms are  known to admit triangular
decompositions  $\U_\A= \U_\A^+ \otimes_{_{\A}}
\U_\A^\circ \otimes_{_{\A}}\U_\A^-,$
and $\fU_\A= \fU_\A^+ \otimes_{_{\A}}
\fU_\A^\circ \otimes_{_{\A}}\fU_\A^-,$
respectively. The latter decompositions induce 
the corresponding decompositions
\begin{equation}\label{triangle_dec}
\U= \U^+ \otimes_\k
\U^\circ \otimes_\k\U^-
,\quad\fU= \fU^+ \otimes_\k
\fU^\circ \otimes_{\k}\fU^-,\quad
\u= \u^+ \otimes_\k \u^\circ \otimes_\k \u^-\,.
\end{equation}

\ab
The subalgebras
$
\tB_q:=\tU_q^+ \otimes_{_{\k(q)}}
\tU_q^\circ \subset \tU_q\,,\,\sB_\A:=\U_\A^+ \otimes_{_{\A}}
\U_\A^\circ\subset\U_\A,$  and
$\fB_\A:=\fU_\A^+ \otimes_{_{\A}}
\fU_\A^\circ\subset\fU_\A,$ as well as
various  specializations like 
\begin{equation}\label{borel_alg}
\sB
:=\U^+ \otimes\U^\circ \subset\U,\quad
\fB:=\fU^+\otimes\fU^\circ\subset\fU,\quad
\bb:=\u^+ \otimes \u^\circ\subset\u,
\end{equation}
will be referred to as Borel
parts of the corresponding algebras.
All of these "Borel parts"  are known to be Hopf subalgebras in $\tU_q$
with coproduct and antipode given by the formulas:
\begin{equation}\label{coprod}
\Delta(E_i)=E_i\otimes 1 +K_i\otimes E_i\,,\,
\Delta(K_i)=K_i\otimes K_i\;,\;
S(E_i)=-K_i^{-1}\cdot E_i\,,\,S(K_i)=K_i^{-1}.
\end{equation}
 
\ab Note that  formulas \eqref{coprod} show that
$\tU^+_q$ is {\it not} a Hopf subalgebra.

\ab Put $\overline{B}_q:=U_q^\circ\otimes_{_{\k(q)}}U_q^-.$
This is a Hopf subalgebra in $U_q$, and
Drinfleld constructed a perfect  pairing:
\beq\label{drin1}
\overline{B}_q \otimes B_q \map\k(q)\,,\quad \bar{b}\times b\mapsto
\langle\bar{b},b\rangle\,.
\eeq
The Drinfeld's pairing  enjoys an
invariance property. To formulate it, one first uses
\eqref{drin1} to define a  "differentiation-action"
of the algebra $\overline{B}_q$ on  $B_q$
by the formula
\beq\label{drin4}
\bar{b}:\ b\mto \partial_{_{\bar{b}}}(b):=\sum
\langle\bar{b},b'_\imath\rangle\cdot b''_\imath,
\quad\text{where}\quad \sum\nolimits_\imath\, b'_\imath\otimes b''_\imath:=\Delta(b).
\eeq
 The differentiation-action
makes  $B_q$ a $\overline{B}_q$-module.
The invariance property states, cf. e.g. \cite{L2}:
$$\langle\bar{x}\cdot\bar{y},z\rangle\,+\,\langle\bar{y},
\partial_{_{\bar{x}}}(z)\rangle
=\epsilon(\bar{x})\langle\bar{y},z\rangle\,,\quad\forall
\bar{x},\bar{y}\in\overline{B}_q\,,\,z\in B_q.
$$

\subsection{Frobenius functor.}\label{ufrob}
Let $G$ be a connected semisimple  group of adjoint type
(with trivial center) such that
$\Lie G=\g$. 
Let $\widetilde G$ be the simply-connected
covering of $G$,
and $\fZ(\widetilde G)$ the center of $\widetilde G$
(a finite abelian group).
Thus, we have a short exact sequence
$$1\too \fZ(\widetilde G)\too\widetilde G\stackrel{\pi}\too G\too 1.$$
The pull-back functor $\pi^*: \rep(G)\to
\rep(\widetilde G)$ identifies a
finite-dimensional $G$-module
with  a finite-dimensional $\widetilde G$-module,  such that
the  group $\fZ(\widetilde G)$ acts trivially on it.

\ab Let 
${\mathcal U}\g$ denote the (classical)
universal enveloping algebra of~$\g$.
 Lusztig introduced a certain completion,
$\barU$, of the algebra ${\mathcal U}\g$ such that
the category $\rep(\barU)$ of finite-dimensional $\barU$-modules
may be identified with the category $\rep(G)$.
In more detail, one has
the canonical  algebra map $\jmath: {\mathcal U}\g\to
\barU$, which induces a functor
$\jmath^*: \rep(\barU)\to \rep({\mathcal U}\g)$.

\ab On the other hand, any  finite-dimensional ${\mathcal U}\g$-module
may be regarded, via the
 exponentiation, as a  $\widetilde G$-module.
The completion $\barU$ has the property that  in
the diagram below the images of the
two imbeddings $\pi^*$ and $\jmath^*$ coincide
$$ 
\xymatrix{
{\rep(G)\;}\ar@{^{(}->}[rr]^<>(0.5){\pi^*}&&
{\;\rep(\widetilde G)=\rep({\mathcal U}\g)\;}&&
{\;\rep(\barU).}\ar@{_{(}->}[ll]_<>(0.5){\jmath^*}
}
$$
Thus, simple objects of the category $\rep(\barU)$ are labelled
by the elements of $\Y^{++}$.

\ab Let  $\,\big\{e_i\,,\,f_i\,,\,h_i\big\}_{i\in I}\,$
denote the standard Chevalley generators of  the Lie algebra
$\g$. Lusztig 
proved  that the assignment:
$E_i^{(\li)}\mapsto e_i\,,\, E_i\mapsto 0,$ and
$F_i^{(\li)}\mapsto f_i\,,\,F_i\mapsto 0\,,\,i\in I,$
can be  extended to a well-defined
algebra homomorphism $\phi: \U \to \barU$, called the
{\it Frobenius} map.
Further, the subalgebra $\u$
is known to be
{\it normal} in $\U$, cf. Definition \ref{normal}. Moreover,
 Lusztig has 
proved that $\Ker\phi= (\u)$, i.e.,
one has an exact sequence of bi-algebras
\begin{equation}\label{frobenius}
0\map ({\u})\map \U\stackrel{\phi}{\longrightarrow}\barU.
\end{equation}

\ab The pull-back
via the Frobenius morphism $\phi:
\U \to \barU$ gives rise to an exact tensor functor
\begin{equation}\label{frob_functor}
\phi^*: \rep(G)=\rep(\barU) \too \Umod,\quad
V\mto \fr V\quad \text{\bf {(Frobenius functor)}}.
\end{equation}

\begin{remark} The reader may have observed that
 the quantum algebras $U_q(\g)\,,\,\U,$ etc.,
that we are using are of 
"adjoint type". The results of the paper can be adapted to
"simply-connected" quantum algebras as well. In that case,
the group $G$ must be taken to be the simply-connected
group with Lie algebra $\g$.
Therefore, $G^\vee$, the Langlands dual group, is of adjoint type.
Hence, the corresponding loop Grassmannian $\Gr$ considered in Part II
becomes disconnected; the group $\pi_0(\Gr)$ of its connected
components is canonically isomorphic to 
$\pi_0(\Gr)\cong \pi_1(G^\vee)=\Hom\bigl(\fZ(G)\,,\,\C^*\bigr),$
the Pontryagin dual of 
the center
 of the  simply-connected
group $G$.
 $\quad\lozenge$\end{remark}

\ab Let $e_\al\in\n$ and $\f_\al\in\overline{\n}$ denote 
root vectors corresponding to a root $\al\in R_+$
(so, for any $i\in I$, we have $e_{\al_i}=e_i,\,f_{\al_i}=f_i$).
Throughout the paper,
 we fix a reduced expression for $w_0\in W$, the element of maximal
length. This puts a normal (total) linear order on the set $R_+$
and, for each $\alpha\in R_+$, gives, via the braid group action
on $\tU_q$ (see [L3] for details), an element $E_\al\in\tU^+_q$ and
$F_\al\in\tU^-_q.$ 
The elements $\{E_\al^\li\,,\,F_\al^\li\}_{\al\in R_+}$
are known to be central in $\fU$,
by [DK, Corollary 3.1].

\begin{definition}\label{cent} Let $\fC$ denote the central
subalgebra in $\fB$ generated by the elements $\{E_\al^\li\}_{\al\in
R_+}$.
\end{definition}

\ab Set  $(\fC):=\fC\eps\cdot\fB=\fB\cdot\fC\eps,$ a
 two-sided ideal in $\fB$. Part (i)  of the following Lemma
is due to De Concini-Kac \cite{DK}, part (ii)
is due to Lusztig, and other statements
can be found in  De Concini-Lyubashenko \cite[~\S3]{DL}
(cf. also \cite[pp. 12-13, 72-88]{BFS} for related results 
on the pairing $\u^-\otimes\u^+\to\k$).
\begin{lemma}\label{ldk}\vi $\fC$ is a Hopf subalgebra in $\fB$.

\vii The projection $\fU\onto \u$ induces, by restriction,
 an exact sequence of bi-algebras:
\[0\too (\fC) \too \fB \stackrel{\phi}\too \bb\too 0\,.\]

\viii\pbox{Drinfeld's pairing \eqref{drin1}
restricts to a well-defined $\A$-bilinear pairing $\U^-_\A \otimes\fU^+_\A \map
\A$; the latter gives, after specialization at $q=\zeta$,
a perfect pairing $\U^-\otimes\fU^+\to\k\,.$}

\iv\pbox{The annihilator of the subspace $Z\sset\fU^+$ with respect to the
pairing in \viii is equal to
the  ideal $(\u^-)\sset\U^-$, i.e., we have $(\u^-)=Z^\perp$.\qed}
\end{lemma}

\ab Parts (iii)-(iv)  of the Lemma combined with
the isomorphism $\U^-/(\u^-)={\mathcal{U}}\overline{\n}$,
give rise to a perfect pairing
\beq\label{drin3}
{\mathcal{U}}\overline{\n}\otimes Z\map\k\,.
\eeq

\subsection{Smooth coinduction.}\label{smooth}
Let $A\subset \tU_q$ be a subalgebra
with  triangular decomposition:
$A=A^+\otimes A^\circ\otimes A^-$,
where $A^{?} := \tU_q^{?}\cap A\,,\,?=+,\circ,-$.
There is a natural algebra map $A^\pm\otimes A^\circ\onto  A^\circ$
given by $a^\pm\otimes a^\circ\mapsto \epsilon(a^\pm)\cdot a^\circ$.
We write $\k_{_{A^\circ}}(\la)$ for the 1-dimensional
$A^\circ$-module
corresponding to an algebra homomorphism $\la: A^\circ\to\k$,
 and let $\k_{_{{A^{^\pm}}\otimes A^{^\circ}}}(\la)$ denote its pull-back
via the projection $A^\pm\otimes A^\circ\onto  A^\circ$.

\ab Given a possibly infinite dimensional $A$-module $M$,
define an $A$-submodule $M\alg\sset M$ as follows:
$$ M\alg := \{m\in M\;|\; \dim(A\cdot m) < \infty \enspace\text{and}\enspace
A^\circ\mbox{\it -action on $A\cdot m$ is
diagonalizable}\}\,.$$

\begin{notation}\label{rep}
Let $\modt(A)$  denote the  abelian 
category of finite dimensional  $A^\circ$-diagonalizable
 $A$-modules. 
\end{notation}

\ab We write 
$\limind\modt(A)$
for the category of all (possibly infinite dimensional)
$A$-modules $M$ such that $M=M\alg$. Clearly,
$\modt(A)\sset \limind\modt(A),$
 and  any object of the category $\limind\modt(A)$ is a direct limit
of its finite dimensional submodules.

\ab Let $A=A^+\otimes A^\circ\otimes A^-$ be an algebra as above,
which is stable under the antipode {(anti-)}
homomorphism $S: A \to A$.
Given a subalgebra $\aa\sset A$, one has
a {\it smooth coinduction} functor 
$$\,\ind_{\aa}^{^{A}}: \modt(\aa) \too
\limind\modt(A),\quad
\ind_{\aa}^{^{A}}\ N:= \left(\Hom_{\aa}(A, N)\right)\alg\,,$$
where $\Hom_{\aa}(A, N)$ is an infinite dimensional vector space
with {\it left} $A$-action given by
$(x\cdot f)(y) = f(S(x)\cdot y)$, see e.g. [APW${}_{2}$].
It is clear that the functor $\ind_{\aa}^{^A}$ is the right adjoint to the
obvious restriction functor ${\res}_{\aa}^{^{A}}\,: \modt(A)\to
\modt(\aa),$ in other words,
for any $M\in \modt(A),\, N \in \modt(\aa)
\,,$ there is a canonical adjunction isomorphism
(Frobenius reciprocity):
\begin{equation}\label{adjun1}
\Hom_{\aa}({\res}_{\aa}^{^{A}}M\,,\,N) \simeq 
\Hom_{_{A}}(M\,,\,\ind_{\aa}^{^{A}} N)\,.
\end{equation}

\subsection{Ind- and pro-objects.}\label{tb}
We view  $\k[G]$, the coordinate ring of the
algebraic group $G$, as a $\g$-module via
the left regular representation. It is clear that
$\k[G]$ is a direct limit of its finite dimensional
$G$-submodules, that is an ind-object in the category
$\rep(\Ug)$. Let 
 $\fr\k[G]\in \limind\Umod$ 
be the corresponding Frobenius pull-back.

\ab Let $T\subset B$ be the maximal torus and the Borel
subgroup corresponding to the Lie algebras $\h\subset \b$,
respectively.
Given $\lambda\in \Y=\Hom(T,\k^\times)$, we 
let ${\mathbf{I}}_\lambda\, := \,\indf_{_T}^{^B}\lambda$ be
the induced $B$-module formed by the
regular algebraic functions on $B$ that
transform via the character $\lambda$ under right translations by $T$.
One can also view ${\mathbf{I}}_\lambda$ as
a locally-finite $\Ub$-module, which is smoothly co-induced up to  $\Ub$
from the character $\lambda: \Uh\to \k$.

\begin{definition}\label{tb1}
We introduce the algebra
 $\,\tb:=\bb\cdot\U^\circ=\u^+\otimes\U^\circ\;\sset\; \U$.
\end{definition}

\ab The algebra $\tb$ is  slightly
larger than $\bb$; it plays the same role as the
group scheme $B_1T\sset B$ plays for the Borel group in
a reductive group $G$ over
an algebraically closed field of positive characteristic, cf. \cite{Ja}.

\ab  For any $\la \in \Y$,
the weight ${l}\la$
clearly defines a one-dimensional representation of $\tb$.
The pull-back 
 via the Frobenius morphism yields the
following
 isomorphism
\beq\label{bI}
\ind_{\tb}^{\B}(\k_\tb({l}\la))\;\simeq\;
\fr({\mathbf{I}}_\lambda)
\eeq

\ab Given a module $M$ over any   Hopf algebra,
there is a well-defined notion of 
 contragredient module $M^\vee$
constructed using the antipode anti-automorphism.
The duality functor  $M\mapsto M^\vee$ 
commutes with the 
 Frobenius functor  $M\mapsto\fr M$,  and sends ind-objects
into pro-objects.

\ab Let $N\sset B$ denote the unipotent radical of $B,$
so that $\Lie N=\n$. 
View 
$\k[N]$ as 
the left regular $\Un$-representation.
We make $\k[N]$ into a $\Ub$-module by
letting the Cartan algebra $\h$ act on  $\k[N]$ via
 the adjoint action.
We have an $\Ub$-module isomorphism $\k[N]\simeq {\mathbf{I}}_0=\k[B/T].$
Thus,
$\fr({\mathbf{I}}_0)$ is a $\sB$-module.

\begin{lemma}\label{indres}
  For any $M\in\Umod$, there is
a natural  isomorphism of $\U$-modules: $
\ind_{\mathsf{u}}^{^{\U}} ({\res}_{\mathsf{u}}^{^{\U}}\,M) \simeq
M\otimes_{_{\k}}\fr\k[G].$
Similarly, for a finite dimensional $\sB$-module $M$,
there is
a natural  isomorphism of $\sB$-modules  $
\ind_\tb^{^{\sB}}({\res}_\tb^{^{\sB}}\,M) \simeq
M\otimes_{_{\k}}\fr\bbI_0.$
\end{lemma}

\noindent
{\sl Proof:} We view elements $x\in \Ug$ as left invariant
differential operators on $G$. One has a perfect
paring:
$\Ug \times \k[G]\to \k$, given by  $(x,\,f) \mapsto
(x f)({\mathbf{1}}_{_G}).$ This pairing induces
  a canonical isomorphism of $\Ug$-modules:
$\k[G] \iso \Hom_{_{\k}}(\Ug,\k)\alg$. Applying the Frobenius functor
$\phi^*$,
one obtains a natural $\U$-module isomorphism
$\,\ind\,\k_{\u}\; \simeq\; \fr\k[G]\,,$
(here and below we use shorthand notation $\ind:=\ind_{\mathsf{u}}^{^{\U}}$
and ${\res}:={\res}_{\mathsf{u}}^{^{\U}}$).

\ab Now, for any $L, M \in \Umod$ and $N\in \umod$, we
have:
$$
\begin{array}{ll}\dis
\Hom_{{{{\sf U}}}}\Bigl(L\,,\, 
\ind\,\bigl(({\res}\,M)\otimes N\bigr)\Bigr) &=
\Hom_{\u}\left({\res}\,L\,,\,({\res}\,M)\otimes N\right)\\
 &=\Hom_{\u}\bigl({\res}(M^\vee\otimes L)\,,\,N\bigr)\\
&=\Hom_{{{{\sf U}}}}\bigl(M^\vee\otimes L\,,\,\ind\,N\bigr)=
\Hom_{{{{\sf U}}}}\bigl(L\,,\, M\otimes\ind\,N\bigr)\;.
\end{array}
$$
Since $L$ is arbitrary, we get a functorial isomorphism
$\ind(N\otimes{\res}\,M) \simeq M\otimes\ind N$.
The first isomorphism of the Lemma
 now follows by setting $N=\k_{\u}$, and using 
the isomorphism
$\ind\,\k_{\u}\, \simeq\, \fr\k[G].$
The second isomorphism is proved similarly. \qed

\subsection{Hopf-adjoint action.} 
Given a Hopf algebra $A$, we always
 write $\Delta$ for the coproduct and
$S$ for the antipode in $A$,
and use  Sweedler notation:
$\Delta(a)=\sum a'_i\otimes a''_i$. The Hopf algebra structure on $A$
makes the category of left $A$-modules, resp. $A$-bimodules,
a monoidal category with respect to the tensor product $\otimes$ (over
$\k$).
  For each $a\in A$, the
map $\,\adh(a) :  \,
m \mapsto \sum a'_i\cdot m\cdot S(a''_r)$ defines
 a {\it Hopf-adjoint} $A$-action on any $A$-bimodule $M$,
 such that the action map:
$\,A\otimes M\otimes A\to M\,,\,a_1\otimes m\otimes a_2\mapsto a_1ma_2,$ is a morphism of
$A$-bimodules. We call an element $m\in M$
{\it central} if $am=ma\,,\,
\forall a\in A,$ and we say that   $m$ is
$\,{\adh}A$-{\it invariant}
if $\,{\adh}a(m)=\epsilon(a)\cdot m,\,
\forall a\in A.$
The following result is well-known, see e.g.
 \cite[Proposition 2.9]{APW}, \cite{Jo}.

\begin{lemma}\label{standard_hopf}
\vi For any $A$-bimodule $M$, 
the assignment: $a \mapsto {\adh}(a)$
gives an algebra homomorphism: $A \to \End_\k M$.

\ab \vii An element $m\in M$ is $\,{\adh}A$-invariant
 if and only if  it is central.
\hfill\qed
\end{lemma}

\begin{proposition}\label{GBbis} \vi$\,$ The Hopf-adjoint 
$\U_\A$-action on $U_q(\g)$ preserves the $\A$-module $\fU_\A$.
Similarly, the Hopf-adjoint $\sB_\A$-action
preserves the $\A$-modules: $\fU^+_\A\sset\fB_\A\sset\tB_q$.
These actions induce, after specialization, an $\U$-module structure
on $\fU$, resp. a
$\sB$-module structure
on $\fB$. 

\ab \vii 
The   subalgebras $\fC\eps\sset \fC\subset\fU^+\sset \fB$,
are all ${\adh}\sB$-stable.

\ab \viii The Hopf-adjoint action of the subalgebra $\bb\subset \sB$ on $\fC$
is trivial.
\end{proposition}

\begin{proof}
In $U_q$,
consider the set  $C=\{E^l_i\,,\,K_i^l-K_i^{-l}\,,\,F_i^l\}_{i\in I}.$
 According to De Concini-Kac \cite{DK},  every element
$c\in C$ projects to a central element
 of the algebra $\fU_\A/\m\cd\fU_\A$,
where $\m=(q-\zeta)$ is the maximal ideal corresponding to our root of
unity. Hence,
Lemma \ref{standard_hopf}(ii) implies that,
 for any $c\in C$, in the algebra   $\fU_\A/\m\cd\fU_\A$ one has
$\adh{c}(u)=0.$ Hence, $\adh{c}(u)\in \m\cd\fU_\A,$ for any $u\in
\fU_\A$. We deduce that the map
$$\adh\bigl({\frac{c}{q-\zeta}}\bigr):\;
 u\mto\frac{1}{q-\zeta}\cd\adh{c}(u)
$$ 
takes the $\A$-algebra $\fU_\A$ into itself. But the 
elements of the form $\frac{c}{q-\zeta}\,,\,c\in C,$ generate
the $\A$-algebra $\U_\A$. Thus, we have proved that
 $\adh\U_\A^{\,}(\fU_\A)\sset \fU_\A$.

\ab Next, from \eqref{coprod}, for any $i\in I,$  we find: 
${\adh}E_i(x)=
E_ix-K_ixK_i^{-1}E_i$,
and ${\adh}K_i(x)=K_ixK_i^{-1}.$ 
We deduce that the subalgebras
$\tU_q^+\subset\tB_q\sset \tU_q$
are stable under the adjoint $\tB_q$-action
on $\tU_q$. Further, a
 straightforward computation yields
that, for any $i\in I,$ one
has:
\begin{equation}\label{GB-preserved}
({\adh}E_i^{\,\,\li})\,E_j\;\subset\; \m\cdot\cO
\;,\quad\text{and}\quad
\bigl({\adh}(K_i^{\,\,\li}-1)\bigr)\,E_j\;
\subset \;\m\cdot\cO\,,
\end{equation}
It follows that the Hopf-adjoint action of the subalgebra
$\cB\subset \tB_q$ preserves the subspaces 
$\fU^+_\A\sset \cO\subset \tB_q$. Specializing at $\zeta$, we
 obtain part (i), as well as $\adh\sB$-invariance of $\fU^+$.
$\hphantom{x}$\hfill\end{proof}

\ab Recall the  complex semisimple group
 $G$ 
with Lie algebra $\g$. Let $B\sset G$, resp.  $\overline{B} \sset G,$
be the two opposite Borel subgroups with Lie algebras
$\b=\h\oplus\n$ and $\overline{\b}=\h\oplus\overline{\n}$, respectively.
We consider the flag manifold $G/B$ with base point $B/B$,
and the `opposite open cell' $\overline{B}\cdot B/B\sset G/B$.
The left $G$-action on $G/B$ induces a Lie algebra map
$\b=\Lie B\mto\text{\it Vector fields on } \overline{B}\cdot B/B$.
This makes the coordinate ring  $\k[\overline{B}\cdot B/B]$
a $\Ub$-module (via Lie derivative).\footnote{Note that since the open cell
 $\overline{B}\cdot B/B\sset G/B$ is not a $B$-stable subset in the
flag variety, the Lie algebra action of $\b$ on  $\k[\overline{B}\cd
B/B]$
cannot be exponentiated to a $B$-action.}

\begin{theorem}\label{GB_thm} There is an algebra isomorphism
$\fC\simeq \k[\overline{B}\cdot B/B]$
that intertwines 
 the Hopf-adjoint  action of the algebra $\Ub=\sB/(\bb)$
on  $\fC$ resulting from Proposition \ref{GBbis}(ii)
and the $\Ub$-action on $\k[\overline{B}\cd B/B]$
described above.
\end{theorem}

\begin{proof}
Let $D(B_q)$ be the Drinfeld double of the Hopf 
$\k(q)$-algebra $B_q$. Thus, $D(B_q)$ is a Hopf algebra
that contains $\overline{B}_q=U_q^\circ\otimes_{_{\k(q)}}U_q^-\sset D(B_q)$ 
and $B_q\sset D(B_q)$ as Hopf subalgebras,
and is isomorphic to $\overline{B}_q\otimes_{_{\k(q)}}B_q$
as a vector space. We combine differentiation-action, see \eqref{drin4} of
the algebra $\overline{B}_q$ on  $B_q$, and the $\adh B_q$-action on
$B_q$
 to obtain a map
$$a:\
 \bigl(\overline{B}_q\otimes_{_{\k(q)}}B_q\bigr)\,\otimes \,B_q\too
B_q\,,\quad
(\bar{b}\otimes b)\,\otimes \,b'\mto \partial_{_{\bar{b}}}\bigl(\adh b(b')\bigr)\,.
$$
One verifies that the commutation relations in
 $D(B_q)$ insure that the $\overline{B}_q$-
and $\adh B_q$-actions on   $B_q$ fit together
to make  the map $a$ a $D(B_q)$-algebra action  on
 $B_q$ such that the multiplication
map $B_q\otimes B_q\map B_q$ is a morphism of 
 $D(B_q)$-modules.

\ab By Proposition \ref{GBbis}(ii), the subspace $U^+_q\sset B_q$
is stable under the $\adh B_q$-action. Further,
the inclusion $\Delta(U_q^+)\sset B_q\otimes U_q^+$
implies that the subspace $U^+_q\sset B_q$
is also stable under the differentiation-action of
$\overline{B}_q$ on $B_q$. Thus, $U_q^+$
is a $D(B_q)$-submodule in  $B_q$.

\ab  According to Drinfeld, one has an algebra isomorphism
$D(B_q)\cong U_q\otimes U_q^\circ.$
Therefore, the $D(B_q)$-action on  $U^+_q$ constructed above gives, by
restriction, to the subalgebra $U_q\sset D(B_q)$ an $U_q$-action on
$U^+_q$. Further, specializing Drinfeld's paring \eqref{drin1}
at  $q=\zeta$, that is, using the pairing
$\U^-\otimes \fU^+\to\k$ and  Proposition \ref{GBbis}(i)-(ii),
we see that 
the action of the $\k(q)$-algebra  $U_q$ on $U_q^+$ defined above
can be specialized at  $q=\zeta$ to give
a well-defined $\U$-action on $\fU^+$.

\ab From Proposition \ref{GBbis}(ii)-(iii) we deduce

\pb{The  subspace $Z\sset \fU^+$ is stable under the $\U$-action on
$\fU^+$, moreover, the subalgebra $\u\sset \U$ acts
trivially (via the augmentation) on $Z$.}

\ab It follows that  the action on $Z$ of the algebra $\U$
factors through $\U/(\u)$. Thus, we have constructed
an action of the  Hopf algebra $\Ug$
(more precisely, of its completion $\barU$) on the algebra $Z$;
in particular, the Lie algebra $\g\sset \Ug$ of primitive elements
acts on $Z$ by derivations.

\ab Observe next that both the algebra $\Ug$ and the space
$Z$ have  natural $\Y$-gradings, and the $\Ug$-action on $Z$
is  clearly compatible with the gradings. Write $\Y^+\sset \Y$ for the
semigroup generated by the positive roots.
Clearly, all the weights occurring in $Z$ belong to
$\Y^+$, while all the  weights occurring in $\mathcal{U}\overline{\n}$
 belong to
$-\Y^+$. It follows that the action of $(\mathcal{U}\overline{\n})\eps$,
the
augmentation
ideal of
$\mathcal{U}\overline{\n}$, on
$Z$ is locally nilpotent, i.e.,
for any $z\in Z$, there exists $k=k(z)$ such that
$(\mathcal{U}\overline{\n})^{\,k}\eps z=0$.

\ab Further, the invariance of the Drinfeld pairing \eqref{drin1}
implies that the
perfect pairing $\mathcal{U}\overline{\n}\otimes Z\to\k$, see \eqref{drin3},
is a morphism of $\mathcal{U}\overline{\n}$-modules.
This gives an
$\mathcal{U}\overline{\n}$-module imbedding $ Z\into 
\Hom_\k(\mathcal{U}\overline{\n},\k)$, where $\Hom_\k(\mathcal{U}\overline{\n},\k)$ 
is viewed as a contregredient representation to
the left regular  representation of 
the algebra $\mathcal{U}\bar{\n}$ on itself.

\ab Now, let $f\in \k[\bar{N}]$, be
a  regular function on the unipotent group
$\bar{N}\sset G$
corresponding to the Lie algebra $\bar{\n}$.
The assignment $\,u\mto\hat{f}(u):= (uf)(1)$, where $u$
runs over the space of left invariant differential operators on the group
$\bar{N}$, 
gives a linear function on $\mathcal{U}\bar{\n}$,
hence
an element $\hat{f}\in \Hom(\mathcal{U}\bar{\n},\k).$ 
 The map $f\mto\hat{f}$ identifies
the coordinate ring $\k[\bar{N}]$ with the
subspace 
$$\{\psi\in\Hom(\mathcal{U}\bar{\n},\k)\enspace|\enspace
\exists k=k(\psi)\enspace\text{such
that}\enspace(\mathcal{U}\bar{\n})^{\,k}\eps\psi=0\}.
$$
But the action of  $(\mathcal{U}\bar{\n})\eps$
on $Z$ being locally-nilpotent, we see that
the image of the imbedding  $Z\into 
\Hom(\mathcal{U}\bar{\n},\k)$
must be contained in the space above, that is,
in  $\k[\bar{N}]$.
This way, we obtain an $\mathcal{U}\bar{\n}$-module imbedding
$Z\into \k[\bar{N}]$. We claim that $Z=\k[\bar{N}]$.
Indeed, the group $\bar{N}$ is isomorphic as an algebraic variety
to a vector space $V$, so that the pairing
$\mathcal{U}\bar{\n} \otimes \k[\bar{N}]\to \k$
may be identified with the canonical pairing
$\sym(V) \otimes \sym(V^*)\map \k$.
Hence, if $Z$ were a proper subspace in $\k[\bar{N}]=\sym(V^*)$ then,
the pairing $\mathcal{U}\bar{\n} \otimes Z\to \k$
could not have been perfect. Thus,  $Z=\k[\bar{N}]$.

\ab Next we use the perfect pairing to identify
$\mathcal{U}\bar{\n}$ with $Z^\dag$, the continuous dual (in the adic topology)
of $Z=\k[\bar{N}]$. The  $\Ug$-module structure on
$Z$ defined above gives rise to an   $\Ug$-module structure on
 $Z^\dag$. It is clear that restricting the  $\Ug$-action to the subalgebra
$\mathcal{U}\bar{\n}\sset\Ug$ we have

\pb{$Z^\dag$ is a rank 1 free $\mathcal{U}\bar{\n}$-module
generated by the element $\epsilon\in Z^\dag$, and}

\pb{The action of the ideal $\Un\eps\sset \Un\sset\Ug$ annihilates
the element $\epsilon\in Z^\dag$.}

Observe that these two properties, combined with the commutation
relations in the algebra $\Ug$, completely determine
the  $\Ug$-module structure on
 $Z^\dag$.

\ab Now, the Lie algebra $\g$ acts on
 $\k[\overline{B}\cdot B/B]$, the coordinate ring of the "big cell".
Let $\k[\overline{B}\cdot B/B]^\dag$ denote the continuous dual
equipped with a natural $\Ug$-module structure.
This latter  $\Ug$-module also 
satisfies the two
properties above. Therefore, there exists
 an $\Ug$-module isomorphism $Z^\dag\cong
\k[\overline{B}\cdot B/B]^\dag$.
Dualizing, we obtain  an $\Ug$-module isomorphism
$Z\cong\k[\overline{B}\cdot B/B]$,
and the Theorem is proved.
\end{proof}
 
\begin{remark}
 Theorem \ref{GB_thm} is closely related to the
results of De Concini-Kac-Procesi on "quantum coadjoint
action". In particular, it was shown in \cite[Theorem 7.6]{DKP}
that the isomorphism
$Z\cong\k[\overline{B}\cdot B/B]$ is a Poisson algebra
isomorphism.
 $\quad\lozenge$\end{remark}

\ab  General properties of Hopf-adjoint actions imply,
by part (i) of Proposition \ref{GBbis}, that the product map:
$\,\,\fB\otimes \fB\to \fB$ is a morphism of
$\sB$-modules. Therefore, we deduce from (ii) that
$(\fC\eps)^2\sset\fC\eps$ is again an  ${\adh}\sB$-stable
subspace. This makes the finite dimensional vector
space  $\fC\eps/(\fC\eps)^2$ into a $\b$-module,
and we have

\begin{corollary}\label{tangent}
\vi There is a $\b$-equivariant vector space
isomorphism $\fC\eps/(\fC\eps)^2\simeq \n$,
where the Lie algebra $\b$ acts on $\n$ via the adjoint action.

\ab  \vii  There is a $\b$-equivariant graded algebra isomorphism
$\Tor^{Z}(\k_{_Z},\k_{_Z})\cong\wedge\hdot\n.$
\end{corollary}

\proof The cotangent space to $G/B$ at the base
point is the space $(\g/\b)^*=\b^\perp\sset\g^*$, which is  $\b$-equivariantly
isomorphic to $\n\sset\g$ via an invariant bilinear form on $\g$.
Now, by Theorem \ref{GB_thm}, one may identify
the open cell $\overline{B}\cd B/B$ with $\Spec Z$. 
The base point goes under this identification,
to the augmentation ideal $Z\eps\in\Spec Z$.
The cotangent space  $T^*\eps(\Spec Z)$ at that
point equals $\fC\eps/(\fC\eps)^2,$ by definition.
This yields part (i) of the Corrolary.

\ab To prove (ii) recall that, for any smooth affine variety $X$
and a point $x\in X,$ one has a canonical
 graded algebra isomorphism
$\Tor^{\k[X]}(\k_x,\k_x)\cong\wedge\hdot(T^*_xX),$
where $\k_x$ denotes the 1-dimensional $\k[X]$-module
corresponding to evaluation at $x$.
Part (ii)  now follows from~(i).
$\hphantom{x}$\hfill\qed

\begin{remark} In this paper we will not
use the isomorphism
$\fC\simeq\k[\overline{B}\cd B/B]$ itself
but only the resulting isomorphism of  Corollary
\ref{tangent}.
 $\quad\lozenge$\end{remark}

\subsection{Cross-product construction.} In the next section
 we will use the following general
construction, see e.g. [Mo]. 
Let
$\aa$ be an associative algebra and $A$ a Hopf algebra.

\begin{proposition}\label{Montgomery}
Let  $A$  act on
$\aa$ in such a way that
 the multiplication map $m: \aa\otimes \aa \to \aa$ is
a morphism of $A$-modules. Then there is a 
natural  associative algebra structure
on the vector space $A\otimes\aa$, to be denoted $A\ltimes\aa$,  such that

\ab \vi$\;\,$ $\aa=1\otimes \aa\,$ and $\,A=A\otimes 1\,$ are subalgebras 
in $A\ltimes\aa$;

\ab \vii$\;\,$ The $\aa$-action on itself by left multiplication, and the
$A$-action on $\aa$ can be combined to give
 a well-defined $A\ltimes\aa$-action on $\aa$.
Furthermore, the multiplication map
$(A\ltimes\aa)\otimes (A\ltimes\aa)\too(A\ltimes\aa)$
is $\adh A$-equivariant with respect to the
tensor product $A$-module structure on $A\ltimes\aa$.

\ab\viii For any Hopf algebra $A$ acting on itself by the Hopf-adjoint action, 
the assignment $a\ltimes a_1\mto a\otimes  (a\cdot a_1)$ gives an algebra
isomorphism $\gamma: A\ltimes A \iso A\otimes A$.
\qed
\end{proposition}

\ab In the special case of a Lie algebra ${\mathfrak{a}}$ acting by derivations on an
associative algebra $A$, the Proposition reduces, 
for $A=\mathcal{U}{\mathfrak{a}}$, to the
very well-known construction of the cross product algebra
$\mathcal{U}{\mathfrak{a}}\ltimes\aa$.
We will apply Proposition \ref{Montgomery}
to the Hopf algebra $A:=\U$ acting
 via the Hopf-adjoint action
on various algebras
described in Proposition \ref{GBbis}. 

\subsection{Cohomology of Hopf algebras.}\label{cohomology}
Given an augmented algebra $\aa$ and a left $\aa$-module
$M$, one defines the cohomology of $\aa$ with coefficients in $M$
as $H\hdot(\aa,M):=\Ext\hdot_{\aa\mmod}({\k_\aa},M).$
In particular, we have $\Ext^0_{\aa\mmod}({\k_\aa},M)=M^\aa
:=\{m\in M\;|\; \aa\eps m=0\},$ is the
space of $\aa$-invariants in $M$.
Further, the space
$H\hdot(\aa,\k):=\Ext\hdot_{\aa\mmod}({\k_\aa},{\k_\aa}),$
has a natural graded algebra structure given by the Yoneda product.
This
 algebra structure is made explicit by identifying 
$\Ext\hdot_{\aa\mmod}({\k_\aa},{\k_\aa})$ with
 the cohomology algebra
of the \dg-algebra $\Hom_{\aa}(P , P )$,
where $P=\,(P^i)\, $ is a projective $\aa$-module resolution of $\k_\aa$.

\ab Let $A$ be an augmented algebra and $\aa\sset A$
a normal  subalgebra, see Definition \ref{normal}.
  For any left $A$-module $M$, the space
$M^\aa\sset M$ of $\aa$-invariants is $A$-stable.
Moreover, the $A$-action on $M^\aa$ descends to the
quotient algebra $A/(\aa)$ and, clearly, we have
$M^A=(M^\aa)^{A/(\aa)}$.  According to general
principles, this gives rise
to a spectral sequence
\beq\label{sp_seq}
  E_2^{p,q}=
H^p\bigl(A/(\aa)\,,\,H^q(\aa, M)\bigr)\quad\Longrightarrow
\quad E_\infty^{p+q}=\gr H^{p+q}(A,M)\,.
\eeq

\ab Below, we will use a special case of the
spectral sequence  where $A$ is a  Hopf algebra and $\aa$ is a normal Hopf
subalgebra. Fix two left $A$-modules $M, N$.
The vector space $\Hom_\k\hdot(M,N)$ has a natural 
 structure of $A$-bimodule,
hence of a left $A$-module, via the Hopf-adjoint action.
Observe further that we may identify the subspace
$\Hom_\aa\hdot(M,N)\sset\Hom_\k\hdot(M,N)$
with the space
$\bigl(\Hom_\k\hdot(M,N)\bigr)^\aa=$
$\{f\in \Hom_\k\hdot(M,N)\;\big|\;a\ccirc f=f\ccirc a\,,\,\forall a\in\aa\}$
of  {\it central} elements of $\Hom_\k(M,N)$,
viewed as  an $\aa$-bimodule. But any central element 
of an $\aa$-bimodule is $\adh\aa$-{\it invariant}, 
by Lemma \ref{standard_hopf}(ii).
Therefore, we deduce that the  Hopf-adjoint
$\aa$-action on $\Hom_\aa\hdot(M,N)$
is trivial. Hence, the $\adh A$-action on  $\Hom_\aa\hdot(M,N)$
descends to the algebra
$A/(\aa)$, and the spectral sequence in \eqref{sp_seq} yields

\begin{lemma}\label{resind} 
  For any left $A$-modules $M,N,$ there is a natural 
Hopf-adjoint $A/(\aa)$-action on $\Ext\hdot_\aa(M,N)$; it
 gives rise to a
spectral sequence
$$  E_2^{p,q}=
H^p\bigl(A/(\aa)\,,\,\Ext^q_\aa(M,N)\bigr)\quad\Longrightarrow
\quad E_\infty^{p+q}=\gr\Ext_A^{p+q}(M,N)\,.\qquad\square
$$
\end{lemma}

\section{Algebraic category equivalences}\label{sec_block}
 The goal of the next three sections is to
 construct
a chain of functors that will provide the following 
equivalences
of triangulated categories (undefined notations will be explained later):
\begin{align*}
\xymatrix{
D^G_{\text{coherent}}(\NN)
\ar[rr]^<>(.5){\text{restriction}}_<>(.5){\sim}&&
D^B_{\text{coherent}}(\n)\ar[rr]^<>(.5){\text{global}}_<>(.5){\text{sections}}&&
D\bigl(\Ub\ltimes\k\hdot[\n]\bigr)}\qquad\qquad\qquad\qquad\quad\hphantom{x}\\
\hphantom{x}\qquad\qquad\quad\xymatrix{
\ar[rr]^<>(.5){\text{Koszul}}_<>(.5){\text{duality}}&&
  D\bigl(\Ub\ltimes(\wedge\hdot\n[1])\bigr)
\ar[rr]_<>(.5){\S\ref{sec_formality}}^<>(.5){\text{formality}}&&
\modB\ar[rr]_<>(.5){\S\ref{sec_induction}}^<>(.5){\text{induction}}&&
D\bigl(\cat\bigr).}
\end{align*}
The composite of the equivalences above will yield  (a non-mixed version of)
a quasi-inverse of the functor $Q$ on the left of diagram
\eqref{sum-up}.

\subsection{Reminder on dg-algebras and dg-modules.}\label{reminder}
Given an algebra $A$, we write $A\bimod$ for the category
of $A$-bimodules, and
${\underline{A}}$ for $A$, viewed  either as a rank one
free left $A$-module, or as an
$A$-bimodule.
Similar notation will be used  below for differential graded
(\dg-)algebras.

\begin{notation}  Let $\,[n]\,$ denote the shift functor in the derived
category, and also the grading shift by $n$ in  a
dg-algebra or a
\dg-module. 
\end{notation}

\ab Given a dg-algebra $A=\bigoplus_{i\in\Z}\,A^i$,
write $\DGMod(A)$ for the homotopy category of all dg-modules
$M=\bigoplus_{i\in\Z}\, M^i$ over $A$ (with
differential $d: M\hdot\to M^{\hhdot+1}$),
 and $D(\DGMod(A))$ for the corresponding derived category.
 Given two objects  $M,N \in \DGMod(A)$
and $i\in\Z$,
we put $\Ext_{_{A}}^i(M,N):=
\Hom_{_{\DGMod(A)}}(M,N[i]).$
The graded space
$\Ext_{_{A}}\hdot(M,M)$
$=\bigoplus_{j\geq
0}\,\Ext_{_{A}}^j(M,M)$
 has a natural algebra structure, given by the Yoneda
product.

\ab An object $M \in \DGMod(A)$ is said to be
{\it projective} if it belongs to  the smallest full
subcategory of $\DGMod(A)$ that contains the rank one
\dg-module $A$,   and which
is closed under taking
mapping-cones and infinite direct sums. 
Any object of
$\DGMod(A)$ is quasi-isomorphic to a projective
object, see [Ke] for a proof. (Instead of projective objects, one can
use {\em semi-free objects} considered e.g.
in \cite[Appendix A,B]{Dr}.)

\ab
Given $M\in \DGMod(A)$, 
choose a quasi-isomorphic projective object $P$.
The graded vector space
$\bigoplus_{n\in\Z}\,\Hom_\k(P, P[n])$
has a natural  algebra structure given by composition.
Commutator with the differential $d\in \Hom_\k(P, P[-1])$
makes this algebra into a \dg-algebra, to be denoted
$\REnd_{_{A}}(M):=
\bigoplus_{n\in\Z}\,\Hom_\k(P, P[n]).$ This \dg-algebra does not depend,
 up to
quasi-isomorphism,
on the choice of  projective representative $P$.

\ab Given a \dg-algebra
morphism
$f: A_1\to A_2$, we let $f_*:
D(\DGMod(A_1))\to D(\DGMod(A_2))$ be
the push-forward
functor
$M\mapsto f_*M:=A_2\lotimes_{A_1} M$,
and
$f^*: D(\DGMod(A_2))\map D(\DGMod(A_1))
$
 the pull-back functor given by change of scalars.
Note that the functor $f^*$ is the right adjoint of $f_*$
and, moreover, if the map $f$ is a \dg-algebra quasi-isomorphism
then the functors $f_*$ and $f^*$ are triangulated equivalences,
quasi-inverse to each other.

\subsection{From coherent sheaves on ${\NN}$ to {$\k[\n]$}-modules.}
\label{from_coh} We say that a linear action of an
 algebraic group ${\mathbf{G}}$ on a (possibly infinite-dimensional)
vector space $M$ is {\em algebraic} if
$M$ is a union of finite-dimensional ${\mathbf{G}}$-stable
subspaces $M_s$ and, for each $s$, the action homomorphism
${\mathbf{G}}\to \op{GL}(M_s)$ is an
 algebraic group homomorphism.

\begin{notation}\label{cohG} Given a 
$\k$-algebra $A$ and an algebraic group ${\mathbf{G}}$, acting {\em algebraically}
 on
$A$ by algebra automorphisms, we let $\Mod^{{\mathbf{G}}}(A)$
denote the abelian category of ${\mathbf{G}}$-equivariant
$A$-modules, i.e.,
$A$-modules $M$ equipped with an algebraic
${\mathbf{G}}$-action such that the action-map
$A\otimes M\to M$ is  ${\mathbf{G}}$-equivariant. If, in addition,
$A$ has a $\Z$-grading preserved by the ${\mathbf{G}}$-action,
we write $\,\Mod^{{\mathbf{G}}\times\Gm\!}(A),$
where  $\Gm$ denotes the multiplicative group, 
for the abelian category of ${\mathbf{G}}$-equivariant
 $\Z$-graded  $A$-modules.

\ab If $A$ is noetherian, we let $\Mod^{{\mathbf{G}}}_f(A)$
be the full subcategory in  $\Mod^{{\mathbf{G}}}(A)$
formed by finitely-generated $A$-modules.
In the case of trivial group ${\mathbf{G}}=\{1\}$
we drop the superscript ${\mathbf{G}}$
and write $\Mod_f(A)$ for the corresponding category.

\ab Given an algebraic ${\mathbf{G}}$-variety $X$, we write
$\Coh^{\mathbf{G}}(X)$ for the abelian category of
${\mathbf{G}}$-equivariant
coherent sheaves on $X$. A quasi-coherent sheaf on $X$ is
said to be ${\mathbf{G}}$-equivariant if it is a direct
limit of its  ${\mathbf{G}}$-equivariant coherent subsheaves.

\ab 
Given, a $\Gm$-equivariant  sheaf $\fcal$, 
write  $z^k\otimes\fcal$ for a $\Gm$-equivariant  sheaf
obtained by twisting $\Gm$-equivariant 
structure on $\fcal$
by the 1-dimensional character $z\mapsto z^k$ of the group
$\Gm$.
\end{notation}

 \ab The  Borel subgroup $B\sset G$  acts on $\n$,
 the nilradical of $\b=\Lie B$,  by  conjugation. 
Further, the multiplicative 
group $\Gm$ acts on $\n$ by dilations:
we let $z\in \Gm$ act via multiplication by $z^2$.
The two actions commute, making
$\n$ a ${B\times\Gm}$-variety.
This gives a ${B\times\Gm}$-action on
the polynomial algebra $\k\hdot[\n]$.
In particular,  the algebra $\k\hdot[\n]$ 
acquires a natural  $\Z\times \Y$-grading:
the $\Z$-component of the grading
is given by {\em twice} the degree of polynomial
(this is consistent with our earlier convention that
$z\in \Gm$ act via multiplication by $z^2$),
and the $\Y$-component  of the grading
is induced from the natural $\Y$-grading on $\n$.

\ab We also consider the  Springer resolution $\NN := G\times_B \,\n$,
which
is a $G\times\Gm$-variety in a natural way.

\ab The closed imbedding $i :\n=\{1\}\times_{_B}\n\into
G\times_{_B}\n=\NN$ gives rise to a natural restriction
functor $i^* : \Coh^G(\NN)\to\Coh^B(\n)$. This functor is an equivalence
of categories whose inverse is
provided by the  induction functor $\Ind: 
\Coh^B(\n)\to \Coh^G(\NN).$
Further, the variety $\n$ being affine, we deduce that
the functor of global sections
yields an  equivalence of abelian 
categories $\Gamma: \Coh^B(\n)\iso \Mod^B_f(\k\hdot[\n]),
\,\fcal\mapsto \Gamma(\n,\fcal).$

\ab We are now going to extend the considerations above to
triangulated categories. In general, let
 ${\mathbf{G}}$ be  an algebraic
group and $X$  a ${\mathbf{G}}\times \Gm$-variety.
We let  $\dgcoh^{{\mathbf{G}}}(X)$
 denote the category whose
objects are diagrams $\fcal_+\leftrightarrows\fcal_-$,
where $\fcal_+,\fcal_-$ are ${\mathbf{G}}\times \Gm$-equivariant quasi-coherent sheaves
on $X$,  such that

\pb{The arrows in the diagram are morphisms of quasi-coherent
sheaves;
specifically, we have  ${\mathbf{G}}\times\Gm$-equivariant
 morphisms
$\partial: \fcal_+ \map z\otimes\fcal_-,$ and $\partial: \fcal_-\map
z\otimes\fcal_+,$ such that $\partial
\ccirc\partial=0$.}

\pb{The  cohomology sheaves (with respect to the
differential $\partial$) ${\mathcal{H}}(\fcal_\pm)$  are coherent.}

\ab 
Similarly, let   $A$ be a dg-algebra,
such that the cohomology
$H\hdot(A)$ is a finitely generated
Noetherian algebra. Let an algebraic group ${\mathbf{G}}$
act  algebraically on $A$ by algebra automorphisms
preserving the grading and commuting with the differential.

\begin{notation}\label{dgmodG} We write $\DGMod^{{\mathbf{G}}}_f(A)$
for a subcategory  in $\Mod^{{\mathbf{G}}\times\Gm}(A)$
formed by differential graded $A$-modules
$M=\bigoplus_{i\in\Z}\, M^i$ (with
differential $d: M\hdot\to M^{\hhdot+1}$)
such that
the ${\mathbf{G}}$-action on $M$ preserves
the grading and commutes with $d$ and, moreover,

\pb{The  cohomology $H\hdot(M)$ is a finitely-generated $H\hdot(A)$-module.}
\end{notation}

\ab
The category  $\dgcoh^{{\mathbf{G}}}(X)$, resp. 
 $\DGMod^{{\mathbf{G}}}_f(A)$,
has a natural structure of  homotopy category,
and we write $\dcoh^{\mathbf{G}}(\n)$, resp. 
$D^{{\mathbf{G}}}_f(A)$,
for the corresponding derived category obtained by localizing at
quasi-isomorphisms. Again, if ${\mathbf{G}}=\{1\}$,
the superscript ${\mathbf{G}}$ will be   dropped.

\begin{notation}\label{SS}
Let $\SS:=\k\hdot[\n]=\sym\hdot(\n^*[-2])$, resp.,  $\La:=\wedge\hdot(\n[1]),$
be the Symmetric algebra
of $\n^*$, resp., the exterior  algebra
of $\n$,  viewed as a differential graded
algebra with zero differential generated by the
 vector space
 $\n^*$ placed in grade degree 2,
 resp.,  by the
 vector space
 $\n$, placed in grade degree $(-1)$.
\end{notation}

\ab Thus,  we have triangulated categories
$\dcoh^{B}(\n),\,
\dcoh^{G}(\NN)$, and also
$D^B_f(\k\hdot[\n])=D^B_f(\SS)$ and  $D^B_f(\La).$
As above, we obtain  equivalences
\begin{equation}\label{res_gamma}
\xymatrix{
\dcoh^{G}(\NN)
\enspace
\ar@<1ex>[rr]^<>(0.5){i^*}&&
\dcoh^{B}(\n)
\ar@<1ex>[ll]^<>(0.5){\RInd}
\ar[rrrr]^<>(0.5){(\fcal_+,\fcal_-)\mto \Gamma(\n,\fcal_+\oplus\fcal_-)}
_<>(0.5){\sim}&&&&
D^B_f(\SS).}
\eeq

\subsection{Koszul duality.}\label{koszul_duality}
Let $\underline{\SS}$ be a free rank 1 graded $\SS$-module
(with generator in degree zero), and 
  $\k_\wedge$  the trivial 
$\La$-module (concentrated in degree zero).
We may view $\uSS$ as an object of $D_f^B(\SS),$
and $\k_\wedge$ as an object of $D_f^B(\La)$.
We recall the well-known  Koszul duality,
cf. \cite{BGG} and \cite{GKM}, \cite{BGS},
between the graded algebras $\SS$ and $\La$.
The  Koszul duality provides
an equivalence:
\begin{align}\label{kappa}
\kappa:   D_f^B(\SS)\iso
D_f^B(\La),\quad\text{such that}
\enspace\kappa(\underline{\SS})=\k_\wedge.
\end{align}

\ab In the {\em non-equivariant} case, this equivalence
amounts, essentially, to a dg-algebra quasi-isomorphism
\begin{equation}\label{SL_kos}
\RHom\hdot_{D_f(\La)}(\k_\wedge,\k_\wedge)\qisto\Ext\hdot_{D_f(\La)}(\k_\wedge,\k_\wedge)
\cong\SS,
\eeq
given by the standard Koszul complex.

\ab The reader should be warned that
 the Koszul duality  we are using here
is slightly different from the one used in \cite{BGG}
in two ways.

\ab First of all,
we consider $\SS$ and $\La$ as dg-algebras, and
our  Koszul duality is an equivalence of
derived categories of the corresponding homotopy categories 
of dg-modules over $\SS$ and $\La$, respectively.
In \cite{BGG},
the authors consider  $\SS$ and $\La$ as {\em plain}
graded algebras, and establish an equivalence between
derived categories of
the {\em abelian} categories of {\em $\Z$-graded} $\SS$--modules
and
$\La$-modules, respectively. In our case, the proof of the equivalence
is very
similar; it is discussed in \cite{GKM} in detail.

\ab Second, in \cite{BGG} no
  $B$-equivariant structure was involved.
However,  
  the
construction of the equivalence given
in \cite{BGG},\cite{GKM}, being canonical,
it extends in a straightforward maner
to the equivariant setting as well.

\subsection{The principal block of $\U$-modules.}\label{W_aff_gr} We keep the notation
introduced in \S2, in particular, we fix a
primitive  root of unity of order $l$.
 Form the semidirect product 
$\waf:=W\ltimes \Y$, cf. (\ref{XY}), to be called
 the affine Weyl group.
Let $\rho=\frac{1}{2}\sum_{i\in I}\,\alpha_i\in \X$
be the half-sum of positive roots.
We define a $\waf$-action on the lattice  $\X$ as follows.
We let an element
$\la\in\Y\subset \waf$ act on $\X$ 
by  translation: $\tau\mto \tau+{l}\la$,
and let the subgroup $W\subset \waf$ act on $\X$ via the
`dot'-action 
centered at $(-\rho)$, i.e.,
for an arbitrary element $w_a=(w\ltimes\lambda)\in\waf$ we put:
$$w_a=w\ltimes\lambda:\;\; \tau\; \mto\; 
w_a\bullet\tau:=w({l}\la+\tau-\rho)+
\rho,\quad\forall \tau\in\X
\,.
$$

\ab 
Recall that  the weight $w_a\bullet 0 \in\X^{++}$
 is {\it dominant}, see (\ref{XY}), if and only if $w_a\in \waf$
is
the minimal element (relative to the standard Bruhat order)
in the corresponding left  coset $W\cdot w_a$.
    From now on we will identify the set of such minimal length
representatives in the left cosets
$W\backslash{\waf}$ with the
 lattice $\Y$ using the natural bijection: 
\beq\label{max}
\Y\stackrel{_\sim}{\longleftrightarrow}W\backslash{\waf}\;\;,\;\;
\lambda\mto\; w_{_{\sf min}}^\lambda:=\,
\mbox{\it minimal element in the coset } W\cdot(1\ltimes\lambda)\;\,
\sset \waf\,.
\eeq
Thus, at the same time we also obtain a bijection
\beq\label{min}
\Y\;\;\iso\;\;
\X^{++}\,\;\cap\,\;\waf\bullet 0 ,\quad
\lambda\;\mto\; w_{_{\sf min}}^\lambda\bullet 0
\eeq

\ab For each $\nu\in \X$ we let $L^\nu$ be the
simple $\U$-module with highest weight $\nu$.
It is known that $L^\nu$ is finite dimensional
if and only if $\nu\in\X^{++}$.

\ab For $\la\in \Y$, let $\nu=w_{_{\sf min}}^\lambda\bullet 0\in
\X^{++}.$ Writing  $w_{_{\sf min}}^\lambda$ in the form
 $w_{_{\sf min}}^\lambda=w\ltimes\la
\in W\ltimes\Y$, we have
$w_{_{\sf min}}^\lambda\bullet 0=w\bullet({l}\la),$
 where $w\in W$
is the unique
element   such that $w\bullet({l}\la) \in\X^{++}$.

\begin{definition}
  For $\la\in\Y,$ let $L_\lambda$ denote the simple $\U$-module with highest weight
$\nu=w_{_{\sf min}}^\lambda\bullet 0$. 
Thus, 
\beq\label{Lsimple}
L_\lambda=L^\nu\quad\text{has highest 
weight}\quad\nu=w({l}\la-\rho)
+\rho
= w\bullet({l}\la).
\eeq
\end{definition}
\ab Note that if $V\in \rep(G)$ is a simple $G$-module with highest weight
$\la\in\Y^{++},$ then $w_{_{\sf min}}^\lambda=1$,
hence, we have $\fr{V}=L_\la$.

\ab In this paper we will be frequently using the following
three partially ordered sets: 
{\renewcommand{\arraystretch}{1.3}
\beq\label{order}
\begin{tabular}{|l|lcl|}
\hline
 $(\Y, \succeq)$ & $\la\succeq\mu$ &$\Longleftrightarrow$ &
 $w_{_{\sf min}}^\la\bullet 0\,-\,
 w_{_{\sf min}}^\mu\bullet 0\enspace
 \mbox{= sum of a number of $\check\alpha_i$'s}$\\
 $(\Y, \rhd)$ & $\lambda \rhd\mu$ &$\Longleftrightarrow$ &
 $\la-\mu\in \Y^{++}=\Y\cap\X^{++}\,$\enspace is dominant.\\
$(\Y, \ge)$ & $\lambda\ge\nu$ &$\Longleftrightarrow$ &
 $\lambda -\nu \enspace
 \mbox{= sum of a number of $\alpha_i$'s}$\\
 \hline
\end{tabular}
\eeq
}

We write $\la\succ\mu$ whenever $\la\succeq\mu$ and $\la\neq\mu$.
This order relation corresponds, as will be
explained in section 6 below, to the closure
relation among  Schubert varieties in the loop Grassmannian
for the Langlands dual group $G^\vee$.

\begin{definition}\label{cat}
Let $\cat\,\sset\,\Umod$ be the `principal block', i.e., the
full subcategory of the abelian category
 $\U\mmod$ formed by finite-dimensional $\U$-modules
$M$ such that all
 simple subquotients  of $M$
are of the form $L_\lambda\,,\,\lambda\in \Y,$
see (\ref{Lsimple}).
\end{definition}

\ab The abelian category $\cat$ is known to have enough
projectives and injectives, and we let $D^b\cat$ denote
the corresponding bounded derived category.

\subsection{Induction.}\label{induction_subsec}
  Given a subalgebra $A^+\otimes A^\circ\sset \U$ and
 an algebra map $\la: A^\circ\to\k$,
we write $\ind_{_{{A^{^\pm}}
\otimes A^{^\circ}}}^\U(\la)$ instead of $\ind_{_{{A^{^\pm}}
\otimes A^{^\circ}}}^\U\left(\k_{_{{A^{^\pm}}
\otimes A^{^\circ}}}(\la)\right)$.

\ab Using injective resolutions as in
\cite{APW}, one defines a {\em derived induction functor}
$\rind_{_\B}^{^\U}$ 
 corresponding to a smooth (co)-induction functor $\ind_{_\B}^{^\U}$.
Let $\tR^i\ind_{_\B}^{^\U}$ denote its $i$-th cohomology functor.

\ab Let $\ell: \waf \to \Z_{\ge 0}$ denote the standard length function on 
$\waf$.
We recall the following (weak)
 version of  Borel-Weil theorem for  quantum groups, proved in [APW].

\begin{lemma}[Borel-Weil theorem]\label{borel_weil}
  For $\lambda\in \Y$, let $w\in W$ be an element of minimal length
such that $w\bullet({l}\la)\in\X^{++}$. Then we have:

\pb{$\tR^{\ell(w)}\ind_{_\B}^{^\U}({l}\la)$
contains $L_\la$ as a simple subquotient with multiplicity one;
Any other simple subquotient of 
$\tR^{\ell(w)}\ind_{_\B}^{^\U}({l}\la)$
is  isomorphic
to $L_\mu$ with $\mu\prec\la$.}

\pb{For any $j\neq \ell(w),$ each simple subquotient
of $\tR^j\ind_{_\B}^{^\U}({l}\la)$ is isomorphic
to $L_\mu$ with $\mu\prec\la$.\qed}
\end{lemma}

  \ab  From Lemma \ref{borel_weil} we immediately deduce

\begin{corollary}\label{generate} \vi
  For any $\la\in \Y$, we have $\RInd_{_{\B}}^{^{\U}}(l\lambda)
\in D^b\cat$.

  \ab \vii
The category $D^b\cat$ is generated,
as a triangulated category, by the family of objects
$\{\RInd_{_{\B}}^{^{\U}}(l\lambda)\}_{\,\lambda\in \Y}$.\qed
\end{corollary}

\begin{definition}\label{modB_cat}
Let $\modB$ be a  triangulated category whose
objects are complexes  of $\Y$-graded $\B$-modules
$M=\{\ldots\to M_i\to M_{i+1}\to\ldots\},\,i\in\Z,
M_i=\bigoplus_{\nu\in\Y}\,M_i(\nu),$
such that

\pb{For any $i\in\Z$, we have $M_i=(M_i)\alg$, see Sect. \ref{smooth}.}

 \pb{We have  $um=\nu(u)\cdot m$ for any $u\in U^\circ, m\in M_i(\nu),
 i\in\Z,
\nu\in\Y$;}

\pb{The total cohomology module,  $H\hdot(M)=\bigoplus_{i\in\Z}\,H^i(M),$ 
has a finite composition series with all subquotients of the form
$\Bbbk_{_{\B}}(l\la)\,,\,\lambda\in\Y$.}
\end{definition}

\begin{remark} 
 Note that the subalgebra $\bb\sset\B$ acts
trivially on the module $\Bbbk_{_{\B}}(l\la)$, for any $\lambda\in\Y$.
\end{remark}

\ab
In section \S\ref{sec_induction}   we will prove 

\begin{theorem}[Induction theorem]\label{proposition_roma1}
 The  functor  $\rind_{_{\B}}^{^{\U}}$ yields
 an equivalence
of triangulated categories $\modB \iso D^b\cat$.
\end{theorem}

\begin{remark}\label{charp}
An analogue of Theorem \ref{proposition_roma1} holds also for the
principal block of complex representations of the algebraic group
$G({\mathbb{F}})$ over an algebraically closed field of characteristic
$p>0$. Our proof of the Theorem  applies  to the
latter
 case as well.
 $\quad\lozenge$\end{remark}

\subsection{Quantum group Formality theorem.}\label{some}
The second main result of
the algebraic part of this paper is the following
theorem that will be deduced from the results of  section \ref{sec_formality}
below. Recall Notation \ref{SS}.

\begin{theorem}[Equivariant Formality]\label{formality4}
There exists a triangulated equivalence
$\bfF:D^B_f(\La)\iso
\modB,$
such that $\bfF(\k_\wedge)=\uhb$, and such that
$\bfF\bigl(\k_\wedge(\la)\otimes M\bigr)\cong
\k_{_\B}(l\la)\otimes\bfF(M),$ for any $\la\in\Y\,,\,M\in D^B_f(\La).$
\end{theorem}

\ab 
The proof of the theorem is based, as has been already indicated
in the Introduction, on 
a much more general result saying
 that an infinite order deformation of an associative 
algebra $\aa$ parametrized by a vector space $V$ yields a homomorphism of dg 
algebras ${\mathsf{def}}: \sym(V[-2])\to \RHom_{\aa\bimod}(\aa,\aa)$. 
This result  will be discussed in 
detail in the forthcoming paper [BG]; here we only 
sketch the main idea.

\ab 
Let $A$ be  an infinite order deformation of $\aa$
parametrized by a vector space $V$, that is, a flat $\k[V]$ algebra such
that
$A/V^*\cd A=\aa$.
We replace $\aa$ by a quasi-isomorphic dg-algebra 
$R$ which is a flat $\k[V]$-algebra such that $H^0(R)=A$. Consider the
tensor product $R^e:=R\otimes_{\k[V]} R\opp$. The tensor category of 
dg-modules over $R^e$ acts naturally on the category of left dg-modules 
over $R$ (which is quasi-equivalent to the derived category of 
$\aa$-modules). It turns out that this action encodes the desired 
homomorphism  ${\mathsf{def}}:\sym(V[-2])\to \RHom_{\aa\bimod}(\aa,\aa)$.
Specifically, we first construct a quasi-isomorphism
$\wedge\hdot(V^*[1])=\op{Tor}\hdot_{\k[V]}(\k,\k)\map R^e$,
where $\wedge\hdot(V^*[1])$ is the exterior algebra viewed 
as a dg-algebra with zero differential.
The homomorphism ${\mathsf{def}}$ is then obtained
from the latter  quasi-isomorphism by  Koszul  duality, cf.
Sect \ref{koszul_duality}.

\ab In \S\ref{sec_formality}, we carry out this argument for the particular case where
$\aa=\bb$, $V=\n^*$, and $A=\B$.
Additional efforts are required  to keep track 
of the adjoint action of $\B$ on $\bb$:
the dg-algebra $R$ comes equipped with an
action of quantized Borel algebra, while the  algebra $\sym(\n^*[-2])$
is acted upon by the classical enveloping algebra
$\Ub$. This difficulty is overcome using the Steinberg 
module. We refer to \S\ref{sec_formality} for full details.

\subsection{Digression: Deformation formality.}\label{digression} This
subsection will not be used elsewhere in the paper. Its sole
purpose is to put Theorem
\ref{formality4} in context. 

\ab Recall that the algebra $\bb$ is a normal Hopf
subalgebra in $\B$, and we have $\B/(\bb)=\Ub$.
 This gives, by Lemma \ref{resind},
a canonical $\Ub$-action on the graded algebra
$H\hdot(\bb,\k_\bb) =\Ext_{\bb\mmod}\hdot (\k_\bb,\k_\bb)$.
One of the main results  of \cite{GK} says 
\begin{proposition}\label{gk}
There is a natural  $\Ub$-equivariant
graded algebra isomorphism
$\sym(\n^*[-2])$
$\iso\Ext_{\bb\mmod}^{2\hhdot}(\k_\bb,\k_\bb)\,.$
Moreover, we have $\Ext_{\bb\mmod}^{\tt{odd}}(\k_\bb,\k_\bb)=0.$
\end{proposition}

\begin{remark}\label{compare_act} In [GK], we  used the standard adjoint (commutator)
action $\ad a: x \mapsto ax-xa$, rather than $\adh\B$ action
on $\bb$, in order to get an  $\Ub$-equivariant structure on
the $\Ext$-group $\Ext_{\bb}\hdot(\k_\bb,\k_\bb)$. 
Although these two  actions  are
 {\it different} they turn out
to induce the {\em same} $\Ub$-equivariant structure
on the cohomology. 

\ab To explain this, we recall that
the  construction of \cite{GK} is based on a certain natural
transgression map $\tau: H^1(Z,\k_{_Z})\to
H^2(\bb,\k_\bb):=\Ext^2_{\bb\mmod}(\k_\bb,\k_\bb),$
see \cite[Corollary 5.2]{GK}.
The  isomorphism   of Proposition \ref{gk}
is 
obtained by extending the map $\tau$ by multiplicativity
to an algebra morphism 
\begin{equation}\label{taut}
\sym(\tau):\ \sym\hdot H^1(Z,\k_{_Z})
\map \Ext_{\bb\mmod}^{2\hhdot}(\k_\bb,\k_\bb),
\end{equation}
and then using the canonical vector space isomorphisms
$H^1(Z,\k_{_Z})=Z\eps/(Z\eps)^2\cong\n^*$
(the last one is due to Corollary \ref{tangent}).

\ab Observe that the algebra $\B$ acts on each
side  of \eqref{taut} in two ways: either via the
commutator action, or via the Hopf-adjoint action.
Furthermore, it has been shown in
\cite[Lemma 2.6]{GK} that the map $\tau$, hence
the isomorphism in  \eqref{taut},
commutes with the commutator action.
It is immediate from definition of  Hopf-adjoint action that
the map $\tau$ automatically commutes
with the $\adh\B$-actions as well.
Now, the point is that although
the commutator action of $\B$ on $Z$
 {\it differs} from the $\adh\B$-action on $Z$,
the induced actions on $Z\eps/(Z\eps)^2\cong\n^*$
coincide, as can be easily seen from
explicit formulas for the two actions.
Thus, isomorphism \eqref{taut} implies that
the two actions on $\Ext_{\bb\mmod}^{2\hhdot}(\k_\bb,\k_\bb)$
are equal, as has been  claimed at the beginning of Remark \ref{compare_act}.

\ab In spite of that, it will be essential for us (in the present paper)
below to   use the Hopf-adjoint action rather than the ordinary
commutator.  The difference between
the two actions becomes important since
these actions agree only on the cohomology
level, but may not agree at the level
of  \dg-algebras. $\quad\lozenge$
\end{remark} 
\smallskip

\ab Using  (an appropriate version of) 
the equivalence $\bfF$ of Theorem  \ref{formality4} 
(more precisely, the equivalence $\fF'$ of Theorem \ref{BG_formality2} below),
and using that $\bfF(\k_{_\La})=\k_\bb$,
one  obtains
 $\Ub$-equivariant
 dg-algebra quasi-isomorphisms, cf. \eqref{SL_kos}:
\begin{equation}\label{equiv_roma}
\REnd\hdot_{\bb}(\k_\bb)\simeq
\REnd\hdot_{\bb}\bigl(\fF'(\k_{_\La})\bigr)
\simeq\REnd\hdot_{D^B_f(\La)}(\k_{_\La})=\SS.
\end{equation}
Our construction of a concrete DG-algebra representing 
the object
$\REnd\hdot_{\bb}(\k_\bb)$ on the left, provides it with
additional $\Ub$-action, so that all the
 quasi-isomorphisms in \eqref{equiv_roma} turn out to be compatible with
 $\Ub$-equivariant structures.

\ab
Comparing with Proposition \ref{gk},
this yields the following result  saying
that the \dg-algebra $\REnd\hdot_{\bb}
(\Bbbk_\bb)$ 
 is $\Ub$-equivariantly formal:

\begin{theorem}\label{formality_thm} The
\dg-algebra $\Ub\ltimes\REnd\hdot_{\bb}
(\Bbbk_\bb)$ is quasi-isomorphic
to the algebra\linebreak $\Ub\ltimes\Ext_{\bb}\hdot(\k_\bb,\k_\bb)
= \Ub\ltimes \sym\hdot(\n^*[-2])$, viewed as
a \dg-algebra with trivial differential.  
\end{theorem}

\ab It is not difficult to show that
 the equivariant formality above is in effect
{\em equivalent} to Theorem \ref{formality4}. To see this,
one has to
recall that the very definition of $\Ub$-action
on $\REnd\hdot_{\bb}(\k_\bb)$ appeals to the
quantum algebra $\B$. Hence, the equivariant
 formality statement must involve in one way or 
the other the  algebra $\B$ as well. Trying to make this
precise leads (as it turns out) inevitably
 to Theorem \ref{formality4}.

\bigskip 
\begin{remark} {\small
 A rough  idea of our approach to the proof of  Theorem \ref{formality_thm} is to
 replace the algebra $\bb$ by 
a `larger' \dg-algebra $\ka,$ which is
 quasi-isomorphic to it. We then  construct
an explicit $\Ub$-equivariant algebra homomorphism
$\SS\to \REnd_{_{\ka}}\hdot(\k_{_{\ka}})$
 with  described properties.
 The main difficulty in proving Theorem \ref{formality_thm} is in
the $\Ub$-equivariance  requirement. Insuring equivariance
crucially involves the existence of the Steinberg
$\U$-module, see \eqref{fF}.
Without the equivariance requirement, the result
reduces, in view of Proposition \ref{gk}, to a special case of the
following simple
\vskip -2mm

\begin{proposition}\label{maxim} Let $\aa$ be a Hopf algebra such that
$\Ext\hdot_{\aa}(\k_\aa,\k_\aa)$ is a free commutative algebra generated by
finitely many elements of even degree.  Then the \dg-algebra
$\REnd\hdot_\aa(\k_{\aa})$ is formal.
\end{proposition}
\vskip -2mm

\proof Let $\,h_1,\ldots,h_n\in \Ext^{2\hhdot}_{\aa}(\k_\aa,\k_\aa)\,$ be a
finite set of homogeneous generators of the cohomology algebra, and
let $P$ be a projective resolution of $\k_{\aa}$.  The Hopf algebra
structure on $\aa$ gives rise to a tensor product on the category of
complexes of $\aa$-modules.  In particular, we may form the complex
$P^{\otimes n}$ which is a projective resolution of
$\k_{\aa}^{\,\otimes n} \simeq\k_{\aa}$ (note that tensor
product of any $\aa$-module and a projective $\aa$-module is again a
projective $\aa$-module).  Thus, the \dg-algebra
$\,\Hom\hdot_\aa(P^{\otimes n}, P^{\otimes n})\,$ represents
$\RHom\hdot_\aa(\Bbbk_\aa,\,\Bbbk_\aa)$.  For
each $i=1,\ldots,n,$ choose $\hat{h}_i\in \Hom\hdot_\aa(P, P)$
representing the class $h_i \in H\hdot(\aa,\k_\aa).$ Then, it is
clear that the element 
$$
{\mathbf{\hat{h}}}_i := \underbrace{\id_{_P}
\otimes\ldots\id_{_P}}_{(i-1)\;{\footnotesize{{times}}}} \,\otimes\,
\hat{h}_i\,\otimes \, \underbrace{\id_{_P}
\otimes\ldots\otimes\id_{_P}}_{(n-i)\;{\footnotesize{times}}}
\enspace\in\enspace\Hom\hdot_\aa(P^{\otimes n}, P^{\otimes n})
$$
also represents the class $h_i \in \Ext\hdot_{\aa}(\k_\aa,\k_\aa).$
Furthermore, for $i\neq j$, the morphisms ${\mathbf{\hat{h}}}_i$ and
${\mathbf{\hat{h}}}_j$ act on different tensor factors, hence commute.
Therefore, the subalgebra generated by the
$\,{\mathbf{\hat{h}}}_1,\ldots,{\mathbf{\hat{h}}}_n\,$ is a
commutative subalgebra in the \dg-algebra
$\,\Hom\hdot_{\aa}(P^{\otimes n}, P^{\otimes n})\,$ which is
formed by cocycles and which maps surjectively onto the cohomology
algebra.  The latter being free, the map is necessarily an
isomorphism, and we are done.  \endproof }
\end{remark}

\ab An analogue of Theorem
\ref{formality_thm} holds for algebraic groups over $\bF$, an algebraically
closed field of finite characteristic.  Specifically, let $G_{_\bF}$
be a connected  reductive group over $\bF$, let
$B_{_\bF} \subset G_{_\bF}$ be a Borel subgroup,  let $B^{(1)}$
denote the first Frobenius kernel of $B_{_\bF}$,
and write $\bF_{_{B^{(1)}}}$ for the trivial $B^{(1)}$-module.
One  can see, going through the proof of Theorem \ref{formality4},
that our argument also may be adapted to prove the following result:

\ab
{\it The
\dg-algebra $\REnd\hdot_{_{B^{(1)}}}(\bF_{_{B^{(1)}}})$
is formal as a \dg-algebra in the category of $B_{_\bF}$-modules.}

\subsection{Equivariance and finiteness conditions.}\label{equiv_fin}
To prove Theorem \ref{formality4} we need to introduce 
several  auxiliarly triangulated categiries.

\ab 
Below, we will be considering various dg-algebras
$A=\oplus_{i\leq 0}\,A_i$ (concentrated in non-positive degrees),
that will come equipped with the following additional data:

\pb{A natural grading by the root
lattice $\Y$ (which is preserved by the differential, as opposed to the
$\Z$-grading): $A=\oplus_{i\in\Z}\,A_i,$
where $
A_i=\oplus_{\mu\in\Y}\,A_i(\mu).
$}

\pb{A differential $\Z\times\Y$-graded subalgebra
$C=\oplus_{i\in\Z}\,C_i\sset A,
\,C_i=\oplus_{\mu\in\Y}\,C_i(\mu)$ such that the cohomology
$H\hdot(C)$ is a finitely generated graded Noetherian algebra.}

\pb{An $\Y$-graded  subalgebra $U\sset A_0$ equipped with
a `triangular decomposition'
$U=U^+\otimes U^\circ$,
such that $U^\circ$ has $\Z\times\Y$-degree zero,
and $U^+(\mu)=0$ unless $\mu\in\Y^+$ is
a sum of positive roots. Further, we require  $U$ to be annihilated by the
differential,
i.e.,  that $d(U)\equiv 0$.}
 
\pb{A Hopf algebra structure on $U$.}

\begin{remark}
Observe that   multiplication
map $m: A\otimes A \to A$ is automatically an $\adh U$-equivariant
map. To see this, consider
the iterated coproduct
$\Delta^3: U\to U^{\otimes 4}$. Given $u\in U$,  write
$\Delta^3(u)=
u^{(1)}\otimes u^{(2)}\otimes 
u^{(3)}\otimes 
u^{(4)}.$
Then, for any $a,\tilde{a}\in A$, we have
$m\bigl(\adh u(a\otimes \tilde{a})\bigr)
=u^{(1)}\cdot a\cdot S(u^{(2)})\cdot u^{(3)}\cdot \tilde{a}\cdot S(u^{(4)}).$
But the axioms of Hopf algebra imply that,
writing $\Delta(u)=u'\otimes u''$ (Sweedler notation), 
in $U^{\otimes 3}$ one has
$u^{(1)}\otimes S(u^{(2)})\cdot u^{(3)}\otimes u^{(4)}$
$=u'\otimes 1\otimes u''$.
Our claim follows from this equation.
$\quad\lozenge$
\end{remark}

\ab
Let  $(A,C,U)$ be a data  as above. In what follows, the
algebra  $U$  will be either classical or quantum
enveloping algebra of a Lie subalgebra in our semisimple Lie algebra
$\g$. So,  any weight $\mu\in\Y$ will
give rise to a natural algebra homomorphism $\mu: U^\circ\to\k$.

\ab Let $M=\bigoplus_{i\in\Z,\nu\in\Y}\,
M_i(\nu)$ be a
 $\Z\times\Y$-graded $A$-module,
 equipped with a differential $d$ such that
$d\bigl( M_i(\nu)\bigr)\sset M_{i+1}(\nu)\,,\,\forall
i\in\Z,\nu\in\Y$. The tensor product of the $\adh U$-action
on $A$ and $U$-action
on $M$ obtained by restricting the $A$-action
make $A\otimes M$ an $U$-module.

\ab We say that the module $M$ is
compatible with  $(A,C,U)$-data if the following holds:

\pb{The action map $A\otimes M\to M$ is an $U$-module
morphism compatible with $\Z\times\Y$-gradings;}

\pb{We have $um=\nu(u)\cdot m$ for any $u\in U^\circ, \nu\in\Y,$ and
$m\in
M_i(\nu)$.}

\ab Let $\DGMod_\Y^{U}(A)$  denote the homotopy category of 
differential $\Z\times\Y$-graded $A$-modules
compatible with  $(A,C,U)$-data.

\ab Further,
write $H:=H\hdot(C).$ We have a natural graded algebra map $H\to H\hdot(A)$,
and also an  algebra map $U\to H^0(A)$, since
$d(U)=0$.

\begin{definition}\label{basic_def} Let $D_\Y^U(A,H)$ denote a full
subcategory in
the triangulated category $D\bigl(\DGMod_\Y^{U}(A)\bigr)$
formed by the objects $M\in D\bigl(\DGMod_\Y^{U}(A)\bigr)$
such that

\pb{The cohomology module $H\hdot(M)$ is finitely generated over $H$;}

\pb{The restriction of the $H\hdot(A)$-action on $H\hdot(M)$ to
the subalgebra $U$ is locally finite, i.e., $\dim U m<\infty,\,\forall m\in
H\hdot(M)$.}
\end{definition}

It is clear that $D_\Y^U(A,H)$ is a triangulated category.
The objects of  $D_\Y^U(A,H)$ may be called
 {\em $U$-equivariant, homologically $H$-finite}, dg-modules over $A$.

\begin{remark}\vi In the case $U=\k$
we will drop the superscript `$U$'  from the notation.

\ab \vii Observe that in the notations
$\DGMod_\Y^{U}(A)$ and $D_\Y^U(A,H)$ the superscript 
`$U$' has different meanings:
according to our definition, the objects of $\DGMod_\Y^{U}(A)=\DGMod_\Y^{U^\circ}(A)$
are required to have
a weight decomposition only with respect to the
subalgebra $U^\circ\sset U$, while
in the  $D_\Y^U(A,H)$-case there is an additional
local finiteness condition for the $U^+$-action
on the {\em cohomology} of $M\in  D_\Y^U(A,H)$.

\ab \viii Let $A$ be an ordinary (noetherian) algebra,
viewed as a dg-algebra concentrated in degree zero and equipped
with zero differential, and let $C=A$. Then we have $H=A$,
and the category
$D_\Y(A,H)=D_\Y(A,A)$ is in this case (an $\Y$-graded version of)
the
full subcategory of the bounded derived category of $A$-modules,
formed by the objects $M$ such that 
the cohomology $H\hdot(M)$ is a finitely
generated $A$-module.
$\quad\lozenge$
\end{remark}

\ab 
The adjoint action of the group $B$ on the algebra $\La=\wedge\hdot(\n[1])$
gives rise to an $\Ub$-action on $\La$. Therefore, we may
perform the cross-product construction of Proposition
\ref{Montgomery}.

\begin{notation}\label{SS2} Let $\AA:=\Ub\ltimes\La$ denote  
the  cross-product algebra, 
viewed
as a \dg-algebra with zero differential and with the
grading given by the natural grading  on $\La=\wedge\hdot(\n[1])\sset\AA$
and such that the subalgebra $\Ub\ltimes\{1\}\sset \AA=\Ub\ltimes\La$ 
is placed in grade degree zero. 
\end{notation}

\ab Applying  Definition \ref{basic_def} to the triple
$A:=\AA=\Ub\ltimes\La,\, C:=\La$, and $U:=\Ub$, one obtains 
a triangulated category $D_\Y^\Ub(\AA,\La)$.
 The algebra $\La$ being finite
dimensional, for any $M\in D_\Y^\Ub(\AA,\La)$ one has
$\dim H\hdot(M)<\infty$. 
 In particular, the action on $H\hdot(M)$ of the
Lie algebra $\n$ is nilpotent, hence, can be exponentiated
to an algebraic action of the corresponding unipotent group.
Combined with the $\Y$-grading, this makes $H\hdot(M)$
a finite-dimensional
algebraic $B$-module.

\ab Now consider $\B$, the quantum Borel subalgebra,
 as a differential
$\Z\times\Y$-graded algebra concentrated in $\Z$-degree zero, and equipped
with zero differential. The algebra $\B$ contains $\bb$ as a
subalgebra.
 Associated to the data
$A=U=\B$ and $C=H\hdot(C)=\bb$, we have
the triangulated category
$D_\Y^\B(\B,\bb)$, see Definition \ref{basic_def}. 
Again, we have $\dim\bb<\infty$, hence, for any $M\in D_\Y^\B(\B,\bb)$,
the cohomology  $H\hdot(M)$ acquires
a natural structure of finite-dimensional
$\B$-module.

\subsection{Comparison of derived categories.}
\label{comparison_cat}
In \S\ref{sec_formality} we will prove the following

\begin{theorem}\label{formality2}
There exists a fully faithful triangulated functor
$\fF: D_\Y^{\Ub}(\AA,\La)\map D^\B_\Y(\B,\bb),$
such that $\fF(\k_{_\AA})=\uhb$, and such that
$\fF\bigl(\k_{_\AA}(\la)\otimes M\bigr)\cong
\k_{_\B}(l\la)\otimes\fF(M),$ for any $\la\in\Y\,,\,M\in D_\Y^{\Ub}(\AA,\La).$
\end{theorem} 

\ab 
In order to deduce from this result  the equivalence
of the Equivariant formality Theorem
\ref{formality4}, we  need  to replace
triangulated categories on each side of the  equivalence
in  Theorem
\ref{formality4}
by  larger categories.
Specifically, we should replace the category
$D^B_f(\La)$, see Sect. \ref{dgmodG}, by the category 
 $\D_\Y^\Ub(\AA,\La)$.
The objects of the former category
are $B$-equivariant dg-modules over $\La$
with {\em algebraic} $B$-action,
while the objects of the 
latter are  dg-modules over $\La$ with $\Ub$-action
which is {\em not} required to be locally-finite,
 hence, cannot be be exponentiated to
a $B$-action, in general. 
More precisely, for any $M\in \DGMod_\Y^\Ub(\AA,\La)$,
the action in $M$ of the Cartan subalgebra of $\b$ is diagonalizable
(according to the $\Y$-grading on $M$), while
the action of the subalgebra $\Un\sset \Ub$ may be arbitrary:
only the induced $\Un$-action on the cohomology of $M$
is algebraic.
 Since
any algebraic $B$-module may be clearly viewed as a $\Lie B$-module,
 we see  that every object of 
$\DGMod^B_f(\La)$ may be also viewed as an object of $\DGMod_\Y^\Ub(\AA,\La)$.
Thus, we have a natural  functor
$\bi_\AA: D^B_f(\La)\to D_\Y^\Ub(\AA,\La)$.

\ab Similarly, we would like to replace the category
$\modB$ in  Theorem
\ref{formality4}
by a larger category $D_\Y^\B(\B,\bb)$ introduced in \S\ref{equiv_fin}.
Again, there is a natural triangulated functor
$\bi_\B: \modB\to D_\Y^\B(\B,\bb)$.

\ab In order to compare
 Theorem
\ref{formality4} with  Theorem
\ref{formality2} we are going to prove

\begin{proposition}\label{big_small} \vi The functor
$\bi_\AA: D^B_f(\La)\to D_\Y^\Ub(\AA,\La)$ is an
equivalence of triangulated categories.

\ab 
\vii The  functor $\bi_\B: \modB\to D_\Y^\B(\B,\bb)$
is fully faithful, i.e., makes $\modB$ a {\sf full}
subcategory in $D_\Y^\B(\B,\bb)$.
\end{proposition}

\ab The proof of  the Proposition exploits the following
`abstract nonsense' result  that  will be  also  used at
several other places.
\begin{lemma}\label{abstract_nonsense}
Let $\bi: \scra\to \scra'$ be an exact functors
between two triangulated categories.
Assume given a set $S$ of objects in $\scra$ such that

\pb{The minimal full triangulated subcategory of  $\scra$,
resp. of  $\scra'$,
containg all the objects $M\in S$,
resp.  all the objects $\bi(M),\, M\in S$, is equal to  $\scra$, resp. 
 to  $\scra'$.}

\pb{For any $M_1,M_2\in S$, the functor $\bi$ induces isomorphisms}
$$\Hom_{\scra}(M_1,M_2[k])\iso \Hom_{\scra'}(\bi(M_1),\bi(M_2)[k]),
\quad\forall k\in\Z.$$
\ab
Then $\bi$ is an equivalence of triangulated categories.\qed
\end{lemma}

\begin{proof}[Proof of Proposition \ref{big_small}.]
To prove (i), we first show that the functor $\bi_\AA$ induces
 isomorphisms
\begin{equation}\label{fullf}
\Ext\hdot_{D^B_f(\La)}(\k_{_\La}(\la),\,
\k_{_\La}(\mu))\iso\Ext\hdot_{D^\Ub_\Y(\AA,\La)}((\k_{_\AA}(\la),\,
\k_{_\AA}(\mu))
\quad\forall \la,\mu\in\Y.
\eeq
We argue as follows. Write $N$ for the unipotent 
radical of $B$ and ${\mathcal{U}}_+\sset \Un$ for the augmentation
ideal. For any $k\geq 0$, let $(\Un/{\mathcal{U}}_+^k)^*$ be an
$\Un$-module dual to the finite-dimensional left
$\Un$-module $\Un/{\mathcal{U}}_+^k$.
It is clear that $\underset{k\to\infty}\limind(\Un/{\mathcal{U}}_+^k)^*$
is an {\em injective} $\Un$-module,
moreover, we have an imbedding $\k_{_\Un}\into 
\underset{k\to\infty}\limind(\Un/{\mathcal{U}}_+^k)^*$.

\ab Observe further that the $\Un$-action on  $(\Un/{\mathcal{U}}_+^k)^*$ 
can be exponentiated to give an
algebraic representation of the group $N$.
By a standard argument, one therefore obtains 
a  resolution $\k_{_\Un}\into I_0\to I_1\to\ldots,$
where each $I_k,\,k=0,1,\ldots,$ is an object
of $\limind\rep(N)$ which is injective as
an $\Un$-module. Now, for any $\la\in\Y$, we may treat
each  $I_k$ as a $\Ub\ltimes\La$-module such
that $\Uh$ acts via the character $\la$ and the algebra
 $\La$ acts trivially.
Since $I_k$ is injective as $\Un$-module, it follows  by 
a standard homological algebra that
the morphisms $\RHom\hdot_{D^B_f(\La)}(I_k,I_l)\to
\RHom\hdot_{D^\Ub_\Y(\AA,\La)}(I_k,I_l)$
are isomorphisms, for any $k,l=0,1,\ldots.$
This implies \eqref{fullf}.

\ab We claim next that the objects of the
form $\{\k_{_\AA}(\la)\}_{\la\in\Y}$
generate $D^\Ub_\Y(\AA,\La)$,
resp. the objects of the
form $\{\k_{_\La}(\la)\}_{\la\in\Y}$
generate $D^B_f(\La),$
as a triangulated category.
This is proved by the standard `devissage',
the key point being that  $\dim H\hdot(M)<\infty$
for any $M\in  D^\Ub_\Y(\AA,\La)$.
In more detail, let $D$ be the smallest triangulated subcategory
in  $D^\Ub_\Y(\AA,\La)$ containing the objects 
 $\{\k_{_\AA}(\la)\}_{\la\in\Y}$.
One then shows by descending induction  on $\dim H\hdot(M)<\infty$,
that,  $M\in  D^\Ub_\Y(\AA,\La)\enspace\Rightarrow
\enspace M\in D$. This is done using 
standard truncation functors $\tau^{\leq j}$,
which take $\AA$-modules into $\AA$-modules since
 the algebra $\AA$ is concentrated in
non-negative degrees.
This proves
our claim for the category  $D^\Ub_\Y(\AA,\La)$;
the proof  for the category $D^B_f(\La)$ is identical.

\ab The proof of part (i) of the Proposition is now
completed, in view of \eqref{fullf}, by Lemma \ref{abstract_nonsense}.
The proof of part (ii) is entirely similar.
\end{proof}

\ab To deduce Theorem \ref{formality4} from 
Theorem \ref{formality2} we use  Proposition \ref{big_small} 
and consider the following diagram
$$\xymatrix{
D^B_f(\La)\ar[rr]^<>(.5){\bi_\AA}_<>(.5){\text{Prop. \ref{big_small}}}&&
D_\Y^\Ub(\AA,\La)\ar[rr]^<>(.5){\fF}_<>(.5){\text{Thm. \ref{formality2}}}
&&D_\Y^\B(\B,\bb)&&
\modB.
\ar[ll]_<>(.5){\bi_\B}^<>(.5){\text{Prop. \ref{big_small}}}
}
$$
The functor $\bi_\AA$ in the diagram is an equivalence,
and the functors $\fF, \bi_\B$ are both fully faithful.
Further, Theorem \ref{formality2} insures that $\fF\ccirc\bi_\AA(\k_{_\AA}(\la))
=\bi_\B(\k_{_\B}(l\la)),$ for any $\la\in\Y$.
Let  $D\sset D_\Y^\B(\B,\bb)$ denote
the full triangulated category  generated by
the objects
$\{\fF\ccirc\bi_\AA(\k_{_\AA}(\la))\}_{\la\in\Y}$.
It follows from Proposition \ref{big_small}(ii)
and Lemma  \ref{abstract_nonsense} that 
this  category is the same as the category  generated by
the objects
$\{\bi_\B(\k_{_\B}(l\la))\}_{\la\in\Y}$;
moreover, our functors induce triangulated
 equivalences $\dis 
\xymatrix{
D^B_f(\La)\ar[r]^<>(0.5){\fF\ccirc\bi_\AA}_<>(0.5){\sim}&
D&\modB.\ar[l]_<>(0.5){\bi_\B}^<>(0.5){\sim}}$
Inverting the equivalence $\bi_\B$, we obtain this way
an equivalence $\bfF: D^B_f(\La)\iso \modB$,
which is by definition the equivalence of Theorem
\ref{formality4}.

\ab Summing-up, we have the following equivalences of triangulated
categories
\begin{align}\label{naive22}
\xymatrix{
\dcoh^{G}(\NN)\ar[r]^<>(0.5){i^*}_<>(.5){\eqref{res_gamma}}&
\dcoh^{B}(\n)\ar[r]^<>(0.5){\Gamma}&
D^B_f(\SS)\ar[r]^<>(.5){\eqref{kappa}}&
D^B_f(\La)}\qquad\qquad\quad\hphantom{x}\\
\hphantom{x}\qquad\qquad\quad\xymatrix{
\ar[rr]_<>(.5){\text{Thm. \ref{formality4}}}^<>(.5){\bfF}&&
\modB
\ar[rr]_<>(.5){\text{Thm. \ref{proposition_roma1}}}^<>(.5){\RInd_\B^\U}
&&D^b\cat.}\nonumber
\end{align}

\ab 
Let $\on(\la)$ be a $G\times\Gm$-equivariant
line bundle on $\NN$ obtained by pull-back from the flag manifold
$G/B$ of a standard $G$-equivariant
line bundle corresponding to the caharacter $\la\in \Y$, see Notation \ref{O(l)}.
Applying the functors $i^*$ and $\Gamma$ in the top row
of  \eqref{naive22}, we clearly get
$\Gamma(\n, i^*\on(\la))\cong \uSS(\la)$, where
$\uSS(\la)$
a rank 1  free $\SS$-module,
viewed as a dg-module over $\Ub\ltimes\SS$ with zero differential and  with
the action of
the subalgebra $\Ub\sset \Ub\ltimes\SS$  being
the natural one on $\SS$ twisted by
the character~$\la$. Further, twisting the $\Gm$-equivariant structure on
$\on(\la)$ by the character $z\mapsto z^k$ 
corresponds to degree shift by $k$ in
the  $\SS$-module.
Thus,  we obtain

\begin{theorem}\label{Psi_equiv}
The  composite functor  $Q': \dcoh^G(\NN)\iso D^b\cat$
in \eqref{naive22} provides an equivalence of triangulated categories
such that $Q'\bigl(z^k\otimes\on(\la)\bigr)=\ind_{_{\B}}^{^{\U}}({l}\la)[k]\,,\,
\forall k\in\Z,\la\in\Y.$ \qed
\end{theorem}

\ab The (non-mixed version of the) functor $Q$ in diagram
\eqref{sum-up} is defined to be the inverse of the equivalence
 $Q'$ in the Theorem above.

\section{Proof of Induction Theorem.}\label{sec_induction}
The goal of this section is to prove Theorem \ref{proposition_roma1}.
\subsection{Intertwining functors.} 
  For every simple {\it affine} root $\alpha\in I\cup\{0\}$, 
let  $s_\alpha\in \waf$ denote  the corresponding simple reflection.
We partition the lattice $\X$  into alcoves of `size' $l$
in such a way that the `base vertex' of the fundamental alcove
is placed at the point $(-\rho)$.
Given $\lambda \in \Y$, let $C_\alpha$ denote 
the unique $\alpha$-wall of the alcove containing $\lambda$, and
let $\lambda^{s_\alpha}$ be
the reflection of $\lambda$ with respect to $C_\alpha$.
The assignment: $\lambda \mapsto \lambda^{s_\alpha}$
extends to a $\waf $-action, that we call the {\it right}
$\waf $-action. When restricted to  points $\nu\in\X$ of
the form $\nu=w_a\bullet 0\in \waf\bullet 0,$ this action
becomes the right multiplication
$\nu=w_a\bullet 0 \mto \nu^{s_\alpha}=
(w_a s_\alpha)\bullet 0$.

\ab It is clear from definition that, for any $\la,\mu\in\Y$
and $w\in\waf$, one has $(\mu+l\la)^w=(\mu^w)+l\la$.

\ab Below, we will use the following
general construction of homological algebra, see e.g. [GM].
Let ${\mathscr{A}},{\mathscr{B}}$ be two abelian categories
with enough projectives, and let
$F_1,F_2: {\mathscr{A}}\to{\mathscr{B}}$ be two exact functors.
Assume in addition that we have a morphism of
functors $\varphi: F_1\Rightarrow F_2$. Then
there is a well-defined {\em mapping-cone functor}
${\mathtt{Cone}}(\varphi):\
D^b({\mathscr{A}})\to D^b({\mathscr{B}})$,
which is a triangulated functor between the corresponding
bounded derived categories.

\ab  Recall now
that there are so-called {\em  reflection functors}
$\Xi_\alpha: \cat\too \cat$ defined
 by composing translation functor `to the wall' $\,C_\alpha,$ see e.g. \cite[\S8]{APW},
with translation functor `out of
the wall' $\,C_\alpha\,$. Translation functors being exact
(as direct summands of  functors of the form $V\otimes_\k(-)$,
for a finite dimensional $\U$-module $V$),
it follows that $\Xi_\alpha$ is an exact functor. Furthermore,
there are canonical ``adjointness'' morphisms $\id\to \Xi_\alpha$, and
 $\Xi_\alpha \to \id$. Applying the above mentioned general construction
of the mapping-cone functor  to the morphism $\id\to \Xi_\alpha$, 
resp. $\Xi_\alpha \to \id$, one obtains a triangulated functor
$\theta_\alpha^+$, 
resp. $\theta_\alpha^-$.
The functors $\theta^\pm_\alpha: D^b\cat\to D^b\cat$
defined in this way are usually referred to as {\it intertwining functors}.

\begin{lemma}\label{roma_equiv1}
${\sf (i)}$ In $D^b\cat$ we have canonical isomorphisms:
\[\theta^+_\alpha \circ \theta^-_\alpha\cong \id\cong 
\theta^-_\alpha\circ \theta^+_\alpha,\quad\forall \alpha\in I\cup\{0\},\]
in particular, the functors $\theta^+_\alpha,\theta^-_\alpha:
D^b\cat\to D^b\cat$ are auto-equivalences.

\ab ${\sf (ii)}$ If $\lambda \in \waf\bullet 0$, and
 $s_\alpha$ is  the reflection with respect to
a simple affine  root $\alpha\in I\cup\{0\}$
such that $\lambda^{s_\alpha} \ge \lambda$,
 then
$\;
\theta^+_\alpha(\RInd_{_{\B}}^{^{\U}}\lambda)
\cong \RInd_{_{\B}}^{^{\U}}(\lambda^{s_\alpha}).$
\end{lemma}

{\em Sketch of Proof.}\,
Part (ii) of the Lemma follows directly from \cite[Theorem 8.3(i)]{APW}.

\ab A statement analogous to part (i) of the Lemma is well-known
in the framework of
the category $\mathcal{O}$ for a complex semisimple Lie algebra,
see \cite{Vo}. Specifically, it is clear from the adjunction
properties that
the functor $\theta^+_\alpha \circ \theta^-_\alpha$,
resp. $\theta^-_\alpha\circ \theta^+_\alpha,$
is quasi-isomorphic to a complex  represented
by the following commutative square
$$
\xymatrix{
\Xi_\alpha\ar[r]\ar[d]&\Id\ar[d]\\
\Xi_\alpha\ccirc\Xi_\alpha\ar[r]&\Xi_\alpha.
}
$$
On the other hand,  the left vertical and low horisontal arrows of the
square form a short exact sequence
$0\to\Xi_\alpha\to\Xi_\alpha\ccirc\Xi_\alpha\to\Xi_\alpha\to 0.$
This is proved similarly to
the corresponding statement for the category $\mathcal{O}$,
 with  all the necessary
ingredients (in our quantum group setting)
being provided by \cite[Theorem 8.3]{APW}.
It follows that the square is quasi-isomorphic to
its upper-right corner,
and we are done.
\qed

\subsection{Beginning of the proof of Theorem \ref{proposition_roma1}.}
 The functor $\rind_{_{\B}}^{^{\U}}$ takes
the set $\{\Bbbk_{_{\B}}({l}\la)\}_{\la\in\Y}$ that generates
 category $\modB$ as a triangulated category to the set $\{\rind_{_{\B}}^{^{\U}}
({l}\la)\}_{\la\in\Y}$ that generates category
$D^b\cat$ as a triangulated category, by Corollary \ref{generate}.
Thus, by Lemma \ref{abstract_nonsense},
in order to prove  Induction theorem we must show that:
{\it For all $\lambda , \mu \in \Y $,
and $i\geq 0$, the canonical morphism, induced by  functoriality
of induction, gives an isomorphism  }
\begin{equation}\label{1}
\Ext^i_{_{\B}}(\Bbbk_{_{\B}}({l}\la)\,,\,
\Bbbk_{_{\B}}({l}\mu))
\;\iso\; \Ext^i_{_{\!\cat\!}}(\rind_{_{\B}}^{^{\U}} ({l}\la)\,,\,
\rind_{_{\B}}^{^{\U}}({l}\mu)).
\end{equation}

\ab This isomorphism will be proved in three steps.

\begin{lemma}\label{Step1}
Both sides in (\ref{1}) have  the same dimension  
$$\dim\Ext^i_{_{\B}}(\Bbbk_{_{\B}}({l}\la)\,,\,
\Bbbk_{_{\B}}({l}\mu))=
\dim\Ext^i_{_{\!\cat\!}}(\rind_{_{\B}}^{^{\U}} ({l}\la)\,,\,
\rind_{_{\B}}^{^{\U}}({l}\mu))\;,\quad\forall i\geq 0.
$$
\end{lemma}

\proof
  From Lemma \ref{borel_weil} (quantum version of 
 Borel-Weil theorem [APW]) we obtain
\begin{equation}\label{BW}
\Rt^0\ind_{_{\B}}^{^{\U}}\Bbbk_{_{\B}}=\Bbbk_{_{\U}} 
\quad\mbox{and}\quad \Rt^i\ind_{_{\B}}^{^{\U}}\Bbbk_{_{\B}}=
0\quad\mbox{if}\quad i>0\,.
\end{equation}
Hence, for any  $\B$-module $M$,  we find   
\[
\Hom_{_\B}(\k_{_\B}\,,\,M) \;\;\stackrel{_{\tt{adjunction}}}{\;=\;}\;\;
\Hom_{_{\U}}(\Bbbk_{_{\U}}\,,\, 
\rind_{\,_\B}^{^{\U}}M)\;\;\stackrel{(\ref{BW})}{\;=\;}\;\;
\Hom_{_{_\U}}(\rind_{_{_\B}}^{^{\U}}\Bbbk_{_{_\B}}
\,,\,\rind_{\,_\B}^{^{\U}}M)\,.\]
This yields  isomorphism (\ref{1}) in the special case
$\lambda=0$, and arbitrary $\mu\in\Y$.

\ab The general case will be reduced to the special case above by means
of translation functors.
Specifically, for any $\lambda,\mu\in \Y$ and $\nu \in \Y^{++}$, we
are going to
establish
an isomorphism  
\begin{equation}\label{roma_eq1}\small{
\RHom_{_{\cat}\!}\left(\rind_{_{\B}}^{^{\U}}({l}\la)\,,\,
\rind_{_{\B}}^{^{\U}}({l}\mu) \right)=
\RHom_{_{\cat}\!}\left(\rind_{_{\B}}^{^{\U}}({l}\lambda+{l}\nu)\,,\,
\rind_{_{\B}}^{^{\U}}({l}\mu+{l}\nu)\right)\,.}
\end{equation}
To prove this isomorphism, we
view the root lattice $\Y$ as the subgroup of $\waf$ formed by translations.
Let
$\nu = s_{\alpha_1}\cdots s_{\alpha _r}\in \waf$
be  a reduced expression of $\nu \in \Y^{++} \subset \waf$.
Using the right $\waf$-action    $\tau \mapsto \tau^w$,
 we 
can write 
$l\nu=(0)^{(s_{\alpha_1}\cdots s_{\alpha_{r}})}.$
Since $\nu$ is
dominant, and
$\,\tau^{(yw)}= (\tau^y)^w\,,\, \forall w,y\in \waf,$ 
 we obtain
$(0)^{(s_{\alpha_1}\cdots s_{\alpha_{r}})}=
\bigl((0)^{(s_{\alpha_1}\cdots s_{\alpha _{r-1}})}\bigr)^{s_{\alpha_{r}}} \ge
  (0)^{(s_{\alpha_1}\cdots s_{\alpha _{r-1}})}\ge
\ldots\ge 0.$ Therefore, 
for any $\la\in\Y$ and $j=1,\ldots,r-1$, we deduce
\[(l\la)^{(s_{\alpha_1}\cdots s_{\alpha _{j+1}})}=
l\la+(0)^{(s_{\alpha_1}\cdots s_{\alpha _{j+1}})}\geq
l\la+(0)^{(s_{\alpha_1}\cdots s_{\alpha _j})}=
(l\la)^{(s_{\alpha_1}\cdots s_{\alpha _j})}.\]
Thus, we obtain
$$l\la+l\nu
= \,({l}\la)^{(s_{\alpha_1}\cdots s_{\alpha_{r}})} \ge
  ({l}\la)^{(s_{\alpha_1}\cdots s_{\alpha _{r-1}})}\ge
\ldots\ge l\la,\quad\forall\la\in\Y.
$$
Hence, Lemma \ref{roma_equiv1}(ii) yields  
 $\,\ind_{_{\B}}^{^{\U}}(l\lambda+l\nu)$
$\simeq\dis
\theta_{\alpha _r}^+\circ \theta_{\alpha_{r-1}}^+ \circ \cdots \circ
\theta_{\alpha_1}^+\Bigl(\ind_{_{\B}}^{^{\U}}({l}\la)\Bigr)\,.$
At this point, isomorphism (\ref{roma_eq1}) follows from
 part (i) of Lemma \ref{roma_equiv1}, saying that the functor
$\theta_\alpha$ is an equivalence of categories.

\ab To complete the proof, fix an arbitrary pair $\la,\mu\in\Y$,
and choose $\nu\in\Y^{++}$ sufficiently large, so that
$\nu-\la\in\Y^{++}$. Then, from \eqref{roma_eq1} we deduce
\begin{align}\label{roma_eq1+}
&\enspace\enspace\;\RHom_{_{\cat}\!}\left(\rind_{_{\B}}^{^{\U}}({l}\la)\,,\,
\rind_{_{\B}}^{^{\U}}({l}\mu) \right)\quad\text{shift by}\;\nu-\la\nonumber\\
&\simeq
\RHom_{_{\cat}\!}\left(\rind_{_{\B}}^{^{\U}}({l}\nu)\,,\,
\rind_{_{\B}}^{^{\U}}({l}(\mu+\nu-\la)) \right)\quad\quad\text{shift by}\;(-\nu)
\\
&\simeq
\RHom_{_{\cat}\!}\left(\rind_{_{\B}}^{^{\U}}(0)\,,\,
\rind_{_{\B}}^{^{\U}}({l}(\mu-\la)) \right)\,.\nonumber
\end{align}
We deduce that, for fixed $i$, the corresponding  groups
$\tR^i\Hom_{_{\cat}\!}$ in \eqref{roma_eq1+} all  have the same dimension.
Moreover, by the special case   $\la=0$ of isomorphism
\eqref{1}, that has been already proved,
this dimension equals $\dim\tR^i\Hom_{_{\B}}(\Bbbk_{_{\B}}(0)\,,\,
\Bbbk_{_{\B}}({l}\mu-l\la))$.
Further, using an obvious isomorphism
$\RHom_{_{\B}}(\Bbbk_{_{\B}}(0)\,,\,
\Bbbk_{_{\B}}({l}\mu-l\la))\simeq \RHom_{_{\B}}\bigl(\Bbbk_{_{\B}}(l\la)\,,\,
\Bbbk_{_{\B}}({l}\mu)\bigr)$ we conclude
\begin{equation}\label{roma_*}
\dim \Ext_{\!_{\B}\!}^i\left(\Bbbk_{_{\B}}({l}\la)\,,\,
\Bbbk_{_{\B}}({l}\mu)\right) = \dim
\Ext_{\!_{\cat}\!}^i\left(\rind_{_{\B}}^{^{\U}}({l}\la) \,,\,
\rind_{_{\B}}^{^{\U}}({l}\mu)\right)\,.
\end{equation}
\ab This completes the proof  of the Lemma.\qed

\begin{remark} Observe that formulas
\eqref{BW} and \eqref{roma_eq1+} actually produce, for any $\la,\mu\in \Y,$
a certain  map of the form required  in \eqref{1}.
Unfortunately, we were unable to show that the map
so constructed is indeed induced
via the functor $\rind_{_{\B}}^{^{\U}}$, by functoriality.
Therefore, below we will use an alternative, more round-about,
approach.
\end{remark}

\subsection{A direct limit construction.}
Let $\overline{\b}=\h\oplus \overline{\n}$ denote
the Borel subalgebra opposite to $\b$, so that
the Chevalley generators $\{f_i\}_{i\in I}$ generate
its nilradical $\overline{\n}$. For any
 $\mu \in \X^{++}$, the simple $\g$-module
$V_\mu$ with highest weight $\mu$ is cyclically
generated over ${\mathcal{U}}\overline{\b}$ by its highest weight
vector, i.e., a nonzero vector annihilated by $\n$. 
Specifically, one has a ${\mathcal{U}}\overline{\n}$-module
isomorphism $\dis V_\mu=
{\mathcal{U}}\overline{\n}/\langle f_i^{\langle\mu,
\check{\alpha}_i\rangle+1}\rangle_{i\in I}.$
We see that, for any $\nu,\mu \in \X^{++}$ such that
$\mu-\nu\in \X^{++}$, there is a unique,
up to nonzero factor, map of
${\mathcal{U}}\overline{\b}$-modules
$\,\Bbbk_{_{{\mathcal{U}}\overline{\b}}}(\nu-\mu)
\otimes V_\mu \too V_\nu$, sending the highest weight line
to the highest weight line. 

\ab Dualizing the construction and using the Cartan
involution on $\Ug$ that interchanges
$\Ub$ and ${\mathcal{U}}\overline{\b}$, we deduce that
there is a unique,
up to nonzero factor, map of  $\Ub$-modules $V_\nu \to V_\mu
\otimes \Bbbk_{_{\Ub}}(\nu-\mu)\,,$
sending the highest weight line
to the highest weight line. For any fixed $\lambda\in \Y\sset \X$
and $\nu,\mu\in\Y$,
the induced maps 
$\,\varepsilon_{\nu,\mu}: V_\nu\otimes \Bbbk_{_{\Ub}}(\lambda
-\nu) \too V_\mu
\otimes \Bbbk_{_{\Ub}}(\lambda
-\mu)\,$
form a direct system with respect to the partial order  $\nu\lhd\mu$ on $\Y$,
cf. \eqref{order}.
We let $
\lim\limits_{\stackrel{\too}{\nu\in\Y^{++}}}
\;\left(V_\nu\big|_{\;\Ub}\,\otimes \Bbbk_{_{\Ub}}(\lambda
-\nu)\right)\,$ denote the resulting
direct limit $\Ub$-module.
This $\Ub$-module is clearly co-free over the
subalgebra ${\mathcal{U}}\n\subset \Ub$,
and is co-generated by a single vector of
weight $\lambda$.

\ab Recall the $\Ub$-module ${\mathbf{I}}_\lambda\, =
\,\indf_{_T}^{^B}\lambda$
introduced in\S \ref{tb}.
It is clear that there is a natural
$\Ub$-module isomorphism
\begin{equation}\label{roma_eq(2-)}
\lim\limits_{\stackrel{\too}{\nu\in\Y^{++}}}
\;\left(V_\nu\big|_{\;\Ub}\,\otimes \Bbbk_{_{\Ub}}(\lambda
-\nu)\right)\;
\simeq\; {\mathbf{I}}_\lambda\, = \,\indf_{_T}^{^B}\lambda\,.
\end{equation}
Applying the Frobenius functor to each side of isomorphism \eqref{roma_eq(2-)}
we obtain, for any $\la \in \Y$,  the
following
 isomorphisms
of $\B$-modules:
\beq\label{roma_eq(2)}
\ind_{\tb}^{\B}({l}\la)\;\simeq\;
\fr({\mathbf{I}}_\lambda)
\;\simeq\;
  \lim\limits_{\stackrel{\too}{\nu\in\Y^{++}}}
\;\left(\fr{V_\nu}\Big|_{\;\B}
\;\otimes \Bbbk_{_{\B}}(l\lambda -l\nu )\right),
\eeq
where 
the subalgebra $\tb\sset\B$ was defined in section \ref{tb}.

\begin{lemma}\label{Step2}
  For any $\lambda,\mu\in \Y$,
 the following canonical morphism, induced by  functoriality
of induction, is injective:
\begin{equation}\label{roma_eq(1)}
\Ext^i_{_{\B}} \left(\ind_{\tb}^{^{\B}}({l}\la)\,,\,
\ind_{\tb}^{^{\B}}({l}\mu)\right)
\;\hookrightarrow
\; \Ext^i_{\!_{\cat}\!} \left(\rind_{\tb}^{^{\U}}({l}\la)\,,\,
\rind_{\tb}^{^{\U}} ({l}\mu)\right)
\end{equation}
\end{lemma}

\proof
Recall that by formula (\ref{roma_eq(2)})
we have   
$\ind_{\tb}^{\B}({l}\lambda)\simeq \fr{\mathbf{I}}_\lambda\,,$
and therefore  
$\,\rind_{\tb}^{^{\U}}({l}\la)$
$=
\rind_{_\B}^{^{\U}}\left(\ind_{\tb}^{\B}({l}\la)\right)\,=\,
\rind_{_\B}^{^{\U}}(\fr{\mathbf{I}}_\lambda)\,.$
To prove injectivity part of Lemma \ref{Step2},
we rewrite the morphism in
(\ref{roma_eq(1)}) using  Frobenius reciprocity  as follows:
$$
\begin{array}{ll}\dis
\!\Ext^i_{\tb}\bigl(\ind_{\tb}^{\B}({l}\la)\,,\!
\k_\tb({l}\mu)\bigr)
&=\;\;\Ext^i_{\tb}(\fr{\mathbf{I}}_\lambda\,,\,
\k_\tb({l}\mu))\\
&
\stackrel{\tilde{\rho}}{\to}\,
\Ext^i_{\tb}(\rind_{_\B}^{^{\U}}(\fr{\mathbf{I}}_\lambda)\,,\,
\k_\tb({l}\mu))\;=\;
\Ext^i_{\tb}\left(\ind_{\tb}^{^{\U}}({l}\la)\,,\,
\k_\tb({l}\mu)\right).
\end{array}
$$
Here the morphism $\tilde{\rho}$ is induced by the canonical
$\B$-module adjunction morphism
 $\,\rho: \rind_{_\B}^{^{\U}}(\fr{\mathbf{I}}_\lambda)
\too\fr{\mathbf{I}}_\lambda,$
restricted to $\tb$.
Injectivity of (\ref{roma_eq(1)}) would follow,
provided we show that   

{\bf{(i)}} the object
$\rind_{_\B}^{^{\U}}(\fr{\mathbf{I}}_\lambda)\in D^b\cat$
is concentrated in degree 0, i.e., is an actual 
$\U$-module and, moreover;

{\bf{(ii)}} the morphism $\rho$ is a surjection, which is split
as a morphism of $\tb$-modules.

\ab To prove (i), we apply  the functor $\rind_{_{\B}}^{^{\U}}$
to  isomorphism (\ref{roma_eq(2)})
and obtain 
\begin{align}\label{roma3}
\rind_{_{\B}}^{^{\U}}(\fr{\mathbf{I}}_\lambda) \; \simeq\;
\rind_{_{\B}}^{^{\U}}\Bigl(\lim\limits_{\stackrel{\too}{\nu\in\Y^{++}}}
&\;\bigl(\fr{V_\nu}\Big|_{\;\B}
\;\otimes \Bbbk_{_{\B}}({l}\lambda -{l}\nu )\bigr)\Bigr)\\
&=\;
\fr{V_\nu}\,\bigotimes\,\Bigl(\lim\limits_{\stackrel{\too}{\nu\in\Y^{++}}}
\;\rind_{_{\B}}^{^{\U}}({l}\lambda -{l}\nu )\Bigr)
\nonumber
\end{align}

\ab By the quantum group analogue of the Kempf vanishing,
see \cite{APW}, for $\lambda\in -\Y^{++}$, the 
object $\rind_{_{\B}}^{^{\U}}({l}\la) \in D^b\cat$
is isomorphic to $\tR^0\ind_{_{\B}}^{^{\U}}({l}\la),$
an actual $\U$-module,
and all other cohomology groups vanish, i.e.,
$R^i\ind_{_{\B}}^{^{\U}}({l}\la)= 0$,
for all $i\neq 0$.
 Hence, formula
(\ref{roma3}) shows that $\rind_{_{\B}}^{^{\U}}\fr{\mathbf{I}}_\lambda$
is an actual $\U$-module, because $\lambda -\nu \in -\Y^{++}$,
for $\nu$ large enough.

\ab To prove property (ii), we use the isomorphism
$\rind_{_{\B}}^{^{\U}}({l}\la)=\tR^0\ind_{_{\B}}^{^{\U}}({l}\la),$
for $\lambda\in -\Y^{++}$. Then, by Frobenius reciprocity one has a canonical
$\B$-module 
projection $\gamma: \rind_{_{\B}}^{^{\U}}({l}\la)=
\ind_{_{\B}}^{^{\U}}({l}\la)\onto 
\Bbbk_{_{\B}}({l}\la)$.
Furthermore, 
there is also an $\u$-module morphism  
$\Bbbk_{_{\u}}({l}\la) \to
\ind_{_{\B}}^{^{\U}}({l}\la)$.
The latter morphism provides a $\bb\, (=\B\cap\u)$-equivariant section of
the projection $\gamma$ 
that, moreover, respects the $\Y$-gradings. Hence the projection  
$\ind_{_{\B}}^{^{\U}}({l}\la)\onto 
\Bbbk_{_{\B}}({l}\la)$
  is split as a morphism
of $\tb$-modules. Using formula 
(\ref{roma3}) we deduce from this, by taking direct limits
as in the previous paragraph,
that the projection $\rho$ is also split
as a morphism of $\tb$-modules.
Therefore, formula (\ref{roma_eq(2)}) implies that
the restriction of $\rho$ to $\tb$ is the
projection to a direct 
summand.
It follows that the map $\tilde{\rho}$ in
(\ref{roma_eq(1)}) is injective.\qed

\begin{lemma}\label{Step3}
The objects of the form  
$\;\ind_{\tb}^{\B}({l}\la)\,,\, \lambda\in \Y\;,$
generate $\modB$
as a triangulated category; moreover, the morphism
(\ref{roma_eq(1)}) is an isomorphism.
\end{lemma}

\proof
The first part of the Lemma is clear,
since the algebra $\Un$ has finite homological dimension
(hence, any $\B$-module in  $\modB$ has a finite resolution
by objects of the form $\ind_{\tb}^{\B}({l}\la)$). Therefore, it remains
to prove that (\ref{roma_eq(1)}) is an isomorphism.
Observe that, for any given $i$, both sides in (\ref{roma_eq(1)})
are finite-dimensional vector spaces. This is so because 
the $Ext$-groups involved are finitely generated graded modules over
the corresponding $Ext$-algebra $\Ext\hdot(\Bbbk,\Bbbk)$,
and the latter is known to be a finitely-generated graded algebra.
Hence, by Lemma \ref{Step2}, we must only show that, for each $i\geq 0$,
 both sides in (\ref{roma_eq(1)})
are of the same dimension.

\ab To prove this we observe that, for any finite dimensional $G$-module
$V$ (viewed as an $\barU$-module), translation functors on $\cat$ commute with 
the functor   $M \mapsto M\otimes \fr{V}$. Hence,
tensoring by $\fr{V}$,  from Lemma \ref{Step1} we deduce  
\begin{equation}\label{roma_**}
\dim \Ext_{\!_{\B}\!}^i\left(\Bbbk_{_{\B}}({l}\la)\,,\,
\Bbbk_{_{\B}}({l}\mu)\otimes \fr{V}\big|_{\;\B}\right) = 
\dim \Ext_{\!_{\cat}\!}^i\left(\rind_{_{\B}}^{^{\U}} ({l}\la)\,,\,
\rind_{_{\B}}^{^{\U}}({l}\mu)\otimes
\fr{V}\right)\,.
\end{equation}
We put $V=V_\nu$, a simple module with highest weight $\nu$.
The equality of dimensions in
(\ref{roma_eq(1)}) follows from equation
(\ref{roma_**}) by taking direct limit as $\nu \to +\infty$
in $\Y^{++}$,
and using formula 
(\ref{roma3}), in the same way  as above. 
 Lemma \ref{Step3} is proved.
\qed

\ab This completes the
proof of Theorem
\ref{proposition_roma1}.

\section{Proof of Quantum group  formality theorem}\label{sec_formality}
\subsection{Constructing an equivariant dg-resolution.} 
\label{resol}
In order
to begin the proof of Theorem \ref{formality2} we
recall the central subalgebra
$Z\sset \fB$, see Definition \ref{cent}.
\begin{lemma}\label{tK}
 There exists a (super)commutative dg-algebra $\tK=\bigoplus_{i\leq
0}\,\tK^i,$ equipped
with an $\Ub$-action, and such that

\pb{The $\Ub$-action on $\tK$ preserves the grading, moreover,
  for each $i$,  there is a direct sum decomposition
$\tK^i=\bigoplus_{\nu\in\Y}\,\tK^i(\nu)$ 
such that $ur=\nu(u)\cdot r\,,\,\forall
u\in {\mathcal{U}}\h\sset\Ub\,,\,r\in \tK^i(\nu).$}

\pb{$\tK^0=Z$, and the graded algebra $\tK$ is a free $\tK^0$-module;}

\pb{$H^0(\tK)=\k$, and $H^i(\tK)=0,$ for all $i\neq 0$.}
\end{lemma}

\proof The argument is quite standard. 
We will construct inductively a sequence of $\Y$-graded
 $\Ub\ltimes Z$-modules
$R^i\,,\,
i=0,-1,-2,\ldots,$ starting with $R^0:=Z$,
and such that each  $R^i$ is free over $Z$.
At every step, we put a differential (of degree $+1$)
on the graded algebra
$\sym\bigl(\bigoplus_{-n\leq i \leq 0}\,R^i\bigr)$,
referred to as an $n$-{\it truncated \dg-algebra}.
We then set $\tK:= \sym\bigl(\bigoplus_{i\leq
0}\,R^i\bigr).$

\ab To do the induction step, assume we have already constructed all
the modules $R^i\,,\,
i=0,-1,\ldots,-n,$  and  differentials $d$ in such a
way that, for the $n$-truncated \dg-algebra we have 
$$
H^j\bigl(\sym(\oplus_{-n\leq i \leq 0}\,R^i)\bigr)
=\begin{cases}
\k &\text{if } j=0\\
0  &\text{if } -n+1\leq j< 0\,.
\end{cases}
$$
\ab Inside the
$n$-truncated algebra, we have  the following $\Ub\ltimes Z$-submodule,
$C^{-n},$ formed by degree $(-n)$-cycles:
$$C^{-n}:=
\Ker\bigl(\sym(\oplus_{-n\leq i \leq 0}\,\tK^i)\stackrel{d}{\too}
\sym(\oplus_{-n\leq i \leq 0}\,R^i)\bigr)\,,
$$
(if $n=-1$ we set $C^{-1}:=Z\eps$).
We can find  an  $\Ub\ltimes Z$-module  surjection $R^{-n-1}\onto C^{-n}$,
such that $R^{-n-1}$ is free as a $Z$-module.
We let $R^{-n-1}$ be the space of degree $(-n-1)$-generators
of our $(n+1)$-truncated algebra, and define the differential on
these new generators to be the map  $R^{-n-1}\onto C^{-n}$.
This completes the induction step. The Lemma is proved.\qed

\ab The \dg-algebra $\tK$ has a natural  augmentation
given by the composite map
\beq\label{tK_aug}
\epsilon_{_\tK}:\;\tK=\oplus_{i\leq 0}\tK^i \too
\tK\big/\left(\oplus_{i<0}\tK^i\right)=
Z\stackrel{\epsilon}\too \k\,.
\eeq
We observe that the last property stated in Lemma \ref{tK} implies that
the map $\epsilon_{_\tK}$ is a quasi-isomorphism.

\ab Next, we form the tensor product  $\tK\otimes_{_Z}\tK$.
This is again a $Z$-free (super-)commutative \dg-algebra concentrated
in non-positive degrees.
We shall see (Lemma \ref{tKtK_lemma}(i))
that  $H\hdot(\tK\otimes_{_Z}\tK)
\cong \La$, where  $\La=\wedge\hdot(\n[1])$ is the exterior
algebra generated by the vector space $\n$ placed in degree $(-1)$.

\ab The \dg-algebra  $\tK\otimes_{_Z}\tK$ acquires 
a natural $\Ub$-action, the tensor product of the $\Ub$-actions on 
both factors.  We form the cross-product dg-algebra
$\Ub\ltimes(\tK\otimes_{_Z}\tK)$, where
the subalgebra $\Ub$ is placed in grade degree zero.
 Lemma \ref{tKtK_lemma} below implies  that
we have
\begin{equation}\label{AA}
H\hdot(\Ub\ltimes(\tK\otimes_{_Z}\tK))=\Ub\ltimes
H\hdot(\tK\otimes_{_Z}\tK)=\Ub\ltimes\La=\AA,
\eeq
cf., Notation \ref{SS2}.
Thus, we  may consider triangulated categories
$  D_\Y^\Ub\bigl(\Ub\ltimes(\tK\otimes_{_Z}\tK), \La\bigr)$
and $D\yh(\AA,\La),$ where in the
latter case $\AA$ and $\La$ are treated  as a \dg-algebras 
with trivial differential.

\ab The $\tK$-bimodule  $\uK$ may be viewed naturally as an
object of  $  D_\Y^\Ub\bigl(\Ub\ltimes(\tK\otimes_{_Z}\tK), \La\bigr)$.

\begin{lemma}\label{tKtK_lemma} \vi $H\hdot(\tK\otimes_{_Z}\tK)
\cong \La\,$ (isomorphism of algebras and  $\Ub$-modules);

\ab \vii The \dg-algebra  $\tK\otimes_{_Z}\tK$ 
is  $\Ub$-equivariantly formal, i.e., there is
a \dg-algebra quasi-isomorphism $\imath: \Ub\ltimes \La
\qisto\Ub\ltimes (\tK\otimes_{_Z}\tK).$

\ab \viii The
induced equivalence $\imath_*:\ D\yh(\AA,\La)
\iso  D\yh\bigl(\Ub\ltimes (\tK\otimes_{_Z}\tK),\La\bigr)$
sends $\k_{_\AA}$ to $\imath_*(\k_{_\AA})=\uK$.
\end{lemma}

\proof By construction, the augmentation 
$\epsilon_{_\tK}:\tK\onto\k$ in \eqref{tK_aug}
gives a free $Z$-algebra resolution of the
trivial $Z$-module $\k_{_Z}$. Thus, by definition
of derived functors, the \dg-algebra  $\tK\otimes_{_Z}\tK$
represents the object $\k_{_Z}\lotimes_{_Z}\k_{_Z}$
in the derived category of \dg-algebras. Therefore,
the cohomology algebra  $H\hdot(\tK\otimes_{_Z}\tK)$ is
isomorphic to the Tor-algebra ${\op{Tor}}^Z\idot(\k_{_Z},\k_{_Z})$.
By  Proposition \ref{GBbis},
we have $Z\simeq \k[\overline{B}\cd B/B].$ 
Hence, by Corollary \ref{tangent}(ii), we obtain
  $\Ub$-equivariant graded algebra isomorphisms
$
H\hdot(\tK\otimes_{_Z}\tK)\cong\Tor^Z(\k_{_Z},\k_{_Z})
\cong 
\La.
$

\ab Thus, to prove part (ii) of the Lemma we must construct an
$\Ub$-equivariant \dg-algebra
quasi-isomorphism $\tK\otimes_{_Z}\tK\qisto\La$.
We first construct such a map  that will only be a
 morphism 
of complexes of $\Ub$-modules (with the algebra structures forgotten).

\ab To this end, we
 use the standard (reduced) bar-resolution
$\left(\cdots \to Z\otimes Z\eps\otimes Z\eps\to\right.$
$\left. Z\otimes Z\eps\to
Z\right)$
$\qisto \k_{_Z},$ 
and replace the trivial $Z$-module $\k_{_Z}$ by a quasi-isomorphic
$\Ub$-equivariant complex of free $Z$-modules. Applying the
functor $\k_{_Z}\otimes_{_Z}(-)\,$ to this resolution
term by term, we represent
the object $\k_{_Z}\lotimes_{_Z}\k_{_Z}$
by the following complex
$$
\op{Bar}\hdot(Z\eps):\;\; \ldots \too Z\eps\otimes 
Z\eps\otimes Z\eps\too Z\eps\otimes Z\eps\too  Z\eps\to \k_{_Z}\,.
$$
Now, given $a\in Z\eps$, let $\bar{a}\in Z\eps/Z^2\eps$
denote its image. It is  well-known (see e.g. \cite{Lo}) that
the assignment
\beq\label{tor3}
a_1\otimes \cdots \otimes a_n \mto \bar{a}_1\wedge \bar{a}_2\wedge\ldots\wedge
\bar{a}_n,\quad
\op{Bar}\hdot(Z\eps)\too \wedge\hdot(Z\eps/Z\eps^2)
\eeq
yields
 an isomorphism of cohomology.
The map \eqref{tor3} is clearly $\Ub$-equivariant, hence,
we obtain a chain of  $\Ub$-equivariant quasi-isomorphisms
\beq\label{tor4}
\tK\otimes_{_Z}\tK\,\cong\,\k_{_Z}\lotimes_{_Z}\k_{_Z}\,\cong\,
\op{Bar}\hdot(Z\eps)\,
\underset{\eqref{tor3}}\qisto\,
\wedge\hdot(Z\eps/Z\eps^2)\,=\,\La\,.
\eeq

\ab We can finally construct a  \dg-{\it algebra}
quasi-isomorphism required in \eqref{tor2} as follows.
Equip the vector space $\n[1]$
 with the trivial
$Z$-action (via the augmentation $Z\to\k$) and with the natural
adjoint $\Ub$-action.
Let $P\hdot$ be an  $\Ub\ltimes Z$-module resolution
of  $\n[1]$ such that each term $P^i$ 
 is free as a $Z$-module.

\ab Restricting the quasi-isomorphism
$\La\qisto\tK\otimes_{_Z}\tK$
in \eqref{tor4} to the subspace $\n\sset\La,$
we get a morphism $\n[1] \to \tK\otimes_{_Z}\tK$
in the triangulated category of $\Ub$-modules.
We can represent this morphism as an actual
 $\Ub$-module map of complexes
$f: P\hdot\too\tK\otimes_{_Z}\tK$,
where $P\hdot$ is the $Z$-free resolution of the
vector space $\n[1]$ (with trivial
$Z$-action) chosen in the previous paragraph.
The morphism $f$
can be uniquely extended,  by multiplicativity,
to a  \dg-algebra morphism
$f_{_{\sf{alg}}}: \La \to\tK\otimes_{_Z}\tK$.
The latter map gives,
by construction, an isomorphism on cohomology:
$\n=\wedge^1\n \iso H^{-1}(\tK\otimes_{_Z}\tK)$.
Therefore, since the cohomology algebra $H\hdot(\tK\otimes_{_Z}\tK)\cong\La$
is freely generated by its first component $\wedge^1\n$, we deduce that
the \dg-algebra morphism
$f_{_{\sf{alg}}}$ induces a graded algebra  isomorphism $
\La \iso
H\hdot(\tK\otimes_{_Z}\tK)$.

\ab Thus, performing the
cross-product construction yields a
graded algebra isomorphism
\beq\label{tor2}
\id_{\Ub}\ltimes f_{_{\sf{alg}}}:\ 
\Ub\ltimes (\tK\otimes_{_Z}\tK)\,\stackrel{\tt{qis}}\simeq\,
\Ub\ltimes \La \,,
\eeq
that induces the isomorphism of cohomology constructed at the beginning of
the proof.\qed

\begin{remark}  The \dg-algebra $\tK\otimes_{Z}\tK$ is likely to be
quasi-isomorphic to the bar complex $\op{Bar}(Z\eps)$, equipped
with the shuffle product algebra
structure, cf. \cite{Lo}.
 $\quad\lozenge$\end{remark}

\subsection{DG-resolution of $\bb$.}\label{dgres_sec} 
  Recall  that the quantum Borel algebra $\fB$ is free over its central
subalgebra
$Z$. We put 
$\tK\bb:= \tK\otimes_{_Z}\fB$. Thus, $\tK\bb=
\bigoplus_{i \leq 0}\,\tK^i\otimes_{_Z}\fB$ is a \dg-algebra
concentrated in non-positive degrees and 
such that its degree zero component is isomorphic to $\fB$.
Further,  the augmentation \eqref{tK_aug} induces
 an algebra map
$$
\tK\bb\,= \,\tK\otimes_{_Z}\fB\,\onto\, \k_{_\tK}\otimes_{_Z}\fB\,=\,
\fB/(Z)\,=\,\bb\,.
$$

\ab We may view the algebra $\bb$ on the right as a \dg-algebra with trivial
differential, concentrated in degree zero. Then, Lemma \ref{tK}
implies that the map above
gives a  quasi-isomorphism $\pi:\tK\bb\qisto \bb$. 

\ab Recall next the  $\adh\sB$-action on $\fB$ (see Proposition \ref{GBbis}),
and view $\tK\otimes_{_Z}\fB$ as a tensor product of $\sB$-modules,
where the $\sB$-action on the first factor is
obtained from 
$\Ub$-action via the Frobenius functor.
Performing the cross-product construction
we obtain a \dg-algebra  quasi-isomorphism
$\pi: \sB\ltimes\tK\bb\qisto 
\sB\ltimes\bb.$

\ab Since $\bb=H\hdot(\ka)\sset
H\hdot\bigl(\sB\ltimes\ka\bigr),$
we may consider the category
$D_\Y^\B\bigl(\sB\ltimes\ka,\bb\bigr)$.
The quasi-isomorphisms constructed above induce the following 
category equivalences
\beq\label{tK_qis}
\pi_*:\ D_\Y(\tK\bb,\bb) \iso D_\Y(\bb,\bb)\;,\quad\text{and}\quad
\pi_*:\ D_\Y^\B\bigl(\sB\ltimes\tK\bb,\bb\bigr)\iso
D_\Y^\B\bigl(\sB\ltimes\bb,\bb\bigr)\,.\eeq

\subsection{Construction of a bi-functor.}\label{bi_functor}
 A key ingredient used to construct the equivalence of Theorem
\ref{formality2} is the following {\it bifunctor}:
\beq\label{bifunctor}
  D_\Y^{\Ub}\bigl(\Ub\ltimes(\tK
\otimes_{_Z}\tK),\La\bigr) \times
  D_\Y^\B(\sB\ltimes\ka,\bb)
\longrightarrow D_\Y^\B(\sB\ltimes\ka,\bb),\quad
M,N \mto M \lotimes_{\tK}N\,. 
\eeq
In this formula, 
the tensor product $M \lotimes_{\tK}N$ is taken with
respect
to the action on $M$ of the {\it second} factor in the algebra $\tK
\otimes_{_Z}\tK$ and with respect to the $\tK$-module structure
on $N$ obtained by restriction to the subalgebra
$\tK\sset\tK\otimes_{_Z}\fB=\tK\bb$.
The object $M \lotimes_{\tK}N$ thus obtained has 
an additional $\tK$-action coming from the
action of the {\it first} factor $\tK\sset \tK
\otimes_{_Z}\tK$ on $M$. It is instructive to 
view the action (on $M$) of the first factor in $\tK
\otimes_{_Z}\tK$ as a {\it left} action,
and of the second factor, as  a {\it right} action,
and thus regard $M$ as an $\tK$-bimodule. 
Then, the left
 $\tK$-action on  $M \lotimes_{\tK}N$
commutes with the  $\fB$-action on  $M \lotimes_{\tK}N$
induced (via the imbedding $\fB\into\tK\bb$)
from the one  on $N$.
Furthermore, the  left $\tK$- and the $\fB$-actions 
agree  on the subalgebra $Z\sset (\tK\otimes1)\cap (1\otimes\fB)$
(where the intersection is taken inside $\ka=\tK\otimes_{_Z}\fB$),
and combine together to make $M \lotimes_{\tK}N$ into a 
$\tK\bb$-module, that is, we put:
$$(r\otimes b)\cdot(m\otimes n)=
\bigl((r\otimes 1)\cdot m\bigr)\otimes b\cdot n\,,\,
\forall m\in M,\,n\in N,\,r\otimes b\in \tK\otimes_{_Z}\fB=\ka.$$

\ab There is also a $\sB$-action on
 $M \lotimes_{\tK}N$ defined as the tensor product of 
$\sB$-action on $M$ induced via the Frobenius functor
from the $\Ub$-action, and the given $\sB$-action on $N$.
It is straightforward to verify that the above actions provide 
 $M \lotimes_{\tK}N$ with a well-defined $\sB\ltimes\ka$-module structure.

\begin{remark} We note that
although the tensor product in $M \lotimes_{\tK}N$ is taken
over the commutative algebra $\tK$, the resulting  $\tK$-action on
 $M \lotimes_{\tK}N$  is {\it not} compatible with other
structures described above, due to the fact that the subalgebra $\tK\sset\sB\ltimes\tK\bb$
is not {\it central}. Thus, it is imperative to use the
`additional' $\tK$-action on
 $M \lotimes_{\tK}N$ (arising from the {\it first} tensor factor
in $\tK
\otimes_{_Z}\tK$) in order to get an
$\sB\ltimes\tK\bb$-module structure on $M \lotimes_{\tK}N$.

\ab The bifunctor \eqref{bifunctor} should be thought of as
`changing' the 
$\ka$-module structure on $N$ via
the $\tK$-bimodule $M$.
 $\quad\lozenge$\end{remark}

\subsection{Main result.}
We now change our point of view and consider the algebra
$\bb$ as a {\it sub}algebra in $\sB$, rather than  a
 {\it quotient} of $\fB$. The imbedding $\bb\into\sB$ gives,
via multiplication in $\sB$, a morphism of $\adh\sB$-modules
$\mult: \sB\ltimes\bb\map\sB$, $x\otimes y\mapsto
xy$. 
By Proposition \ref{Montgomery}, the map $\mult$ is in effect an algebra
map. The induced direct image functor 
$\mult_*: D_\Y^{\sB}(\sB\ltimes\bb,\bb) \too D_\Y^\B(\sB,\bb)$
is given by the (derived) tensor product functor
 $M\mto \bigl((\sB\ltimes\bb)/\Ker \mult\bigr)\lotimes_{_{\sB\ltimes\bb}}M$.
Thus, we can define the following composite \dg-algebra map,
and the corresponding  direct image functor
\begin{align}\label{alpha}
 \alpha:\ \sB\ltimes\ka
\underset{\pi}\qisto 
&
\sB\ltimes\bb\stackrel{\mult}\longrightarrow\sB,
\\ &
\alpha_*:\   D_\Y^\B(\sB\ltimes\ka,\bb)
\overset{\eqref{tK_qis}}\iso   
D_\Y^{\sB}(\sB\ltimes\bb,\bb) \stackrel{\mult_*}\longrightarrow  D_\Y^\B(\sB,\bb).  
\nonumber
\end{align}

\ab Next, we consider the  projection
to the first factor $\beta:\
\sB\ltimes\ka
\longrightarrow\sB\,,\,
b\ltimes r\mto b\cdot\epsilon_{_\ka}(r),$ 
where the augmentation $\epsilon_{_\ka}$
is given by the tensor product $\epsilon_{_\ka}:=
\epsilon_{_\tK}\otimes\epsilon_{_\fB}:\
 \ka=\tK\otimes_{_Z}\fB\map \k_{_\tK}\otimes\k_{_\fB}=\k$.
By Proposition \ref{Montgomery}, 
the map $\beta$ is an algebra morphism that gives rise to a pull-back functor
$ \beta^*$. Thus, we obtain the following diagram of algebras and
functors:
\begin{equation}\label{curve}
\xymatrix{
&\sB\ltimes\ka
\ar[dl]^{\beta}\ar[dr]_{\al}&
&&  D_\Y^\B(\sB\ltimes\ka,\bb)
\ar[dr]_{\al_*}&\\
\sB&&\sB& 
  D_\Y^\B(\sB,\bb)\ar[ur]_{\beta^*}
&&  D_\Y^\B(\sB,\bb)
}
\end{equation}

\begin{definition}\label{St} Let $\St$ be  the {\sl Steinberg
$\U$-module}, twisted by an appropriate
1-dimensional character  $\sB\to\k$
in such a way that its highest weight becomes equal to zero.
Specifically,
in the notation of sect. \ref{W_aff_gr} below, we define
$\St$ as the following $\B$-module
$\St:=L^{(l-1)\rho}\otimes\k_{_\B}((1-l)\rho)$. 
\end{definition}

\ab It is known, cf. \cite{AP}, that
the module $\St$ becomes, when restricted to the subalgebra $\u^+\sset\B$,
a rank one free $\u^+$-module.  We regard $\St$ as an object
of $D_\Y^\B(\sB,\bb)$.
A key property of this object  exploited below, cf. \eqref{St_property},
is that, in  $D_\Y^\B(\sB,\bb)$, one has a (quasi)-isomorphism
$$\alpha_*\beta^*\St\simeq\k.$$

\ab
We are now in a position to combine the above constructions 
in order to 
define a functor $\fF:\
 D_\Y^{\Ub}\bigl(\Ub\ltimes (\tK\otimes_{_Z}\tK),\La\bigr)
\longrightarrow 
D_\Y^\B(\B,\bb)$. To this end, recall the notation
$\AA:=\Ub\ltimes\La$, see \eqref{SS2}.
We  use the bifunctor \eqref{bifunctor},
and the two functors in diagram \eqref{curve},
to introduce the following composite functor:
\begin{align}\label{fF}
\fF:\;D_\Y^{\Ub}\bigl(\AA,\La\bigr)\overset{\imath_*}{\iso}
  &D_\Y^{\Ub}\bigl(\Ub\ltimes(\tK\otimes_{_{Z}}\tK),\La\bigr)
{\xymatrix{\ar[rr]^<>(.5){(-)\lotimes_{_\tK} \beta^*\St}
_<>(.5){\eqref{bifunctor}}&&
 D_\Y^\B(\sB\ltimes\ka,\bb)}}\stackrel{\alpha_*}\too D_\Y^\B(\B,\bb),\nonumber\\
&
M\, \mto\, \fF(M):=\al_*\bigl(M\lo_{_{\tK}}\beta^*\St\bigr)
=
\mult_*\pi_*\bigl(\imath_*M\lo_{_{\tK}}\,\beta^*\St\bigr)\,.
\end{align}

\ab Here is a more precise version of Theorem \ref{formality2}.
 \begin{theorem}\label{formality3} We have $\fF(\k_{_{\AA}})=\k_{_\B}$,
and $\fF\bigl(\k_{_{\AA}}(\la)\otimes M\bigr)=\k_{_\B}(l\la)\otimes \fF(M),$
for any $\la\in\Y$ and $M\in D_\Y^{\Ub}(\AA,\La).$
Moreover,
the functor  $\fF$ induces, for any $\la,\mu\in\Y,$ natural isomorphisms
\begin{align*}
\Ext\hdot_{_{D_\Y^{\Ub}(\AA,\La)}}
\bigl(\k_{_{\AA}}(\la),\,\k_{_{\AA}}(\mu)\bigr)
\iso
\Ext\hdot_{_{D_\Y^\B(\B,\bb)}}\bigl(\fF(\k_{_{\AA}}(\la)),
\,\fF(\k_{_{\AA}}(\mu))\bigr)=
\Ext\hdot_{_{D_\Y^\B(\B,\bb)}}(\k_{_\B}(\la),\k_{_\B}(l\mu)).
\end{align*}
\end{theorem} 

\ab The rest of this section is mainly devoted to the proof of the
Theorem.

\subsection{Comparison of functors.} The functor \eqref{fF} has a `non-equivariant analogue',
obtained by forgetting the Hopf-adjoint actions. Specifically,
we form the composite map $\alpha': \bb\ltimes\ka\stackrel\pi\too
\bb\ltimes\bb\stackrel\mult\too\bb$,
 where  $\pi: \ka\qisto\bb$ is the quasi-isomorphism constructed in
\S\ref{dgres_sec}.
Let $\alpha'_*:\ D_\Y(\bb\ltimes\ka,\bb)\map D_\Y(\bb,\bb)$
denote the corresponding functor.
Also, we have the following non-equivariant counterpart of
the bifunctor \eqref{bifunctor}:
\beq\label{bifunctor2}
  D_\Y(\tK
\otimes_{_Z}\tK,\La) \times
  D_\Y(\bb\ltimes\ka,\bb)
\longrightarrow D_\Y(\bb\ltimes\ka,\bb),\quad
M,N \mto M \lotimes_{\tK}N\,. 
\eeq

We use the
equivalence $\imath_*$ induced  by a
non-equivariant analogue  of Lemma \ref{tKtK_lemma} (cf. also \eqref{tor2}),
to obtain the following 
functor
\begin{align}\label{functor_F}
\bff:\; 
  D_\Y(\La,\La)\overset{\imath_*}{\iso}
  &D_\Y(\tK\otimes_{_{Z}}\tK,\La)
\;\underset{\eqref{bifunctor2}}{\stackrel{(-)\lotimes_{_\tK}\k}{\tooo}}
\;    D_\Y(\bb\ltimes\ka,\bb)\stackrel{\alpha'_*}\too   D_\Y(\bb,\bb),\nonumber\\
&
M\, \mto\, \bff(M):=\alpha'_*\bigl(M\lo_{_{\tK}}\,\k\bigr)
=
\mult_*\pi_*\bigl(\imath_*M\lo_{_{\tK}}\,\k\bigr).
\end{align}

\ab In order to prove Theorem \ref{formality3}, we will need
to relate the functors $\fF$ and $\bff$.
To this end, given an algebra $A$ 
and a
subalgebra $\aa\sset A$, we let $\res^A_\aa: A\mmod\to\aa\mmod$
denote the obvious restriction functor, and use similar notation for 
derived categories.

\begin{lemma} \label{functors}
There is an isomorphism of  functors: 
$
\res_{\bb}^{\sB}\circ\fF\;\simeq\;\bff\circ 
\res_{\La}^{\AA}$, in other words, the following diagram commutes:
$$
\xymatrix{
  D_\Y^{\Ub}(\AA,\La)\ar[rr]^<>(.5){\res_{\La}^{\AA}}
\ar[d]_<>(.5){\fF}&&  D_\Y(\La,\La)
\ar[d]_<>(.5){\bff}\\
  D_\Y^\B(\sB,\bb)\ar[rr]^<>(.5){\res_{\bb}^{\sB}}&&  D_\Y(\bb,\bb)\;.
}
$$
\end{lemma}

\proof We restrict the algebra morphisms $\al$ and $\beta$ in
diagram \eqref{curve} to the subalgebra
$\bb\ltimes\ka\sset\B\ltimes\ka$,
and consider the following diagram:
\begin{equation}\label{curve2}
\xymatrix{
&&&&\bb\ltimes\ka=\bb\ltimes(\tK\otimes_{_Z}\fB)
\ar[dllll]_<>(.5){\beta|_{\bb\ltimes\ka}}
\ar[d]_<>(.5){\Id_\bb\ltimes\pi}\ar[drrr]^<>(.5){\al|_{\bb\ltimes\ka}}&&&\\
\bb&&\bb\otimes\bb\ar[ll]^<>(.5){\Id\otimes\epsilon_{_\bb}}
&&\bb\ltimes\bb\ar[ll]_<>(.5){\sim}^{\gamma}\ar[rrr]_<>(.5){\mult}&&&\bb
}
\end{equation} 
In this diagram, the map $\pi: \ka\to\bb$ is the quasi-isomorphism
of \S\ref{dgres_sec}, and the map $\mult: \bb\ltimes\bb\to\bb$ (given
by multiplication in the algebra $\bb$) is an algebra
morphism, by Proposition \ref{Montgomery}.
Thus, the  right triangle in diagram \eqref{curve2} commutes,
by definition of the map $\al$.

\ab In the left triangle of diagram \eqref{curve2},
we have the map
$$\beta|_{\bb\ltimes\ka}:\ \bb\ltimes\ka=\bb\ltimes (\tK\otimes_{_Z}\fB)\too
\bb,\quad
b\ltimes(r\otimes \tilde{b})\mto b\cdot\epsilon_{_\tK}(r)\cdot\pi(\tilde{b})\,.
$$
Further, the map $\gamma$ in  \eqref{curve2} is the algebra isomorphism
of Proposition \ref{Montgomery}(iii), which is given by:
$b\ltimes b'\mto b\otimes bb'$. Thus, we see that 
the triangle on the left of diagram \eqref{curve2} commutes also.

\ab Computing the inverse of $\gamma$, we get
$\gamma^{-1}(b\otimes b')=b\ltimes (S(b)\cdot b')$.
Therefore, we find:
$\mult\ccirc\gamma^{-1}(b\otimes b')= \epsilon_\bb(b)\cdot b'$.
Further, observe that since the $\adh\bb$-action on $Z$ is trivial, 
we have an algebra isomorphism
$\bb\ltimes(\tK\otimes_{_Z}\fB)\simeq
\tK\otimes_{_Z}(\bb\ltimes\fB)$. 
We introduce the composite quasi-isomorphism
\begin{align*}
\theta:= \gamma\ccirc(\Id_\bb\ltimes\pi):\;\;
\bb\ltimes\ka= &\tK\otimes_{_Z}(\bb\ltimes\fB)\qisto\bb\otimes\bb,
\\ &
r\otimes (b\ltimes\tilde{b})
\mto \epsilon_{_\tK}(r)\cdot\bigl(b\otimes(b\cdot \pi(\tilde{b}))\bigr).
\end{align*}
\ab Thus,
we can rewrite  diagram \eqref{curve2} in the following more symmetric form:
\beq\label{pr2} 
\xymatrix{
&&\tK\otimes_{_Z}(\bb\ltimes\fB)\ar[dll]_<>(.5){\beta|_{\bb\ltimes\ka}}
\ar[d]_<>(.5){\theta}\ar[drr]^<>(.5){\al|_{\bb\ltimes\ka}}&&\\
\bb&&\bb\otimes\bb\ar[ll]^<>(.5){\Id\otimes\epsilon_{_\bb}}
\ar[rr]_<>(.5){\epsilon_{_\bb}\otimes\Id}
&&\bb
}
\eeq

\ab Next, let $\ste:=\res^{\sB}_\bb\St$ be  the Steinberg module $\St$ viewed
as an object of $  D_\Y(\bb,\bb)$.
Since $\St|_{\u^+}\simeq \u^+$,
we deduce that: $(\epsilon_{_\bb})_*(\ste)=\k$, where
$(\epsilon_{_\bb})_*$
is the direct image functor
corresponding to the augmentation $\epsilon_{_\bb}: \bb\to\k$.
On the other hand, one may also view the $\bb$-module $\ste$
as a $\fB$-module via the projection $\fB\onto
\fB/(Z)=\bb$. Then,
using the left triangle in diagram \eqref{pr2}, we get:
\begin{equation}\label{St_property}\res_{\bb\ltimes\ka}^{\sB\ltimes\ka}
(\beta^*\St)=\theta^*\ccirc(\Id\otimes\epsilon_{_\bb})^*\St'=
\theta^*(\St'\otimes\k)=\St'\otimes\k\,. 
\eeq
Therefore,
for any $\tilde{M}\in  D_\Y^{\Ub}\bigl(\Ub\ltimes(\tK\otimes_{_Z}\tK),\La\bigr)$,
applying  bifunctors \eqref{bifunctor}, \eqref{bifunctor2}, we obtain
\beq\label{cas}
\theta_*\bigl(\res_{\bb\ltimes\ka}^{\sB\ltimes\ka}(\tilde{M}\lo_{_{\tK}}\beta^*\St)\bigr)=
\theta_*\bigl((\res_{\tK\otimes_{_Z}\tK}^{\Ub\ltimes(\tK\otimes_{_Z}\tK)}\tilde{M})
\lo_{_{\tK}}(\St'\otimes\k)\bigr).
\eeq

\ab Now,  let $\tilde{M}:=\imath_*M,$ for some
$M\in   D_\Y^{\Ub}(\AA,\La).$ Then,
using the right triangle in \eqref{pr2}
and the definition of the functor $\bff$, see \eqref{functor_F}, we obtain
\begin{align}
\res_{\bb}^{\sB}\fF(M)&=
\res_{\bb}^{\sB}\circ\al_*\bigl(M\lo_{_{\tK}}\beta^*\St\bigr)
\enspace\text{by \eqref{cas}}\nonumber\\
&=(\epsilon_{_\bb}\otimes\Id)_*\circ\theta_*\bigl(M\lo_{_{\tK}}\beta^*\St\bigr)
\nonumber\\
&=
(\epsilon_{_\bb}\otimes\Id)_*\circ\theta_*\bigl(\tilde{M}\lo_{_{\tK}}
(\ste\otimes\tK)\bigr)\nonumber\\
&=(\epsilon_{_\bb}\otimes\Id)_*
\bigl(\ste\otimes
\bff(\res_{\La}^{\AA}M)\bigr)\nonumber\\
&=\bigl((\epsilon_{_\bb})_*(\ste)\bigr)\otimes 
\bff(\res_{\La}^{\AA}M)=
\bff(\res_{\La}^{\AA}M)\,,\nonumber
\end{align}
and the Lemma is proved.\qed

\subsection{`Deformation' morphism.} 
Below, we will use  a
 well-known result of Gerstenhaber saying that,
 for any algebra $A$, the graded algebra
$\Ext_{_{A\bimod}}\hdot(\underline{A},\underline{A})$
 is  always commutative.
We also remind the reader that 
the category $A\mmod$ may be viewed as a {\it module category}
over the category $A\bimod,$ of
$A$-bimodules. This gives, for any
$M\in A\mmod$, a canonical
graded algebra morphism, to be referred to as {\it evaluation} at $M$:
 
\beq\label{module_category} 
\ev_{_M}:\ \Ext_{A\bimod}\hdot({\underline{A}},
{\underline{A}})\too\Ext_{A\mmod}\hdot(M,M)\,.
\eeq

\ab Now,
recall that the quantum Borel algebra $\fB$ is a free module over
its central subalgebra $Z$. Let
$\epsilon\in \Spec Z$ denote  the `base-point' corresponding to the
augmentation ideal $Z\eps\sset Z$.
We will view $\fB$ as a flat family
of (non-commutative) algebras over the smooth base $\Spec Z$
whose   fiber over the  base-point 
 is the algebra
$\bb=\k\eps\otimes_{_Z}\fB\cong\fB/(Z)$. 
Otherwise put, the algebra  $\fB$ is a multi-parameter
deformation of $\bb$. 
By the classical
work of Gerstenhaber, such a deformation gives
a linear map $T_\epsilon(\Spec Z)\to
\Ext_{_{\bb\bimod}}^2(\uaa,\uaa)$, 
where $T_\epsilon(\Spec Z)$ denotes the tangent space at the point $\epsilon$.
By commutativity of the algebra
$\Ext_{_{\bb\bimod}}\hdot(\uaa,\uaa),$
 the linear map above can be uniquely extended,
by multiplicativity, to a degree doubling
 algebra morphism
\beq\label{deform_map}
{\mathsf {deform}}:\;\; \sym\hdot\bigl(T_\epsilon(\Spec Z)\bigr)\too
\Ext_{_{\bb\bimod}}^{2\hhdot}(\uaa,\uaa).
\eeq

\ab Next,  we would like to take the $\adh\B$-action on 
$\bb$, cf. Proposition \ref{GBbis}, into
considerations. The $\adh\B$-action induces, for each $j\geq 0$,
a  $\sB$-action on $\Ext_{_{\bb\bimod}}^j(\uaa,\uaa)$, that
makes $\Ext_{_{\bb\bimod}}\hdot(\uaa,\uaa)$
a graded $\sB$-algebra.
Enhancing Gerstenhaber's construction to the equivariant setting, one 
finds that the algebra map \eqref{deform_map} 
is actually a morphism of $\sB$-modules. 

\ab It turns out that the $\sB$-action on each side of 
\eqref{deform_map} descends  to the algebra
$\sB/(\bb)$,
which is isomorphic to $\Ub$ via the Frobenius map.
 This follows, for the RHS of \eqref{deform_map}, 
from the general result saying that the
 Hopf-adjoint action of any Hopf algebra $\aa$ on 
$\Ext\hdot_{\aa\bimod}(\underline{\aa},\underline{\aa})$ is trivial, see
\S \ref{cohomology}.
  For the LHS, we use
 Corollary \ref{tangent} saying that
there is a canonical $\adh\sB$-equivariant isomorphism of vector spaces
$T_\epsilon(\Spec Z)\simeq Z\eps/Z\eps^2\simeq \n^*$.
Thus, the morphism in \eqref{deform_map} becomes the following
  $\Ub$-equivariant graded algebra morphism
\beq\label{def_map_bimod}
{\mathsf {deform}}:\;\; \sym\hdot(\n^*[-2])\too
\Ext_{_{\bb\bimod}}\hdot(\uaa,\uaa)\,.
\eeq

\subsection{General deformation formality theorem.} Our proof of Theorem
\ref{formality3}
 is based on a much more general Theorem \ref{BG_formality2} below,
 proved in
\cite{BG}.

\ab To explain the setting of \cite{BG}, let $\bb$ be  (temporarily)
 an {\it arbitrary} associative
algebra and
 $\fB$  an {\it arbitrary} flat deformation of $\bb$ over a smooth base
$\Spec Z$. Choose
$\tK$, a $Z$-free \dg-resolution
  of the trivial $Z$-module $\k_{_Z}$
(as in Lemma \ref{tK}, but with $\Ub$-action ignored),
corresponding to the base point $\epsilon\in \Spec Z$.
Put $\n:= T\eps Z$ and $\La:=\wedge\hdot(\n^*[1]).$
We form the \dg-algebras $\ka:=\tK\otimes_{_Z}\fB$,
and $\tK\otimes_{_Z}\tK$. Then we 
 establish, as we have done in\S \ref{resol},
  \dg-algebra  quasi-isomorphisms
$\imath: \
\La\qisto\tK\otimes_{_Z}\tK$,
cf. \eqref{tor2},
and $\pi: \ka\qisto\bb$. 
Thus formula \eqref{functor_F} gives, in our general situation,
a well-defined functor $\bff:   D_\Y(\La,\La)\too   D_\Y(\bb,\bb)$.

\ab 
In \cite{BG}, we  prove the following result.
\begin{theorem}
\label{BG_formality2}
We have $\bff(\k_\La) =\k_\bb$.
Furthermore, the
induced 
map $\bff_*$ makes the following diagram commute 
$$
\xymatrix{
\sym(\n^*[-2])\ar[r]_<>(.5){\sim}^<>(.5){\eqref{SL_kos}}
\ar[d]_<>(.5){{\mathsf {deform}}}^<>(.5){\eqref{def_map_bimod}}
&\Ext\hdot_{_{D_\Y(\La,\La)}}(\k_\La,\k_\La)
\ar[rr]^<>(.5){\bff_*}_<>(.5){\sim}
&&\Ext\hdot_{_{  D_\Y(\bb,\bb)}}\left(\bff(\k_\La)\,,\,\bff(\k_\La)\right)
\ar@{=}[d]\\
\Ext\hdot _{\bb\bimod}(\uaa,\uaa)
\ar[rrr]^<>(.5){\ev}_<>(.5){\eqref{module_category}}&&&
\Ext\hdot _{D_\Y(\bb,\bb)}(\k_\bb,\k_\bb).
}
$$
\end{theorem}

\subsection{Proof of Theorem  \ref{formality3}.} Recall that $\AA:=\Ub\ltimes\La$.
We first prove that
\beq\label{FtoF}
\fF(\k_{_{\AA}}(\la))=\k_{_\B}(l\la),\quad\forall\la\in\Y\,.
\eeq
To this end, we use
 Lemma~\ref{functors} and Theorem \ref{BG_formality2} to
 deduce that the \dg-module $\res_\bb^\sB\fF(\k_{_{\AA}})\in   D_\Y(\bb,\bb)$
is 
 quasi-isomorphic
 to the trivial  module $\k$. In particular, it has 
a single non-zero cohomology group: $H^0\bigl(\fF(\k_{_{\AA}})\bigr)\simeq \k,$
no matter whether it is considered as a $\sB$-module, or as a $\bb$-module.
But the action of  the augmentation ideal
$\,(\U^+)\eps\,$ of the 
(sub)algebra $\U^+\sset\B$
on the cohomology of any object of the category $D_\Y(\B,\bb)$ is necessarily
nilpotent. Hence the subalgebra  $\U^+$
acts trivially (that is, via the augmentation)
on the 1-dimensional vector space $H^0\bigl(\fF(\k_{_{\AA}})\bigr)\simeq \k.$
Furthermore, since the module $\St$ has been normalized
so that its highest weight is equal to zero, it immediately follows
that the  1-dimensional  space  $H^0\bigl(\fF(\k_{_{\AA}}(\la))\bigr)$ 
has weight $l\la$ 
with respect to the  $\U^\circ$-action.
Thus, we have a $\B\ltimes\bb$-module isomorphism
$H^0\bigl(\fF(\k_{_{\AA}}(\la))\bigr)\simeq \k_{_{\B}}(l\la)$. We conclude that
the object $\fF(\k_{_{\AA}}(\la))$
is quasi-isomorphic to $\k_{_{\B}}(l\la)\in
D_\Y(\B,\bb)$, and \eqref{FtoF} is proved.

\ab To complete the proof the Theorem, we must  show that the functor $\fF$ induces
a graded algebra isomorphism
\begin{align}\label{ugu1}
\Ext\hdot _{  D_\Y^{\Ub}(\AA,\La)}(\k_{_{\AA}}(\la),
\k_{_{\AA}}(\mu))\iso 
\Ext\hdot_{D_\Y^\B(\B,\bb)}&\bigl(\fF(\k_{_{\AA}}(\la))\,,
\,\fF(\k_{_{\AA}}(\mu))\bigr)\\
&=\Ext\hdot _{D_\Y^\B(\B,\bb)}(\k_{_\B}(l\la),\k_{_\B}(l\mu))\,.
\nonumber
\end{align}
To compare the Ext-groups on the LHS and on the RHS, we use the spectral
sequence provided by
Lemma \ref{resind}. Specifically, since
$\AA/(\La)=\Ub$ and $\B/(\bb)=\Ub$
 we have the following 
two spectral sequences, see \eqref{sp_seq}:
\beq\label{sp_seq2}
\xymatrix{
H^p\bigl(\Ub\,,\,\Ext_{  D_\Y(\La,\La)}^q(\k_{_{\AA}}(\la),
\k_{_{\AA}}(\mu))\bigr)=
E_2^{p,q}\;\ar[d]^<>(.5){\bff_*}\ar@{=>}[r]
&\;
E^{p+q}_\infty=\gr\Ext_{  D_\Y^{\Ub}(\AA,\La)}^{p+q}(\k_{_{\AA}}(\la),
\k_{_{\AA}}(\mu))\ar[d]^<>(.5){\fF_*}\\
H^p\bigl(\Ub\,,\,\Ext_{  D_\Y(\bb,\bb)}^q(\k_{_\B}(l\la),\k_{_\B}(l\mu))\bigr)=
E_2^{p,q}\;\ar@{=>}[r]
&\;
E^{p+q}_\infty=\gr\Ext^{p+q}_{D_\Y^\B(\B,\bb)}(\k_{_\B}(l\la),\k_{_\B}(l\mu))\,.
}
\eeq

\ab The vertical arrow  on the left of the diagram is induced by
the map $\bff_*: \Ext_{  D_\Y(\La,\La)}\hdot(\k_{_{\La}},
\k_{_{\La}})$
$\to \Ext_{  D_\Y(\bb,\bb)}\hdot
\bigl(\bff(\k_{_{\La}})\,,\,\bff(\k_{_{\La}})\bigr)$. Theorem \ref{BG_formality2}
implies that this map may be identified with the composite map
$${\mathsf {deform}}\ccirc\ev_{_{\k_\bb}}
:\; \Ext_{  D_\Y(\La,\La)}\hdot(\k_{_{\La}},\k_{_{\La}})\cong
\sym(\n^*[-2])\too \Ext_{  D_\Y(\bb,\bb)}\hdot
\bigl(\bff(\k_{_{\La}})\,,\,\bff(\k_{_{\La}})\bigr)\,.
$$
But
the latter map is exactly the map that was used 
in \cite{GK}  to construct the
$\Ub$-equivariant isomorphism in
Proposition \ref{gk}. 
Thus, we conclude that the functor $\bff_*$
induces an isomorphism between the $E_2$-terms of the
two  spectral sequences in \eqref{sp_seq2}.

\ab Further, the  vertical map between the $E_\infty$-terms 
on the right of diagram \eqref{sp_seq2}
is  induced by the functor $\fF$.
This map coincides, by Lemma \ref{functors},
with the map induced by the isomorphism
 between the $E_2$-terms of the spectral sequences.
Hence, it is itself an isomorphism. It follows that
 morphism \eqref{ugu1} is an isomorphism.
\qed

\bigskip
\pagebreak[3]
\centerline{\bf PART $\,$II$\,$: $\;$ Geometry.}
\bigskip

\ab Throughout Part II (with the exception of \S9) we let $\k=\C$.

\section{The loop Grassmannian and the Principal nilpotent.}
\label{sec_loop}


\ab In this section we recall a  connection, discovered
in [G2],
between the cohomology of a loop Grassmannian, to be introduced below,
and the principal nilpotent element in the Lie algebra $\g$.

\ab Let $D^b(X)$ be the bounded derived category
of constructible complexes on an algebraic variety $X$,
cf. \cite{BBD}.
Given an algebraic group ${\mathbf{G}}$ and a
${\mathbf{G}}$-action on $X$, we let
$D^b_{_{\mathbf{G}}}(X)$
denote the $\mathbf{G}$-equivariant
 bounded derived category  
 on $X$; see
[BL] for more information on the equivariant derived category.
We write $D^b_{_{\mathbf{G}\mbox{-}{\sf{mon}}}}(X)$,
 for the full subcategory of
$D^b(X)$ formed by "$\mathbf{G}$-monodromic" complexes,
that is,  formed  by complexes whose cohomology sheaves are locally
constant along $\mathbf{G}$-orbits.

\ab We let $\pe_{_{\mathbf{G}}}(X)\sset D^b_{_{\mathbf{G}}}(X)$,
resp., $\pe_{_{\mathbf{G}\mbox{-}{\sf{mon}}}}(X)\sset 
D^b_{_{\mathbf{G}\mbox{-}{\sf{mon}}}}(X)$ stand for the abelian category
of $\mathbf{G}$-equivariant, resp.
$\mathbf{G}$-monodromic, perverse sheaves on $X$.

\subsection{The loop group.} Let $\Gd$ be a 
complex connected
semisimple group with  maximal torus 
$\Td=\C^*\otimes_{_\Z}\Y$, which is  dual to $(G, T)$ in the sense
of Langlands. Thus,  $\Gd$ is a simply-connected group
 such that the root system of $(\Gd, \Td)$
is dual to that of $(\g,\h)$. 
Let
$\gd=\Lie\Gd$ be the Lie algebra
of $\Gd$. The Lie algebra of the maximal torus $\Td\sset \Gd$
gives a distinguished Cartan subalgebra:
$\Lie\Td =\hd=\C\otimes_{_\Z}\Y=\h^*$ in $\gd$.

\ab Let $\K=\C((z))$ be the field of formal
Laurent power series, and $\oo=\C[[z]] \sset\K$ its ring of integers,
that is, the ring of
formal power series regular at $z=0$. 
Write $\Gd(\K)$, resp. $\Gd(\oo)$,
for the set of  $\K$-rational, resp. $\oo$-rational,
points of $\Gd$. The coset space $\Gr := \Gd(\K)/\Gd(\oo)$ is called
the {\it loop Grassmannian}. It
has the natural structure of an ind-scheme, more precisely, $\Gr$ is
a direct limit of a sequence
of $\Gd(\oo)$-stable
projective varieties of increasing dimension, see
 e.g. [BD], [Ga], [G2] or [L1], such that
the action of $\Gd(\oo)$ on any such variety factors through
a finite-dimensional quotient of  $\Gd(\oo)$.

\ab An Iwasawa decomposition for $\Gd(\K)$,
 see [G2], [PS], implies that the loop Grassmannian is isomorphic,
as a topological space, to  the space
 of based loops into a compact
form of the complex group $\Gd$.
 It follows that $\Gr$ is an $H$-space,
hence, the cohomology $H\hdot(\Gr, \C)$ has the natural structure
of a graded commutative and cocommutative Hopf algebra.
Further, the group $G^\vee$ being simply-connected, we deduce that
the loop  Grassmannian
 $\Gr$ is connected. 

\subsection{Cohomology of the loop Grassmannian.}
We recall that 
the group $\Gd(\oo)$ is homotopy equivalent to $\Gd$,
hence for the $\Gd(\oo)$-equivariant cohomology of a point
we have
\beq\label{point}
H\hdot\GO(pt)=H\hdot_{_{G^\vee}}(pt) =
H\hdot(BG^\vee)=\C[\hd]^W=\C[\h^*]^W=\C[\g^*]^G=(\sym\g)^G\,,
\eeq
where $BG^\vee$ stands for the classifying space of the group $G^\vee$.
Thus
$H\hdot\GO(\Gr)$,
the $\Gd(\oo)$-equivariant cohomology of the Grassmannian,
has a natural
$H\hdot\GO(pt)$-module structure, hence, a $\C[\h^*]^W$-module structure.

\ab Next, we introduce the notation $\g^x\sset\g$ for the Lie algebra
of the isotropy group of an element  $x\in\g^*$ under the
coadjoint action (thus, if one uses the identification
$\g\simeq \g^*$ provided by an invariant form,
then $\g^x$ becomes the centralizer of $x$ in $\g$).
Let $\greg\sset\g^*$
be the Zariski open dense subset of regular (not necessarily
semisimple) elements in $\g^*$. 
The family of spaces $\,\{\g^x\,,\,x\in \greg\}
,$ gives an $\Ad G$-equivariant vector bundle
on $\greg$. We let $\z$ denote the corresponding vector bundle on the
adjoint quotient space $\greg/\Ad G$. 
The latter space is isomorphic, due to Kostant [Ko],
to $\h^*/W$, hence is a smooth affine variety. We will often regard $\z$ as a 
vector
bundle on $\h^*/W$ via the Kostant isomorphism. The fibers of $\z$ are abelian
Lie subalgebras in $\g$, and we let ${\mathcal{U}}\z$ denote the
vector bundle on $\greg/\Ad G$ with fibers ${\mathcal{U}}(\g^x)\,,\,x\in \greg$. 
Let  $\Gamma(\h^*/W\,,\,{\mathcal{U}}\z)$ be the
 commutative algebra of 
global regular sections of ${\mathcal{U}}\z$.

\ab 
Fix a
principal $\sll$-triple $\langle \s,\e,\f \rangle \subset \g$,
such that $\s=\sum_{\alpha\in R_+}\,\check\alpha\in \h,$
and such that the principal
nilpotent $\e \in \g$ is a linear combination of  simple
root vectors
 with non-zero coefficients.
The  Lie algebra
${\gee}$  is an abelian Lie subalgebra in $\g$ of dimension
$\rk\g$. The ad-action of $\s$ puts a grading on $\g$,
and we endow ${\gee}$ and ${\mathcal{U}}({\gee})$,
the enveloping algebra of ${\gee}$, with
 induced gradings. 

\ab  The 
natural $\C^*$-action on $\h^*$ by dilations
makes $\h^*/W$ a $\C^*$-variety. Moreover, the grading on $\g$
considered above, gives a $\C^*$-action on $\g$, hence,
 makes $\z$ a $\C^*$-equivariant vector bundle
on $\h^*/W$. Thus, $\Gamma(\h^*/W\,,\,{\mathcal{U}}\z)$
acquires a grading compatible with the algebra structure.
In \cite{G2} we have proved the following

\begin{lemma}[\cite{G2}]\label{H_G(Gr)1}
There is a natural  graded  Hopf algebra isomorphism
$\,\varphi_{_{\mathcal{U}}}:
H\hdot\GO(\Gr)$
$\simeq
\Gamma(\h^*/W\,,\,{\mathcal{U}}\z).$
\qed
\end{lemma}

Observe that
the fiber over $0\in\h^*/W$ of the vector bundle
$\z$ clearly identifies with $\g^e$. Hence, from Lemma 
\ref{H_G(Gr)1} we get

\begin{corollary}\label{H(Gr)}{{\sf[G2, Proposition 1.7.2]}}$\,$
There is a natural graded Hopf algebra isomorphism $\varphi:
H\hdot(\Gr, \C) 
\stackrel{\sim}{\longrightarrow}
{{\mathcal{U}}}(\g^e).$
\qed
\end{corollary} 

\subsection{Geometric Satake equivalence.}\label{satake_sec}
 Let $D^b(\Gr)$ denote the bounded derived category of
constructible complexes on $\Gr$, to be understood as a direct limit
of the
corresponding  bounded derived categories on finite dimensional projective
subvarieties that exhaust $\Gr$.  
One  similarly defines $\pe(\Gr)\sset D^b(\Gr)$, the abelian category 
of  perverse sheaves.

\begin{definition} Let $\PO$
be the (full) abelian subcategory in $\pe(\Gr)$ formed by
 semisimple\footnote{Any $\Gd(\oo)$-equivariant 
perverse sheaf  on $\Gr$ is, in effect, automatically semisimple,
cf. e.g. \cite{MV}.}
 $\Gd(\oo)$-equivariant 
perverse sheaves on $\Gr$.
\end{definition}

\ab   For any $\L \in \PO$, there is a standard convolution functor: 
$D^b(\Gr)\to
D^b(\Gr),$ $\,\M\mapsto\M\star\L:=a_*(\M\tilde{\boxtimes}\L),$
where
$a: \Gd(\K)\times_{_{G^\vee(\oo)}}\Gr\too\Gr$ is the action-map,
and $\M\tilde{\boxtimes}\L$ stands for a twisted version of
external tensor product, see [MV] or [G2] for more details.
A fundamental result due to Gaitsgory says that this  functor
takes perverse sheaves 
 into  perverse sheaves, that is, we have the following
\begin{theorem}[\cite{Ga}]\label{denis} The convolution 
gives  an exact bifunctor
$$\pe(\Gr)\times\PO \too \pe(\Gr)\,,\quad \M,\L\mto\M\star\L.\qquad\Box$$
\end{theorem}

\begin{remark} In the special case,  where $\M\in \pe(\Gr)$ is a perverse
sheaf which is constant along the Schubert cell stratification
(by Iwahori orbits) of the loop Grassmannian, cf. \S\ref{eq_perv_sub}
below, the Theorem above has been first conjectured
in \cite[p.22]{G2}, and proved by Lusztig \cite{L4} shortly
after that. This special case of Theorem \ref{denis} 
is the only case that will be actually used in the present paper.
$\quad\lozenge$
\end{remark}
 
\ab   For each $g\in \Gd(\K)$, the double coset $\Gd(\oo)\cdot g\cdot
\Gd(\oo)$ contains an element $\lambda \in \Hom(\C^*,\Td)$,
viewed as a loop in $\Gd$. Moreover, such an element is unique
up to the action of $W$, the Weyl group. This gives
a parametrization of $\Gd(\oo)$-orbits in 
$\Gr$ by dominant (co)weights 
$\lambda \in \Hom(\C^*,\Td)^{++}=\Y^{++}$. 
We write $\Gd(\oo)\cdot\lambda$ for the 
$\Gd(\oo)$-orbit  corresponding to a dominant (co)weight
$\lambda$. The closure, $\overline{\Gd(\oo)\cdot\lambda}\sset \Gr$ is known
to be a finite dimensional projective variety, singular in general.
Let $IC_\la$ denote the intersection complex on 
$\overline{\Gd(\oo)\cdot\lambda}$
corresponding to the constant sheaf on $\Gd(\oo)\cdot\lambda$ 
(extended by zero on $\Gr\smallsetminus
\overline{\Gd(\oo)\cdot\lambda}$, and normalized to be a perverse sheaf).
The $IC_\la\,,\, \la\in\Y^{++}$, are the simple objects of the
category
$\PO$. Theorem \ref{denis}
puts on $\PO$  the structure of a monoidal category,
via the  convolution product.

\ab Recall the tensor category $\rep(G)$
 of finite dimensional rational
representations of  $G$.   For each $\la\in\Y^{++}$,
let $V_\la\in \rep(G)$ denote an irreducible representation with 
highest weight $\la$. The proof of the following fundamental 
result, inspired by Lusztig \cite{L1},
can be found in
\cite[Theorem 1.4.1]{G2} (following an idea of Drinfeld);
a more geometric proof (involving a different
commutativity constraint, also suggested by Drinfeld)
has been found later in [MV]; the most conceptual
argument is given in [Ga].

\begin{theorem}\label{tensor}
There is an equivalence $\p: \rep(G)\iso\PO,$ of  monoidal categories 
which sends  $V_\lambda$ to $ IC_\lambda$, for any $\la\in\Y^{++}$. \qed
\end{theorem}

\subsection{Fiber functors.}
 We will need a more elaborate `equivariant' version
of Theorem \ref{tensor}, established in [G2]. To formulate it,
identify $\Td$ with the subgroup in $\Gd(\K)$ formed
by constant loops into $\Td$. Thus, any
object $A\in \PO$ may be regarded as a
$\Td$-equivariant perverse sheaf,
hence there are well defined $\Td$-equivariant (hyper)-cohomology
groups,
$H_{\Td}\hdot(\Gr,A)$.
Given $s\in \Lie \Td$, we write
$H_s(A)=H_{\Td}\hdot(\Gr,A)\big|_s$ 
for the $\Td$-equivariant (hyper)-cohomology
of $A$ specialized at  $s$, viewed
as a  point in $\Spec H\hdot_{\Td}(pt)$.
  For $s=0$, we have   $H_s(A)=H\hdot(\Gr,A)$, is the ordinary cohomology of
$A$
(due to  the collapse of the spectral sequence for  equivariant
intersection cohomology).

\ab Observe that,  for $s$ regular, the $s$-fixed
point set in $\Gr$ is the  lattice $\Y$, viewed as a discrete subset in
$\Gr$ via
the  natural imbedding
$i: \Y=\Hom(\C^*,\Td)\into \Gr$.
  For each $\la\in \Y$,
let $i_\la: \{\la\}\into\Gr$ denote the corresponding one
point imbedding. 
By the Localisation theorem in equivariant cohomology,
 the map: $H_s(i^!A)\to H_s(A)$, 
 induced by the adjunction morphism:
$i_!i^!A\to A$ yields, see [G2, (3.6.1)], 
 the following direct sum 

{\bf Fixed point decomposition:}
\vskip -3pt
\begin{equation}\label{fp_decomposition}
H_s(A)= \bplus_{\la\in\Y}\;H_s(i_\la^!A),\quad \forall A\in \PO\,.
\end{equation}

\ab Recall
the principal $\sll$-triple $\langle \s,\e,\f \rangle \subset \g$.
Observe that the element $\s+e$ is $\Ad G$-conjugate to $\s$, hence is a regular
semisimple element in $\g$. Thus,
  $\fh:=\g^{\s+e}$
is a Cartan subalgebra. Furthermore, the fiber of the
vector bundle ${\mathcal{U}}\z$ over the 
$\Ad G$-conjugacy class $\Ad G(\s+e)=\Ad G(\s)
\subset \greg$ gets identified with
${\mathcal{U}}(\g^{\s+e})={\mathcal{U}}\fh$.

\begin{proposition}[\cite{G2}]\label{main_perv_prop}
  \vi$\;$   For any $s\in \Lie\Td$,
the assignment $A \mapsto H_s(A)$ gives a fiber
functor on the tensor category $\PO$.

\ab \vii$\;$ There is an isomorphism of the functor $H_s(-)$ on
$\PO$ with the forgetful functor on $\rep(G)$
(i.e.,  a system of isomorphisms $\varphi_{_V}:
 H_s(\p V)\iso V\,,
\,\forall V\in \rep(G)$, compatible with morphisms in
$\rep(G)$ and with the tensor structure) such that:

\ab \parbox[t]{140mm}{{\sl For any
$u\in H_s(\Gr)$,
the natural action of $u$ on the hyper-cohomology $H_s(\p V)$ corresponds,
via $\varphi_{_V}$ and the isomorphism $\varphi_{_{\mathcal{U}}}:
H_s(\Gr) 
\stackrel{\sim}{\longrightarrow}
{{\mathcal{U}}}(\fh)$ of Lemma  \ref{H_G(Gr)1},
to the natural
action of $\varphi_{_{\mathcal{U}}}(u)
\in {{\mathcal{U}}}(\g^{\s+e})$ in the  $G$-module $V$.}
\quad\qed}
\end{proposition}

\subsection{Equivariant and Brylinski's filtrations.}
The standard grading on the equivariant cohomology
$H_{\Td}\hdot(\Gr,A)$
induces, after specialization at a point
 $s\in \Spec H_{\Td}\hdot(\Gr)=\check{\tg}$, a canonical increasing filtration,
$W\idot H_s(\Gr,A),$ on the  specialized
equivariant
cohomology. Furthermore, the collapse of the spectral sequence
for equivariant intersection cohomology yields a natural isomorphism
$H\hdot(\Gr,A)\cong H\hdot_{\{o\}}(A)$, where the LHS stands for
the (non-equivariant)   cohomology of $A$
and the RHS  stands for the  specialization 
of  equivariant   cohomology of $A$
at the zero point: $o\in \check{\tg}=\Spec H_{\Td}\hdot(\Gr)$. 
On the other hand, for any
 $s\in \check{\tg}=\Spec H_{\Td}\hdot(\Gr)$, one has 
a  canonical
graded space isomorphism
$\gr^W\idot H_s(A) \simeq H\hdot_{\{o\}}(A)$,
by  the definition of filtration $W\idot$.
Thus, composing the two isomorphisms we obtain, for any 
 $s\in \Spec H_{\Td}\hdot(\Gr)$, a
natural
graded space isomorphism
$\gr^W\idot H_s(A) \simeq H\hdot(\Gr,A)$.

\ab From now on we will make a particular choice of 
the point $s \in \Spec H_{\Td}\hdot(\Gr,A)=\check{\tg}.$
Specifically, we let $s=\sum_{\alpha\in R_+}\,\alpha
\in\check{\tg}
$ be the element "dual", in a sense, to 
 $\s=\sum_{\alpha\in R_+}\,\check\alpha\in\tg$.
Further, the eigenvalues of
the $\s$-action in any 
 finite dimensional $G$-module $V$ are known to be integral.

\begin{definition}[Brylinski filtration]\label{bry} We define an increasing 
filtration $W\idot V$ on
 $V\in \rep(G)$ by letting $W_kV$ be the direct sum of all
eigenspaces of $\s$ with eigenvalues $\leq k$.
\end{definition}
 
\ab
Further, given a finite dimensional $\g$-module $V$,
and a weight $\mu\in \fh^*$, write $V(\mu)$ for the
corresponding weight space of $V$
(with respect to the Cartan subalgebra $\fh=\g^{e+t}$,
 not~$\h$).

\begin{theorem}[\cite{G2}, Thm. 5.3.1]\label{main_perv}
 If $s=\sum_{\alpha\in R_+}\,{\check\alpha}$ then 
the isomorphisms $\varphi_{_V}:$\linebreak
 $H_s(\p V)\iso V$ (of Proposition \ref{main_perv_prop}(ii))
 can be chosen so that, in addition to 
claims of Proposition
\ref{main_perv_prop}, one has:

\pb{The canonical filtration
$W\idot H_s(\p V)$ goes,  under the isomorphism
$\varphi_{_V}$, to the filtration~$W\idot V$;}

\pb{The 
fixed point decomposition (\ref{fp_decomposition})
corresponds,  under the isomorphism
$\varphi_{_V}$, to the weight decomposition:
$V=\bigoplus_{\mu\in\Y}\;V(\mu)$
with respect to the Cartan subalgebra $\fh.\quad$\qed}
\end{theorem}

\ab To replace equivariant cohomology by the ordinary
cohomology in the Theorem above, note first that,
 for a $G$-module $V$ and for
$\mu=0$ we have $V(0)= V^{\fh}$
(note that the weights of any finite-dimensional
$G$-module  belong to the
root lattice $\Y$).
The filtration $W\idot V$ induces by restriction a filtration
on $V^{\fh}$, and R. Brylinski [Br] proved
\begin{proposition}\label{lim} For any $V\in\rep(G)$ there is
a canonical graded space isomorphism
$\gr^W\idot(V^\fh)\simeq V^{\g^e}.$\quad\qed
\end{proposition}

\ab In particular, for $V=\g$, the adjoint representation,
the Proposition yields a canonical graded space isomorphism
 $\gr(\fh)\simeq\g^e$ (which has been constructed earlier by Kostant),
hence, a graded algebra isomorphism
$\gr({{\mathcal{U}}}\fh)\simeq{{\mathcal{U}}}(\g^e)$.
Thus, passing to associated graded objects 
in Theorem \ref{main_perv}, and using the canonical isomorphisms:
$\gr^W\idot  H_s(\Gr)\simeq H\hdot(\Gr,\C)$ and
$\gr^W\idot   H_s(A)\simeq H\hdot(\Gr,A)\,,\,\forall A\in \PO,$
yields

\begin{corollary}[\cite{G2}, Theorems 1.6.3, 1.7.6]\label{ord}
The isomorphism of functors $\varphi_{_V}:
 H_s(\p V)\iso V$, of Theorem \ref{main_perv},
gives an  isomorphism of tensor functors $\gr(\varphi_{_V}):$
\linebreak
$H\hdot(\Gr, \p V)\iso V\,,
\,\forall V\in \rep(G)$,  such that

\pb{The grading on $H\hdot(\Gr, \p V)$
goes,  under the isomorphism
$\varphi_{_V}$, to the grading on $V$ by the eigenvalues of the
$t$-action.}

\pb{For any $u\in H\hdot(\Gr, \C)$,
the natural action of $u$ on the hyper-cohomology $H\hdot(\Gr,\p V)$
 corresponds,
via $\gr(\varphi_{_V})$ and the isomorphism $\varphi_{_{\mathcal{U}}}:
H\hdot(\Gr,\C) 
\stackrel{\sim}{\longrightarrow}
{{\mathcal{U}}}(\g^e)$ of Corollary \ref{H_G(Gr)1},
to the natural
action of $\varphi_{_{\mathcal{U}}}(u)
\in {{\mathcal{U}}}(\g^e)$ in the  $G$-module $V$.\quad\qed}
\end{corollary}
\ab We have used here an obvious canonical identification
of $V$, viewed as graded space, with $\gr^W\idot {V}$, the associated
graded space corresponding to the filtration $W\idot V$.

\section{Self-extensions of the Regular sheaf}\label{sec_extensions}

\ab Let $\N\subset \g$ be the nilpotent variety in $\g$,
and $e\in\N$ a fixed regular element. The results of \cite{G2} outlined
in the previous section allow to `see' the element $e$,
as well as
its centralizer $\g^e$, in terms of perverse sheaves
on the loop Grassmannian $\Gr=G^\vee(\K)/G^\vee(\oo)$. One of the goals
 of this section is to show how to
reconstruct the whole nilpotent variety $\N$,
not just the principal nilpotent conjugacy class $
\Ad G\cdot e\subset\N$,
from the category of perverse sheaves on $\Gr$.

\subsection{The regular perverse sheaf $\R$.}\label{regular_sheaf}
 Let
 $\hu:=\C[G]^\vee$ be the continuous dual of the Hopf algebra $\C[G]$,
viewed as a topological algebra with respect to
the topology induced by the augmentation. Thus,
 $\hu$ is a topological Hopf algebra
equipped with a canonical continuous perfect Hopf
pairing $\hu \times \C[G] \to \C$.
The pairing yields, cf. the proof
of Lemma \ref{indres},  a  canonical Hopf algebra imbedding
$\jmath: \mathcal{U}\g\into \hu$, thus identifies 
$\hu$ with a
completion of the enveloping algebra $\Ug$.
Similarly to the situation considered in \S\ref{ufrob},
we have a natural equivalence $\rep(\hu)\cong\rep(G)$,
so that the (isomorphism classes of)
simple objects of  $\rep(\hu)$ are labelled by the set $\Y^{++}$.
We will view the left regular representation of the
algebra $\hu$ as a projective (pro-)object in the category $\rep(\hu)$.

\ab Let $\R:=\p(\C[G])$ be  the ind-object of the category
$\PO$
corresponding to the regular representation $\C[G]$,
viewed as an ind-object of the category
$\rep(G)$. Then, $\R^\vee=\p(\C[G]^\vee)=\p(\hu)$,
is the dual  pro-object. Explicitly,
applying the functor $\p$ to the 
 $G$-bimodule direct sum decomposition of the regular representation
on the left (below), 
and using Theorem \ref{tensor}(i) we deduce:
$$\C[G]\;=\;\bigoplus_{\lambda \in\Y^{++}}\;
V_\lambda\otimes_{_\C} V^*_\lambda,\quad
\R = \bigoplus_{\lambda \in\Y^{++}}\;
IC_\lambda\otimes_{_\C} V^*_\lambda
,\quad
\R^\vee=\prod_{\lambda \in\Y^{++}}\;IC_\lambda\otimes_{_\C}
V^*_\lambda
\,.$$
\ab Observe that right translation by an element
$g\in G$ gives a morphism $R_g: \C[G]\to\C[G]$ of {\em left}
$G$-modules. Hence, applying the functor $\p(-)$ we
get, for any $g\in G$, a morphism $R_g: \R\to\R$,
that corresponds to the $g$-action on the factor
$V_\la^*$ in the above decomposition
$\R=\bigoplus\, IC_\lambda\otimes_{_\C} V^*_\lambda.$
The collection of morphisms $R_g\,,\,g\in G,$ satisfies an obvious 
associativity. Therefore,
for any objects $M,N\in D^b(\Gr),$ these morphisms induce,
by functoriality, a $G$-action on the graded vector space
$\Ext\hdot_{_{D^b(\Gr)}}(M\,,\,N\star\R).$
This is the $G$-module structure on the various Ext-groups
that will be considered below.

\subsection{Two Ext-algebras.}\label{Ext1} The first $\Ext$-algebra that we
are going to consider is the space
$\,\Ext_{_{D^b(\Gr)}}\hdot(\R,\R)\,,$ equipped with
Yoneda product. More explicitly, write the ind-object
$\R$ as
$\R=\lim\limits_{\longrightarrow}\,
\R_\alpha$, and accordingly write the dual pro-object as
$\R^\vee=\lim\limits_{\longleftarrow}\,
\R_\alpha^\vee$. Then, we have  
\begin{align}\label{ind_ext}
\Ext\hdot_{_{D^b(\Gr)}}(\R\,,\,\R) &:= \;
\lim\limits_{{\stackrel{\longleftarrow}{\beta}}}
\lim\limits_{{\stackrel{\too}{\alpha}}}
\,\Ext_{_{D^b(\Gr)}}(\R_\beta\,,\,\R_\alpha)\\
&\;\simeq\;
\lim\limits_{\stackrel{\longleftarrow}{\beta}}
\lim\limits_{{\stackrel{\too}{\alpha}}}
\,\Ext_{_{D^b(\Gr)}}\left((\R_\alpha)^\vee\,,\,(\R_\beta)^\vee\right) =:
\Ext\hdot_{_{D^b(\Gr)}}(\R^\vee\,,\,\R^\vee)
\;,\nonumber
\end{align}
where we have used the canonical isomorphism
$\Ext\hdot(\L,\,\M) =\Ext\hdot(\M^\vee,\,\L^\vee)$.

\ab To define the second algebra observe first that
multiplication of functions makes $\C[G]$ a ringed ind-object
of the category $\rep(G)$, that is, the product map
induces the following morphisms
\begin{equation}\label{product_map_m}
m_{_{\C[G]}}:\
\C[G]\otimes\C[G] \too\C[G],\quad\text{resp.},
\quad\mt=\p(m_{_{\C[G]}}):\
\R\star\R\too\R,
\eeq
in the categories $\limind\rep(G)$ and
$\limind\PO$, respectively.

\ab Write $\one=IC_0=\p\C$
for the sky-scraper sheaf 
(corresponding to the trivial 1-dimensional $G$-module)
at the base point of $\Gr$,
and let $\R[i]$ denote  the shift of $\R$ in the derived category.

\ab We define an
 associative graded algebra structure on the vector
space
$\,\Ext\hdot_{_{D^b(\Gr)}}(\one\,,\,\R)
\;:=\;
\bigoplus\nolimits_{\lambda \in\Y^{++}}\;
\Ext_{_{D^b(\Gr)}}\hdot
\left(\one\,,\,
 IC_\lambda\right)\otimes_{_\C} V^*_\lambda
\,$
 as follows. 
Let $x\in \Ext^i_{_{D^b(\Gr)}}(\one,\,\R)= \Hom_{_{D^b(\Gr)}}(\one,\,\R[i])$.
Taking convolution of the identity morphism $\id_\R:\R\to\R$
with $x$, viewed as a  "derived morphism", gives a morphism
$\R=\R\star \one \stackrel{\star x}{\longrightarrow}\R\star \R[i]$.
Given $y \in \Ext^j_{_{D^b(\Gr)}}(\one,\,\R)$, we define
$y\cdot x\in \Ext^{j+i}_{_{D^b(\Gr)}}(\one,\,\R)$ to be the 
composite:
\beq\label{composite}
y\cdot x\;
:\;\;\one\,\stackrel{y}{\longrightarrow}\,\R[j]=(\R\star \one)[j]
\,\stackrel{\star x}{\longrightarrow}\,\R\star \R[i+j]
\,\stackrel{\mt}{\longrightarrow}\,\R[i+j]\,.
\eeq
Similarly, for any ${\mathcal M}\in \PI$, the following maps
make $\Ext\hdot_{_{D^b(\Gr)}}(\one,\,\M\star\R)$  a graded
$\Ext\hdot_{_{D^b(\Gr)}}(\one,\,\R)$-module:
$\;
\one\,\longrightarrow\,\M\star\R[j]
\,\stackrel{\star x}{\longrightarrow}\,\M\star\R\star \R[i+j]
\,\stackrel{\mt}{\longrightarrow}\,$ $\M\star\R[i+j].$

\ab Note  that the complex $\one$ 
is the unambiguously determined direct summand of $\R$,
and the corresponding projection: $\R\,\twoheadrightarrow
\,\one $ induces an imbedding
$\varepsilon_{_{geom}}:\,\Ext\hdot(\one,\,\R)
\,\hookrightarrow\,
\Ext\hdot(\R,\,\R)$. It is easy to check that this  imbedding
becomes an algebra homomorphism, provided the Ext-group on the
right is equipped with the Yoneda product, and the 
Ext-group on the left is equipped   with the product defined in (\ref{composite}).
Dually, there is a 
 map $\varepsilon_{_{geom}}:
\Ext\hdot(\R^\vee,\,\one)
\;\hookrightarrow\;
\Ext\hdot(\R^\vee\,,\,\R^\vee)$  induced by
the imbedding: $\one\,\hookrightarrow\,
\R^\vee$.

\subsection{Main result.}
 The adjoint $\g$-action on $\N$ makes the coordinate ring
$\C[\N]$
a locally finite $\Ug$-module, hence an $\hu$-module.
Let 
$\hu{\st} \C[\N]$ be 
the corresponding cross-product algebra. 
There is an obvious algebra imbedding
$\varepsilon_{_{alg}}:\,\C[\N]\,\hookrightarrow\,
\hu{\st} \C[\N]\,.$
We put a grading
on $\hu{\st} \C[\N]$ by taking the natural grading on
the subalgebra $\C[\N] \sset \hu{\st} \C[\N]$,
and by placing $\hu$ in grade degree zero. 

\begin{theorem}\label{main1} 
There are  natural  $G$-equivariant
graded algebra isomorphisms ${\Psi}$ and  $\psi$
making the following
diagram commute:
\pagebreak[3]
$$
\Ext\hdot_{_{D^b(\Gr)}}(\R^\vee,\,\one)=
\Ext\hdot_{_{D^b(\Gr)}}(\one,\,\R)
\hskip 1mm\stackrel{\varepsilon_{_{geom}}}{\longrightarrow}\hskip 1mm
\Ext\hdot_{_{D^b(\Gr)}}(\R\,,\,\R)=
\Ext\hdot_{_{D^b(\Gr)}}(\R^\vee\,,\,\R^\vee)$$
$\hphantom{x}\hskip 48mm 
\big|\!\!\downarrow\psi \hskip 21mm\qquad
\qquad\big|\!\!\downarrow\Psi$
$$
\C[\N]\hskip 10mm\stackrel{\varepsilon_{_{alg}}}{\longrightarrow}
\hskip 10mm\hu{\st}\C[\N]\;.
$$
\end{theorem}

\ab The rest of this section is devoted to the proof of the Theorem.

\subsection{Some general results.}
According to the well-known results of
Kostant, the centralizer in $G$ (an adjoint
group) of the principal nilpotent element $e$ is  a connected unipotent subgroup
$G^e\sset G$
 with Lie algebra $\g^e$.
Furthermore, the $\Ad{G}$-conjugacy class
of $e$ is known to be the open dense subset  $\Oe \sset \N$
formed
by the regular nilpotents. Thus,  we get: $\Oe =G^e\backslash G$,
where  $\Ad{G}$-action is viewed as a
right
action.
Moreover, Kostant has shown in [Ko] that the natural imbedding: $\Oe
\hookrightarrow \N$ induces an isomorphism of the rings of regular
functions.
Thus we obtain a chain of  natural algebra isomorphisms,
where the superscripts stand for invariants under left translation
\beq\label{normality}
\C[G]^{\gee}\;=\;
\C[G]^{G^e}\;=\;\C[G^e\backslash G]\;=\;
\C[\Oe]\;=\;\C[\N]\,.
\eeq

\ab Next we consider
 the
vector space
 $\Hom_{_{{\mathsf{cont}}}}(\C[G],\,\C[G])$
 of $\C$-linear {\it continuous} maps $\C[G]\to\C[G]$,
that is a pro-ind limit
of finite dimensional $\Hom_{_{\C}}$-spaces,
defined the same way as we have earlier
defined the $\Ext$-spaces between ind-objects, see \eqref{ind_ext}. The space
 $\Hom_{_{{\mathsf{cont}}}}(\C[G],\,\C[G])$
has a natural  $G$-action by conjugation; furthermore, 
it has 
 the structure of topological algebra via composition. We have
\begin{lemma}\label{topology} There is a natural
 topological algebra isomorphism 
$$
\bigl(\Hom_{_{{\mathsf{cont}}}}(\C[G]\,,\,\C[G])\bigr)^{\gee}
\;\simeq\;\hu{\st}\C[\N]\,.$$
\end{lemma}

\proof Write $\dd(G)$ for the algebra of regular
algebraic differential operators on the group $G$. The action of
any such differential
operator $u\in \dd(G)$ clearly gives a {\it continuous}
$\C$-linear map
$u: \C[G] \to\C[G]$. This way one obtains a
$G$-equivariant imbedding:
$\dd(G)\;\hookrightarrow\;\Hom_{_{{\mathsf{cont}}}}(\C[G],\,\C[G])$
with dense image. 
View
the algebra $\Ug$ as  left invariant differential operators
on $G$, and the algebra $\C[G]$ as multiplication-operators.
Then the algebra $\dd(G)$ is  isomorphic to the
cross-product:
$\dd(G)\simeq \Ug{\st}\C[G]$. One can show
that the composition:
$\,\Ug {\st}\C[G] \iso\dd(G)\;
\hookrightarrow\;\Hom_{_{{\mathsf{cont}}}}(\C[G],\,\C[G])$
extends by continuity to a $G$-equivariant
{\it topological algebra} isomorphism
$\hu{\st}\C[G] \iso
\Hom_{_{{\mathsf{cont}}}}(\C[G],\,\C[G])$. 
Observe that the group $G$ acts {\it trivially}
on the space $\Ug$ formed by left-invariant differential operators,
and acts on $\C[G]$ via left translations.
Hence, taking $G^e$-invariants
on each side of the isomorphism above,  we obtain
\begin{align*}
\Bigl(\Hom_{_{{\mathsf{cont}}}}(\C[G]\,,\,\C[G])\Bigr)^{\g^e}\;&=\;
\bigl(\hu{\st}\C[G]\bigr)^{\g^e}\;=\;
\hu{\st}\bigl(\C[G]^{\g^e}\bigr)\stackrel{\eqref{normality}}{\longeq}
\hu{\st}\C[\N]\,.\qquad\qed
\end{align*}

\begin{notation}\label{Hom^*}
Given two graded objects $L_1, L_2$, we set
$\Hom\hdot(L_1, L_2) := \bigoplus_i\;\Hom^i(L_1, L_2)$,
where $\Hom^i(L_1, L_2)$ stands for the space of morphisms
shifting the grading by $i$.
\end{notation}

\ab 
Next, we remind the reader
that,
for any $\L \in D^b(\Gr)$, the hyper-cohomology 
group $H\hdot(\Gr,\L)$ has a natural structure 
of  graded $H\hdot(\Gr,\C)$-module. Further, for any
$\L_1, \L_2 \in D^b(\Gr)$, there is a functorial
linear map of graded vector spaces:
\beq\label{G1map}
\Ext\hdot_{_{\!D^b(\Gr)}}(\L_1, \L_2)\;
\longrightarrow\;
\Hom\hdot_{_{H(\Gr)}}\left(H\hdot(\Gr,\L_1)\,,\,
H\hdot(\Gr,\L_2)\right)\,.
\eeq
\ab
The main result of [G1], in the special case of a $\C^*$-action on $\Gr$
implies 
\begin{proposition}\label{GG1} If $\L_1, \L_2\in D^b(\Gr)$
are semisimple perverse sheaves
constructible relative to a Bialinicki-Birula stratification
of $\Gr$ (cf. \cite{G1} for more details),
Then the map (\ref{G1map}) is an isomorphism.\hfill\qed
\end{proposition}

\ab We note that
since any $G^\vee({\mathcal O})$-equivariant perverse sheaf on $\Gr$ is constructible
relative to the Schubert stratification, cf. [PS] or \S8 below,
the map (\ref{G1map}) is an isomorphism
for any $\L_1, \L_2\in \PO$.

\ab  The action of the semisimple element
$t\in\h$ (from the principal $\sll$-triple)
puts a grading on the underlying vector space of any representation
$V\in \rep(G)$.  Further, by Corollary \ref{H(Gr)}
 we have an algebra isomorphism
$H\hdot(\Gr)\simeq \cU({\gee})$.
Hence, 
from Corollary 
\ref{ord} and Proposition \ref{GG1}
we deduce
\begin{corollary}\label{invariants}
  For any
$V_1,\,V_2\in\rep(G)\,,$ there is a functorial
isomorphism
$$
\Ext\hdot_{D^b(\Gr)}(\p V_1\,,\,\p V_2)
\,\;\;\iso\,\;\;\Hom\hdot_{_{\cU({\gee})}}(V_1\,,\,V_2)
\;=\;\Bigr(\Hom\hdot_{_{\C}}(V_1\,,\,V_2)\Bigr)^{\gee}.\quad\square
$$
\end{corollary}

\subsection{Proof of Theorem \ref{main1}:} 
Let $r\in \Hom_{_{\!D^b(\Gr)}}(\C_{_\Gr},\one)=\Ext^0_{_{\!D^b(\Gr)}}(\C_{_\Gr},\one)$
be the natural restriction morphism. Composing with $r$ yields a canonical
 map
\begin{equation}\label{GrG1}
\Ext\hdot_{_{\!D^b(\Gr)}}(\one\,,\,\L)\too
\Ext\hdot_{_{\!D^b(\Gr)}}(\C_{_\Gr},\L)=H\hdot(\Gr,\L),\quad \forall \L\in D^b(\Gr)\,.
\eeq
This map can also be identified
with
the map of Corollary \ref{invariants}, 
for $\one=\p \C$ and $\L=\p V_2.$ Hence, 
 Corollary \ref{invariants}  says that \eqref{GrG1} gives an isomorphism
$$
\Ext\hdot_{_{\!D^b(\Gr)}}(\one\,,\,\L)\iso
\Hom\hdot_{_{H(\Gr)}}\left(\C\,,\,H\hdot(\Gr,\L)\right)
\;=\;H\hdot(\Gr,\L)^{{\gee}},\quad \forall \L\in\PO.
$$
The isomorphism above holds, in particular,  for $\L:=\R$, an ind-object of 
$\PO$.
We use this isomorphism to define a linear isomorphism
$\psi$
as the following composite:
\begin{align}\label{H(u)_computation}
\xymatrix{
\Ext\hdot_{_{D^b(\Gr)}}(\one\,,\,\R)
\ar[rr]^<>(0.5){\text{Cor. \ref{invariants}}}_<>(0.5){\sim}&&
H\hdot(\Gr,\,\R)^{{\gee}}
\ar[rrr]^<>(0.5){\text{Cor. \ref{ord} for}\;\,\R=\p\C[G]}
_<>(0.5){\sim}&&&
\C[G]^{G^e}\ar@{=}[r]^<>(0.5){\eqref{normality}}&\C[\N].}
\end{align}
\ab 
To complete the proof of the theorem, it suffices to verify
that the chain of isomorphisms (\ref{H(u)_computation})
transports the above defined algebra structure on
$\Ext\hdot(\one,\R)$ to the standard algebra structure on
 $\C[G]^{G^e}=\C[\Oe]\simeq\C[\N]$. To this end, we start
with the
 canonical identification
 $\C[G]=H\hdot(\Gr,\R).$ Since $H\hdot(-)$ is a fiber functor on
$\PO$, we also have
$H\hdot(\Gr, \R\star\R)=\C[G]\otimes\C[G]$.
By construction, the natural imbedding $\C[G]^{{G^e}}\into \C[G]$
corresponds, via the identification
$\Ext\hdot(\one,\R) = \C[G]^{{G^e}}$ in \eqref{H(u)_computation},
to the morphism \eqref{GrG1}.
Similarly, writing $G^e_{_{\textsf{diag}}}
\sset{G^e}\times{G^e}$ for the diagonal, we may identify 
$\bigl(\C[G]\otimes\C[G]\bigr)^{G^e_{_{\textsf{diag}}}}\sset
\C[G]\otimes\C[G]$ with 
$\Ext\hdot(\one,\R\star\R)\sset H\hdot(\Gr,\,\R\star\R).$

\ab Observe further that
the multiplication  $m_{_{\C[G]}}$ in the coordinate ring $\C[G]$
is recovered from the morphism
$\mt: \R\star\R \to \R$, cf. \eqref{product_map_m},
 as the induced morphism of hyper-cohomology:
\beq\label{mult2}
H\hdot(\Gr, \R) \otimes H\hdot(\Gr, \R)
\;{=}\;H\hdot(\Gr\times \Gr, \R\boxtimes\R)
\iso
H\hdot(\Gr, \R\star\R)\;\stackrel{\mt}{\longrightarrow}\;
H\hdot(\Gr, \R).
\eeq 
Restricting the map $m_{_{\C[G]}}$, resp., $\mt$,
to $G^e$-invariants yields the corresponding map in the
top, resp., bottom,
row of the following  diagram:
$$
\xymatrix{
\C[G]^{G^e}\otimes\C[G]^{G^e}\enspace 
\ar@{^{(}->}[r]^<>(0.5){i}\ar@{=}[d]&
\bigl(\C[G]\otimes\C[G]\bigr)^{G^e_{_{\textsf{diag}}}}
\ar[r]^<>(0.5){m_{_{\C[G]}}}\ar@{=}[d]&\C[G]^{G^e}\ar@{=}[d]
\\
\Ext\hdot(\one,\R) \otimes\Ext\hdot(\one,\R)
\ar[r]^<>(0.5){\eqref{composite}}&\Ext^{\hhdot}(\one,\R\star\R)
\ar[r]^<>(0.5){\mt}&\Ext\hdot(\one,\R).
}
$$
As has been explained, the identifications we have made
insure that the diagram commutes. Further, it is clear
that the composite map in the top row of the diagram is
the multiplication map
$\C[G]^{G^e}\otimes\C[G]^{G^e}\to\C[G]^{G^e}$.
It follows that the latter map
corresponds, geometrically, to
the composite map in the bottom row.
This proves that the map $\psi$ of the Theorem is an
algebra
isomorphism.

\ab We now similarly construct the map $\Psi$. By Corollary \ref{invariants},
we have the
following isomorphisms
\beq\label{main_proof}
\Ext\hdot_{_{D^b(\Gr)}}(\R\,,\,\R)\;=\;
\Ext\hdot_{D^b(\Gr)}\bigl(\p(\C[G])\,,\,\p(\C[G])\bigr)
\;
=\;\bigl(\Hom_{_{{\mathsf{cont}}}}(\C[G]\,,\,\C[G])\bigr)^{\gee}.
\eeq
Hence, from Lemma \ref{topology} we deduce
$$
\Ext\hdot_{_{D^b(\Gr)}}(\R,\,\R)\;=\;
\left(\Hom\hdot_{_{{\mathsf{cont}}}}(\C[G]\,,\,\C[G])\right)^{\gee}\;=\;
\;=\;\hu{\st}\C[\N]\,.
$$
Compatibility of the maps that we have constructed above with algebra, resp.
module, structures is
verified in a similar way,
as we did for $\Ext\hdot(\one,\,\R)$. We leave details to the reader.
\qed

%

\subsection{Equivariant version.} Write
 $\Ext\hdot\GO$ for
the Ext-groups in  $D^b\GO(\Gr)$,
$\Gd(\oo)$-equivariant derived category of constructible complexes
in the sense of [BL]. These Ext-groups  have
canonical $H\hdot\GO(pt)$-module structure, cf. \eqref{point}.
Clearly, $\PO$ is an abelian subcategory of 
$D^b\GO(\Gr)$, hence $\R$ is an ind-object of $D^b\GO(\Gr)$. 
As in (\ref{composite}), the map $\mt: \R\star\R\to \R$ makes
$\Ext\hdot\GO(\one\,,\,\R)$
into a graded  algebra, and we have
the following
equivariant analogue of Theorem \ref{main1}.

\begin{theorem}\label{Ext_G}
There are  natural  algebra isomorphisms ${\Psi}$ and  ${\psi}$
making the following
diagram commute:
\pagebreak[3]
$$
\Ext\hdot\GO(\one,\,\R)
\hskip 1mm\stackrel{\varepsilon_{_{geom}}}{\longrightarrow}\hskip 1mm
\Ext\hdot\GO(\R\,,\,\R)$$
$\hphantom{x}\hskip 48mm 
\big|\!\!\downarrow\psi \hskip 21mm\qquad
\qquad\big|\!\!\downarrow\Psi$
$$
\C[\g^*]\hskip 10mm\stackrel{\varepsilon_{_{alg}}}{\longrightarrow}
\hskip 10mm\hu{\st}\C[\g^*]\;.
$$
\end{theorem}

\begin{remark}
It is tempting to use a kind of {\it delocalized}
equivariant cohomology, see [BBM], in order to be able to replace the algebra
$\hu{\st} \C[\g]$ above (where
we identify $\g$ with $\g^*$ via an invariant form) by the algebra
$\hu{\st} \C[G]$. Recall
that both 
$\hu$ and $\C[G]$ have natural structures of Hopf algebras,
topologically dual to each other, i.e.,
$\hu$ is an inverse limit, while $\C[G]$ is a direct limit.
The algebra $\hu{\st} \C[G]$
is the Drinfeld's {\it double} Hopf algebra. 
 $\quad\lozenge$\end{remark}

\ab In order to begin the proof of Theorem \ref{Ext_G}
we need to introduce more notation. Given any 
vector space $E$ we will write $\underline{E}$ for the 
quasi-coherent sheaf (trivial
vector bundle) on $\h^*/W$ with geometric fiber $E$.
Thus, $\Gamma(\h^*/W,\,\underline{E}) = \C[\h^*]^W \otimes E$.

\ab The result below is an equivariant version of the last
statement of Corollary \ref{ord}, applied to
the regular perverse sheaf $\R$. 

\begin{lemma}\label{H_G(Gr)}
There is a natural $\C[\h^*]^W$-module isomorphism
$$H\hdot\GO(\Gr\,,\,\R)\;\simeq\;\C[\h^*]^W \otimes_{_{\C}} \C[G] =
\Gamma(\h^*/W\,,\,\underline{\C[G]})\;;$$
The canonical $H\hdot\GO(\Gr)$-module structure on
$H\hdot\GO(\Gr\,,\,\R)$ corresponds, under the isomorphism 
of Lemma \ref{H_G(Gr)1},
to the natural left $\Gamma(\h^*/W\,,\,{\mathcal{U}}\z)$-action on 
$\Gamma(\h^*/W\,,\,\underline{\C[G]})$.
\qed
\end{lemma}

\ab With the results of [G2] mentioned in \S6, the
proof of the Lemma is straightforward and will be omitted.\hfill\qed
\medskip

{\sc{Proof of Theorem \ref{Ext_G}.}}\quad
 We will only establish the isomorphism $\psi$;
once it is understood, the
construction
of $\Psi$ is entirely similar to that in the proof of Theorem \ref{main1}.

\ab Using
Lemma \ref{H_G(Gr)1},
we  perform the following
calculation, similar to (\ref{H(u)_computation}):

\begin{align}\label{equivariant_computation}
\Ext\hdot\GO(\one\,,\,\R) &=
\Hom\hdot_{_{H\GO(\Gr)}}\left(H\hdot\GO(\one)\,,\,H\hdot\GO(\Gr,\,\R)\right)
\quad\mbox{(by [G1])}\nonumber\\
&=\Hom\hdot_{_{H\GO(\Gr)}}\left(H\hdot\GO(pt)\,,
\,H\hdot\GO(\Gr,\,\R)\right)\nonumber\\
\mbox{(by Lemma \ref{H_G(Gr)})}\quad
&=\Hom\hdot_{_{\Gamma(\h^*/W\,,\,{\mathcal{UG}})}}\left(\C[\h^*]^W\,,\,\C[\h^*]^W 
\otimes_{_{\C}} \C[G]\right)\\
&=\Gamma\Bigl(\h^*/W\,,\,{\mathcal H}\mbox{\it om}_{_{{\mathcal{UG}}}}(\underline{\C}
\,,\,\underline{\C[G]})\Bigr)\nonumber\\
&=\Gamma\left(\h^*/W\,,\,\,\underline{\C[G]}^{^{{\mathcal{UG}}}}\right)
=\Gamma(\greg/\Ad G\,,\,\,\underline{\C[G]}^{^{{\mathcal{UG}}}})\;.\nonumber
\end{align}

The last expression may be identified naturally with the algebra of regular
functions on the total space of the canonical fibration $p :\greg \to\greg/\Ad G$,
because for any $x\in\greg$, we have  
$p^{-1}(p(x)) \simeq G^x\backslash G$. We see that the algebra in question
equals $\C[\greg]$. Since the complement $\g^*\smallsetminus \greg$
is known to have codimension $\geq 2$ in $\g^*$, we conclude that
$\C[\greg]=\C[\g^*]$. The result is proved.\qed

\subsection{A fiber
 functor on perverse sheaves.}\label{fiber_fun1} In this subsection, we will make a link
of our results with  general Tannakian formalism.

\ab Let $\scrc$ be an abelian category, which is a (right)
{\em module category} over the tensor category $\rep(G)$.
This means that we are given an exact bifunctor
$\scrc \times\rep(G)\too\scrc\,,\,(M,V) \mto M\otimes V,$
satisfying a natural associativity constraint:
$(M\otimes V)\otimes V'\cong M\otimes(V\otimes V').$

\ab Let  $\scrc$ be a module category over $\rep(G)$.
View $\C[G]$ as an ind-object of $\rep(G)$ and,
given $M\in \scrc$, form an ind-object $M\otimes \C[G]$
in $\scrc$. Now, fix $M,N\in\scrc$.
As we have explained  (in a special case) at the end of section 
\ref{regular_sheaf},
the action of $G$ on $\C[G]$ by {\em right} translations
gives rise to a natural $G$-module structure
on the vector space $\Ext\hdot_\scrc(M\,,\,N\otimes\C[G])$.

\ab Now, given  $L\in\scrc$ and $V\in\rep(G)$, 
we may apply the above construction to the objects
$N=L$ and $N=L\otimes V,$ respectively. The proof of the following
`abstract nonsense' result is left for the reader.

\begin{lemma}\label{nonsense} For any $M,L\in\scrc$, there is a natural functorial
$G$-module isomorphism
$$\Ext\hdot_\scrc\bigl(M\,,\,(L\otimes V)\otimes\C[G]\bigr)\simeq
V\otimes \Ext\hdot_\scrc\bigl(M\,,\,L\otimes\C[G]\bigr)\,
\quad\forall V\in\rep(G).\qquad\Box
$$
\end{lemma}

\ab By  Gaitsgory theorem \ref{denis}, convolution of perverse sheaves gives
an exact bifunctor
$\pe(\Gr)\times \PO\too\pe(\Gr)\,,\,\M,\A\mapsto \M\star \A .$
This way, the category $\pe(\Gr)$ becomes a module category
over $\PO$. 
Transporting the module structure
by means of  Satake equivalence $\p:\rep(G)\iso \PO$,
we thus make $\pe(\Gr)$ a module category over
$\rep(G)$.
Applying Lemma \ref{nonsense} to the module category $\pe(\Gr)$
 we get, for any  $\M,\L\in \pe(\Gr),$ the following  natural $G$-module isomorphism
\beq\label{nonsense2}
\Ext\hdot_{D^b(\Gr)}\bigl(\M\,,\,\L\star\p V\star \R\bigr)
\simeq
V\otimes \Ext\hdot_{D^b(\Gr)}(\M,\L\star\R)\,,\quad
\forall V\in\rep(G).
\eeq
\begin{remark}
Having in mind further applications
of isomorphism \eqref{nonsense2}, we have written the Ext-groups in the triangulated
category $D^b(\Gr)$ rather than in the abelian category $\pe(\Gr)$. These
two Ext-groups are known (due to Beilinson) to be actually the same.
The Beilinson's result is not however absolutely necessary for
justifying
\eqref{nonsense2}; indeed, the formal nature of the setup of Lemma \ref{nonsense}
makes it applicable, in effect, to a  more general
setup of a triangulated `module category'.$\quad\lozenge$
\end{remark}

\ab Recall now the graded algebra $\Ext_{D^b(\Gr)}\hdot(\one,\R)$,
that comes  equipped with a natural $G$-action.
The isomorphism $\psi: \Ext_{D^b(\Gr)}\hdot(\one,\R)
\iso \C[\N],$ of Theorem \ref{main1}, induces
an equivalence of categories
$$\psi_*:\ 
\op{Mod}^{G\times\C^*\!}\bigl(\Ext_{D^b(\Gr)}\hdot(\one,\R)\bigr)
\iso \op{Mod}^{G\times\C^*\!}\bigl(\C[\N]\bigr)=\Coh^{G\times\C^*}(\N),
$$
where the category of
$G\times\C^*$-equivariant coherent sheaves
on $\N$ is identified  naturally with $\op{Mod}^{G\times\C^*\!}\bigl(\C[\N]\bigr),$
see \ref{cohG}.

\ab 
Next we observe that, for any $\L\in\pe(\Gr)$ one can define,
by modifying appropriately formula \eqref{composite},
a natural pairing 
$$\Ext_{D^b(\Gr)}^i(\one,\,\L\star \R)\; \otimes\;
\Ext_{D^b(\Gr)}^j(\one,\R)\too
\Ext_{D^b(\Gr)}^{i+j}(\one,\,\L\star\R).
$$
This pairing makes the space
$\Ext_{D^b(\Gr)}\hdot(\one,\,\L\star \R)$ a graded
$\Ext_{D^b(\Gr)}\hdot(\one,\R)$-module.
We also have the $G$-action on the Ext-groups
involved.
Therefore, we obtain an exact functor
\beq\label{fib_fun1}
\pe(\Gr)\too
\op{Mod}^{G\times\C^*\!}\bigl(\Ext_{D^b(\Gr)}\hdot(\one,\R)\bigr),\quad
\L\mto\Ext_{D^b(\Gr)}\hdot(\one,\,\L\star \R).
\eeq

\ab We now use the category equivalence  $\psi_*$, mentioned above,
together with the Satake equivalence, and form 
 the following composite  functor:
\beq\label{fib_fun2}
{\mathbb{S}}: \ \rep(G)\stackrel{\p}\too
\pe(\Gr)\stackrel{\eqref{fib_fun1}}\too
\op{Mod}^{G\times\C^*\!}\bigl(\Ext_{D^b(\Gr)}\hdot(\one,\R)\bigr)
\stackrel{\psi_*}\too \Coh^{G\times\C^*}(\N).
\eeq

\ab The canonical isomorphism in \eqref{nonsense2} then translates
into the following result.

\begin{proposition}\label{bpsi} The functor ${\mathbb{S}}$
is isomorphic to the functor $V\mto V\otimes \oo_\N,$
that assigns to any $V\in \rep(G)$ the free
$\oo_\N$-sheaf  with fiber $V$ (and equipped with the tensor product
$G$-equivariant structure).\qed
\end{proposition}

\section{Wakimoto sheaves}\label{sec_wakimoto}

\subsection{} We have shown in the previous section that the coordinate
ring $\C[\N]$, together with the $\Ad G$-action on it, can be
reconstructed as an $\Ext$-group between certain
perverse sheaves on the loop Grassmannian.

\ab The goal of this section is to give a similar construction
for $\NN=G\times_{\!_B} \n$, the Springer resolution of $\N$,
see e.g.  [CG, ch.3] for a survey.
 One obstacle for doing so is that $\NN$ is {\it not} an
affine
variety, and it is known that $\C[\NN]=\C[\N]$. Therefore,
the variety $\NN$ is not determined by the ring of its global regular
functions. Instead, we will consider the affine cone over a
kind of Pl\"ucker imbedding of
$\NN$. Our main result shows how to reconstruct
 the (multi)-homogeneous coordinate
ring of that cone as an $\Ext$-algebra of perverse
sheaves on the
loop Grassmannian.

\subsection{The affine flag manifold.}\label{baf}  Let $\Td\sset \Bd$ denote
the
maximal torus and the Borel 
subgroup in $\Gd$ corresponding
to our choice of positive roots. We 
write $\;\I=\{f\in\Gd(\oo)\;|\; f(0)\in \Bd\}$ for the
corresponding Iwahori (=affine Borel) subgroup in $\Gd(\K)$.
We let $\baf=\Gd(\K)/\I\,,$
be the {\it affine Flag variety}; it has a natural
 ind-scheme structure. Since $\I\subset \Gd(\oo)$, the projection
$\vp: \Gd(\K)/\I\onto\Gd(\K)/\Gd(\oo)$ gives
a smooth and proper morphism of ind-schemes
$\vp: \baf \onto \Gr$, whose fiber is isomorphic
to the finite dimensional flag manifold $\Gd/\Bd$.

\ab The left $G(\K)$-action on
$\baf$ gives rise to the following {\sf{convolution diagram}}:
\beq\label{baft_conv}
G^\vee(\K)\times_{_\I}\baf\stackrel{a}\too\baf,\quad (g, x) \mto gx.
\eeq
Given $\A,\M\in D^b_{_\I}(\baf)$, 
one defines, using  $\I$-equivariance of
$\M$, (cf.  e.g., [G2], [MV]  or [Ga]), an object $\A\tilde{\boxtimes} 
{{\mathcal{M}}}\in D^b\bigl(G^\vee(\K)\times_{_\I}\baf\bigr)$.
The assignment $(\A,\M)\mto \A\star \M:=a_*(\A\tilde{\boxtimes} 
{{\mathcal{M}}})$ gives the standard convolution
 bifunctor
$\star: D^b_{_\I}(\baf)\times D^b_{_\I}(\baf)\too D^b_{_\I}(\baf).$

\ab  Each $\I$-orbit on either $\baf$ or $\Gr$ is isomorphic to a finite dimensional
vector space $\C^n$. Moreover, each $\I$-orbit on 
 $\Gr$,  contains a unique
coset  $\lambda\cdot\Gd(\oo)/\Gd(\oo)$, where
$\lambda \in \Hom(\C^*,\Td)$ is viewed
as an element of $\Gd(\K)$. This way one gets a natural bijection
between 
the set of $\I$-orbits in $\Gr$ and the set $\Hom(\C^*,\Td)=\q$.
We will write $\Gr_\lambda$ for the $\I$-orbit
corresponding to a weight
$\lambda\in \q$. Such an $\I$-orbit, $\Gr_\lambda$, is open dense in
the $\Gd(\oo)$-orbit of $\lambda\in \Gr,$ if and only if 
$\lambda \in \q^{++}$ is a dominant weight.
Similarly, each  $\I$-orbit in $\baf$ contains a unique
coset $y\cdot\lambda\cdot\I/\I\,,$ where
$y\cdot\lambda\in \waf.$
We let 
$j_w: \FI_w \into \baf$ denote the corresponding 
$\I$-orbit imbedding. The orbits form a
stratification: $\baf = \coprod_{w\in \waf}\;\FI_w$.

\begin{notation}
  For any $w\in \waf$ we set: $\,{{\mathcal{M}}}_w
:=(j_w)_!\C_w[\dim \FI_w]\,$ and
$\,{{\mathcal{M}}}_w^\vee:=(j_w)_*\C_w[\dim \FI_w],\,$ where
$\C_w$ stands for the constant sheaf on the cell $\FI_w$.
\end{notation}
\ab The imbedding $j_w$ being affine, it is known that
${{\mathcal{M}}}_w\,,\,{{\mathcal{M}}}_w^\vee \in\BI.$ 
It is known further that, for any $y,w\in\waf,$  such that
$\ell(y)+\ell(w)=\ell(y\cdot w)$, one has:
\begin{equation}\label{conv_additive}
{{\mathcal{M}}}_y \starb {{\mathcal{M}}}_w =
{{\mathcal{M}}}_{yw},\quad
{{\mathcal{M}}}^\vee_y \starb {{\mathcal{M}}}^\vee_w=
{{\mathcal{M}}}^\vee_{yw}\,,\quad\mbox{and}\quad
{{\mathcal{M}}}^\vee_w\starb {{\mathcal{M}}}_{w^{-1}} =
{{\mathcal{M}}}_e\,,
\end{equation}
where $e\in \waf$ denotes the identity.

\ab
The following result will play a crucial role in our approach.

\begin{proposition}\label{wakimoto}{{\sf{(I. Mirkovi\'c)}}}$\;$
  For any $w,y\in \waf$ we have   \;\vi
${{\mathcal{M}}}^\vee_y\star{{\mathcal{M}}}_w\in \BI.$ 

\ab  \vii $\supp\bigl({{\mathcal{M}}}^\vee_y\star{{\mathcal{M}}}_w\bigr)
=\overline{\baf}_{yw}$, moreover: 
$\bigl({{\mathcal{M}}}^\vee_y\star{{\mathcal{M}}}_w\bigr)\Big|_{\baf_{yw}}=
\C_{yw}[-\dim\baf_{yw}].$
\end{proposition}

\proof To prove \vi$\;$, fix any $w\in\waf$.
According to the definition of convolution, for any 
complex ${{\mathcal{M}}}$ supported in $\overline{\FI}_w$,
and any $\A\in \BI$, we
have  
$\A\starb {{\mathcal{M}}} = \overline{a}_*(\A
\tilde{\boxtimes} {{\mathcal{M}}})$. Here
$\overline{a}_*=\overline{a}_!$ is the direct image with respect to a 
natural {\it proper} action-morphism
$\,\overline{a} : \Gd(\K) \times_{\I} \overline{\FI}_w \too \baf\,,$
and "$\,\tilde{\boxtimes}\,$" stands for a twisted version of
 external tensor product that has been already mentioned earlier.

\ab  In the special case
${{\mathcal{M}}}=(j_w)_!\C_w[\dim \FI_w]$, the same result
is obtained by replacing $\overline{a}$ by the {\it non-proper}
map $\,a : \Gd(\K) \times_{\I} {\FI}_w \too \baf\,,$
and taking $a_!(\A\tilde{\boxtimes} \C_w[\dim \FI_w])$.
But the morphism $a$ being {\it affine},  by [BBD] one has:
$\,a_!\left(D^{\leq 0}(\Gd(\K) \times_{\I} {\FI}_w\right)\, \subset\,
D^{\leq 0}(\baf)$. Hence, $\A\starb {{\mathcal{M}}}_w\in
D^{\leq 0}(\baf)\,,$ for any $\A\in \BI$.
Dually, one obtains:
${{\mathcal{M}}}_y^\vee\starb \A\in
D^{\geq 0}(\baf)\,,$ for any $\A\in \BI$.
It follows that:
${{\mathcal{M}}}_y^\vee\starb {{\mathcal{M}}}_w\in
D^{\leq 0}(\baf)\cap D^{\geq 0}(\baf)$.
Thus, ${{\mathcal{M}}}_y^\vee\starb {{\mathcal{M}}}_w$
is a perverse sheaf, and part \vi$\;$ is proved.

\ab  To prove (ii), 
let $K_{_\I}(\baf)$ denote the Grothendieck group of
the  category $\DI(\baf)$. The 
 classes 
$\,\big\{[{{\mathcal{M}}}_{x}]\big\}_{x\in \waf}\,$
form a natural $\Z$-basis of $K_{_\I}(\baf).$
The
convolution on 
$\DI(\baf)$ makes $K_{_\I}(\baf)$ into an associative 
ring, and the assignment: $x \mapsto [{{\mathcal{M}}}_{x}]$
is well-known (due to Beilinson-Bernstein, Lusztig, MacPherson, and others)
 to yield a ring isomorphism
$\Z[\waf]\iso K_{_\I}(\baf)$, where $\Z[\waf]$ denotes the group
algebra of the group $\waf$. Further,
in $K_{_\I}(\baf)$,  one has:
$[{{\mathcal{M}}}_{x}]=[{{\mathcal{M}}}^\vee_{x}].$
Hence the class 
$\big[{{\mathcal{M}}}_y^\vee\starb {{\mathcal{M}}}_w\big]$
corresponds,  under the 
isomorphism above, to the element $y\cdot w\in \waf\subset
\Z[\waf],$ in other words, in
$K_{_\I}(\baf)$ we have an equation:
$\big[{{\mathcal{M}}}_y^\vee\starb {{\mathcal{M}}}_w\big]
=[{{\mathcal{M}}}_{yw}]$.
But since $K_{_\I}(\baf)=K(\PI)$ may be identified with the
 Grothendieck group of
the  abelian category $\PI$, and since
${{\mathcal{M}}}_y^\vee\starb 
{{\mathcal{M}}}_w$ is a perverse sheaf by (i), the latter
equation 
in $K(\PI)$ yields:
$\supp\bigl({{\mathcal{M}}}_y^\vee\starb {{\mathcal{M}}}_w\bigr)
=\supp{{\mathcal{M}}}_{yw}=
\overline{\baf}_{yw}.$
\qed

\subsection{Wakimoto sheaves.} Given any $\lambda \in \q$, one can find a pair of dominant
weights $\mu,\nu \in \q^{++}$ such that $\lambda=\mu-\nu$.
Imitating Bernstein's well-known
construction of a large commutative
subalgebra in the affine Hecke algebra, see e.g. [CG, ch.7],
I. Mirkovi\'c introduced the following objects 
\begin{equation}\label{wak_def}
\W_\lambda\; := {{\mathcal{M}}}^\vee_\mu\starb 
{{\mathcal{M}}}_{-\nu}\,,
\end{equation}
 where $\mu,\nu$ are viewed as elements
of $\Y\subset \waf$.
Now, a standard argument due to Bernstein implies
that the above definition of $\W_\lambda$ is independent of the choice of  presentation:
$\lambda=\mu-\nu$. Indeed, if $\mu',\nu' \in \q^{++}$ is another pair
such that
$\lambda=\mu'-\nu'$, then in $\waf$ we have   $\ell(\mu+\nu')=
\ell(\mu) + \ell(\nu')$, $\ell(\mu'+\nu)=
\ell(\mu') + \ell(\nu)$, and $\mu+\nu'=\mu'+\nu$.
Independence of the presentation now follows from formulas
(\ref{conv_additive}). Further, Proposition \ref{wakimoto}(i)
 guarantees that $\W_\lambda\in \PI$. These objects
 are called {\it Wakimoto sheaves}.

\begin{corollary}\label{wakimoto2}{{\sf{(I. Mirkovi\'c)}}}$\;$
\vi $\W_\lambda\starb\W_\mu= \W_{\lambda+\mu}\,,$ for any $\lambda,
 \mu\in \q$.

\ab \vii $\supp(\vp_*\W_\lambda) = \overline{\Gr}_\lambda;$
If $w\in W$ is an element of minimal length such that\linebreak
$w(\lambda)\in\q^{++},$ then:
 $(\vp_*\W_\lambda)\Big|_{\Gr_\lambda}=
\C_\lambda[-\dim\Gr_\lambda-\ell(w)]$.
\end{corollary}

\proof  Part (i) is a straightforward application of formulas
(\ref{conv_additive}). 
To prove (ii), choose dominant weights $\mu, \nu\in \q^{++}$, such that
$\lambda=\mu-\nu$.
Observe that, for any $w\in W$,
in $\waf$ we have   $\ell(-\nu)=\ell\bigl((-\nu)\cdot w^{-1}\bigr)
+\ell(w)$.
Hence, equations (\ref{conv_additive}) yield:
${{\mathcal{M}}}_{-\nu}={{\mathcal{M}}}_{(-\nu)\cdot w^{-1}}
\star
{{\mathcal{M}}}_w$.
Thus, since $\lambda=\mu-\nu$, we find:
\begin{equation}\label{wak5}
\W_\lambda\;=\;{{\mathcal{M}}}^\vee_\mu\starb 
{{\mathcal{M}}}_{-\nu}\;=\;{{\mathcal{M}}}^\vee_\mu\star
{{\mathcal{M}}}_{(-\nu)\cdot w^{-1}}
\star
{{\mathcal{M}}}_w\,.
\end{equation}

We now let $w\in W$ be as in part (ii) of the Corollary.
Then, the element
$\lambda\cdot w^{-1}$ is minimal
in the right coset $\lambda\cdot W\subset \waf,$ 
and we put $x=\lambda\cdot w^{-1}$. Geometrically,
this means that
 the Bruhat cell $\baf_x=\I\cd x \sset\baf$ is closed in
$\vp^{-1}\bigl(\vp(\baf_x)\bigr)$.
Further, let
 $A={{\mathcal{M}}}^\vee_\mu\star
{{\mathcal{M}}}_{(-\nu)\cdot w^{-1}}.$
Part (ii) of  Proposition \ref{wakimoto}
implies that $\supp A=\overline{\baf}_x$
and, moreover, $A\Big|_{\baf_x}=
\C_{\baf_x}[-\dim\baf_x].$ It follows that
$\supp\vp_*(A\star
{{\mathcal{M}}}_w)=\vp(\supp{A})=\overline{\Gr}_\lambda,$
and that 
$\bigl(\vp_*({A}\star{{\mathcal{M}}}_w)\bigr)\Big|_{\Gr_\lambda}$
$=
\C_{_{\Gr_\lambda}}[-\dim\Gr_\lambda-\ell(w)]$.
Part (ii) of the Corollary now follows from formula
(\ref{wak5}).
\qed

\subsection{An Ext-composition.}\label{E(J)}
Recall the setup of equation \eqref{conv_additive}.
In addition to the convolution-product functor
$\starb:$ $\DI({\mathcal{B}}) \times 
\DI({\mathcal{B}}) \too \DI({\mathcal{B}}),\,$
there is also a well-defined convolution functor
$\starb:\DI({\mathcal{B}}) \times   \DI(\Gr)
\too \DI(\Gr)$ arising from the convolution-diagram
$\Gd(\K)\times_{_\I}\Gr\too\Gr\,,\,(g,x)\mapsto gx$. Moreover,
for any $A\in \DI({\mathcal{B}})$ and $M\in \DI(\Gr),$
one has $A\starb M=(\vp_* A)\star M$.
In particular, for any $\mu\in \q^{++}$ and $\la \in \q,$
there is a well-defined object $\W_\la \starb
IC_\mu \in \DI({\mathcal{B}})$.

\ab   For any $M\in \DI(\Gr),$  and  $\la\in \q^{++},$
one may form an ind-object $\W_\la\starb M\starb\R \in
\limind\DI(\Gr)$.
Therefore, given two objects $M_1,M_2\in \DI(\Gr),$
we  may  consider the following
$\Z\times \q$-graded vector space
\begin{align}\label{J}
\EE(M_1,M_2)=\bigoplus\nolimits_{\la \in
\q^{++}\;}&\EE\hdot(M_1,M_2)_\la\,,\\
&
\EE\hdot(M_1,M_2)_\la:=
\Ext_{_{D^b(\Gr)}}\hdot (M_1\,,\,\W_\la\starb M_2\starb\R)\,.
\nonumber
\end{align}

\ab Next we introduce, for any three objects
$M_1,M_2,M_3\in \DI(\Gr),$ a  "composition-type"
pairing of $\EE$-groups, similar to the one defined in \S\ref{Ext1} on the
vector space $\Ext\hdot(\one, \R)$. Specifically,
given any $\mu,\la\in \q^{++}$ and $i,j\in\Z$, we define
a pairing
\beq\label{Ext2}
\EE^i(M_1,M_2)_\la\otimes\EE^j(M_2,M_3)_\mu\too
\EE^{i+j}(M_1,M_3)_{\la+\mu}
\eeq
as follows.
Let 
$\dis x\in\Ext_{_{D^b(\Gr)}}^i(M_2\,,\,\W_\la\starb M_3\starb\R)\,,\quad
y\in\Ext_{_{D^b(\Gr)}}^j(M_1\,,\,\W_\mu\starb M_2\starb\R).
$
We may view $x$ as a morphism
$x: M_2\to\W_\la\starb M_3\starb\R[i]$. 
Applying convolution, we get a morphism
\begin{align*} \W_\mu\star x\star\R:\;\;
\W_\mu\star M_2\star\R\too&\W_\mu\star (\W_\la\starb
M_3\starb\R[i])\star\R
\,=\,
(\W_\mu\star \W_\la)\starb M_3\star (\R\starb \R)[i].
\end{align*}
Note that $\W_\mu\star \W_\la\cong\W_{\la+\mu}$,
furthermore, the ring-object structure on $\R$
yields a morphism $\mt: \R\star \R\to\R$.
We now define $y\cdot x\in \EE^{i+j}(M_1,M_3)_{\la+\mu}
=\Ext_{_{D^b(\Gr)}}^{i+j}(M_1\,,\,\W_{\la+\mu}\starb
M_3\starb\R)$
as the following composite morphism
\begin{align*}
\xymatrix{{y\cdot x:\;
M_1}\ar[r]^<>(0.5){y}&\W_\mu\starb M_2\starb\R[j]\ar[rr]^<>(0.5){\W_\mu\star^{\,}
x^{\,}\star\R}&&(\W_\mu\star \W_\la)\starb M_2\star (\R\starb \R)[i][j]=}\\
\xymatrix{&&&=\W_{\la+\mu}\starb M_2\star (\R\starb \R)[i+j]
\ar[r]^<>(0.5){\mt}&\W_{\la+\mu}\starb
M_3\star\R[i+j].}
\end{align*}
The product $(x,y)\mapsto y\cdot x$ thus defined is associative
in a natural way.

\subsection{Homogeneous coordinate ring of $\NN$ as an Ext-algebra.}
\label{homogen}
  Let $\BL$ be the (finite dimensional) flag variety for the group
$G$, i.e., the variety of
all Borel subalgebras in $\g$. 
We write $\pi: \NN=T^*\BL \to \BL$ for the cotangent bundle on
$\BL$.
Recall that $G$ is of adjoint type,
hence, $\Hom(B,\C^\times)=\Hom(T,\C^\times)=\Y$.

\ab   For each $\la\in \Y,$ we
write $\oo_{_\BL}(\la)=G\times_{_B}\C(\la)$ for the standard $G$-equivariant line bundle
on $\BL$ induced from $\la$, the latter being viewed as a character $B\to\C^\times$.

\begin{notation}\label{O(l)}
  For any  $\la\in \Y,$ we put: $\on(\la) = \pi^*\oo_{_\BL}(\la)$.
\end{notation}

\ab The natural $\C^*$-action on $T^*\BL$ by dilations commutes with the
$G$-action, making $T^*\BL$ into a $\C^*\times G$-variety.
Clearly, $\on(\la)$ is 
a $\C^*\times G$-equivariant line bundle on $T^*\BL$.
The
$\C^*$-structure on $\on(\la)$ gives a $\Z$-grading  on
$\Gamma\hdot(\NN, \on(\la))$, the space of global sections. Further,
the obvious canonical isomorphism $\on(\la) \otimes \on(\mu)
\iso \on(\la+\mu)$ induces, for any $\la,\mu\in \Y,$
a bilinear pairing of the spaces of global sections:
$$\Gamma\hdot(\NN, \on(\la)) \otimes \Gamma\hdot(\NN, \on(\mu))
\too \Gamma\hdot(\NN, \on(\la+\mu))\,.$$
These pairings make
$\,\bigoplus_{\la \in \q^{++}}\;\Gamma\hdot(\NN, \on(\la))\,$
a $\Z\times\q^{++}$-graded algebra.

\ab Now, the construction of \S\ref{E(J)}, applied to the
special case $M_1=M_2=\one$, yields a $\Z\times\q^{++}$-graded algebra
 $\EE(\one,\one)=\bigoplus_{\la \in \q^{++}\;}
\Ext_{_{D^b(\Gr)}}\hdot (\one\,,\,\W_\la\starb\R).$
This algebra comes equipped with a $G$-action, as has been explained
at the end of \S \ref{regular_sheaf}.
It is known further that the Ext-group above has
no odd degree components, due to a standard parity-type vanishing
result for IC-sheaves on the loop Grassmannian (that is, for the affine
Kazhdan-Lusztig polynomials).

\ab 
The main result of this section is the
following

\begin{theorem}\label{ranee} There is a canonical $G$-equivariant 
$(\Z\times\q^{++})$-graded
algebra isomorphism
$$\bigoplus\nolimits_{\la \in \q^{++}\;}
\Ext_{_{D^b(\Gr)}}^{2\hhdot}(\one\,,\,\W_\la\starb\R)
\enspace\simeq\enspace
\bigoplus\nolimits_{\la \in \q^{++}}\;\Gamma\hdot(\NN\,,\, \on(\la))\,.$$
\end{theorem}

\ab To prove the Theorem recall the
principal $\sll$-triple $\langle \s,\e,\f \rangle \subset \g$,
 and the  Cartan
subalgebra $\fh=\g^{\s+e}$,
introduced in \S6. Associated to this $\sll$-triple,
is the Brylinski
 filtration, $W\idot V$, cf. Definition \ref{bry}, on any $G$-representation $V\in\rep(G)$.

\ab Our proof of Theorem \ref{ranee} exploits 
two totally different geometric interpretations of Brylinski
 filtration.
The first one is in terms of Intersection cohomology,
due to [G2], and the second one is in terms of 
equivariant line
bundles on the flag manifold for $G$, due to [Br]. 
The compatibility of the two interpretations, which is crucial
for the proof below, is yet another manifestation of Langlands duality.

\ab We begin by reviewing the main construction of [Br].
\subsection{Brylinski filtration in terms of Springer resolution (after \cite{Br}).}
Write  $H\subset G$ for the maximal torus
corresponding to the Cartan subalgebra
$\fh=\g^{t+e}.$ Given a  weight $\la\in \Y$
let 
 $e^{\la} : H\to \C^*$ denote the corresponding homomorphism.
Put
\begin{equation}\label{ind_weight} 
\C[G](-\la) = \{f\in \C[G]\enspace\Big|\enspace
f(h\cdot g)= e^{-\la}(h)\cdot f(g)\;\,,\,\;\forall g\in G\,,\,
h\in H\}\,.
\end{equation}
Observe that the $G$-action  on $\C[G]$ on the 
{\it right} by the formula 
$(R_gf)(y)=f(y\cdot g^{-1})$ makes the space $
\C[G](-\la)$ above into a {\it left} $G$-module,
isomorphic (by means of the
anti-involution: $y\mapsto y^{-1}$ on $G$) to the
 induced (from  character $e^\la$) representation 
$$
\indf_H^Ge^\la=\{f\in \C[G]\enspace\Big|\enspace
f(g\cdot h)= e^\la(h)\cdot f(g)\}.
$$

\ab Consider the  Borel subgroup   $B$ corresponding
to the Lie algebra $\Lie B=\n+\h$.
Since $e,t\in \Lie B$,
we have $H\subset B$. We identify the flag manifold $\BL$
with $G/B$.
The projection: $\Ad g(t+e) \mapsto gB/B$ makes
the conjugacy class
 $\Ad G(t+e)\subset \g$  an affine bundle over 
$\BL$ relative to the underlying vector bundle $\pi: T^*\BL\to
\BL$,
the cotangent bundle.

\ab Observe next that
the map $g \mapsto \Ad g(t+e)$, descends to
 a $G$-equivariant isomorphism $G/H\iso$
$ \Ad G(t+e).$ 
We  view the space $\indf_H^Ge^\la$, cf (\ref{ind_formula}),
as the space of global regular sections of the $G$-equivariant 
line bundle on the conjugacy class $\Ad G(t+e)$,
corresponding to the character $e^\la$. 
 Following 
[Br,\S5], consider the natural filtration
$F\idot^{^{_{\mathsf{fib}}}}$
on
$\indf_H^Ge^\la$ by {\it fiber degree}, see [Br, Theorem 4.4].
Write
$\gr^{^{_{\mathsf{fib}}}}\idot\bigl(\indf_H^Ge^\la\bigr)$
for the corresponding associated graded space.
By [Br, Theorem 5.5] we have a $\Z$-graded $G$-module isomorphism
\begin{equation}\label{ind_BR}
\gr^{^{_{\mathsf{fib}}}}\idot\bigl(\indf_H^Ge^\la\bigr)
\;=\;\Gamma\bigl(T^*\BL\,,\, \pi^*\oo(\la)\bigr) =
\Gamma\hdot\bigl(\NN, \on(\la)\bigr)\,.
\end{equation}
\ab The main idea of Brylinski is that the fiber degree filtration on
the line bundle corresponds to a `principal' filtration on 
(the $\la$-weight space of) the regular representation $\C[G]$,
see [Br, Theorem 5.8], whose $k$-th term is defined  by the formula:
$\Ker(e^{k+1})\cap \C[G]=$
$\{\varphi\in \C[G]\enspace\Big|\enspace
e^{k+1}(\varphi)=0\},$ where
$e(\varphi)$ denotes the  infinitesimal $e$-action on $\varphi$ on the
right.

\ab In more detail, for any $V\in\rep(G)$, we consider the
 weight decomposition $V=\bigoplus_{\mu\in\Y}\; V(\mu),$
 with respect 
to the Cartan subalgebra $\fh=\g^{t+e}$, and also
the Brylinski filtration $W\idot V$, cf. Definition \ref{bry}.
We will use the notation $\,W_kV(\mu):=V(\mu)\cap W_kV$,
and $\gr^W_kV(\mu):=W_kV(\mu)/W_{k-1}V(\mu).$
Write $\hht(\la):=\sum_{i\in I}\,\langle\la,\check{\alpha}_i\rangle$
for the {\em height} of $\la\in\Y$.

\ab
By [G2, Lemma 5.5.1] we know that,
for any $V\in \rep(G)$, $\la\in \Y^{++}$ and $k\in \Z$, one has:
$\Ker(e^{k+1})\cap V(-\la)= W_{k+\hht(\la)}V
\cap V(-\la).$ 
This equation, applied for $V=\C[G]$,
 and Brylinski's Theorem [Br, Thm. 5.8] cited above imply that
the graded space on the LHS of formula
(\ref{ind_BR}) coincides with 
the graded space 
$\gr^W_{\!\hhdot+\hht(\la)}\bigl(\indf_H^Ge^\la\bigr)$.
Thus, for any $\la\in \q^{++}$,
we obtain a chain of isomorphisms
\begin{equation}\label{BR_fin} 
\gr^W_{\!\hhdot+\hht(\la)}\bigl(\indf_H^Ge^\la\bigr)
\,=\,
\gr^{^{_{\mathsf{fib}}}}_{\!\hhdot+\hht(\la)}\bigl(\indf_H^Ge^\la\bigr)
\,=\,
\Gamma\hdot\bigl(\NN, \on(\la)\bigr) 
.
\end{equation}

\subsection{Proof of  Theorem \ref{ranee}.}
\ab We first fix $\la \in \q^{++}$, and construct a
$\Z\times\q$-graded
vector space isomorphism $
\Ext_{_{D^b(\Gr)}}\hdot(\one\,,\,\W_\la\starb\R)
\simeq\Gamma\hdot(\NN\,,\, \on(\la)).$
We will then verify compatibility of the constructed isomorphisms
for different $\la$'s
with graded algebra structures on both sides
of the isomorphism of Theorem \ref{ranee}.

\begin{notation}
 Given $\la\in\Y$, 
let  $j_\la: \Gr_\la \into \Gr$ be the Bruhat cell imbedding,
and $\,\Delta_\la:=(j_\la)_!\C_{_{\Gr_\la}}[\dim\Gr_\la],$
the corresponding `standard'
perverse sheaf on $\Gr$.
\end{notation}

\ab We will use below that, for $\la\in \q^{++},$ 
the restriction of the projection
$\vp:\baf\onto\Gr$ to the Bruhat cell $\baf_{-\la}$ is a fibration
$\baf_{-\la} \to \Gr_{-\la}$ 
with fiber $\C^{\ell(w)}$,
where $w\in W$ is the element of minimal length such that
$w(-\la)\in \q^{++}$. It follows that
$\vp_*{\mathcal{M}}_{-\la}=
\Delta_{-\la}[-\ell(w)].$ Further, for $\la\in \q^{++}$ we have by definition,
$\W_{-\la}={\mathcal{M}}_{-\la}$, hence,
$\W_{-\la}\starb \one=\vp_*\W_{-\la}=\vp_*{\mathcal{M}}_{-\la}=
\Delta_{-\la}[-\ell(w)].$

\ab Note that, since convolving with $(\W_{-\la}\star\W_\la)$
acts as the identity functor on $\DI(\Gr)$,
the functor $\W_{-\la}\starb(-)$ is an equivalence.
Therefore, we find:
\begin{align*}
\Ext_{_{D^b(\Gr)}}\hdot(\one\,,\,\W_\la\starb\R)\;&=
\Ext_{_{D^b(\Gr)}}\hdot(\W_{-\la}\starb \one\,,\,
\W_{-\la}\star\W_\la\starb\R)\nonumber\\
&=\Ext_{_{D^b(\Gr)}}\hdot(\W_{-\la}\starb \one\,,\,\R)
=\Ext_{_{D^b(\Gr)}}\hdot(\vp_*\W_{-\la}\,,\,\R)\\
&=\Ext_{_{D^b(\Gr)}}\hdot(\vp_*{\mathcal{M}}_{-\la}\;
\,,\,\R) =\Ext_{_{D^b(\Gr)}}\hdot\bigl(\Delta_{-\la}[-\ell(w)]
\,,\,\R\bigr).
\end{align*}

\ab Further, view an element $\nu\in \q$ as a point 
in $\Gr$, the center of the cell $\Gr_{\nu}$,
and write $i_{\nu}: \{\nu\}\into \Gr$ for the corresponding imbedding.
Then, for any object  $N\in D^b\PI$, we clearly have  
$\Ext_{_{D^b(\Gr)}}\hdot(\Delta_\nu, N)=
H\hdot\bigl(i_\nu^!N[\dim\Gr_{\nu}]\bigr)
=H^{\hhdot+\dim\Gr_{\nu}}(i_\nu^!N)$. 
Hence, for $\la\in\Y^{++}$, we get:
$\Ext_{_{D^b(\Gr)}}\hdot\bigl(\Delta_{-\la}[-\ell(w)], N)
=H^{\bullet-\ell(w)+\dim\Gr_{-\la}}(i_{-\la}^!N)$.
Combining all the observations above,
we obtain a canonical isomorphism
\begin{equation}\label{lhs_thm}
\Ext_{_{D^b(\Gr)}}\hdot(\one\,,\,\W_\la\starb\R)
\;=\;
H^{\hhdot+\hht(\la)}(i_{-\la}^!\R)\enspace,\enspace
\hht(\la):=\mbox{height}(\la)\,,\,\forall \la\in \q^{++},
\end{equation} 
where
we have used that if $\la$ is dominant, and $w\in W$ is an element of
minimal length such that $-w(\la)\in\Y^{++},$ then
$\dim\Gr_{-\la}=\ell(w)+\hht(\la)$.

\ab We will express the RHS of (\ref{lhs_thm})
in representation theoretic terms.
To this end, we recall first that
there is a canonical
filtration, $W\idot H_s(M)$,
on the specialized  equivariant cohomology
of any object $M \in D^b_{\Td}(\Gr)$.
Now, let $\mu\in\Y\subset \Gr$ be
a $\Td$-fixed point, and $M=i_\mu^!A$,
for some $A\in \PO$.
By [G2, Proposition 5.6.2],
 the filtration $W\idot$
is strictly compatible with the fixed point decomposition
(\ref{fp_decomposition}), i.e.,
we have
\begin{equation}\label{strict}
W_kH_s(i_\mu^!A) = H_s(i_\mu^!A)\, \cap\, W_kH_s(A)
,\quad \forall A\in \PO\,,\,k\in\Z\,.
\end{equation}
 Further, for any $V\in\rep(G)$, the
fixed point decomposition on $H_s(\p V)$ corresponds,
by  Theorem \ref{main_perv},
to the weight decomposition on $V=\bigoplus_{\mu\in\Y}\; V(\mu),$
 with respect 
to the Cartan subalgebra $\fh=\g^{t+e}$ and, moreover,
the filtration $W\idot H_s(\p V)$
corresponds to Brylinski filtration $W\idot V$.
Hence, (\ref{strict}) implies that, for any $k\in\Z$,
 the 
subspace $W_kH_s\bigl(i_\mu^!(\p V)\bigr) \subset
H_s(\p V)$ corresponds to the subspace
$V(\mu)\cap W_kV \subset V$.

\ab It is known, cf. e.g. \cite{G1}
that, for any $A\in \PO$, the complex
$i_\mu^!A$ is {\em pure} in the sense of [BBD].
Hence, the spectral sequence for
equivariant hyper-cohomology of $i_\mu^!A$
 collapses.  It follows, see  \cite[Proposition 5.6.2]{G2},
that, for any $A\in \PO$,  one has a canonical
isomorphism $H\hdot(i_\mu^!A)
=\gr^W\idot\bigl(H_s(i_\mu^!A)\bigr)$.
Thus, we obtain canonical isomorphisms
\begin{equation}\label{key_compare}
H^k\bigl(i_\mu^!(\p V)\bigr) =
\gr^W_kH_s\Bigl(i_\mu^!(\p V)\Bigr) =\gr^W_kV(\mu)\,.
\end{equation}
\ab Since $\R=\p \C[G]$, the isomorphism (\ref{key_compare})
yields:
$
H^k(i_{-\la}^!\R)=
\gr^W_k\bigl(\C[G](-\la)\bigr),$ where $\C[G](-\la)$ is
the $(-\la)$-weight space  of
the {\it left} regular representation of $G$.
Thus, the considerations above and formulas (\ref{lhs_thm}) and
\eqref{BR_fin}
yield isomorphisms 
\begin{equation}\label{ind_formula}
\Ext_{_{D^b(\Gr)}}\hdot(\one\,,\,\W_\la\starb\R)
\;\simeq\;
\gr^W_{\!\hhdot+\hht(\la)}\bigl(\indf_H^Ge^\la\bigr)
\;\simeq\;
\Gamma\hdot\bigl(\NN, \on(\la)\bigr) 
,\quad\forall \la\in \q^{++}\,.
\end{equation}
The composite isomorphism
 gives the isomorphism claimed in Theorem \ref{ranee}.

\ab To complete the proof, we must show that, for any
$\la,\mu \in \q^{++},$ the isomorphisms
(\ref{ind_formula}) and (\ref{BR_fin})  transport the product map,
see \eqref{Ext2}:
\begin{equation}\label{ext_product_W}
\Ext_{_{D^b(\Gr)}}\hdot(\one\,,\,\W_\la\starb\R)
\otimes \Ext_{_{D^b(\Gr)}}\hdot(\one\,,\,\R\starg\W_\mu)
\too
\Ext_{_{D^b(\Gr)}}\hdot(\one\,,\,\R\starg\W_{\la+\mu})
\end{equation}
to the natural product-pairing:
$\Gamma\hdot\bigl(\NN, \on(\la)\bigr)
\otimes \Gamma\hdot\bigl(\NN, \on(\mu)\bigr)
\too
\Gamma\hdot\bigl(\NN, \on(\la+\mu)\bigr).$

\ab To this end, consider  the action-map
$a: \Gd(\K)\times_{_{G^\vee(\oo)}}\Gr\too\Gr$, see Sect. \ref{satake_sec}.
For any $\chi\in\Y$,
set
$\,\Gr^2_\chi:= a^{-1}(\{\chi\})\subset \Gd(\K)\times_{_{G^\vee(\oo)}}\Gr,$
and
write $j_\chi :\Gr^2_\chi\into \Gd(\K)\times_{_{G^\vee(\oo)}}\Gr$
for the imbedding. By base change, we have  
$i_\chi^!\ccirc a_*=a_*\ccirc j_\chi^!.$
Further, recall
the  morphism
$\mt: \R\star\R=a_*(\R\tilde{\boxtimes} \R)
\to \R$ corresponding to the algebra structure on
$\C[G]$, see \eqref{product_map_m}. We form the composite morphism
\begin{equation}\label{comp_mor}
a_*\ccirc j_\chi^!(\R\tilde{\boxtimes} \R)=
i_\chi^!\ccirc a_*(\R\tilde{\boxtimes} \R)
\stackrel{\mt}{\too}
i_\chi^!\R
\end{equation}
Now, given $\nu,\eta\in\Y$, we have an obvious
imbedding $i_\nu\boxtimes i_\eta: \{\nu\}\times
\{\eta\}\into \Gr^2_{\nu+\eta}$.
Composing the maps of cohomology induced by the latter imbedding
and by morphism (\ref{comp_mor})
for $\chi=\nu+\eta$, one obtains  natural maps

\begin{align}\label{comp_mor2}
H\hdot(i_\nu^!\R)\otimes H\hdot(i_\eta^!\R)=
H\hdot(i_\nu^!\R \boxtimes i_\eta^!\R) & \to
 H\hdot\bigl(j_{\nu+\eta}^!(\R\tilde{\boxtimes} \R)\bigr)\\
&=H\hdot\bigl(a_*\ccirc j_{\nu+\eta}^!(\R\tilde{\boxtimes} \R)\bigr)
\stackrel{(\ref{comp_mor})}{\too}H\hdot(i_{\nu+\eta}^!\R)\,.
\nonumber
\end{align}

\ab It is a matter of routine argument  involving base change, cf. e.g. proof
of [G2,~Proposition 3.6.2], to
show that  the following diagram commutes
{\small
$$
\diagram
\Ext_{_{D^b(\Gr)}}\hdot(\one\,,\,\W_\la\starb\R)
\otimes \Ext_{_{D^b(\Gr)}}\hdot(\one\,,\,\R\starg\W_\mu)
\dto^{(\ref{ext_product_W})}\rdouble_<>(.5){(\ref{lhs_thm})}&
H^{\hhdot+\hht(\la)}(i_{-\la}^!\R)\otimes 
H^{\hhdot+\hht(\mu)}(i_{-\mu}^!\R)\dto^{(\ref{comp_mor2})}\\
\Ext_{_{D^b(\Gr)}}\hdot(\one\,,\,\R\starg\W_{\la+\mu})
\rdouble_<>(.5){(\ref{lhs_thm})}&
H^{\hhdot+\hht(\la+\mu)}(i_{-\la-\mu}^!\R)\,.
\enddiagram
$$}

\ab
Thus, completing the proof of the Theorem amounts
to verifying that: {\em For any $\nu,\eta$,
and  $V=\C[G]$,
the  diagram below (arising from
 (\ref{ind_weight}) by means of
isomorphisms \eqref{key_compare}) commutes}
{\small
$$
\diagram
H\hdot(i_\nu^!\R)\otimes
H\hdot(i_\eta^!\R) \dto^{(\ref{comp_mor2})}\rdouble
&\gr\idot^WH_s(i_\nu^!\R)\otimes
\gr\idot^WH_s(i_\eta^!\R)\dto^{(\ref{comp_mor2})}\rdouble
&\gr\idot^W\bigl(\C[G](\nu)\bigr)\otimes
\gr\idot^W\bigl(\C[G](\eta)\bigr)\dto^{{\mathtt{mult}}}\\
H\hdot(i_{\nu+\eta}^!\R)\;\rdouble&
\;\gr\idot^WH_s(i_{\nu+\eta}^!\R)\;\rdouble&
\;\gr\idot^W\bigl(\C[G](\nu+\eta)\bigr)\,.
\enddiagram
$$}
The square on the left of the last diagram commutes because
$H_s(-)$ is a tensor functor on $\PO$, and
the canonical filtration $W\idot$ on the specialized
equivariant cohomology is compatible with convolution.
The square on the right commutes, because the filtration
on $\g$-representations by the eigenvalues of the $t$-action
is obviously compatible with tensor products.
This completes the proof of the Theorem. \quad\qed

\begin{remark} Write: $\widetilde{\g}=G\times_{\!_B}\b$
for  Grothendieck's
simultaneous resolution, see e.g., [CG, ch.3].
Similarly to Theorem \ref{main1},  one also has an equivariant version
of Theorem \ref{ranee}, saying that there is a canonical $\Z\times\q^{++}$-graded
algebra isomorphism
$$\bigoplus\nolimits_{\la \in \q^{++}\;}
\Ext_{_{\DI(\Gr)}}\hdot(\one\,,\,\W_\la\starb\R)
\enspace\simeq\enspace
\bigoplus\nolimits_{\la \in \q^{++}}\;\Gamma\hdot(\widetilde{\g}\,, \,
{{\mathcal{O}}}_{_{\widetilde{\g}}}(\la))\,.\qquad\lozenge$$
\end{remark}
\subsection{Another fiber functor on $\PO$.}\label{another}
Following \S\ref{homogen}, we  consider a
$\Z\times\Y^{++}$-graded algebra
$$\CA :=
\bigoplus\nolimits_{\la \in \q^{++}}\;\Gamma\hdot(\NN\,,\, \on(\la)).
$$
The action of $G$ on $\NN$ induces a 
natural $G$-action on $\CA$ by graded algebra automorphisms.

\ab By construction, the algebra $\CA$
may be viewed as a multi-graded homogeneous coordinate ring
of the variety $\NN$. In particular, any object
$M\in \mcac$ gives rise to a $G\times\C^*$-equivariant coherent sheaf
$\ff(M)\in\Coh^{G\times\C^*\!}(\NN)$. 
The assignment $M\mto \ff(M)$ gives an exact
functor $\ff: \mcac\too\Coh^{G\times\C^*\!}(\NN)$.

\ab Next, recall the notation \eqref{J}
and observe that, for any $M\in \D^b_{_\J}(\Gr)$, 
the space $\EE(\one, M)$ has a natural
graded $\EE(\one,\one)$-algebra structure,
via the pairing defined in~\eqref{Ext2}.

\ab Now, the graded algebra
isomorphism $\EE(\one,\one)\cong \CA$, of
Theorem  \ref{ranee},
gives rise to  an equivalence
$\Mod^{G\times\C^*\!}\bigl(\EE(\one,\one)\bigr)
\iso \Mod^{G\times\C^*\!}\bigl(\CA\bigr).$ We consider
 the  composite functor $\widehat{\mathbb{S}}: V\mto \ff\bigl(\EE(\one, \p{V})\bigr)$;
explicitly, our functor is  the following composite:
\begin{align}\label{bE}
\widehat{\mathbb{S}}:\ \rep(G)\stackrel{\p}\iso \PO\into  D^b_{_\J}(\Gr)&\too 
\op{Mod}^{G\times\C^*\!}\bigl(\EE(\one,\one)\bigr)\\
&\iso \Mod^{G\times\C^*\!}\bigl(\CA\bigr)
\stackrel{\ff}\too\Coh^{G\times\C^*\!}(\NN).\nonumber
\end{align}

\ab Very similar to Proposition \ref{bpsi}, using the isomorphism in
\eqref{nonsense2},
one proves

\begin{proposition}\label{bpsi2} The functor $\widehat{\mathbb{S}}$
is isomorphic to the functor
$V\mapsto V\otimes \on.$\qed
\end{proposition}

\ab We will not go into more details here, since a much more elaborate
version of this result will be proved in  section \ref{sec_equivalence}.

\subsection{Monodromic sheaves and extended affine flag manifold.}
At several occasions, we will need to extend various
constructions involving convolution of $\I$-equivariant
sheaves to the  larger category of $\I$-monodromic sheaves.
In particular, we would like to
extend the bifunctor $M_1,M_2\mto \EE(M_1,M_2),$
from the category $D^b_{_\I}(\Gr)$ to the category
$D^b_{_\J}(\Gr)$. The construction of such an `extension'
is analogous to a similar construction in the framework of
 $\mathcal{D}$-modules
on the (finite dimensional) flag manifold, given in \cite[\S5]{BlG}
(see esp. \cite[(5.6), (5.9.1)-(5.9.2)]{BlG}).
We proceed to more details.

\ab In $G^\vee(\K)$, we consider the following subgroup
 $\I_1:=\{\gamma\in G^\vee(\oo)\;|\;\gamma(0)=1\}.$
Thus, $\I_1$ is a pro-unipotent algebraic group
that may be thought of as the (pro)-unipotent radical of the Iwahori
group $\I$. We have $\I=T^\vee\cdot\I_1.$ The coset space
$\baft:=G^\vee(\K)/\I_1$ has a natural ind-scheme structure,
and will be called the {\em extended affine flag manifold}.
The torus $T^\vee$ normalizes the group $\I_1$. Hence, there is a
natural $T^\vee$-action on $\baft$ on the right, making the
 canonical projection $\pi: \baft=G^\vee(\K)/\I_1\too G^\vee(\K)/\I=\baf$
a principal $T^\vee$-bundle.

\ab Let $X$ be an (ind-) $\I$-variety.
We recall that   the notions of being
$\I_1$-monodromic and  of being $\I_1$-equivariant are 
known to be
equivalent, since the group
$\I_1$ is  pro-unipotent.
In particular,  we have 
$D^b_{_\J}(X)\subseteq D^b_{_{\I_1\mbox{-}{\sf{mon}}}}(X)
=  D^b_{_{\I_1}}(X).$

\ab Assume now  $X$ is  an (ind-) $G^\vee(\K)$-variety.
The  $G(\K)$-action on $X$ gives rise to the following
{\em convolution diagram}:
 $G^\vee(\K)\times_{_{\I_1}} X \map X\,,\,(g, x) \mapsto gx,$ cf.,
\eqref{baft_conv}.
For any $\A\in D^b_{_{\J}}(G(\K)/\I_1)=D^b_{_{\J}}(\baft)$
and $\M\in D^b_{_\J}(X)\subseteq D^b_{_{\I_1}}(X)$, there is a well-defined object
${\A}\tilde{\boxtimes} \M\in D^b\bigl(G^\vee(\K)\times_{_{\I_1}}X\bigr)$.
We put $\A\star \M:=a_*(\A\tilde{\boxtimes} 
{{\mathcal{M}}})$,  an $\I$-monodromic complex on $X$.
This way, one defines a convolution bifunctor
$\star: D^b_{_\J}(\baft)\times D^b_{_\J}(X)\too D^b_{_\J}(X).$
This applies, in particular, for $X=\baf$ and $X=\Gr$.

\ab
  For each $w\in\waf,$ we put
$\FIt_w:=\pi^{-1}(\FI_w)$, a $T^\vee$-bundle over the cell $\FI_w$,
which is (non-canonically) isomorphic to the trivial $T^\vee$-bundle
$T^\vee\times\FI_w\to \FI_w.$
Write $\widetilde{j}_w: \FIt_w\into \baft$ for the
imbedding.

\ab Let ${\mathscr{E}}$ be a "universal" pro-unipotent local
system on $T^\vee$; specifically, if $\Pi$ denotes the
group algebra of the fundamental group of the torus  $T^\vee$
with respect to a base point $*\in T^\vee$,
then the fiber of the  local
system ${\mathscr{E}}$ at $*\in T^\vee$ equals the completion of
$\Pi$ at the augmentation ideal. Given $w\in \waf$,
let ${\mathscr{E}}_w$ denote the
pull-back of ${\mathscr{E}}$  via the  projection $\FIt_w\to T^\vee$
arising from
a (non-canonically) chosen $T^\vee$-bundle trivialization $\FIt_w\simeq
T^\vee\times\FI_w.$
We put $\,\Mt_w
:=(\widetilde{j}_w)_!{\mathscr{E}}_w[\dim \FIt_w],\,$ and
$\,\Mt_w^\vee:=(\widetilde{j}_w)_*{\mathscr{E}}_w[\dim \FI_w].$
These are pro-objects in $\pe_{_\J}(\baft)$.

\ab Let  $w\in\waf$.
It is straightforward to verify that, for any 
$M\in
D^b_{_\J}(\baf),$
resp., $M\in
D^b_{_\J}(\Gr)$, the pro-objects $\Mt_w\star M\,,\,\Mt_w^\vee\star M$
are indeed actual objects of $D^b_{_\J}(\baf),$ respectively,  of
$D^b_{_\J}(\Gr).$
Moreover, if
 $M\in
D^b_{_\I}(\baf),$
resp., $M\in
D^b_{_\I}(\Gr)$, then one has canonical isomorphisms
\beq\label{reduce_gr}
\Mt_w\star M\simeq
{{\mathcal{M}}}_w\star M\quad\text{and}\quad
\Mt_w^\vee\star M\simeq
{{\mathcal{M}}}_w^\vee\star M.
\eeq
\begin{remark} We warn the reader that, for $M\in
D^b_{_{\J}}(\Gr)$ and $\la\in\Y$, the convolution $\W_\la\star M$ is
{\em not} defined,
in general.
$\quad\lozenge$ 
\end{remark}
\ab We are now able to extend the constructions of \S\ref{E(J)}, and
 assign $\EE$-groups to objects of the category
$D^b_{_{\J}}(\Gr)\supset D^b_{_{\I}}(\Gr),$ as follows. 
 Given  $\la\in\Y$, choose
 $\mu,\nu \in \q^{++},$ such that $\lambda=\mu-\nu$,
and
consider the objects $\Mt^\vee_\mu\,,\,\Mt_{-\nu}\in
D^b_{_{\J}}(\baft).$
Although the convolution of two objects of
the category $D^b_{_{\J}}(\baft)$ has not been defined,
 we can define an object
$$\blangle \W_\la\star M\brangle  \;:=\; \Mt^\vee_\mu\star (\Mt_{-\nu}\star M)\in D^b_{_{\J}}(\Gr),
\quad\text{for any}\quad M\in D^b_{_{\J}}(\Gr),
$$
that involves only 
the convolution  $\star: D^b_{_{\J}}(\baft)\times
D^b_{_{\J}}(\Gr)\too D^b_{_{\J}}(\Gr).$
One verifies that the object $\blangle\W_\la\star M\brangle$ thus defined is independent
of the choice of the presentation  $\lambda=\mu-\nu$.
Moreover, it follows readily from \eqref{reduce_gr} that
the functor $\blangle \W_\la\star(-)\brangle : D^b_{_{\J}}(\Gr)\too D^b_{_{\J}}(\Gr),$
when restricted to the subcategory
 $D^b_{_{\I}}(\Gr)\sset  D^b_{_{\J}}(\Gr)$,
agrees with the functor $\W_\la\star(-)$
given by the genuine convolution.
The functor   $\blangle \W_\la\star(-)\brangle$ also has all the other
expected properties, in particular,
we have ${\boldsymbol{\big{\langle}}}
 \W_\mu\star\blangle\W_\la\star M\brangle {\boldsymbol{\big{\rangle}}}
=\blangle \W_{\mu+\la}\star M\brangle .$

\ab Now, for any $M_1,M_2\in  D^b_{_{\J}}(\Gr)$,
we put 
$$\EE(M_1,M_2):=\bigoplus\nolimits_{\la \in
\q^{++}\;}
\Ext_{_{D^b(\Gr)}}\hdot \bigl(M_1\,,\,\blangle \W_\la\starb
M_2\brangle \starb\R\bigr),
$$
where the  ind-object $\blangle \W_\la\starb
M_2\brangle \starb\R$  is  well-defined
 since $\R$ is clearly an $\I$-equivariant ind-sheaf on $\Gr$.

\ab From now on, we will make no distinction between the
functors $\W_\la\star(-)$ and $\blangle \W_\la\starb(-)\brangle $
and, abusing the notation, write simply  $\W_\la\star(-)$.

\section{Geometric Equivalence theorems}\label{sec_equivalence}
\ab In this section (only) we will be working over the ground field
$\k=\overline{{\mathbb{Q}}_\ell}$, an algebraic closure of
the field of $\ell$-adic numbers. We write $\Gm$ for the
multiplicative group, viewed as a 1-dimensional algebraic
group (a torus) over $\k$.

\subsection{Equivalence theorem for perverse sheaves.}\label{eq_perv_sub}
 For any $\lambda \in \q$,
the closure, $\overline{\Gr}_\lambda=\overline{\I\cd\la}\sset \Gr$ is known
to be a finite dimensional projective variety, 
an affine Schubert variety.
Let $IC_\lambda=IC(\overline{\Gr}_\lambda)$ denote the corresponding
intersection
cohomology complex, the Deligne-Goresky-MacPherson extension of
the constant sheaf on ${\Gr}_\lambda$ (shifted by $\dim {\Gr}_\lambda$).
  For $\la\in\Y^{++}$, we have
$\overline{\Gr}_\lambda=\overline{\Gd(\oo)\cd\la},$
hence the notation above  agrees with the one used
in \S\S6-8.

\begin{notation} 
We write $\Perv(\overline{\Gr}_\la):=\Perv_{_\J}(\overline{\Gr}_\la)$
 for the abelian category 
of  $\I$-monodromic $\ell$-adic perverse sheaves on
$\overline{\Gr}_\la$,
and let $\Perv(\Gr):=
\underset{\la\in\Y^{++}}\limind \Perv_{_\J}(\overline{\Gr}_\la),$
be a direct limit of these categories.
\end{notation}

\ab Clearly, $\Perv(\Gr)$ is an abelian subcategory in $D^b(\Gr)$.
The Ext-groups in the  categories $\Perv(\Gr)$ and $D^b(\Gr)$
turn out to be the same.
Specifically, one has
\begin{proposition}[\cite{BGS}, Corollary 3.3.2]\label{BGS}
There is
 a canonical isomorphism
$$
\Ext\hdot_{_{\Perv(\Gr)}}(\L_1, \L_2) \iso
\Ext\hdot_{_{D^b(\Gr)}}(\L_1, \L_2),\quad
\forall\,\L_1,\, \L_2 \in \Perv(\Gr).\qquad\Box
$$
\end{proposition}

\ab 
\ab Recall that, for each $\la\in\Y^{++},$
the  finite-dimensional projective variety $\overline{\Gr}_\la$ admits 
 a (finite) algebraic stratification by
Schubert cells. Therefore, 
each category $\Perv(\overline{\Gr}_\la)$ has finitely many
 simple
objects, and the simple
objects of  the
category $\Perv(\Gr)$ are parametrized by  elements of the
root lattice $\Hom(\Gm,\Td)=\q$.
  For each $\la\in\Y$, one also has the standard
perverse sheaf $\Delta_\la=(j_\la)_!\k_{_{\Gr_\la}}[-\dim\Gr_\la]\in
\Perv(\Gr),$ and the costandard perverse sheaf 
$\nabla_\la=(j_\la)_*\k_{_{\Gr_\la}}[-\dim\Gr_\la]\in
\Perv(\Gr)$, where $j_\la:\Gr_\la\into \Gr$ denotes the
Schubert cell imbedding.

\begin{notation}\label{bW_def} Given $\la\in\Y,$ we put
$\bW_\la:= \W_\la\star \one=\varpi_*(\W_\la)\in \Perv(\Gr)$,
where $\varpi:\baf\to\Gr$ is the standard projection,
cf. \S\ref{baf}.
\end{notation} 

\ab If $\la\in\Y^{++}$, resp., $\la\in-\Y^{++},$ then from the
definition of Wakimoto sheaves we get $\bW_\la=\nabla_\la$, resp., $
\bW_\la=\Delta_\la$.

\ab Recall the category
$D^G_{\text{coherent}}(\NN)$ introduced in \S\ref{reminder},
and invertible coherent sheaves $\on(\la)$ on the Springer resolution 
$\NN$, see Notation \ref{O(l)}.
Later in this section we 
are going to prove the following  result that yields 
(a non-mixed counterpart of) the
equivalence $P$ on the right of diagram~\eqref{sum-up}.

\begin{theorem}\label{P_equiv0}
There is an equivalence of  triangulated categories
$P': D^b\Perv(\Gr)\iso$
$D_{\text{coherent}}^G(\NN),$
such that $P'(\one)=\on,$ and, for any $M\in
D^b\Perv(\Gr),$ one has
$$
P'(\W_\la\star M\star \p{V})=V\otimes\on(\la)\otimes P'(M),
\enspace\text{for all}\quad\la\in\Y\,,\,V\in\rep(G).
$$
\end{theorem}

\subsection{Mixed  categories.}\label{mixed_categories} An   abelian
$\k$-category $\scrc\mix$ is called a {\em mixed} (abelian) category,
cf. \cite[Definition 4.1.1]{BGS}, provided  a map $w: \op{Irr}(\scrc)\to\Z$ (called
{\em weight}) is given, such that for any two simple objects
$M,N\in \op{Irr}(\scrc\mix)$ with $w(M)\leq w(N),$
one has $\Ext^1_{\scrc\mix}(M,N)=0$.

\ab Let  $\scrc\mix$ be a mixed category with  degree one 
{\em Tate twist}, i.e., with  an auto-equivalence $M\to M\langle1\rangle$
of  $\scrc\mix,$ such that $w(M\langle1\rangle)=w(M)+1$
for any $M\in  \op{Irr}(\scrc\mix)$.
Every object $M\in \scrc\mix$ comes
equipped with a canonical increasing `weight' filtration $\Wm\idot M$,
cf. \cite[Lemma~4.1.2]{BGS}; the Tate twist
shifts the weight filtration by 
1. 

\begin{remark} In \cite[Definition 4.1.1]{BGS},
an additional requirement that  $\scrc\mix$ is
an artinian category is included in the definition of a mixed category.
This requirement may be replaced, without affecting the theory,
 by the following two weaker conditions:

\pb{For any $M\in \scrc\mix$ and $i\in\Z$, the object $\Wm_iM/\Wm_{-i}M$ has
finite length, and}

\pb{For any $M,N\in \scrc\mix$ and $i\geq 0$, we have
$\dim\Ext^i_{\scrc\mix}(M,N)<\infty$.}

\ab These two conditions hold for all mixed categories considered below.
 $\quad\lozenge$\end{remark}

\ab We recall the notion of a `mixed version' of an
abelian category, alternatively called a `grading'
on the category, see \cite[\S4.3]{BGS}. Let $\scrc$ be an abelian $\k$-category,
 $\scrc\mix$ a  mixed  abelian category,
and $v:\scrc\mix\to\scrc$ an exact faithful functor.

\begin{definition}[\cite{BGS}, Definition 4.3.1]\label{grading}
 The pair  $(\scrc\mix,v)$ is said to be a
{\em mixed version} of  $\scrc$ (with forgetful functor $v$)
 if the following holds:

\pb{The functor $v$ sends semisimple objects into semisimple objects,
and any simple object of $\scrc$ is isomorphic to one of the form
$v(M)\,,\,M\in\op{Irr}(\scrc\mix)$;}

\pb{There is a natural isomorphism $\varepsilon: v(M)\iso
v(M\langle1\rangle)\,,\,
\forall M\in \scrc\mix$;}

\pb{For any $M,N\in \scrc\mix,$ the functor $v$ induces an isomorphism}
\beq\label{mix_ext}
\bigoplus\nolimits_{n\in\Z}\;\Ext\hdot_{\scrc\mix}(M,N\langle n\rangle)\iso
\Ext\hdot_{\scrc}\bigl(v(M)\,,\,v(N)\bigr).
\eeq
\end{definition}

\ab For example, 
 the category  $\Coh^{G\times\Gm}(\NN)$
has a natural structure of mixed abelian category.
It is a mixed version of 
 the abelian
category $\Coh^G(\NN)$,
with the functor
$v: \Coh^{G\times\Gm}(\NN)\to \Coh^G(\NN)$ forgetting the
$\Gm$-equivariant structure. 

\ab Similarly,
the abelian category $\agr$, see Notation \ref{cohG},
is a mixed abelian category.
The set  $\op{Irr}(\agr)$ of (isomorphism classes of)
 simple
objects of this category consists 
1-dimensional graded modules:
$\k_\wedge(\la)\langle i\rangle\,,\,\la\in\Y,i\in\Z,$ (here $\langle i\rangle$
indicates that the vector space $\k_\wedge(\la)$ is placed in $\Z$-grade degree $i$).
The weight function on $\op{Irr}(\agr)$
is given by
$w: \k_\wedge(\la)\langle i\rangle\mto i$.
Further, the natural functor $\agr\too\ag$ forgetting the $\Z$-grading
makes the category $\agr$ a mixed version of the category $\ag$, see Notation \ref{cohG}.

\ab We will  also use  the notion of a mixed {\em triangulated}
category, for which we refer the reader to [BGSh].
The derived category of a mixed abelian category is a
mixed  triangulated category.
Given  a triangulated category $D$, a mixed triangulated category
$D\mix$, and a triangulated functor $v: D\mix\to D$,
one says that the pair $(D\mix, v)$ is a mixed version of $D$ if
conditions similar to those of Definition \ref{grading} hold.

\ab For example, we have a natural functor
 $\Coh^{G\times\Gm}(\NN)\too
\dcoh^G(\NN)$, which extends to
a triangulated  functor 
${\mathbf{f}}: D^b\Coh^{G\times\Gm}(\NN)\too
\dcoh^G(\NN)$.
We claim that the latter functor 
 makes the category
$D^b\Coh^{G\times\Gm}(\NN)$ a mixed  version of $\dcoh^G(\NN)$.
To see this, we first apply the equivalences of 
diagram \eqref{res_gamma}.
This way,
we reduce the claim to the
statement that a similarly defined functor
$D^b\Mod^{B\times\Gm}_f(\SS)\to D^B_f(\SS)$
 makes the category
$D^b\Mod^{B\times\Gm}_f(\SS)$ 
 a mixed  version of $D^B_f(\SS)$.
We leave  the proof of this latter statement to the reader
(one can either compare the $\Ext$-groups in the  two categories
directly, or else
use Koszul duality  to eventually reduce the
 comparison of $\Ext$'s to Proposition \ref{big_small}(i)).

\ab This way, Theorem \ref{Psi_equiv} may
be suggestively reinterpreted as  follows

\begin{corollary}
\label{Psi_equiv2} The pair
$(D^b\Coh^{G\times\Gm}(\NN)\,,\,v),$
where  the functor `$v$' is defined as a composite
$$v:\;D^b\Coh^{G\times\Gm}(\NN)\stackrel{{\mathbf{f}}}\too
\dcoh^G(\NN)\stackrel{Q'}\iso D^b\cat\,,
$$
is a mixed version of the triangulated category $D^b\cat$.\qed
\end{corollary}

\subsection{Mixed  perverse sheaves.} In the rest 
of this section we will be working with algebraic varieties
over  $\fff$, an algebraic closure of a finite field $\fff_q$, where $q$
is prime to $\ell$.
Thus, the loop Grassmannian $\Gr$ is viewed as an ind-scheme defined
over $\fff_q$. Let $D^b(\Gr)$ be the
constructible derived category of (compactly  supported)
$\ell$-adic sheaves on $\Gr$. Following Deligne \cite{De},
we also consider $D^b_{\tt{mixed}}(\Gr),$ the category 
of {\em mixed} $\ell$-adic constructible complexes on $\Gr$, cf. \cite{BBD}
for more references.
We  have an obvious forgetful functor $v: D^b_{\tt{mixed}}(\Gr)
\map D^b(\Gr).$

\ab For any $M,N\in D^b_{\tt{mixed}}(\Gr)$,  the vector space
 $\Ext\hdot_{_{D^b(\Gr)}}(vM , vN)$ comes equipped
with
a natural (geometric) Frobenius
action, preserving the Ext-degree. 
Below, we will be abusing the notation slightly and
 will be writing $\Ext\hdot_{_{D^b(\Gr)}}(M,N)$
 instead of $\Ext\hdot_{_{D^b(\Gr)}}(vM , vN)$; this cannot lead to 
a confusion since the category we are dealing
with at any particular moment
is
always indicated in the subscript of the Ext-group in 
question.

\ab
Deligne  results from \cite{De} imply
\begin{proposition}\label{weil}
  For any $M_1,M_2\in D^b_{\tt{mixed}}(\overline{\Gr}_\la),$
all the absolute values of eigenvalues of the (geometric) Frobenius action
on the spaces $\Ext_{D^b(\Gr)}\hdot(M_1,M_2)$ are integral
powers of $q^{1/2}$.\qed
\end{proposition}

\ab 
According to \cite[p. 519]{BGS},
there is a well-defined  abelian  subcategory
 $\Perv\mix(\overline{\Gr}_\la)\sset D^b_{\tt{mixed}}(\overline{\Gr}_\la)$
(it is the
category denoted by $\widetilde{\mathcal{P}}$ in \cite[above Theorem 4.4.4]{BGS}),
 which is a mixed abelian category.
Moreover, by \cite[Lemmas 4.4.1(2), 4.4.6, 4.4.8]{BGS}, one has
\begin{proposition}\label{bgs_p}   For each  $\la\in\Y^{++},$ 
the following holds:

\ab \vi  The
forgetful functor
$v: D^b_{\tt{mixed}}(\overline{\Gr}_\la)\map D^b(\overline{\Gr}_\la)$
restricts to a  functor
  $v:\Perv\mix(\overline{\Gr}_\la)$ $\map \Perv(\overline{\Gr}_\la),$
that makes the category $\Perv\mix(\overline{\Gr}_\la)$
a mixed version of $\Perv(\overline{\Gr}_\la)$.

\ab \vii
The 
category 
 $\Perv\mix(\overline{\Gr}_\la)$ has enough projectives. Moreover,
for any projective $P\in \Perv\mix(\overline{\Gr}_\la),$
the object $\,vP\,$ is  projective in $\Perv(\overline{\Gr}_\la)$.
\qed
\end{proposition}

\ab It was further shown in \cite{BGS} that, for each $\mu\leq \la$,  the standard
perverse sheaf $\Delta_\mu$ and the costandard perverse sheaf 
$\nabla_\mu$ admit canonical lifts to $\Perv\mix(\overline{\Gr}_\la)$.
Abusing the notation, we will denote these lifts by  $\Delta_\mu$ and 
$\nabla_\mu$ again.

\ab 
Below, we will make use of the following
\begin{proposition}\label{wakimoto3} \vi Any projective $P\in \Perv\mix(\overline{\Gr}_\la)$
has a $\Delta$-flag. 

\ab\vii   For any $\la,\mu\in \Y$ we have
${\mathcal{M}}_\la^\vee\star \Delta_\mu \in\Perv(\Gr)$,
is a perverse sheaf.
\end{proposition}
\begin{proof} Part (i) is standard, cf. e.g. \cite[Theorem 3.2.1]{BGS}.
The proof of part (ii) is entirely similar to the proof
of Proposition \ref{wakimoto}(i), and is left to the reader.
\end{proof}

\ab Next, we put
$\dis\Perv\mix(\Gr):=\underset{\la\in\Y^{++}}\limind\Perv(\overline{\Gr}_\la).$
  From Proposition \ref{bgs_p} we deduce that   $\Perv\mix(\Gr)$  is
a mixed abelian category, moreover, it is a mixed 
 version of the category $\Perv(\Gr)$.
In particular, the simple objects of $\Perv\mix(\Gr)$ are parametrized by  elements of
$\Z\times\q$. 

\ab The category $\Perv\mix(\Gr)$ does not {\em a priori}
have enough projectives (although, it turns out,
 {\em a posteriori} that it does).
Nonetheless, given a simple
object $L\in \Perv\mix(\Gr)$ and any
  $\la\in\Y^{++},$ such that $\supp L\sset \overline{\Gr}_\la)$, we can use
Proposition \ref{bgs_p}(ii) to find an indecomposable
projective cover $P^\la\onto L,$ in the 
category 
 $\Perv\mix(\overline{\Gr}_\la)$.
The uniqueness of projective covers yields, for any
pair $\la < \mu$ (in $\Y^{++}$), a morphism
$P^\mu\to P^\la$. We conclude that the objects $\{P^\la\}_{\la\in\Y^{++}}\,$
form an inverse system, and we put
$P(L):= \limproj P^\la$. This
is a  pro-object in the category $\Perv\mix(\Gr)$,
which is a projective cover of $L$. If $L$ is pure of
weight $w(L)=n$, then the weight filtration $\Wm\idot  P(L)$
is well-defined, moreover, we have $\Wm_n P(L)=P(L)$.
This way one proves
\begin{lemma}\label{cover}  Any object of  $M\in \Perv\mix(\Gr)$ is
a quotient of a projective pro-object $P$ such that

\ab\vi $\Wm_iP/\Wm_{-i}P\in \Perv\mix(\Gr)$, for any $i\in\Z$, and

\ab\vii If $\Wm_nM=M$, then we have $\Wm_nP=P.$\qed
\end{lemma}

\begin{remark} Note that in the situation of the Lemma one
typically has $\Wm_i P\neq 0$ for all $i\ll 0,$
in general.
 $\quad\lozenge$\end{remark}

\ab An important role in the construction below will be played
by the following

\begin{lemma}\label{EE2}
  For any projective objects $P_1,P_2\in \limproj\Perv\mix(\Gr)$, and
$\la,\mu\in\Y^{++},$
we have 
$$
\Ext^n_{D^b(\Gr)}\bigl(P_1\,,\,\W_\la\star P_2\star
IC_\mu\bigr)=0\quad\text{for any}\quad n\neq 0.
$$
\end{lemma}
\begin{proof} Fix $\la$ and $\mu$ as above.
We claim first that
 $\W_\la\star vP_2\star IC_\mu$ is
a perverse sheaf.
To see this, we note that for $\la\in\Y^{++},$
we have $\W_\la={\mathcal{M}}^\vee_\la$.
Further, by Proposition \ref{wakimoto3}(i),
the projective $P_2$ admits a  $\Delta$-flag.
Hence, we are reduced to proving that
${\mathcal{M}}^\vee_\la\star \Delta_\nu\star IC_\mu$ is
a perverse sheaf, for any $\nu\in\Y$.
But, part (ii) of  Proposition \ref{wakimoto3}
says that ${\mathcal{M}}^\vee_\la\star \Delta_\nu
\in \Perv(\Gr)$. Hence, $\W_\la\star vP_2=
{\mathcal{M}}^\vee_\la\star vP_2\in\Perv(\Gr)$,
by Proposition \ref{wakimoto3}(ii), and therefore
$(\W_\la\star vP_2)\star  IC_\mu$
is a perverse sheaf by Gaitsgory's Theorem \ref{denis}.

\ab To complete the proof, we use Proposition \ref{BGS} to
obtain
$$\Ext^n_{D^b(\Gr)}\bigl(P_1\,,\,\W_\la\star P_2\star
IC_\mu\bigr)=\Ext^n_{\Perv(\Gr)}\bigl(P_1\,,\,\W_\la\star P_2\star
IC_\mu\bigr).
$$
Now, since $vP_1$ is a projective, the Ext-group
on the RHS above vanishes for all $n\neq 0$, and we are done.
\end{proof}

\begin{definition}\label{DD}
Let $\DD\mix(\Gr),$ resp., $\DD(\Gr),$
be the full subcategory in the
homotopy  category of complexes in $\limproj \Perv\mix(\Gr),$ resp., in 
$\limproj \Perv(\Gr)$,
whose objects are  complexes $C\hdot=(\ldots\to C^i\to C^{i+1}\to\ldots)$
such that

\pb{$vC^i$  is a projective pro-object in
$\Perv(\Gr)$, for any $i\in\Z$, and $C^i=0$ for $i\gg 0$;}

\pb{$H^i(C\hdot)=0$ for $i\ll 0,$ moreover,
$vH^i(C\hdot)\in\Perv(\Gr)$,
for any $i\in\Z$.}

\ab In the $\DD(\Gr)$-case, one has to replace $vC\hdot$ by $C\hdot$
in the two conditions above.
\end{definition}

\ab 
  From Lemma \ref{cover} one derives the following
\begin{corollary}\label{DD_cor} \vi
The natural functor $\Theta\mix: \DD\mix(\Gr)\map D^b\Perv\mix(\Gr),$
resp. the functor $\Theta: \DD(\Gr)\map D^b\Perv(\Gr),$
is an equivalence.  \qed
\end{corollary}
\subsection{From perverse sheaves on $\Gr$ to 
coherent sheaves on $\NN$.} The rest of this section is mostly devoted to constructing
the equivalence $P$  on the right of diagram \eqref{sum-up}.
Our strategy will be as follows. 

\ab First, we consider the
$\Z\times\Y^{++}$-graded algebra
$\CA$, introduced  in \S\ref{another},
and the abelian category $\mca$. We
will construct a functor
$\DD\mix(\Gr) \map D^b\mca$.
Then, we will  invert  the equivalence of Corollary \ref{DD_cor}, 
 and form the following composite functor
\begin{align}\label{comp_P}
P:\ D^b\Perv\mix(\Gr)\stackrel{(\Theta\mix)^{^{-1}}}\iso
\DD\mix(\Gr) \map D^b\mca\stackrel{\ff}\iso
D^b\Coh^{G\times\Gm}(\NN),
\end{align}
where  $\ff: \mca\map\Coh^{G\times\Gm}(\NN)$
is the natural functor that has been considered in  \S\ref{another}.

\ab To proceed further, we need the following

\begin{definition}\label{thin_def} An object $M=\bigoplus_{(i,\la)\in\Z\times\Y\;}M^i_\la
\in \op{Mod}^{G\times\Gm}(\CA)$ is said to be 
{\em thin} if there exists $\mu\in\Y$ such that
$M^i_\la=0$ for all $\la\in \mu+\Y^{++}$.
\end{definition}

\ab Thin objects  form a Serre subcategory 
in $\op{Mod}^{G\times\Gm}(\CA)$,
which will be denoted $\mcaf$.
The  functor $\ff$ sends any object of  the subcategory $\mcaf$ to zero,
and, for any $M\in \mca$, one clearly has
$$\bigoplus_{\la \in
\q^{++}}\;\Gamma\bigl(\NN\,,\,\ff(M)\otimes
\on(\la)\bigr)\;\cong\; M\quad\text{in}\quad
\mca\big/\mcaf.
$$
Furthermore, an  equivariant version of the classical result of
Serre implies that the functor $\ff$ induces
an equivalence
$$\ff:\ \mca\big/\mcaf\;\iso\; \Coh^{G\times\Gm}(\NN).
$$

\ab Our goal in this section is to prove the following result
which is, essentially, a mixed analogue of Theorem \ref{P_equiv0}.

\begin{theorem}\label{phi_thm}
The  functor $P$ in \eqref{comp_P}
is an equivalence of triangulated categories,
such that $P(M\star \p{V})=V\otimes P(M)\,,\,\forall V\in\rep(G),$ moreover,
for any $\la\in\pm\Y^{++},$ we have
$P(\bW_\la)=\on(\la)$.
\end{theorem}

\begin{remark} The reason we restrict,
in the last statement of the Theorem,
 to the case of $\la\in\pm\Y^{++}$
 is that we do not know,
for  general $\la\in\Y$, whether the object $\bW_\la\in\Perv(\Gr)$ admits
a natural lifting to $\Perv\mix(\Gr)$
(this would follow from `standard conjectures' saying that the
property that the Frobenius action be semisimple is preserved
under direct and inverse image functors).  In the special
case of $\la\in\Y^{++},$ resp. $\la\in-\Y^{++},$
we have $\bW_\la=\nabla_\la,$ resp. $\bW_\la=\Delta_\la,$
and such a lifting is then afforded by \cite{BGS}.
This is the lifting that we are using in  the Theorem above.
$\quad\lozenge$\end{remark}

\subsection{An Ext-formality result.}
We view the sky-scrapper sheaf $\one$
as a simple object (of weight zero)  in $\Pmix$. 
Using Lemma \ref{cover}, we construct
a projective resolution $\ldots\to  \pone^{-1}\to \pone^0 \onto \one$,
where each $\pone^i$ is a pro-object in  $\Pmix$.
Moreover, since
$\Pmix$ is known by \cite{BGS} to be a {\em Koszul category,} 
 one may choose the resolution in such a way
that $\Wm_i\pone^i=\pone^i$, for all $i=0,-1,-2,\ldots$.
Thus, the direct sum $\pone:=\oplus_{i\leq 0}\,\pone^i$
is a well-defined  pro-object in   $\Pmix$.
The differential in the resolution makes $\pone$ a dg-object,
equiped with a quasi-isomorphism $\pone\qisto\one$. 
Of course, $\pone$ is a projective, hence,
 $v\pone$ is a projective pro-object in $\Perv(\Gr)$,
by Proposition \ref{bgs_p}(ii).

\ab
Recall that for any pair of objects $M_1,M_2\in D^b_\J(\Gr),$
we have defined
a vector space
$\EE(M_1,M_2)$, see \eqref{J}. 
Applying the construction to $M_1=v\pone^i$ and 
$M_2=v\pone^j\,,\,i,j=0,-1,-2,\ldots,$
 we thus get a  $\Z$-graded algebra
$\EE\hdot(\pone,\pone):=
\bigoplus_{n\in\Z}\;\EE^n(\pone,\pone),$
where
\beq\label{EE1}
\EE^n(\pone,\pone)
:=\!\bigoplus_{\left\{
(i,j)\in\Z^2\;|\;i-j=n\,,\,
\la\in\Y^{++}\right\}}\!
\Ext_{_{D^b(\Gr)}}^0(\pone^i\,,\,\W_\la\starb \pone^j\starb\R).
\eeq

Notice, that the compatibility of formulas \eqref{EE1} and
\eqref{J} is insured by Lemma \ref{EE2}, which yields 
$$\Ext_{_{D^b(\Gr)}}^k(\pone^i\,,\,\W_\la\starb \pone^j\starb\R)=0
\quad\text{for all}\quad k\neq 0.
$$
Further, the commutator with the differential $d: P\hdot\to P^{\hhdot+1}$
induces a differential $d: \EE^n(\pone,\pone)\to
\EE^{n+1}(\pone,\pone),$ thus makes $\EE(\pone,\pone)$
a dg-algebra.

\ab Recall next that by Proposition \ref{weil},
 all the absolute values of egenvalues of the (geometric) Frobenius action
on  $\EE(\pone,\pone)$ are integral
powers of $q^{1/2}$. These  integral
powers give rise to an additional $\Z$-grading 
$\EE^i(\pone,\pone)=\bigoplus_{j\in\Z}\;\EE^i_{ j}(\pone,\pone)$
which is preserved by the differential
and is compatible with the algebra structure.
Thus, 
the object $\EE(\pone,\pone)$ becomes a differential {\em bi-graded}
algebra.
We write $H\bigl(\EE(\pone,\pone)\bigr)=\bigoplus H^i_{ j}\bigl(\EE(\pone,\pone)\bigr)$
for the corresponding cohomology algebra
considered as a differential bi-graded algebra with induced bi-grading,
and with zero differential.

\begin{proposition}\label{EEformal} \vi   For any $i\neq j$, we have
$H^i_{ j}\bigl(\EE(\pone,\pone)\bigr)=0.$
Futhermore,
there is a canonical
graded algebra isomorphism
$$\bigoplus_{i\in\Z}\;H^i_i\bigl(\EE(\pone,\pone)\bigr)\;\simeq\;
\EE(\one,\one)\;\left(\,=
\bigoplus_{i\in\Z}\;\,\Bigl(\bigoplus_{\la\in\Y^{++}\;}
\Ext_{_{D^b(\Gr)}}^i(\one,\W_\la\star\R)\Bigr)\right).
$$

\ab \vii The dg-algebra $\EE(\pone,\pone)$
is {\sl formal},
that is, there exists a {\sl{bi-graded}} algebra quasi-isomorphism
$\sigma: \EE(\pone,\pone)\qisto\EE(\one,\one).$
\end{proposition}

\ab To prove the Proposition, we will use some standard
yoga from \cite{De}.

\ab Let $A=\oplus_i\;A^i$
be a \dg-algebra equipped with an additional `weight' grading
$A^i=\bigoplus_{k}\;A^i_{ k}$ which is preserved by the differential,
i.e., is such that $d: A^i_{ k} \to A^{i+1}_{ k}.$
Thus, the `weight' grading on $A$ induces a grading
$H^i(A)=\oplus_k\;H^i_{ k}(A)$,
on each  cohomology group. 
The following formality criterion is (implicitly)
contained in \cite[5.3.1]{De}
(cf. also the  proof of \cite[Corollary 5.3.7]{De}).

\begin{lemma}\label{deligne} 
If $H^i_{ k}(A)=0$ for all  $i\neq k$, then
the \dg-algebra $A$ is formal, i.e., is quasi-isomorphic to
$\bigl(H(A)\,,\,d=0\bigr)$ as a bigraded algebra.\qed
\end{lemma}

\begin{proof}[Proof of Proposition \ref{EEformal}.]
  From isomorphisms \eqref{lhs_thm}
we find
$$
\EE\hdot(\one,\one)=\bigoplus\nolimits_{\la\in\Y^{++}\;}
\Ext_{_{D^b(\Gr)}}\hdot(\one,\W_\la\star\R)
=\bigoplus\nolimits_{\la\in\Y^{++}\;}H^{\hhdot+\hht(\la)}(i_{-\la}^!\R).
$$
By the pointwise purity of the intersection cohomology
of affine Schubert varieties (proved at the end of \cite{KL1})
the space on the right is pure,
that is, for any $j$, all  the eigenvalues  of the   
geometric Frobenius
action on $H^j(i_{-\la}^!\R)$ have  absolute value $q^{j/2}$.
 Hence, we obtain
$\EE^i_j(\one,\one)=0$ unless~$i=j$.

\ab
To complete the proof of part (i), we observe that 
$\pone\qisto\one$ is a projective resolution
and there are no non-zero Ext's between the objects
$\pone^i$ and $\W_\la\star
\pone^j\star\R$ (by Lemma \ref{EE2}).
It follows that the cohomology groups of the
dg-object $\Hom(\pone\,,\,\W_\la\star
\pone\star\R)$ are by definition,
the Ext-groups in the category $\Perv(\Gr)$.
Further, Proposition \ref{BGS} says that we may
replace Ext-groups in the category $\Perv(\Gr)$ by
those in the ambient category $D^b(\Gr)$.
The isomorphism of part (i) follows.

\ab
We deduce from the isomorphism and from the
first paragraph of the proof that
$H^i_{ j}\bigl(\EE(\pone,\pone)\bigr)=0$
for any $i\neq j$. This completes the
proof of part (i).
Part (ii) is  an immediate consequence of Lemma \ref{deligne}.
\end{proof}

\subsection{Re-grading functor.} For any $M\in\Perv(\Gr)$,
the group $G$ acts
naturally on the algebra $\EE(M,M)$ by graded algebra
automorphisms, as explained in \S\ref{regular_sheaf}.
So, we may consider the category $\op{Mod}^G(\EE(M,M))$
of $G$-equivariant algebraic $\EE(M,M)$-modules,
see Notation \ref{cohG}.
In particular, let $\Co(\EE(\pone,\pone))$ be 
the homotopy category of $G$-equivariant algebraic 
differential {\em bi-graded} finitely-generated
$\EE(\pone,\pone)$-modules (the group $G$ preserves each bi-graded
component and the differential 
acts on an object $K\in
\Co(\EE(\pone,\pone))$
as follows $d: K^i_j \to K^{i+1}_j$).
The following is a variation of \cite{BGG}.

\ab  Let $K=(K^i_j)\in \Co(\EE(\one,\one)\bigr)$.
Then, for any element $a\in \EE^n(\one,\one)$,
the $a$-action sends $K^i_j$ to $K^{i+n}_{j+n}$,
since $\EE^n(\one,\one)=\EE_n^n(\one,\one),$ by Proposition \ref{EEformal}.
Therefore, for each integer $m\in\Z$,
the subspace $\GG^m(K):= \bigoplus_{i\in\Z}\,K^{i}_{i+m}\sset
\bigoplus_{i,j\in\Z}\,K^{i}_j$ is  $\EE(\one,\one)$-stable.
We put a $\Z$-grading  $\GG^m(K)=\bigoplus_{i\in\Z}\,\GG^m(K)_i$ 
on this subspace by
 $\GG^m(K)_i:=K_{i}^{i+m}.$
Further, the differential $d: K_{i}^{i+m}\to K_{i}^{i+m+1}$, on $K$,
gives  rise to $\EE(\one,\one)$-module
morphisms $\GG^m(K)\to\GG^{m+1}(K)$ which preserve
the above defined gradings.
This way, an object $K=(K^i_j)\in
\Co(\EE(\one,\one)\bigr)$
gives rise to a complex
$\GG(K)=\bigl(\ldots\to\GG^m(K)\to\GG^{m+1}(K)\to\ldots\bigr),$
of graded $\EE(\one,\one)$-modules.
The assignment $K\mto \GG(K)$ thus defined
yields a functor
$$\GG:\ \Co(\EE(\one,\one)\bigr)\too
D^b\op{Mod}^{G\times\Gm}(\EE(\one,\one)).
$$

\subsection{Construction of the functors $P$ and $\tP$.}\label{functor_Phi}
We begin with constructing a functor
\beq\label{psi_fun1}
\Psi:\ \DD\mix(\Gr) \too \Co(\EE(\pone,\pone))
\eeq
as follows.

\ab View an object of $\DD\mix(\Gr)$
represented by a
 complex $\ldots\to C^i\to C^{i+1}\to\ldots$
as a dg-sheaf $C=\bigoplus_i\,C^i\in\limproj\Perv\mix(\Gr)$.
To  such a $C$, we associate
a graded vector space
$
\EE(\pone,C)=\bigoplus_{n\in\Z\;}\EE^n(\pone,C)$, where
$$
\EE^n(\pone,C)
:=\!\bigoplus_{\left\{
(i,j)\in\Z^2\;|\;j-i=n\,,\,
\la\in\Y^{++}\right\}}\!
\Ext_{_{D^b(\Gr)}}^0(\pone^i\,,\,\W_\la\starb C^j\starb\R).
$$
(By Lemma \ref{EE2} we know that
$\Ext_{_{D^b(\Gr)}}^k(\pone^i\,,\,\W_\la\starb C^j\starb\R)=0$,
for all $k\neq 0,$ since $C^j$ is a projective).
Further, the formula $u\mto d_C\ccirc u- u\ccirc d_\pone$,
where $d_C,d_\pone$ are the differentials on $C$ and on $\pone$,
respectively, induces a differential $\EE\hdot(\pone,C)\to
\EE^{\hhdot+1}(\pone,C)$.
 There is also  a natural
$\EE(\pone,\pone)$-module structure on the space  $\EE(\pone,C)$
coming from
the pairing \eqref{Ext2}.
Finally,  the 
geometric Frobenius action makes $\EE(\pone,C)$
a  bigraded
vector space.
Summarizing, the space $\EE(\pone,C)$ has a natural differential bigraded
$\EE(\pone,\pone)$-module structure.
The assignment $\Psi: C\mto \EE(\pone,C)$ thus obtained
gives the desired functor $\Psi$ in \eqref{psi_fun1}.

\ab Further, we compose $\Psi$ with the pull-back functor induced
by the quasi-isomorphism $\sigma$ of  Proposition \ref{EEformal}(ii).
In view of part (i) of  Proposition \ref{EEformal} we thus
 obtain a functor

\beq\label{functor_Psi}
 \DD\mix(\Gr) \stackrel{\Psi}\too \Co(\EE(\pone,\pone))
\stackrel{\sigma^*}\iso
\Co(\EE(\one,\one)\bigr).
\eeq

\ab Next, we exploit Theorem \ref{ranee} that provides a
$G$-equivariant
$(\Z\times\q^{++})$-graded
algebra isomorphism, cf. \S\ref{another},
\beq\label{phi1}
\CA\cong\EE(\one,\one)\; \left(=
\bigoplus\nolimits_{\la \in \q^{++}\;}
\Ext_{_{D^b(\Gr)}}^{2\hhdot}(\one\,,\,\W_\la\starb\R)\right).
\eeq
This isomorphism induces a category equivalence
$\tau:\ \Mod^{G\times\Gm}\bigl(\EE(\one,\one)\bigr)\iso$
$\mca$.
We  define a functor $\Phi$, cf. \eqref{comp_P},
as the following composite functor: 
\begin{align}\label{Phi}
\Phi:\ \DD\mix(\Gr) \stackrel{\Psi}\too
&\Co(\EE(\pone,\pone))
\stackrel{\sigma^*}\iso
\Co(\EE(\one,\one)\bigr)\\
&\stackrel{\GG}\too
D^b\op{Mod}^{G\times\Gm}(\EE(\one,\one))
\stackrel{\tau}\iso D^b\mca.\nonumber
\end{align}
This is the functor that gives the middle arrow in diagram
\eqref{comp_P}.

\ab Finally, we define  $P:= \ff\ccirc\Phi\ccirc(\Theta\mix)^{^{-1}},
$ the functor used in the statement of Theorem
\ref{phi_thm}; explicitly, this is the following composite functor:
\begin{align}\label{Pfinal}
P:\ D^b&\Perv\mix(\Gr)\stackrel{(\Theta\mix)^{^{-1}}}{\stackrel{_\sim}{\tooo}}
\DD\mix(\Gr)  \stackrel{\Psi}\too
\Co(\EE(\pone,\pone))\stackrel{\sigma^*}\iso
\Co(\EE(\one,\one)\bigr)
\\
&\quad
\stackrel{\GG}\too
D^b\op{Mod}^{G\times\Gm}(\EE(\one,\one))
\stackrel{\eqref{phi1}}\iso D^b\mca\stackrel{\ff}\iso
D^b\Coh^{G\times\Gm}(\NN).
\nonumber
\end{align}

\ab   Recall next that we have introduced in Definition \ref{DD}
two homotopy categories, $\DD\mix(\Gr)$ and $\DD(\Gr)$.
So far, we have only worked with the category $\DD\mix(\Gr)$,
since that category serves, by Corollary \ref{DD_cor}, as a replacement of
$D^b\Perv\mix(\Gr)$, the category that we are interested in.

\ab 
Now, however, it will be convenient for us to start working with
the category $\DD(\Gr)$ instead. We observe that  the construction
of the functor $\Phi$ given in
 \S\ref{functor_Phi} applies with obvious modifications,
such as replacing double gradings by single gradings,
to produce a functor
\beq\label{Phi1}
\Phi': \ \DD(\Gr)\too D_f^G(\CA).
\eeq
\begin{remark} As opposed to the construction of $\Phi$,
in  the construction of $\Phi'$
the step involving the "re-grading" functor $\GG$ 
should be skipped. Note also that we still use
(as we may) the formality statement in Proposition \ref{EEformal}(ii),
since we may exploit the mixed structure on $\pone$.
$\quad\lozenge$
\end{remark}

\ab 
Further, inverting the equivalence
$\DD(\Gr)\cong D^b\Perv(\Gr)$ of
Corollary \ref{DD_cor} and mimicking \eqref{comp_P}, cf. also \eqref{Pfinal},
we define  the following composite functor
\begin{align}\label{comp_P2}
\tP: D^b\Perv(\Gr)&\stackrel{\Theta^{-1}}\iso
\DD(\Gr)\stackrel{\Phi'}\map D_f^G(\CA)\stackrel{\ff}\map
 D_{\text{coherent}}^G(\NN).
\nonumber
\end{align}

\subsection{Properties of the functor $\tP$.} 
An advantage of considering the ``non-mixed'' setting
 is that, for any $\la\in\Y,$ there
is a well-defined functor  $\W_\la\star(-)$ on $D^b\Perv(\Gr)$.

\begin{proposition}\label{W_comm} We have $\tP(\one)=\on$; Furthermore,
for any $M\in D^b\Perv(\Gr),$ there is a natural isomorphism
$$\tP(\W_\mu\star M\star \p{V})=V\otimes\on(\mu)\otimes\tP(M)\,,\quad
\forall\mu\in\Y\,,\,V\in\rep(G).
$$
\end{proposition}

\begin{proof} The isomorphism $\tP(\one)=\on$ is immediate from the definition.

\ab To prove the second statement, we fix
$\mu\in\Y$ and an object  of $\DD(\Gr)$ represented
by a single pro-object $M\in \Perv(\Gr)$,
such that $vM$ is projective in $\Perv(\Gr)$.
To compute $\tP(M)$, we have to consider
the $\EE(\pone,\pone)$-module 
$\EE(\pone,M)=\bigoplus_{\la\in\Y^{++}\,}\EE(\pone,M)_\la,$ 
where $\EE(\pone,M)_\la=\Ext_{_{D^b(\Gr)}}^0(\pone,$
$\,\W_\la\star
M\starb\R).$

\ab Now, let $V\in\rep(G).$ To  compute $\tP(\W_\mu\star M\star \p{V}),$
we 
must replace the object $\W_\mu\star M\star \p{V}$ by
a projective resolution $C\qisto \W_\mu\star M\star \p{V},$ 
viewed as a dg-object $C=\oplus_i\,C^i$,
and set
\beq\label{EE5}
\EE(\pone,C)=\bigoplus_{\la\in\Y^{++}\,}\EE(\pone,C)_\la,
\enspace\text{where}\enspace
\EE(\pone,C)_\la=\bigoplus_{(i,j)\in\Z^2\;}
\Ext_{_{D^b(\Gr)}}^0(\pone^i\,,\,\W_\la\starb C^j\starb\R)
\eeq
Thus,  $\EE(\pone,C)$ is a dg-module
 over the dg-algebra $\EE(\pone,\pone)$.
The differential on $\EE(\pone,C)$ is equal to
$d=d_{_\pone}+d_{_{C}},$ a sum of two
(anti-)commuting differentials, the first being
induced from the differential on
$\pone$, and the second  from  the differential on
$C$. Each of the two  differentials clearly
preserves the weight decomposition on the left of
\eqref{EE5}.

\ab In order to compare the objects $\EE(\pone,M)$
and $\EE(\pone,C)$, we now
 prove the following
\begin{claim}\label{claimW}
  For any $\la\in\Y^{++}$ such that
$\la+\mu\in\Y^{++}$,
the dg-vector space $\left(\EE(\pone,C)_\la\,,\,d\right)$ is canonically quasi-isomorphic
to the dg-vector space $\left(V\otimes\EE(\pone,\W_\mu\star
M)_\la\,,\,d_{_\pone}\right).$
\end{claim}

\begin{proof}[Proof of Claim.] 
We observe first that 
the functor $\RHom_{_{D^b(\Gr)}}(\pone\,,\,\W_\la\star(-)\starb\R)$
applied
to the complex $C$, gives rise to a standard spectral
sequence
\beq\label{EE6} E_2=\Ext_{_{D^b(\Gr)}}\hdot(\pone\,,\,\W_\la\star
H\hdot(C)\starb\R)
\quad\Longrightarrow\quad
H\hdot\left(\Ext_{_{D^b(\Gr)}}\hdot(\pone\,,\,\W_\la\star C
\starb\R)\,,\,d_{_{C}}\right).
\eeq
Since $C\qisto \W_\mu\star M\star \p{V}$, is a resolution,
we get $H\hdot(C)=\W_\mu\star M\star \p{V}$. 
Therefore, the $E_2$-term on the left of \eqref{EE6} equals
$\Ext_{_{D^b(\Gr)}}\hdot(\pone\,,\,\W_\la\star \W_\mu\star
M\star\p{V}\star\R)$. This last Ext-group
is canonically isomorphic
to $V\otimes\Ext_{_{D^b(\Gr)}}\hdot(\pone\,,\,\W_\la\star \W_\mu\star
M\star\R)$, due to \eqref{nonsense2}.

\ab 
A key point is, that our  assumption: $\la+\mu\in\Y^{++}$
(combined with the fact that $M$, being  projective,
has a $\Delta$-flag)
implies, by Proposition
\ref{wakimoto3}, that $\W_{\la+\mu}\star M$ is a  {\em perverse sheaf}.
Hence, $\W_\la\star \W_\mu\star M\star\R=
(\W_{\la+\mu}\star M)\star\R$ is also a perverse sheaf,
by Gaitsgory's theorem.
Further, since $\pone$ is a projective, we obtain
$$\Ext_{_{D^b(\Gr)}}^n(\pone\,,\,\W_\la\star \W_\mu\star
M\star\R)=\Ext_{_{\Perv(\Gr)}}^n(\pone\,,\,\W_{\la+\mu}\star
M\star\R)=0\quad\text{for all}\enspace n\neq 0.
$$
Thus, for the $E_2$-term in \eqref{EE6}, we obtain
$$E_2=V\otimes\Ext_{_{D^b(\Gr)}}^0(\pone\,,\,\W_\la\star \W_\mu\star
M\star\R)=V\otimes\EE(\pone,\W_\mu\star M)_\la.
$$
This implies that
the spectral sequence in \eqref{EE6} degenerates, i.e., reduces to
 the following long exact sequence
$$\ldots\stackrel{d_{_{C}}}\too
\EE(\pone,C^2)_\la\stackrel{d_{_{C}}}\too\EE(\pone,C^1)_\la
\stackrel{d_{_{C}}}\too\EE(\pone,C^0)_\la\stackrel{d_{_{C}}}\too
V\otimes\EE(\pone,\W_\mu\star M)_\la
\too 0.
$$
The long exact sequence yields a canonical
quasi-isomorphism
$\bigl(\EE(\pone,C)_\la\,,\,d_{_{C}}\bigr)\qisto
\bigl(V\otimes\EE(\pone,\W_\mu\star M)_\la\,,\,d_{_{\pone}}\bigr).$
Claim \ref{claimW} now follows from another standard spectral sequence,
the one for a bicomplex, in which $\EE(\pone,C)_\la$
is viewed as a bicomplex 
with two differentials, $d_{_\pone}$ and $d_{_{C}}$.
\end{proof}

\ab Next, recall that we are given a weight $\mu\in\Y^{++}$, and  put 
\begin{align*}
&\EE(\pone,C)^\circ:=
\bigoplus\nolimits_{\{\la\in\Y^{++}\;|\;\la+\mu\in \Y^{++}\}}\,\EE(\pone,C)_\la,
\quad\text{and}\\
&\EE(\pone\,,\,\W_\mu\star M)^\circ:=
\bigoplus\nolimits_{\{\la\in\Y^{++}\;|\;\la+\mu\in \Y^{++}\}}\,
\EE(\pone\,,\,\W_\mu\star M)_\la.
\end{align*}
It is clear that $\EE(\pone,C)^\circ$ is a  $\EE(\pone,\pone)$-submodule
in  $\EE(\pone,C)$,
and $\EE(\pone\,,\,\W_\mu\star M)^\circ$  is a  $\EE(\pone,\pone)$-submodule
in  $\EE(\pone\,,\,\W_\mu\star M).$
By Claim \ref{claimW}, we have a quasi-isomorphism
$\EE(\pone,C)^\circ\qisto V\otimes\EE(\pone\,,\,\W_\mu\star M)^\circ$.
Further, the quotients 
$\EE(\pone,C)/\EE(\pone,C)^\circ$ and
$\EE(\pone\,,\,\W_\mu\star M)/\EE(\pone\,,\,\W_\mu\star M)^\circ,$
are 
both {\em thin} $\EE(\pone,\pone)$-modules,  by Definition
\ref{thin_def}.
Thus, 
we have established the following quasi- isomorphism
$$
\EE(\pone,C) \simeq V\otimes\EE(\pone,\W_\mu\star M)
\quad\text{in}\enspace D^b\left(\mca\big/\mcaf\right).
$$

\ab To complete the proof of the Proposition, assuming
that $\la+\mu\in\Y^{++}$, we compute
$\EE(\pone,\W_\mu\star M)_\la=\EE(\pone,\W_\la\star\W_\mu\star M)=
\EE(\pone, M)_{\la+\mu}$.
Thus, the graded space $\EE(\pone,\W_\mu\star M)$ $=
\bigoplus_{\la\in\Y^{++}\;}\EE(\pone,\W_\mu\star M)_\la$
is isomorphic, up to a thin subspace,
to the space $\EE(\pone, M)=\bigoplus_{\la\in\Y^{++}\;}\EE(\pone,M)_\la,$
with the $\Y$-grading being shifted by $\mu$.
But shifting by  $\mu$ is the same as
tensoring by $\k(\mu)$.
Therefore, in $D^G_{\text{coherent}}(\NN)$, we have
\begin{align*}
\ff\left(\EE(\pone,C)\right)\, \simeq\,
\ff\left( V\otimes\EE(\pone\,,\,\W_\mu\star M)\right)\,
 &\simeq\, V\otimes\ff\left(\k(\mu)\otimes\EE(\pone,M)\right)\\
&=\, V\otimes\on(\mu)\,\otimes\,\ff\left(\EE(\pone,M)\right),
\end{align*}
and the Proposition is proved.
\end{proof}

\begin{proposition}\label{key_ext}
  For any $\la\in\Y^{++}$ and $\L\in \PO,$
 the functor
$\tP$ induces an isomorphism
$$\Ext\hdot_{_{D^b(\Gr)}}(\one,\, \W_\la\star\L)\iso 
\Ext\hdot_{_{D^G_{\text{coherent}}(\NN)}}\bigl(\tP(\one)\,,\,\tP(\W_\la\star\L)\bigr).$$
\end{proposition}
\begin{proof} We  may write $\L=\p{V}$, for some $V\in\rep(G).$ Then,
using formulas \eqref{lhs_thm} and \eqref{key_compare},
we obtain
\beq\label{ExtLHS}
\Ext\hdot_{_{D^b(\Gr)}}(\one,\, \W_\la\star\p{V})=
H^{\hhdot+\hht(\la)}(i_{-\la}^!\p{V})=\gr^W_{k+\hht(\la)}V(-\la).
\eeq

On the other hand,  Proposition \ref{W_comm} yields
 $\tP(\W_\la\star \p{V})=V\otimes\on(\la)\otimes\tP(\one)=
V\otimes\on(\la),$  and also $\tP(\one)=\on$.
Since $\la\in\Y^{++}$, we have
$$\Ext^k_{D_{\text{coherent}}(\NN)}\bigl(\on,\,\on(\la)\bigr)=
H^k(\NN, \on(\la))=0\quad\forall k>0.
$$
Hence, for $\L=\p{V}$, we find
\begin{align}\label{EED}
\Ext\hdot_{_{D^G_{\text{coherent}}(\NN)}}\bigl(\tP(\one)\,,\,\tP(\W_\la\star
\p{V})\bigr)&=\Ext\hdot_{_{D^G_{\text{coherent}}(\NN)}}\bigl(\on\,,\,
V\otimes\on(\la)\bigr)\nonumber\\
&=\Hom_G\left(\k\,,\,V\otimes 
\Ext\hdot_{_{D_{\text{coherent}}(\NN)}}\bigl(\on,\,\on(\la)\bigr)\right)\nonumber\\
&=\left(V\otimes\Gamma\hdot\bigl(\NN,\on(\la)\bigr)\right)^G\\
\text{(by second equality in \eqref{BR_fin})}\quad&=\left(V\otimes
\bigl(\gr^{^{_{\mathsf{fib}}}}_{\!\hhdot+\hht(\la)}(\indf_H^Ge^\la)\bigr)\right)^G
\nonumber\\
\text{(by first equality in \eqref{BR_fin})}\quad&=
\gr^W_{\!\hhdot+\hht(\la)}\left(\bigl(V
\otimes(\indf_H^Ge^\la)\bigr)^G\right).
\nonumber
\end{align}

\ab Further, using Frobenius reciprocity, we obtain
$\bigl(V\otimes(\indf_H^Ge^\la)\bigr)^G=V(-\la)$.
Hence, the last line in \eqref{EED} equals
$\gr^W_{\!\hhdot+\hht(\la)}V(-\la)$,
which is exactly the RHS of \eqref{ExtLHS}.

\ab We leave to the reader to check that
the isomorphism between the LHS and RHS  of \eqref{ExtLHS}
that we have constructed above is the same as the one given by the map in
\eqref{ExtLHS}.  The Proposition is proved.
\end{proof}

\begin{lemma}\label{vanish_lemma} Let $\la\in\Y$ be such
that $\la\not\geq 0$, i.e., $\la$ does not
belong to the semi-group in $\Y$ generated by positive
roots. Then, we have  

\vi $\;\RHom_{_{D^b(\Gr)}}(\one,\bW_\la)=0$;
\quad and \quad 
\vii $\;\RHom_{_{D_{\op{coherent}}^G(\NN)}}\bigl(\on\,,\,\on(\la)\bigr)=0$.
\end{lemma} 

\begin{proof} To prove (i) recall that, for any $\mu,\nu\in\Y$,
we have   $\W_\mu\starb \bW_\nu=\bW_{\mu+\nu}$.
It follows, since $(\W_{-\mu}\star\W_\mu)\starb(-)$
is the identity functor
that, for any $\mu\in\Y$,
we have   $\Ext\hdot_{_{D^b(\Gr)}}(\one,\bW_\la)=
\Ext\hdot_{_{D^b(\Gr)}}(\bW_{\mu},\bW_{\mu+\la})$.
We choose $\mu$ to be anti-dominant and such that
$\mu+\la$ is also anti-dominant. Then we know that
$\bW_\mu= \Delta_\mu$ and $\bW_{\mu+\la}=\Delta_{\mu+\la}$.
Hence, writing $j_\mu: \Gr_\mu\into \Gr$ for the imbedding,
 we get
\begin{align*}
\RHom\hdot_{_{D^b(\Gr)}}(\one,\bW_\la)&=
\RHom\hdot_{_{D^b(\Gr)}}(\bW_\mu,\bW_{\mu+\la})\\
&=
\RHom\hdot_{_{D^b(\Gr)}}(\Delta_{\mu},\Delta_{\mu+\la})=
j_\mu^!\Delta_{\mu+\la},\enspace\forall\mu\ll 0.
\end{align*}
Now, the condition  $\la\not\geq 0$ implies that $\Gr_\mu\not\subset
\overline{\Gr}_{\mu+\la},$
for all
sufficiently anti-dominant $\mu$.
This forces $j_\mu^!\Delta_{\mu+\la}=0$, and part (i) is proved.

\ab \ab To prove (ii),  we first use
 the chain of equivalences  in the first
line of \eqref{naive22} to reduce the Ext-vanishing in the
category $\dcoh^G(\NN)$ by a similar  Ext-vanishing in the
category $D^B_f(\La)$. Further, 
recall the algebra
$\AA=\Ub\ltimes\La$ 
and the triangulated category
$D^\Ub_\Y(\AA,\La),$ see   Notation \ref{cohG} and definitions following it.
By Proposition \ref{big_small},
the Ext-groups are unaffected if category $D^B_f(\La)$
is replaced 
by $D^\Ub_\Y(\AA,\La)$. This way, we obtain natural isomorphisms
\beq\label{Ext_vanish1}
\Ext\hdot_{_{{\dcoh^{G}(\NN)}}}\bigl(\on\,,\,\on(\la)\bigr)=
\Ext\hdot_{_{D^B_f(\n)}}\bigl(\imath^*\on\,,\,\imath^*\on(\la)\bigr)=
\Ext\hdot_{_{D^\Ub_\Y(\AA,\La)}}\bigl(\k_{_\AA}\,,\,
\k_{_\AA}(\la)\bigr).
\eeq
Thus, we are reduced to showing that, for all $i\in \Z$,
one has: $\Ext^i_{_{D^\Ub_\Y(\AA,\La)}}\bigl(\k_{_\AA}\,,\,
\k_{_\AA}(\la)\bigr)=0$. 

\ab To this end, we use the   standard spectral sequence
\eqref{sp_seq} for the cohomology of a semidirect product:
$$
\Ext^p\bigl(\k_{_\Ub}\,,\,
\Ext^q_{_{\La}}(\k_{_\AA}\,,\,
\k_{_\AA}(\la))\bigr)=E_2^{p,q}
\quad\Longrightarrow\quad E_\infty^{p,q}=
\gr\Ext^{p+q}_{_{D^\Ub_\Y(\AA,\La)}}\bigl(\k_{_\AA}\,,\,
\k_{_\AA}(\la)\bigr)\,,
$$
where the group $\Ext^p\bigl(\k_{_\Ub},-)$ on the left
is computed in the category of $\Y$-graded
$\Ub$-modules such that the action of the Cartan subalgebra
$\mathfrak{t}\sset \b$ is compatible with the $\Y$-grading.
By (an appropriate version of)
formula \eqref{SL_kos} we have  $\Ext\hdot_{_{\La}}(\k_{_\La}\,,\,
\k_{_\La}(\la))\bigr)=\uSS(\la),$
furthermore, the above spectral sequence collapses.
Thus we get
$$
\gr\Ext^p_{_{D^\Ub_\Y(\AA,\La)}}\bigl(\k_{_\AA}\,,\,
\k_{_\AA}(\la)\bigr)=
\big[H^p(\n\,,\,\SS(\la))\big]^{(0)},
$$
where we write $[\ldots]^{(0)}$ for the zero-weight component
with respect to the $\Y$-grading.

\ab Recall now that  the Lie algebra cohomology
 on the RHS above are given by  the cohomology
of the Koszul complex $\wedge\hdot\n^*\otimes \SS(\la)$.
Observe that any weight in
$\wedge\hdot\n^*\otimes\SS=\wedge\hdot\n^*\otimes\sym(\n^*[-2])$
is clearly a sum of {\em negative} roots.
Therefore, the  zero-weight component of
the cohomology in the RHS above vanishes, since
$\la$ is not a sum of positive roots, by our assumptions.
Thus, $\gr\Ext^p_{_{D^\Ub_\Y(\AA,\La)}}\bigl(\k_{_\AA}\,,\,
\k_{_\AA}(\la)\bigr)=0,$ for all $p\in \Z$,
and (ii) is proved.
\end{proof}

\ab The proof of the next Lemma 
is easy and will be left to the reader.
\begin{lemma}\label{generate2}The smallest triangulated subcategory
in $\DD(\Gr)$ that contains the following set of objects
$$ X:=\big\{\bW_\la\,,\,\forall \la\not\geq 0;\quad
\W_\nu\star \L\,,\,\forall \nu\in\Y^{++}\,,\,\L\in\PO\big\}
$$
coincides with the category $\DD(\Gr)$ itself.\qed
\end{lemma}

\subsection{Proof of Theorems \ref{P_equiv0} and \ref{phi_thm}.}
We claim first that, for any $M\in D^b\Perv(\Gr),$
the functor $\tP$ induces an isomorphism
\beq\label{one_ok}
\Ext\hdot_{_{D^b(\Gr)}}(\one,M)\iso 
\Ext\hdot_{_{\dcoh^G(\NN)}}\bigl(\tP(\one)\,,\,\tP(M)\bigr).
\eeq
Due to Lemma \ref{generate2}, it suffices to prove
\eqref{one_ok} for all objects $M\in X$.
  For  the objects $M=\overline{\W}_\la\,,\,\la\not\geq 0,$
equation \eqref{one_ok} is insured by Lemma
\ref{vanish_lemma}.   For  the objects $M=\W_\la\star\L,$
where $\la\in\Y^{++}$ and $\L\in\PO,$ 
equation \eqref{one_ok} follows from Proposition
\ref{key_ext}. Thus, \eqref{one_ok} is proved.

\ab Now, let $\mu\in\Y$ and $M\in  D^b\Perv(\Gr).$
  From Proposition \ref{W_comm}, using that $\tp(\one)=\on,$ 
 we obtain a natural commutative diagram
$$
\xymatrix{
\Ext\hdot_{_{D^b(\Gr)}}(\overline{\W}_\mu\,,\,M)\ar[rr]^<>(0.5){\tP}
\ar@{=>}[d]^<>(0.5){\W_{-\mu}\star(-)}&&
\Ext\hdot_{_{\dcoh^G(\NN)}}\bigl(\tP(\overline{\W}_\mu)\,,\,\tP(M)\bigr)
\ar@{=>}[d]^<>(0.5){\on(-\mu)\otimes(-)}\\
\Ext\hdot_{_{D^b(\Gr)}}(\one\,,\,\W_{-\mu}\star M)\ar[rr]^<>(0.5){\tP}
&&
\Ext\hdot_{_{\dcoh^G(\NN)}}\bigl(\tP(\one)\,,\,\tP(\W_{-\mu}\star
M)\bigr)
}
$$
The vertical maps in the diagram are isomorphisms
since the functors $\W_{-\mu}\star(-)$,
resp. $\on(-\mu)\otimes(-)$, are equivalences
of derived categories.
The map in the bottom line of this diagram is already known to be an
isomorphism, by \eqref{one_ok}. Hence,
the map in the top line is also an isomorphism.
But, the set of objects $\,\{\overline{\W}_\mu\}_{\mu\in\Y}\,$
clearly generates the category $D^b\Perv(\Gr).$
Therefore, we deduce that, for any $N,M\in D^b\Perv(\Gr),$
the functor $\tP$ induces an isomorphism
$$
\Ext\hdot_{_{D^b(\Gr)}}(N,M)\iso
\Ext\hdot_{_{\dcoh^G(\NN)}}\bigl(\tP(N)\,,\,\tP(M)\bigr).
$$
Thus, we have proved that the functor $\tP$ is fully faithful.

\ab To prove that $P$ is  fully faithful,
we recall that the category $\Perv\mix(\Gr)$,
resp., $D^b\Coh^{G\times\Gm}(\NN)$, is a mixed version of
$\Perv(\Gr)$,
resp., of  $\dcoh^G(\NN)$. It is clear from the construction
of the functors $P$ and $\tP$ that one has
an isomorphism of functors $\tP\ccirc v=v\ccirc P$.
Thus, using \eqref{mix_ext}, we obtain,
for any $M,N\in D^b\Perv\mix(\Gr),$ a natural commutative square
$$
\xymatrix{
{{\bigoplus_{n\in\Z}}\;\Ext\hdot_{D^b\Perv\mix(\Gr)}(M,N\langle n\rangle)}
\ar[rr]^<>(0.5){\sim}_<>(0.5){v}\ar[d]^<>(0.5){P}&&
\Ext\hdot_{D^b\Perv(\Gr)}(M,N)\ar@{=>}[d]^<>(0.5){\tP}\\
{{\bigoplus_{n\in\Z}}\;
\Ext\hdot_{D^b\Coh^{G\times\Gm}(\NN)}\bigl(P(M)\,,\,P(N)\langle
n\rangle\bigr)}
\ar[rr]^<>(0.5){\sim}_<>(0.5){v}&&
\Ext\hdot_{\dcoh^G(\NN)}\bigl(\tP(M)\,,\,\tP(N)\bigr)
}
$$
We have already proved that the vertical map on the right 
is an isomorphism.
It follows from the diagram, that the vertical map on the 
left is also an isomorphism.

\ab The
set of objects $\,\{\on(\mu)=\tP(\overline{\W}_\mu)\}_{\mu\in\Y}\,$
clearly generates the category $\dcoh^G(\NN)$.
Thus, from Lemma \ref{abstract_nonsense} we deduce,
using
Proposition \ref{W_comm}, that the functor
 $\tP$ is an equivalence. 
That completes the proof of Theorem  \ref{P_equiv0}.

\ab Finally,  any simple object of $\Coh^{G}(\NN)$ has the form
$v(\mathcal{F})$, where $\mathcal{F}$ is a  simple object of
$\Coh^{G\times\Gm}(\NN)$. We deduce, using that  $\tP$
 is an equivalence and applying Lemma \ref{abstract_nonsense}
once again, that  $P$ is also an equivalence.
Theorem \ref{phi_thm} is proved.\qed 

\subsection{Equivalence of abelian categories.} We now combine
Theorem \ref{phi_thm} and Theorem \ref{Psi_equiv} together, and compose
the inverse of the equivalence $Q': \dcoh^G(\NN) \iso$
$ D^b\cat$
with  the inverse of the equivalence $P': D^b\Perv(\Gr)\iso
\dcoh^G(\NN)$.
This way, we obtain the following composite equivalence of
triangulated categories
\beq\label{ppp}
\xymatrix{
\Upsilon: \; D^b\cat\;\ar[rr]^<>(0.5){(Q')^{-1}}&&\;\dcoh^G(\NN)\;
\ar[rr]^<>(0.5){(P')^{-1}}&&\; D^b\Perv(\Gr).
}
\eeq

\ab The triangulated
categories $D^b\cat$ and $D^b\Perv(\Gr)$, each has a natural
$t$-structure, with cores $\cat$ and $\Perv(\Gr)=\Perv_{\J}(\Gr)$, respectively.
Below, we will prove the following result.

\begin{theorem}\label{add1}
The equivalence $\Upsilon:  D^b\cat\iso
D^b\Perv(\Gr)$
 respects the t-structures, hence
induces an equivalence
of abelian categories $\p:\ \Perv(\Gr)\iso \cat,$
 such that
$\,\p L_\lambda =IC_\lambda\,,$ for any $\lambda\in \Y$.
\end{theorem}

\begin{remark} 
 Recall that functor $\phi^*:\rep(G)\to\cat$,
$ M\mapsto\fr M,$  see (\ref{frobenius}),
identifies $\rep(G)$ with
 the full subcategory $\phi^*(\rep(G))\sset \cat$.
Further, it follows from the properties of the functors $P'$ and $Q'$
proved in the previous sections that, for any $V\in\rep(G),$ the
functor $\Upsilon$ in \eqref{ppp} sends the $\U$-module
$\fr{V}$ to $\p{V}$, the perverse sheaf corresponding to
$V$ via the Satake equivalence.
Therefore, the functor $\Upsilon$ 
maps  the subcategory $\phi^*(\rep(G))$
into the subcategory
$\PO\sset\Perv_{_\J}(\Gr)$. Thus, the  restriction of
 $\Upsilon$ to the abelian category $\cat$
may be regarded, in view of Theorem \ref{add1}, as  a natural `extension'
of the functor $\p: \rep(G)\map \Perv(\Gr)$ to
the larger category $\cat$, i.e.,
 one has the following commutative diagram:
\beq\label{add2}
\xymatrix{
{\rep(G)\enspace}\ar@{^{(}->}[rrr]^<>(.5){\text{\tt Frobenius}\enspace\phi^*}
\ar@{=>}[d]_<>(.5){{\text{\footnotesize{\tt Satake}}}\atop{\text{
\footnotesize{\tt  equivalence}}}}^<>(.5){\p}
&&&\enspace\cat
\ar@{=>}[d]^<>(.5){{\text{\footnotesize{\tt  
equivalence}}}\atop{{\Upsilon=P^{-1}\circ Q,\;{\text{\footnotesize{\tt see}}\; (\ref{abelian_eq})}
}}}\\
{\PO\enspace} \ar@{^{(}->}[rrr]^<>(.5){\text{\tt inclusion}}&&&{\enspace\Perv_{_\J}(\Gr)}
}
\eeq
  For this reason, we use the notation $\p$ for the
functor $\Upsilon\big|_{\cat}$.$\quad\lozenge$.
\end{remark}

\ab To prove Theorem \ref{add1}, we will use 
 filtrations on triangulated categories $D^b\cat$ and  $D^b\Perv(\Gr)$,
defined 
as follows.   For any $\la\in\q^{++}$,
let $D^b_{\!_{\leq\la}}\cat$ be the smallest full
triangulated subcategory of  $D^b\cat$
that contains all the simple objects
$L_\mu\in \cat$ with $\mu\leq \la$.
Clearly, we have   $D^b_{\!_{\leq\la}}\cat \subset
D^b_{\leq\nu}\cat $, whenever $\la\leq \nu$.
Moreover, the category 
$D^b_{\!_{\leq\la}}\cat/D^b_{\!_{<\la}}\cat$ is semisimple
and is formed by direct sums of objects of the type
$L_\la[k]\,,\,k\in\Z$.

\ab Similarly, let $D^b_{\!_{\leq\la}}\Perv$
 be the  full
triangulated subcategory of  $D^b\Perv(\Gr)$
formed  by the objects supported on
$\overline{\Gr}_\la$. 

\begin{lemma}\label{perv_equiv_lemma}   For any $\la\in\q^{++}$,
we have

\vi The functor $\Upsilon$ induces an equivalence
 $D^b_{\!_{\leq\la}}\Perv\iso
D^b_{\!_{\leq\la}}\cat$. Moreover,

\vii The induced functor: 
$$
D^b_{\!_{\leq\la}}\Perv\big/D^b_{\!_{<\la}}\Perv
\too
D^b_{\!_{\leq\la}}\cat\big/D^b_{\!_{<\la}}\cat
$$
sends the class of $IC_\la$ to the class of $L_\la$.
\end{lemma}

\proof  We know that, for any $\la\in\q^{++}$,
the functor $\Upsilon$ sends, by construction, the object
$\rind_\B^{^\U}({l}\la)$ into $\overline{\W}_\la$.
Further, fix $\la\in\Y$ and let $w\in W$ be the element of minimal length
such that the weight $w({l}\la)$ is dominant.
Then by
Lemma \ref{borel_weil} we have  
$\tR^{\ell(w)}\ind_\B^{^\U}({l}\la)
=\text{Weyl}_{w({l}\la)},$ is a Weyl module,
and for any $j\neq \ell(w),$ each simple subquotient
of $\tR^j\Ind_\B^{^\U}({l}\la)$ is isomorphic
to $L_\mu$ with $\mu<\la$.
It follows that the category
$D^b_{\!_{\leq\la}}\cat$ is generated by the objects
$\,\big\{\rind_\B^{^\U}({l}\nu)\big\}_{\nu\leq\la}.$
Hence, by Corollary
 \ref{wakimoto2}(ii), we get:
$\supp\Upsilon\bigl(\rind_\B^{^\U}({l}\nu)\bigr)=
\supp\W_\nu\subset \overline{\Gr}_\la$.
Thus, the functor $\Upsilon$ takes
 $D^b_{\!_{\leq\la}}\cat$ into $D^b_{\!_{\leq\la}}\Perv$.
Similarly, the functor $\Upsilon^{-1}$ takes
$D^b_{\!_{\leq\la}}\Perv$ into  $D^b_{\!_{\leq\la}}\cat$,
and part (i) of the Lemma follows.

\ab 
To prove (ii), we use Lemma \ref{borel_weil}, that has been already
exploited
above. The
 Lemma  says in particular that
in $D^b_{\!_{\leq\la}}\cat/D^b_{\!_{<\la}}\cat$ one has
an isomorphism 
$\rind_\B^{^\U}({l}\la)\simeq L_{w({l}\la)}[-\ell(w)]$.
Similarly, in $D^b_{\!_{\leq\la}}\Perv/D^b_{\!_{<\la}}\Perv$  one has
an isomorphism $\W_\la\simeq IC_\la[-\ell(w)]$,
by  Corollary
 \ref{wakimoto2}(ii). Statement (ii) now follows from
(i) and the equation
$\Upsilon\bigl(\rind_\B^{^\U}({l}\la)\bigr)=\W_\la$.~\qed

\ab We will need  the following result
on the gluing of $t$-structures, proved in [BBD].

\begin{lemma}\label{bbd}
Let $\,(D, D_{\!_{\leq\la}})\,$ be a triangulated
category equipped with filtration such that each quotient category
$D_{\!_{\leq\la}}/D_{\!_{<\la}}$ is a semisimple category generated by one
object. Then there exists a unique $t$-structure on $D$
compatible with given $t$-structures on each  category
$D_{\!_{\leq\la}}/D_{\!_{<\la}}$.\qed
\end{lemma}

\begin{proof}[Proof of Theorem  \ref{add1}.]
 Lemma \ref{perv_equiv_lemma}(i) implies that the functor $\Upsilon$
is an equivalence of  triangulated
categories which is, moreover, compatible 
 with filtrations.
Thus, the part of Theorem  \ref{add1} concerning
the $t$-structures follows from Lemma \ref{bbd}.

\ab Alternatively, the same thing can be proved as follows.
The  $t$-structure on $D=D^b\Perv(\Gr)$ is
characterized by the property that,
for any $\la\in\q$, the objects
$\,\Delta_\la:=(j_\la)_!\C_{_{\Gr_\la}}[-\dim\Gr_\la]$
and $\,\nabla_\la:=(j_\la)_*\C_{_{\Gr_\la}}[-\dim\Gr_\la],$
where $j_\la: \Gr_\la \into \Gr$ is the imbedding of the
Bruhat cell $\Gr_\la$,
belong to the core  of the  $t$-structure.
But the object $\Delta_\la$
is completely determined by the following conditions
formulated entirely in terms of the filtration
$D^b_{\!_{\leq\la}}\Perv$, and the  $t$-structures
on $D^b_{\!_{\leq\la}}\Perv/D^b_{\!_{<\la}}\Perv$:

$\bullet\quad$ $\Delta_\la\in D^b_{\!_{\leq\la}}\Perv$;

$\bullet\quad$ $\Hom_{_{D^b\Perv(\Gr)}}\bigl(\Delta_\la\,,\,D^b_{\!_{<\la}}\Perv\bigr)=0$;

$\bullet\quad$ $\Delta_\la=IC_\la\;\mbox{\sl mod}\;D^b_{\!_{<\la}}\Perv$.

\ab Similarly, 
there is a dual characterisation of the $\nabla_\la$'s.
The above characterizations yield the  compatibility 
of the $t$-structures. They
imply also that
$\p L_\lambda =IC_\lambda,$ for any $\la\in\Y.$
This completes the proof of Theorem  \ref{add1}.
\end{proof}

\ab Recall that
convolution of perverse sheaves gives an exact bi-functor
$\star:\Perv_{_\J}(\Gr)  \times \PO \to \Perv_{_\J}(\Gr)$. On the other hand, note that
if $V\in \rep(G)$ and $M\in \cat$,  then
$M \otimes_{_{\C}}\fr V \in \cat$. Thus the tensor product of
$\U$-modules gives an exact bi-functor
$\otimes: \cat \times \rep(G) \to \cat$.

\ab From Theorem \ref{add1} one derives the following result.

\begin{proposition} \label{add3}
 The above bi-functors $\;\otimes:  \cat \times \rep(G)
\to \cat\,,$ and $\enspace$ $\star: \Perv(\Gr)  \times \PO
 \to \Perv(\Gr)$ correspond to each other
under the equivalences of Theorems \ref{add1} and  \ref{tensor},
that is, for any $V\in \rep(G)$ and $M\in \cat$, we have:

\ab\vi Functorial
isomorphism $\p(M \otimes\fr V) \simeq \p M \star \p(\fr V)$.

\ab\vii Functorial  vector space isomorphism
$\,
H\hdot(\Gr\,,\,\p(M \otimes\fr V))\;\simeq\;
V\otimes H\hdot(\Gr\,,\,\p M).$
\end{proposition}

\subsection{Duality and convolution.}
Recall that $\Gr=G^\vee(\K)/G^\vee(\oo)$. 
By an Iwasawa decomposition,
we may identify $\Gr$ with a based loop group $\Omega$, cf. [PS].
The inversion map $\sigma: g\mapsto g^{-1}$ on that based loop group
induces an auto-equivalence: $\M\mapsto\sigma^*\M$, on $D^b(\Gr)$.
Let also: $\M\mapsto{\mathbf{D}}\M$ denote the Verdier duality
functor on $D^b(\Gr)$. The above two functors 
obviously commute,
and we will write: $\,\M \mto \M^\vee := \sigma^*({\mathbf{D}}\M)=
{\mathbf{D}}(\sigma^*\M)\,$
for the composition, an involutive {\it contravariant}
auto-equivalence on $D^b(\Gr)$.

\ab The functor: $\M \mapsto \M^\vee$
corresponds,
via the equivalence of Theorem \ref{tensor},
to the duality: $V\mapsto V^\vee$ on the
category $\rep(G)$, i.e.,
there is a functor isomorphism $\p(\fr V^\vee)=\bigl(\p(\fr M)\bigr)^\vee$.
In particular, for any $\lambda\in\Y^{++}$, we have  
$(IC_\lambda)^\vee = IC_{-w_\circ(\lambda)}$,
where $w_\circ$ is the longest element in the Weyl group $W$.

\ab One has the  standard isomorphism
\beq\label{standard}
\Hom\left(\L_1 \star {\mathcal M}, \,{\mathcal L}_2\right)\;\simeq\;
\Hom\left(\L_1 , \,{\mathcal L}_2 \star ({\mathcal M}^\vee)\right)
,\qquad\forall\L_1 \,,\,{\mathcal L}_2\,,\,\M\in D^b(\Gr)\,.
\eeq

\section{Quantum group cohomology and the loop Grassmannian}\label{Gr}
\ab The goal of this section is to prove the conjecture
given in [GK, \S4.3]. We put $\k=\C$.

\subsection{Equivariant coherent sheaves on $\N$.}
The group  $G\times\Gm$ acts naturally on the nil-cone $\N\sset\g$
(the group $G$ acts by conjugation and the group $\Gm$ by dilations).
Therefore the  coordinate ring $\k\hdot[\N]$
 comes equipped
 with a natural grading,
and  with a $\g$-module structure.

\ab We consider  the abelian category 
 of $G$-equivariant
$\Z$-graded finitely generated $\k[\N]$-modules. 
It  may be identified,
since   the nil-cone $\N$ is an affine variety, with
$\cohgm$, the category
of $G\times\Gm$-equivariant coherent sheaves on $\N$.

\ab Let $e\in\N$ be a fixed
principal nilpotent element,
and  $\a\sset \g\oplus\k$ be the Lie algebra of the isotropy group
of the point $e\in\N$ in $G\times\Gm$. Given a $G\times\Gm$-equivariant
coherent sheaf $\F$ on $\N$, let $\F|_e$ denote its geometric fiber
at $e$.
Thus, $\F|_e$ is a finite dimensional vector space with a natural
$\a$-action.
Clearly, ${\gee}={\gee}\oplus\{0\}
\sset\a$, hence the space $\F|_e$ comes equipped with 
a ${\gee}$-module structure.

\ab Further, we choose an ${\mathfrak{s}\mathfrak{l}}_2$-triple for $e$, and let
$\s$ denote the corresponding semisimple element in the triple.
Then the element $\hat\s:=(\s, -2)\in \g\oplus\k$
belongs to $\a$,
since $\ad\s(e)=2 e$. The action of the semisimple
element $\hat\s$ (or rather of the
one-parameter multiplicative subgroup $\Gm$ 
generated by $\hat\s$) puts a $\Z$-grading on 
$\F|_e$. The action-map $\gee\otimes \F|_e\too \F|_e$
preserves the gradings.

\ab Recall that $G^e$ is the connected unipotent group with Lie algebra $\gee$.
Any finite-dimensional $\gee$-module has a canonical
$G^e$-module structure, by exponentiation.
We write $\emod$ for the abelian category
of finite-dimensional {\it $\Z$-graded} 
${{G^e}}$-modules  equipped with a grading compatible 
with $\Lie G^e$-action.
 Thus, in the notations above, for any $\F\in\cohgm$,
  one has: $\F|_e\; \in \;\emod$.

\subsection{Quantum group cohomology and the principal nilpotent.}
According to Proposition \ref{GBbis}, we have a
  Hopf-adjoint action  of the quantum algebra $\U$ on $\fU$.
The action preserves the kernel of the Frobenius
homomorphism $\phi: \fU\to\U$, hence induces  an $\U$-action on
$\u$, the image of the  Frobenius
homomorphism. This gives rise
to an $\U$-action on the cohomology $H\hdot(\u,\k_{\u})$.
As we have already mentioned in \S\ref{cohomology},
the  Hopf-adjoint action  of any algebra  on its
own cohomology is  trivial. Thus,
the  Hopf-adjoint action  of $\U$ on  $H\hdot(\u,\k_{\u})$
descends to the quotient
algebra $\U/(\u)$, which is isomorphic to $\Ug$ via the
Frobenius homomorphism. This makes  $H\hdot(\u,\k_{\u})$
a $\g$-module, and we have

\begin{theorem}[\cite{GK}]\label{giku}
There is a natural $\g$-equivariant graded algebra 
isomorphism 
$\,H^{2\hhdot}(\u,\k_{\u}) \simeq \k\hdot[\N]$. Moreover,
all odd cohomology groups vanish: $H^{\tt{odd}}(\u,\k_{\u})=0.$
\end{theorem}

\ab   For any $\u$-module $M$, the cohomology
$H\hdot(\u, M)=
\Ext\hdot_{\u}(\k_{\u}, M)$ has a natural 
graded $H\hdot(\u, \k_{\u})$-module structure
via the Yoneda product: 
$\,\Ext^i_{\u}(\k_{\u}, \k_{\u})\;\times\;
\Ext^j_{\u}(\k_{\u}, M)\;\longrightarrow\;$
$\Ext^{i+j}_{\u}(\k_{\u}, M)\,.$ 
One can show that $H\hdot(\u, M)$ is a finitely generated
$H\hdot(\u, \k_{\u})$-module provided $M$ is  finite dimensional over $\k$.

\ab Assume now that $M$ is a 
$\U$-module, and view it as an $\u$-module, by restriction.
Then, there is a natural $\U$-module structure on each cohomology group
$H^i(\u, M)$, that descends to
$\Ug$,
see [GK, \S5.2]. If $M$ is, in addition,  finite dimensional then
$H\hdot(\u, M)$ is finitely generated over $H\hdot(\u, \k_{\u})\simeq \k[\N]$,
hence each cohomology group
$H^i(\u, M)$ is  finite dimensional over $\k$. It follows that the
$\g$-action on $H^i(\u, M)$ may be exponentiated canonically
to an algebraic action of the group $G$. Thus,
we may (and will) view  $H\hdot(\u, M)$ as an object of the
category $\cohgm$. This way we obtain, following [GK, \S4.1],
a functor
\beq\label{functorF}
F: \;\;\rep(\U)\;
\,\longrightarrow\, \cohgm,\quad M\mto H\hdot(\u, \res^\U_\u M)\;.
\eeq
Taking the geometric fiber at $e$, and restricting
attention to $\U$-modules that belong to the principal block,
 we thus obtain the functor
$$
F_e:\,\cat \longrightarrow\,\emod\;,\;
M \mto F(M)\big|_e=H\hdot(\u,\,M)\big|_e\,.
$$

\ab The Theorem below is the main result of this section;
it has been
conjectured in  [GK,~\S4.3].

\begin{theorem}\label{functor_equiv}
There is a canonical isomorphism of the following two functors
$\cat \longrightarrow\,\emod$:
$$M\mto H\hdot(\Gr, \p M)\quad\text{and}\quad M\mto  F_e(M)\,.
$$
\end{theorem}

\ab The rest of this section is devoted to the proof of this Theorem,
which requires some preparations.

\subsection{Induction.} We are going to relate cohomology over
the algebras $\u$ and $\U$.
To this end, consider the smooth-coinduction functor
$\ind=\ind_{\mathsf{u}}^{^{\U}}\,.$ By [APW${}_{2}$, Theorem 2.9.1]
we have

\begin{theorem}\label{APK} \vi Each of the categories
$\Umod$ and $\umod$, has  enough projectives. Moreover,
in each category,
projective and injective objects coincide. 

\ab \vii Furthermore, the functor 
${\res}\,: \Umod\to\umod$ takes projectives into projectives,
and the functor $\ind: \umod \to \limind\Umod$ is exact.\qed
\end{theorem}

\ab Next, given $M_{_{\sf U}}
\in \Umod$,  choose a projective resolution of $M_{_{\sf U}}$, and apply
the adjunction isomorphism \eqref{adjun1} to that resolution
 term by term. This way, the result of Andersen-Polo-Wen 
above yields, for any $N_\u\in\umod$,
a canonical graded space isomorphism
\beq\label{adjunction2}
\Ext\hdot_{\u}({\res}\,M_{_{\!{{\sf U}}}}\,,\,N_\u) \simeq 
\Ext\hdot_{_{\!{{\U\mmod}}}}
(M_{_{\!{{\sf U}}}}\,,\,\ind\,
N_\u)
\,.
\eeq

\ab We can now express the functor $M\mto H\hdot(\u,\,M)$
in terms of cohomology of $\U$ rather than~$\u$.

\begin{lemma}\label{first_reduction}  For any $M\in\Umod$, there 
is a natural  isomorphism $
H\hdot(\u,\, M)=\Ext\hdot_{_{\!{{\U\mmod}}}}
\Bigl(\k_{_{{{\sf U}}}}
\,,\,M\otimes_{_{\k}}\fr\k[G]\Bigr)\;.$
\end{lemma}

\noindent
{\sl Proof:} Write $\Ind:=\Ind^\U_\u$ and $\res:=\res^\U_\u$.
  For any $M\in\Umod$, we calculate:
\begin{align*}
H\hdot(\u,\, M)& =
\Ext\hdot_{\u}\bigl({\res}\,\k_{_{{{\sf U}}}}\,,
\,{\res}\,M\bigr)\\
\mbox{ by (\ref{adjunction2})}\quad &=\Ext\hdot_{_{\!{\U\mmod}}}
\Bigl(\k_{_{{{\sf U}}}}\,,\,\ind ({\res}\,M)\Bigr)\\
\mbox{by Lemma \ref{indres}}\quad 
&=\Ext\hdot_{_{\!{\U\mmod}}}
\Bigl(\k_{_{{{\sf U}}}}
\,,\,M\otimes \fr\k[G]\Bigr)\;.\quad\square
\end{align*}

\subsection{Quantum group cohomology via the loop Grassmannian.}
We are going to re-interpret the functor $F: \Umod\to \cohgm\,,\,
M\longmapsto H\hdot(\u, M)$ geometrically, in terms of 
perverse sheaves on $\Gr$.  

\ab Recall the setup of Theorem
\ref{main1}. We will prove the following

\begin{proposition}\label{geometry} {\sf {(i)}} There is a natural 
graded algebra isomorphism
$$H\hdot(\u,\,\k_{_{\u}})\; \simeq \;\Ext\hdot_{_{D^b(\Gr)}}(\one,\,\R)\;.$$

{\sf {(ii)}}   For any $M\in \cat$, there is 
a graded module isomorphism compatible with the algebra isomorphism of
part
(i):
$$H\hdot(\u,\,M)\quad\simeq\quad \Ext\hdot_{_{D^b(\Gr)}}(\one,\,\M\star\R)
=\Ext\hdot_{_{D^b(\Gr)}}(\R^\vee\,,\,\M)\,.
$$
\end{proposition}

\begin{remark}
It would be very interesting to give an algebraic interpretation of 
the space $\Ext\hdot\GO(\one\,,\,\R)$ in terms of the quantum group $\u$, that 
is
to find an equivariant analogue of the  isomorphism of
part
(i) of the Proposition above. 
 $\quad\lozenge$\end{remark}

{\sl Proof of Proposition \ref{geometry}.}
The category $\cat$
is known [APW] to be a direct summand of the category $\Umod$.
Hence, the $\Ext$-groups in these two categories coincide and,
 for any $M\in \cat$, we compute
\begin{align}\label{third_reduction}
H\hdot(\u, M)
& =\;
\Ext\hdot_{_{\!{{\sf U}}{\tiny -}{{\mathtt {mod}}}}}
\Bigl(\k_{_{{{\sf U}}}}
\,,\,M\otimes_{_{\k}}\fr\k[G]\Bigr)\quad(\mbox{by Lemma \ref{indres}})\nonumber\\
&=\Ext\hdot_{_{\cat}}
\Bigl(\k_{_{\sf U}}
\,,\,M\otimes_{_{\k}}\fr\k[G]\Bigr)\nonumber\\
\mbox{(by Theorem \ref{add1})}\qquad &=\;
\Ext\hdot_{_{\PJ}}
\Bigl(\p(\k_{_{\sf U}})\,,\,
\p(M\otimes_{_{\k}}\fr\k[G])\Bigr)
\\
\mbox{(by Proposition \ref{add3})}\quad
&=\;\Ext\hdot_{_{\PJ}}
\Bigl(\one\,,\,
\p M\star \p(\fr\k[G])\Bigr)\nonumber\\
&=\;\Ext\hdot_{_{\PJ}}
\left(\one\,,\,
\p M\star \R\right)\;\,\stackrel{(\ref{BGS})}{=}\;\,
\Ext\hdot_{_{D^b(\Gr)}}\left(\one\,,\,
\p M\star \R\right).\nonumber
\end{align}

This proves the isomorphism in (ii) which, in
 the special case $M=\k_{_{\sf U}}$, reduces to (i). 

\ab To see that the isomorphism in (i) is an 
{\it algebra} isomorphism we put
$M:=\k$ in \eqref{third_reduction}.
Recall that the cup product on $H\hdot(\u,\,\k_{\u})\,(\cong
\Ext\hdot_{_{\!\U\mmod}}(\k_{_{{{\sf 
U}}}},\,\fr\k[G]))$
is known, see e.g. [GK, \S5.1], to be induced by the algebra structure on $\k[G]$.
It follows that the cup-product on 
$\Ext\hdot_{_{\!\U\mmod}}(\k_{_{{{\sf 
U}}}},\,\fr\k[G])$
can be defined in a way very similar to the one we have defined a
product on $\Ext\hdot_{_{D^b(\Gr)}}(\one,\,\R)$. Specifically, 
given $x\in \Ext^i_{_{\!\U\mmod}}
(\k_{_{{{\sf U}}}},\,\fr\k[G])$, we view it as a 
"derived morphism" $x: \k_{_{{{\sf U}}}} \to \fr\k[G][i]$.
Tensoring (over $\k$) with this "derived morphism" gives, for any
$y\in \Ext^j_{_{\!\U\mmod}}
(\k_{_{{{\sf U}}}},\,\fr\k[G])$, a composition:
$$
y\cdot x:\;\k_{_{{{\sf U}}}}\,\stackrel{y}{\longrightarrow}
\,\fr\k[G][j]=(\fr\k[G]\otimes_{_{\k}}\k_{_{{{\sf U}}}})[j]
\,\stackrel{\otimes x}{\longrightarrow}\,\fr(\k[G]\otimes_{_{\k}}\k[G])[i+j]
\;\stackrel{\tiny product}{\longrightarrow}\;\fr\k[G][i+j]\,.
$$
Comparison with formula (\ref{composite}) yields compatibility of the
algebra structures on the LHS and RHS of the isomorphism in part (i).
 The rightmost isomorphism in (ii)
is due to (\ref{standard}). The Proposition is proved. \qed

\subsection{} 
 Given any finite-dimensional
graded 
${{G^e}}$-module $E$, we 
 form the algebraic
$G$-equivariant vector bundle $\,\indf^{^G}_{_{G^e}} E 
:= E \times_{_{G^e}} G$
on $G^e\backslash G =\Oe$.
 It is clear from definitions that the space of global 
regular sections of the corresponding locally free sheaf on $\Oe$
is given by the formula:
$\,\Gamma(\Oe\,,\,\indf^{^G}_{_{G^e}} E) = 
(E\otimes_{_{\k}} \k[G])^{G^e}\,.\,$
Further, let $j: \Oe\hookrightarrow\N$ denote
 the open
imbedding, and $j_*(\indf^{^G}_{_{G^e}} E)$ denote
the direct image, which is a coherent sheaf on $\N$
 since the complement
$\N\smallsetminus\Oe$ has codimension $\geq 2$ in $\N$.

\begin{proposition}\label{F} \vi There is a morphism between the
following two functors $\Umod$ $\too\cohgm$(= category of
$G\times\Gm$-equivariant $\k[\N]$-modules):
\beq\label{sheaf_morphism} 
F(M)= \Ext\hdot_{_{\!D^b(\Gr)}}(\one,\,\p M \star \R)
\;
\longrightarrow\;
\Gamma\left(\N\,,\, j_*\bigl(\indf^{^G}_{_{G^e}}H\hdot(\Gr,\,\p M)\bigr)\right)\,.
\eeq
This isomorphism  is compatible, via the isomorphism
of Proposition \ref{geometry}, with the algebra actions 
on each side.

\ab \vii
  For any semisimple object
$M\in \cat$ the morphism  (\ref{sheaf_morphism}) is
an isomorphism.
\end{proposition}
 
\proof Let $E$ be a finite dimensional $G^e$-module,
and ${\mathcal E}=j_*(\indf^{^G}_{_{G^e}} E)$
the corresponding coherent sheaf on $\N$.
We observe that by definition one has a
canonical isomorphism of the spaces of global sections:
$\Gamma(\Oe\,,\,\indf^{^G}_{_{G^e}} E)\, =\,\Gamma(\N\,,\,{\mathcal E}).$
Therefore, since $\k[\Oe] \simeq\k[\N]$,
one can alternatively {\it define}
the sheaf ${\mathcal E}$ as the sheaf on $\N$ whose space
of global sections is the vector space 
$\Gamma(\Oe\,,\,\indf^{^G}_{_{G^e}} E)=(E\otimes_{_\k}\k[G])^{G^e}$,
viewed as a $\k[\N]$-module.

\ab Recall next that, for any $M\in \cat$, the hyper-cohomology
$H\hdot(\Gr,$ 
$\p M)$ has a natural graded ${{G^e}}$-module structure.
Thus, we may view this space as a finite dimensional $G^e$-module, and 
perform the
construction of the previous paragraph. We get
\begin{align}\label{fourth_reduction}
\Gamma\left(\N\,,\,j_*\indf^{^G}_{_{G^e}}H\hdot(\Gr,\,\p M)\right)
& =\Gamma\left(\Oe\,,\,\indf^{^G}_{_{G^e}}H\hdot(\Gr,\,\p M)\right)\nonumber\\
& =\Bigl(\k[G]\otimes_{_{\k}}H\hdot(\Gr,\p M)\Bigr)^{G^e}\nonumber\\
(\text{since }\;\R=\p(\k[G]))\quad
& =\Bigl(H\hdot(\Gr,\,\R)\otimes_{_{\k}}
H\hdot(\Gr,\p M)\Bigr)^{\gee}\nonumber\\
\mbox{(by Proposition \ref{add3})}\quad
& = \Bigl(H\hdot(\Gr\,,\,\R\star \p M)\Bigr)^{\gee}\\
\mbox{(by Theorem \ref{add1})}\quad
& = \Hom\hdot_{_{H(\Gr)}}
\Bigl(H\hdot(\Gr,\,\one)\,,\,H\hdot(\Gr,\,\R\star \p M)\Bigr)\;.\nonumber
\end{align}
Notice that the  expression in the bottom line coincides with a special
case of the RHS of
formula (\ref{G1map}). Hence, the formula gives a canonical
map: 
$$\Ext\hdot_{_{\!D^b(\Gr)}}\!(\one,\,\p M \star \R)
\;
\longrightarrow\;\Hom\hdot_{_{H(\Gr)}}\Bigl(H\hdot(\one)\,,\,
H\hdot(\R\star \p M)\Bigr) = 
\Gamma\Bigl(\N\,,\,j_*\indf^{^G}_{_{G^e}}H\hdot(\p M)\Bigr)\,.
$$
It is seen easily that this map intertwines
the $\Ext\hdot_{_{\!D^b(\Gr)}}\!(\one,\,\R)$-module structure on
the RHS with the $\k[\N]$-module structure on
the LHS (using the algebra isomorphism
$\Ext\hdot_{_{\!D^b(\Gr)}}\!(\one,\,\R)\simeq\k[\N]$,
see \S7).
This yields part (i) of the Proposition.

\ab To prove part (ii) observe that, for 
any simple module
$M\in \cat$, the perverse
sheaf $\p M \star \R$ is semisimple, by
 Proposition \ref{add3} and  Decomposition theorem [BBD].
 The result now follows
 from the chain of isomorphisms in
(\ref{fourth_reduction}), and Proposition \ref{GG1}
applied to the semisimple object $\L_2=\p M \star \R$.
\qed

\subsection{Proof of Theorem \ref{functor_equiv}.} 
  Formula (\ref{sheaf_morphism}) gives, by adjunction,
a morphism: $j^*F(M)
\;
\longrightarrow\;
\indf^{^G}_{_{G^e}}H\hdot(\Gr,\,\p M)\,.$
This is a morphism of $\Ad G$-equivariant
coherent sheaves on the orbit $\Oe$.
Giving such a morphism is equivalent to giving 
a $G^e$-equivariant linear map between
the geometric fibers at $e$ of the corresponding
(locally free) sheaves. Thus, we obtain a canonical
morphism of functors
 $\vartheta: F_e(M)= F(M)\big|_e\;
\longrightarrow\;H\hdot(\Gr,\,\p M)\,.$

\ab It remains to prove that, for any
$M\in\cat$, the morphism $\vartheta$
is an isomorphism.   For $M$ simple this
is insured by Lemma \ref{F}.
In the general case, choose
a Jordan-H\"older filtration:
$0=M_0\,\subset\,M_1\,\subset\ldots\subset\,M_n=M$
with simple subquotients $M_i/M_{i-1}\,,\, i=1,\ldots,n\,,$
and put $\gr M=\bigoplus\;M_i/M_{i-1}$.
We have the standard convergent spectral
sequence: $\,H\hdot(\u,\,\gr M) = E_2\quad
\Longrightarrow\quad E_\infty= \gr H\hdot(\u,\,M)\,.$
The restriction functor of the sheaf $F(M)=H\hdot(\u,\,M)$ to the point
$e$ being {\it exact} (since the $G$-orbit of $e$ is open in $\N$,
and $F(M)$ is $G$-equivariant), the spectral sequence above induces
a  spectral sequence of the fibers:
$\,F_e(\gr M)\;\Longrightarrow\;\gr F_e(M)\,.$

\ab Similar arguments apply to the functor:
$M\mto H\hdot(\Gr,\, \p M)$, and yield a
spectral sequence: $\,E_2=H\hdot\left(\Gr,\, \p(\gr M)\right)
\;\Longrightarrow\;{\gr}H\hdot(\Gr,\, \p M)$.
Now the canonical morphism $\vartheta$ gives 
a morphism of the two spectral sequences. This morphism 
induces an isomorphism between the $E_2$-terms, due to
Lemma \ref{F}. It follows that $\vartheta$ is itself
 an isomorphism.
(Equivalently, instead of the spectral sequence argument
above, one
can use that, for each $i$, the map 
$\vartheta$ gives a morphism of two long exact sequences
of derived functors corresponding to the short
exact sequence: $0\to M_{i-1}\to M_i\to
M_i/M_{i-1}\to 0\,.$ The result then follows
by induction on $i$, by means of the 
{\it {five-Lemma}}.)
\qed

\footnotesize{

}
\vskip 1cm

\footnotesize{
{\bf S.A.}: $\,$  Department of Mathematics,
Yale  University, 10 Hillhouse Ave., New Haven, CT 06520, USA\\
\hphantom{x}\ab\, {\tt serguei.arkhipov@yale.edu}}

\footnotesize{
{\bf R.B.}: Department of Mathematics,
 Northwestern University,
Evanston, IL
60208, USA;\\ 
\hphantom{x}\ab\, {\tt bezrukav@math.northwestern.edu}}

\footnotesize{
{\bf V.G.}: Department of Mathematics, University of Chicago, 
Chicago, IL
60637, USA;\\ 
\hphantom{x}\ab\, {\tt ginzburg@math.uchicago.edu}}

\end{document}